\documentclass[10pt,english,a4paper]{article}

\usepackage{amssymb}
\usepackage{amsfonts}
\usepackage{graphicx}
\usepackage{amstext}

\usepackage{amsmath,amscd, amsthm}
\usepackage{mathrsfs}

\usepackage{import}

\usepackage{euscript}
\usepackage[latin9]{inputenc}
\usepackage{geometry}
\usepackage{pdfsync}

\usepackage{hyperref}
\usepackage{cleveref}

\usepackage[all]{xy}

\usepackage[title,titletoc]{appendix}

\usepackage{stmaryrd}

\usepackage{tikz} 
\usetikzlibrary{arrows,decorations.pathmorphing,backgrounds,positioning,fit,petri,shapes.misc,decorations.markings,shapes,shapes.multipart,calc,cd}

\tikzset{
	plain trans/.style = {circle, minimum size = 14, inner sep=2, draw, font=\scriptsize}}
\tikzset{
	plain inv trans/.style = {circle, minimum size = 18, inner sep=1, draw, font=\tiny}}
\tikzset{
	plain dual trans/.style = {circle, minimum size = 14, inner sep=2, draw, font=\scriptsize}}
\tikzset{
	plain inv dual trans/.style = {circle, minimum size = 18, inner sep=1, draw, font=\tiny}}

\tikzset{trans/.style={circle, minimum size = 14, inner sep=1, draw, font=\scriptsize,
	path picture={%
	\filldraw[anchor=center,black]
	($(path picture bounding box.south east)!0.15!(path picture bounding box.north east)$) to [out=-135, in=0]	(path picture bounding box.south) to [out=180, in=-45]
	($(path picture bounding box.south west)!0.15!(path picture bounding box.north west)$) to cycle;
      }}}
      
\tikzset{inv trans/.style={circle, minimum size = 14, inner sep=1, draw, font=\tiny,
	path picture={%
	\filldraw[anchor=center,black]
	($(path picture bounding box.south east)!0.15!(path picture bounding box.north east)$) to [out=-135, in=0]	(path picture bounding box.south) to [out=180, in=-45]
	($(path picture bounding box.south west)!0.15!(path picture bounding box.north west)$) to cycle;
      }}}      

\tikzset{dual trans/.style={circle, minimum size = 14, inner sep=1, draw, font=\scriptsize,
	path picture={%
	\filldraw[anchor=center,black]
	($(path picture bounding box.north east)!0.15!(path picture bounding box.south east)$) to [out=135, in=0]	(path picture bounding box.north) to [out=180, in=45]
	($(path picture bounding box.north west)!0.15!(path picture bounding box.south west)$) to cycle;
      }}}      

\tikzset{inv dual trans/.style={circle, minimum size = 14, inner sep=1, draw, font=\tiny,
	path picture={%
	\filldraw[anchor=center,black]
	($(path picture bounding box.north east)!0.15!(path picture bounding box.south east)$) to [out=135, in=0]	(path picture bounding box.north) to [out=180, in=45]
	($(path picture bounding box.north west)!0.15!(path picture bounding box.south west)$) to cycle;
      }}}      

\tikzset{
	func/.style = {font=\scriptsize},
	in func/.style = {font=\scriptsize, above, text opacity=0},
	in func visible/.style = {font=\scriptsize, above},
	out func/.style = {font=\scriptsize, below, text opacity=0},	
	out func visible/.style = {font=\scriptsize, below},	
	mid func/.style = {font=\scriptsize, inner sep=1, fill=white} 
}

\tikzset{
	cat/.style = {font=\scriptsize}}
\tikzset{
	frame/.style = {double,rounded corners}}

\tikzset{
	paste/.style = {dashed,opacity=0.5}}

\tikzset{
    dir/.style={},
    dir end/.style={},    
    dual/.style={},
    eq/.style={inner sep=0pt},
    plain/.style={} 
}

\newcommand{\N}{\mathbb{N}}

\newcommand{\Q}{\mathbb{Q}}
\newcommand{\R}{\mathbb{R}}
\renewcommand{\P}{\mathbb{P}}

\newcommand{\Hom}{\operatorname{Hom}}
\newcommand{\Map}{\operatorname{Map}}

\newcommand{\Fun}{\operatorname{Fun}}

\newcommand{\im}{\operatorname{Im}}

\newcommand{\iso}{\overset{\sim}{\longrightarrow}}
\renewcommand{\=}{\cong}

\newcommand{\Aut}{\operatorname{Aut}}
\newcommand{\BAut}{\operatorname{BAut}}

\newcommand{\rank}{\operatorname{rank}}

\newcommand{\colim}{\operatorname{\underrightarrow{colim}}}
\newcommand{\hocolim}{\operatorname{\underrightarrow{hocolim}}}

\newcommand{\xto}{\xrightarrow}
\newcommand{\xfrom}{\xleftarrow}

\newcommand{\bnd}{\partial}

\newcommand{\set}[1]{\left\{{#1}\right\}}

\newcommand{\diag}[1]{\vcenter{\xymatrix{#1}}}

\newcommand{\M}{\mathcal{M}}

\newcommand{\E}{\mathcal{E}}
\newcommand{\EE}{\mathbb{E}}

\newcommand{\C}{\mathcal{C}}
\newcommand{\D}{\mathcal{D}}

\newcommand{\U}{\underline{U}}

\mathchardef\mhyphen="2D 

\newcommand{\Sp}{\mathbf{Sp}}
\newcommand{\Ss}{\mathcal{S}}
\newcommand{\GTop}{\mathbf{Top}^G}

\newcommand{\vect}{\mathbf{Vect}}

\newcommand{\adj}{\leftrightarrows}

\newcommand\noloc{%
  \nobreak
  \mspace{6mu plus 1mu}
  {:}
  \nonscript\mkern-\thinmuskip
  \mathpunct{}
  \mspace{2mu}
}

\newcommand{\Cat}{\mathcal{C}\mathrm{at}}
\newcommand{\ul}[1]{\underline{#1}}

\newcommand{\ulTopG}{\ul{\mathbf{Top}}^G}

\newcommand{\OG}{\mathcal{O}_G}
\newcommand{\OGop}{{\OG^{op}}}
\newcommand{\GFin}{\ul{\mathbf{Fin}}^G}
\newcommand{\mfld}{\mathbf{Mfld}}
\newcommand{\ulmfld}{\ul{\mathbf{Mfld}}}
\newcommand{\Gmfld}{\ul{\mfld}^G}
\newcommand{\GmfldD}{\ul{\mfld}^{G,\sqcup}}
\newcommand{\disk}{\mathbf{Disk}}

\newcommand{\Gdisk}{\ul{\disk}^G}
\newcommand{\GdiskD}{\ul{\disk}^{G,\sqcup}}
\newcommand{\Sing}{\mathbf{Sing}}
\newcommand{\NSing}{\mathbf{N}}
\newcommand{\orb}{\operatorname{Orbit}}

\newcommand{\Rep}{\mathbf{Rep}}
\newcommand{\ulRep}{\ul{\mathbf{Rep}}}
\newcommand{\OGTop}{\OG\mhyphen\mathbf{Top}}
\newcommand{\EmbOG}{Emb^{\OG}}
\newcommand{\OGmfld}{\OG\mhyphen\mfld}

\newcommand{\Fin}{\mathbf{Fin}}
\newcommand{\OGFin}{\OG\mhyphen\Fin}
\newcommand{\mfldGD}{\OG\mhyphen\Fin\mhyphen\mfld}
\newcommand{\ulFin}{\ul{\Fin}}
\newcommand{\ulFun}{\ul{\Fun}}
\newcommand{\ultimes}{\ul{\times}}

\newcommand{\Arr}{\ul{\mathbf{Arr}}}

\newcommand{\Moore}{{\operatorname{Moore}}}

\newcommand{\Conf}{\mathbf{Conf}}

\newtheorem{thm}{Theorem}[subsection]
\newtheorem{prop}[thm]{Proposition}
\newtheorem{cor}[thm]{Corollary}
\newtheorem{lem}[thm]{Lemma}
\newtheorem{mydef}[thm]{Definition}

\theoremstyle{definition}
\newtheorem{observation}{Observation}
\newtheorem{rem}[thm]{Remark}
\newtheorem{ex}[thm]{Example}
\newtheorem{notation}[thm]{Notation}
\newtheorem{construction}[thm]{Construction}

\newcommand{\cof}{\hookrightarrow}

\newcommand{\trivcof}{\stackrel{\sim}{\hookrightarrow}}

\newcommand{\fib}{\twoheadrightarrow}

\newcommand{\trivfib}{\stackrel{\sim}{\twoheadrightarrow}}

\newcommand{\pullbackcorner}[1][dr]{\save*!/#1-1.2pc/#1:(1,-1)@^{|-}\restore}

\newcommand{\pullbackcornerleft}[1][dl]{\save*!/#1-1.2pc/#1:(1,-1)@^{|-}\restore}

\newcommand{\pbcorner}{\arrow[dr, phantom, "\ulcorner" description, very near start]}

\newcommand{\rightthreearrows}{%
        \mathrel{\vcenter{\mathsurround0pt
                \ialign{##\crcr
                        \noalign{\nointerlineskip}$\rightarrow$\crcr
                        \noalign{\nointerlineskip}$\rightarrow$\crcr
                        \noalign{\nointerlineskip}$\rightarrow$\crcr
                }%
        }}%
}

\newcommand{\SimplicialDiagram}[2]{
  \cdots \rightthreearrows {#1} \rightrightarrows {#2}
}

\newcommand{\Dbun}{\mathbb{D}}
\newcommand{\cDbun}{\overline{\mathbb{D}}}
\newcommand{\Sbun}{\mathbb{S}}
\newcommand{\Abun}{\mathbb{A}}

\begin{document}

\title{Genuine equivariant factorization homology}
\author{Asaf Horev}
\maketitle

\begin{abstract}
  We construct a genuine $G$-equivariant extension of factorization homology for $G$ a finite group, assigning a genuine $G$-spectrum to a manifold with $G$-action. 
  We show that $G$-factorization homology is compatible with Hill-Hopkins-Ravenel norms and satisfies equivariant $\otimes$-excision.
  Following Ayala-Francis we prove an axiomatic characterization of genuine $G$-factorization homology.
  Applications include a description of real topological Hochschild homology and relative topological Hochschild homology of $C_n$-rings using genuine $G$-factorization homology.
\end{abstract}

\tableofcontents

\section{Introduction}
Factorization homology, introduced by Lurie under the name topological chiral homology (\cite{Lurie_Cobordism_Survey}, \cite{HA}), is an invariant of an $\EE_n$-algebra and a framed $n$-dimensional manifold.
The factorization homology of a framed $n$-dimensional manifold $M$ with coefficients in an $\EE_n$-ring spectrum $A$ is a spectrum denoted $\int_M A$.
If $M$ admits an action of a finite group $G$ then $\int_M A$ admits an $G$-action by functoriality.
However, this action is defined only up to coherent homotopy, as $\int_M A$ is defined by an $\infty$-categorical colimit. 
A fundamental observation of equivariant homotopy theory is that such a ``naive'' action does not determine the homotopy type of the fixed points.
In particular the action of $G$ on $\int_M A$ does not define a \emph{genuine} $G$-spectrum structure on $\int_M A$.

The first goal of this paper is to construct and study such a genuine equivariant extension of factorization homology for a fixed finite group $G$.
We draw on two points of view in order to explain the expected properties of genuine equivariant factorization homology.

\paragraph{Factorization homology as a tensor product}
First, according to \cite[rem. 4.1.19]{Lurie_Cobordism_Survey} one can intuitively think of $\int_M A$ as a continuous tensor product $\otimes_{x\in M} A$ indexed by the points of $M$.
One should have this intuition in mind when considering the behavior of factorization homology with respect to disjoint unions\footnote{We distinguish between disjoint unions and coproducts since disjoint union is not the categorical coproduct in the category of $n$-dimensional manifolds and open embeddings which we consider below.}, namely 
\begin{align}
  \int_{M_1 \sqcup M_2} A \simeq \int_{M_1} A \otimes \int_{M_2} A .
  \label{eq:FH_disj_tensor}
\end{align}

In order to generalize this behavior to genuine $G$-spectra we now recall the interaction of the smash product with the group action. 
If $X$ is a genuine $H$-spectrum for $H<G$ a subgroup then the smash product $ \otimes_{|G/H|} X$ of $G/H$ copies of $X$ has a naive $G$-action, induced by the combining the action of $H$ on $X$ with the action of $G$ on the indexing set $G/H$.
Hill-Hopkins-Ravenel \cite{HHR} extended this naive $G$-action to a genuine $G$-spectrum, $N_H^G(X)$, the Hill-Hopkins-Ravenel norm of $X$.
More generally they define smash products indexed by finite $G$-sets as the smash product of Hill-Hopkins-Ravenel norms (see \cite[app. A.3]{HHR}).
Let $U$ be a finite $G$-set, given by a coproduct of orbits $U= \coprod\limits_{i\in I} G/H_i$ with stabilizers $H_i <G$.
The $U$-indexed smash product of a family $ X_\bullet = \set{X_i}_{i\in I}$, where each $X_i$ is a genuine $H_i$-spectrum, 
is the genuine $G$-spectrum given by the smash product of the norms $ \otimes_U X_\bullet = \otimes_{i\in I} N_{H_i}^G(X_i)$.
The indexed smash product interacts with smash products and norms as follows.
\begin{itemize}
  \item 
  The indexed smash product takes disjoint unions to smash products: 
  if $U', U''$ are of finite $G$-sets then the indexed smash product along $U' \coprod U''$ is equivalent to smash product of the indexed products, $\otimes_{U' \coprod U''} X_\bullet \simeq ( \otimes_{U'} X_\bullet) \otimes ( \otimes_{U''} X_\bullet )$.
  \item 
  The indexed smash product takes topological inductions to norms:
  given a subgroup $H<G$ and a finite $H$-set $U$, denote the quotient $G \times_H U$ by $\coprod_{G/H} U$.
  The left action of $G$ on the first coordinate makes $\coprod_{G/H} U$ a $G$-set which we call the \emph{topological induction} of $U$ from $H$ to $G$.
  The norm of an indexed product is given by an indexed product along the topological induction, \( \otimes_{\coprod_{G/H} U} X_\bullet \simeq N_H^G( \otimes_U X_\bullet) \).
\end{itemize}
Note that stating the second property required us to consider tensor products indexed by finite $H$-sets for all $H<G$.

Interpreting genuine equivariant factorization homology as a tensor product indexed by a $G$-manifold $M$, one would expect a similar behavior. 
To state it, we consider the genuine factorization homology of $H$-manifolds for all subgroups $H<G$.
Namely, for any subgroup $H<G$ and $H$-manifold $M$ we expect genuine equivariant factorization homology to assign a genuine $H$-spectrum $\int_M A\in \Sp^H$,
which interacts with smash products and norms as follows. 
\begin{itemize}
  \item 
    Genuine equivariant factorization homology takes disjoint unions to smash products:
    if $M',M''$ are $n$-dimensional $G$-manifolds then the genuine equivariant factorization homology along $M' \sqcup M''$ is equivalent to the smash products of the genuine equivariant factorization homologies along $M'$ and $M''$, 
    \begin{align*}
      \int_{M' \sqcup M'' } A \simeq ( \int_{M'} A ) \otimes ( \int_{M''} A ),
    \end{align*}
    as genuine $G$-spectra.
  \item 
    Genuine equivariant factorization homology takes topological inductions to norms:
    given a subgroup $H<G$ and an $n$-dimensional $H$-manifold $M$, denote the topological induction $G \times_H M$ by $\sqcup_{G/H} M$,
    with left $G$-action induced by acting on the first coordinate.
    The norm of genuine equivariant factorization homology along $M$ is equivalent to genuine equivalent factorization homology along the topological induction $\sqcup_{G/H} M$, 
    \begin{align}
      \label{eq:G_FH_G_SM_compatibilities}
      \int_{\sqcup_{G/H} M} A \simeq N_H^G( \int_{M} A ) ,
    \end{align}
    as genuine $G$-spectra.
\end{itemize}

\paragraph{Factorization homology as a homology theory} A second point of view on factorization homology is given by Ayala-Francis \cite{AF}, where factorization homology is considered as a homology theory of $n$-dimensional manifolds.
Ayala-Francis start from the observation that factorization homology is functorial with respect to open embeddings of framed $n$-dimensional manifolds.
Let $\mfld_n^{fr}$ be the $\infty$-category of framed $n$-dimensional manifolds and framed open embeddings, and let $\C$ be a cocomplete symmetric monoidal $\infty$-category.
Fixing an $\EE_n$-algebra $A$ in $\C$,
Ayala-Francis consider factorization homology $M \mapsto \int_M A$ as a functor of $\infty$-categories $\int_- A \colon \mfld_n^{fr} \to \C$.
Factorization homology extends to a symmetric monoidal functor $\int_- A \colon \mfld_n^{fr,\sqcup} \to \C^\otimes$ with respect to disjoint union of manifolds (a functorial version of \cref{eq:FH_disj_tensor}) under mild conditions\footnote{Namely that the tensor product in $\C$ distributes over sifted colimits.} on $\C$.

Taking the view that excision is the characterizing property of a homology theory, Ayala-Francis define a homology theory for manifolds as a symmetric monoidal functor $\mfld_n^{fr,\sqcup} \to \C^\otimes$ satisfying $\otimes$-excision, and show that $\int_- A$ is indeed such a homology theory for manifolds. 

Furthermore, they show that the Eilenberg-Steenrod characterization of generalized homology theories admits the following generalization.
Let $\mathcal{H}(\mfld_n^{fr}, \C) \subseteq \Fun^\otimes(\mfld_n^{fr}, \C)$ be the full subcategory of symmetric monoidal functors satisfying $\otimes$-excision.
\begin{thm}[Ayala-Francis]
  There is an equivalence of $\infty$-categories
  \begin{align*} 
    \int \colon Alg_{\EE_n}(\C) \iso \mathcal{H}(\mfld_n^{fr},\C), \quad A \mapsto (\int_- A \colon\mfld_n^{fr} \to \C )
  \end{align*}
  sending an $\EE_n$-algebra $A$
  to factorization homology with coefficients in $A$. 
\end{thm}
In fact, this theorem holds in greater generality, replacing framed manifolds with $B$-framed manifolds and $\EE_n$-algebras with $B$-framed $n$-disk-algebras.
The second goal of this paper is to provide such an axiomatic characterization of genuine equivariant factorization homology (see \cref{thm:G_FH_axiomatic_characterization}).

\paragraph{Framed $G$-manifolds}
We now describe $V$-framed $G$-manifolds, which serve as the geometric inputs of genuine $G$-factorization homology theories.
The notion of $V$-framed $G$-manifolds has already been studied by \cite{Weelinck}, though our construction differs from his.

Fix a finite group $G$ and $n\in \N$.
In what follows a \emph{$G$-manifold} is an $n$-dimensional smooth manifold $M$ with a smooth action of $G$.
We organize $G$-manifolds and $G$-equivariant smooth open embeddings using a topological category $\mfld^G$, which we consider as an $\infty$-category by taking its coherent nerve.

Recall that a framing of $M$ is trivialization of its tangent bundle, i.e an isomorphism of tangent bundles $TM\=M \times \R^n$. 
In order to define a framing of $G$-manifolds we consider $TM$ as $G$-vector bundle, with $G$-action induced from the smooth action of $G$ on $M$ by taking differentials.
Fix a real $n$-dimensional $G$-representation $V$.
A \emph{$V$-framing} of $M$ as an isomorphism of $G$-vector bundles $TM \= M \times V$ over $M$.
The $\infty$-category of $\mfld^G$ can be enhanced to an $\infty$-category $\mfld^{G,V-fr}$ of $V$-framed $G$-manifolds.

In fact, we consider genuine equivariant factorization homology theories of $G$-manifolds with more general tangential structures (see \cref{def:G_FH_as_colimit}).
These tangential structures include unframed $G$-manifolds, equivariant orientations in the sense of \cite{CMW_equivariant_orientation_theory} and manifolds with a free $G$-action (see \cref{sec:G_framed_mfld}).

We plan to compare this notion of an equivariant tangential structure with the one introduced by \cite[sec. 2.2]{Weelinck} in future work.

\paragraph{Equivariant factorization homology as a single functor of $\infty$-categories}
Viewing factorization homology as a homology theory suggests a natural generalization to $G$-manifolds. 
Namely, define a $G$-factorization homology theory as a symmetric monoidal functor 
\[ \mfld^{G,V-fr} \to \C \]
satisfying $\otimes$-excision.  
This is essentially the approach taken by Weelinck in \cite{Weelinck}, which leads to a natural generalization of the axiomatic characterization of factorization homology discussed above.
In particular, taking $\C=\Sp^G$ to be the $\infty$-category of genuine $G$-spectra produces invariants of $G$-manifolds valued in genuine $G$-spectra.

However, this is not the approach we take in this paper, for two reasons.
First, we are looking for an \emph{extension} of factorization homology to genuine $G$-spectra. 
If $M$ is a $G$-manifold and $F\colon\mfld^{G,V-fr} \to \Sp^G$ is a $G$-factorization homology theory in the sense of \cite{Weelinck} then the underlying spectrum of $F(M)\in \Sp^G$ need not agree with the factorization homology of $M$.
Second, using a single functor $\mfld^{G,V-fr} \to \C$ to encode a $G$-factorization homology theory prevents us from expressing its expected compatibility with norms described in \cref{eq:G_FH_G_SM_compatibilities}.

Our emphasis on the compatibly of equivariant factorization homology with norms implies that our notion 
an equivariant disk algebra, serving as a coefficient system for equivariant factorization homology, is different from the one introduced in \cite{Weelinck}.
For a specific example, compare \cite[ex. 1.3]{Weelinck} with the description of $\EE_\sigma$-algebras in \cref{sec:Real_THH}.
A detailed comparison of these two notions will appear in future work.

\paragraph{Parametrized $\infty$-categories.}
In order to express both the functoriality of genuine equivariant factorization homology with respect to equivariant embeddings and the compatibilities of \cref{eq:G_FH_G_SM_compatibilities} we view genuine factorization homology as a collection of symmetric monoidal functors 
\begin{align*}
  \forall H<G: \quad \int_- A\colon \mfld^{H,V-fr} \to \Sp^H 
\end{align*}
from the $\infty$-category of $V$-framed $H$-manifolds\footnote{Here we consider $V$ as an $H$-representation by restricting the $G$-action to $H<G$.}
to the category of genuine $H$-spectra, coherently compatible with restrictions and topological inductions.
\footnote{In particular, $\int_- A$ defines a natural transformation between two functors from $\OGop$ to symmetric monoidal $\infty$-categories. However, such natural transformation does not capture the compatibility of norms with topological inductions.}

To make this coherent compatibilities precise we use the theory of parametrized $\infty$-categories, developed by Barwick-Dotto-Glasman-Nardin-Shah in \cite{Expose1,Expose2,Nardin_thesis,Parametrized_algebra,Expose4}.
Informally, a $G$-$\infty$-category is a diagram of $\infty$-categories $\OGop \to \Cat_\infty$ indexed contravariantly by the orbits of $G$.
A $G$-symmetric monoidal structure encodes a symmetric monoidal structure on each of the $\infty$-categories in the diagram together with norm functors and all their expected compatibilities.
In \cref{sec:G_cats_background} we review parametrized $\infty$-category theory in more detail. 

In particular, we use the $G$-$\infty$-category $\ul\Sp^G$ of genuine $G$-spectra constructed in \cite{Nardin_thesis}.
As a $G$-$\infty$-category $\ul\Sp^G$ encodes the $\infty$-categories $\Sp^H$ for all subgroups $H<G$ and the restriction functors relating them.
The $G$-symmetric monoidal structure on $\ul\Sp^G$ encodes smash products and Hill-Hopkins-Ravenel norms. 
Nardin gives an axiomatic characterization of this $G$-symmetric monoidal, see \cite[cor. 3.28]{Nardin_thesis}.
This characterization allows us to work with the Hill-Hopkins-Ravenel norms at a formal level, avoiding the original point set definition of \cite{HHR}.

\paragraph{$\EE_V$-algebras and $V$-framed disks.}
Genuine equivariant factorization homology is an invariant of a geometric input, a $V$-framed $G$-manifold (described above), 
and of an algebraic input, an $\EE_V$-algebra.
We now briefly describe this algebraic structure.

Conceptually, factorization homology is constructed by gluing local data, given by a coefficient system.
Such a coefficient system is an algebraic structure indexed by the local geometry of manifolds: 
an $n$-disk algebra in the case of factorization homology of $n$-dimensional manifolds
and an $\EE_n$-algebra in the case of factorization homology of framed $n$-dimensional manifolds.

Similarly, the structure of an $\EE_V$-algebra is determined by the local structure of $V$-framed $G$-manifolds.
Let $M$ be a $V$-framed $G$-manifold and $x\in M$ a point with stabilizer $H<G$, then $H$ acts linearly on the tangent space $T_x M$, and the $H$-representation $T_x M$ is isomorphic to $V$ (with the action restricted to $H<G$).
\footnote{To see this, pull the $V$-framing $TM \cong M \times V$ along $\set{x} \to M$.}
It follows that $x\in M$ has an $H$-equivariant neighborhood isomorphic to an open disk in $V$.
\footnote{Choose a $G$-equivariant Riemannian metric on $M$ use the fact that the exponential map $T_x M \dashrightarrow M$ is $H$-invariant.}
Therefore the orbit of $x$ (considered as a $G$-submanifold of codimension $0$) has a $G$-tubular neighborhood isomorphic to the topological induction $\coprod_{G/H} V = G \times_H V$.

Let $\D_V$ be the $G$-operad of little $V$-disks, and $\EE_V$ its genuine operadic nerve (see \cite{Bonventre}).
We define $\EE_V$-algebras in $\ul\Sp^G$ as maps of $G$-$\infty$-operads
\[
  \EE_V \to \ul\Sp^G .
\]
Informally, 
an $\EE_V$-algebra $A$ in $\ul\Sp^G$ 
assigns to $V$ a genuine $G$-spectrum $A$ (the ``underlying $G$-spectrum'' of $A$).
The algebraic structure on $A$ is indexed by $H$-equivariant embeddings\footnote{compatible with the $G$-framing} for $H<G$.
To an $H$-embedding $V \sqcup V \cof V$ the algebra $A$ assigns a map of genuine $H$-spectra $A \otimes A \to A$ (a ``multiplication map''),
and to a $H$-embedding $ \sqcup_{H/K} V \cof V$ the algebra $A$ assigns a map $N_K^H(A) \to A$ (a ``multiplicative norm map'') from the Hill-Hopkins-Ravenel norm of $A$. 
All of these maps are coherently compatible with smash products, restrictions of the group action and with each other. 
We use $G$-$\infty$-category theory, and specifically $G$-symmetric monoidal structures, to handle these coherent compatibilities. 

The $G$-$\infty$-operad $\EE_V$ is closely related to $\ulmfld^{G,V-fr}$, as we now explain.
Let $\ul\disk^{G,V-fr}$ be the full $G$-$\infty$-subcategory of $\ulmfld^{G,V-fr}$ generated from the $G$-manifold $V$ by restricting the group action, disjoint unions and topological induction (see \cref{sec:G_manifolds} for details).
By construction, the $G$-symmetric monoidal structure of $\ulmfld^{G,V-fr}$ induces a $G$-symmetric monoidal structure on $\ul\disk^{G,V-fr}$.
In \cref{sec:DV_and_Gdisk_V_fr} we show that $\ul\disk^{G,V-fr}$ is equivalent to the $G$-symmetric monoidal envelope of $\EE_V$.
In particular, an $\EE_V$-algebra in $\ul\Sp^G$ corresponds to an essentially unique $G$-symmetric monoidal functor 
\[
  \ul\disk^{G,V-fr} \to \ul\Sp^G .
\]
We call such functors \emph{$V$-framed $G$-disk algebras} in $\ul\Sp^G$. 

\paragraph{Genuine equivariant factorization homology}
We encode the functors $\int_- A\colon \mfld^{H,V-fr} \to \Sp^H$ as a single $G$-symmetric monoidal $G$-functor $\ulmfld^{G,V-fr} \to \ul\Sp^G$ from the $G$-$\infty$-category of $V$-framed $G$-manifolds to the $G$-$\infty$-category of genuine $G$-spectra.
Given an $\EE_V$-algebra $A$ in $\ul\Sp^G$, let $A\colon \ul\disk^{G,V-fr} \to \ul\Sp^G$ denote the corresponding $V$-framed $G$-disk algebra.
We construct genuine $G$-factorization homology 
\[
  \int_- A\colon \ulmfld^{G,V-fr} \to \ul\Sp^G
\]
as the $G$-left Kan extension of $A$ along the inclusion $\ul\disk^{G,V-fr} \cof \ulmfld^{G,V-fr}$.
By work of Shah \cite{Expose2} the genuine $G$-spectrum $\int_M A$ has an explicit description as a $G$-colimit indexed by ``little disks in $M$'', see \cref{def:G_FH_as_colimit} and \cref{G_FH_functor}.
This construction is indeed a homology theory of $G$-manifolds, as it extends to a $G$-symmetric monoidal functor satisfying $G$-$\otimes$-excision. 

Genuine equivariant factorization homology satisfies the following extension of the Ayala-Francis axiomatic characterization.
\begin{thm} 
  Let \( \mathcal{H}(\ulmfld^{G,V-fr}, \ul\Sp^G) \subset \Fun_G^\otimes(\ulmfld^{G,V-fr},\ul\Sp^G) \) be the full subcategory of the $\infty$-category of $G$-symmetric monoidal $G$-functors $\ulmfld^{G,V-fr} \to \ul\Sp^G$ which satisfy $G$-$\otimes$-excision and respect sequential unions.
  Then there is an equivalence of $\infty$-categories
  \begin{align} 
    \int \colon Alg_{\EE_V}(\ul\Sp^G) \iso \mathcal{H}(\ulmfld^{G,V-fr},\ul\Sp^G), \quad A \mapsto (\int_- A \colon\ulmfld^{G,V-fr} \to \ul\Sp^G )
  \end{align}
  sending an $\EE_V$-algebra $A$
  to $G$-equivariant factorization homology with coefficients in $A$.  
\end{thm}
The above result holds in greater generality.
First, $V$-framed $G$-manifolds can be replaced with $G$-manifolds with more general equivariant tangential structures\footnote{This requires replacing $V$-framed $G$-disk algebras with more general $G$-disk algebras.}.
Second, the $G$-$\infty$-category of genuine $G$-spectra can be replaced with any presentable $G$-symmetric monoidal $G$-$\infty$-category $\ul\C$\footnote{Conjecturally, this condition can  be weakened to distributivity of the tensor product over parametrized sifted colimits}.
The general statement is given in \cref{thm:G_FH_axiomatic_characterization}, which is the main result of this paper. 

\paragraph{Applications}
As an application of \cref{thm:G_FH_axiomatic_characterization}, we describe two variants of topological Hochschild homology using genuine $G$-factorization homology.

In \cref{sec:Real_THH} we show that the real topological Hochschild homology spectrum of Hesselholt-Madsen \cite{Hesselholt_Madsen_real_alg_K_theory} is equivalent to genuine $C_2$-factorization homology over $S^1$.  
\begin{prop}[\cref{prop:Real_THH_as_GFH}]
  For $A$ an $\EE_\sigma$-algebra
  in $\ul\Sp^{C_2}$ 
  there is an equivalence of genuine $C_2$-spectra
  \begin{align*}
    \int_{S^1} A \simeq A \otimes_{N_e^{C_2} A } A .
  \end{align*}
  where $C_2$ acts on $S^1$ by reflection.
\end{prop}
By a theorem of (\cite{Real_THH}) it follows that for $A$ a flat ring spectrum with anti-involution 
there is as an equivalence of genuine $C_2$-spectra
\[ 
  \int_{S^1} A \simeq THR(A) ,
\]
where $THR(A)$ is the real topological Hochschild homology of $A$, see \cref{rem:THR_as_GFH}.

In \cref{sec:Twisted_THH} we show that the ``twisted'' topological Hochschild homology of a genuine $C_n$-ring spectrum of \cite[sec. 8]{TC_via_the_norm} is equivalent to the geometric fixed points of $C_n$-factorization homology over $S^1$. 
\begin{prop}[\cref{prop:twisted_THH_via_GFH}]
  Let $A$ be an $\EE_1$-ring spectrum in $\Sp_{C_n}$,
  and $C_n \curvearrowright S^1$ be the standard action. 
  Then there exists an equivalence of spectra 
  \begin{align*}
    \left( \int_{S^1} A \right)^{\Phi C_n} \simeq  THH(A;A^\tau) . 
  \end{align*}
  In particular,  $THH(A;A^\tau)$ admits a natural circle action.
\end{prop}
This circle action is equivalent to the circle action on the nerve of the ``twisted cyclic bar construction'' of \cite[sec. 8]{TC_via_the_norm}, which gives an alternative description of the relative norm of \cite[def. 8.2]{TC_via_the_norm}.

\paragraph{Construction of $V$-framed $G$-disk algebras.}
Above we gave a rough description of a $V$-framed $G$-disk algebra as encoding multiplication maps, multiplicative norm maps and their coherent compatibilities.
Unwinding these compatibilities implied by \cref{def:G_disk_alg} is usually a non-trivial task (especially when $\dim V \geq 2$), and so it is inadvisable to construct a $V$-framed $G$-disk algebra by specifying multiplication maps, multiplicative norm maps and associated coherence data.
It is therefore desirable to have some general mechanisms for constructing $V$-framed $G$-disk algebras.

For example, one would expect to be able to construct a $V$-framed $G$-disk algebra from an algebra over the $G$-operad $\D_V$.
Such a construction would provide many examples of coefficients for genuine equivalent factorization homology of $V$-framed $G$-manifolds. 
More generally, it would be reassuring to have a ``rectification'' result showing that classical algebras over $\D_V$ form a model for the $\infty$-category of $V$-framed $G$-disk algebras, in the style of \cite[thm. 7.10]{PavlovScholbach}.

We leave such constructions for future work.

\paragraph{What about compact Lie groups?}
It is natural to want to extend genuine equivariant factorization homology from finite groups to compact Lie groups. 
There are two different points in which one encounters complications. 

First, we prove \cref{thm:G_FH_axiomatic_characterization} (the axiomatic characterization of genuine $G$-factorization homology) inductively using equivariant handle bundle decompositions (see \cite{Wasserman}).
We produce these decompositions using equivariant Morse theory, which is more complicated over a compact Lie group. 
Choosing an invariant Morse function gives rise to a handle bundle decomposition, 
where each handle bundle is an equivariant disk bundle over a critical orbit. 
However, for a compact Lie group of positive dimension these handle bundles can be non trivial, since critical orbits are submanifolds of possibly positive dimension.

Second, and more fundamental, is the lack of good $G$-$\infty$-category theory for a compact Lie group $G$. 
The source of the problem is the lack of multiplicative norms for subgroups $H<G$ of non-finite index.
In order to understand the significance of this fact for genuine $G$-factorization homology, consider $M=\mathbb{C}$ the complex plane with the standard action of the circle group $S^1=\mathbb{C}^\times$.
The unit circle is an $S^1$-orbit in $\mathbb{C}$, with $S^1$-tubular neighborhood given by the open annulus.
The embedding of the open annulus in $\mathbb{C}$ should induce a ``multiplication norm map'' \( \otimes_{S^1} A \to A \) of genuine $S^1$-spectra, where the tensor product is indexed over the free orbit $S^1$.
However, we do not have a good definition for the domain of this map as a genuine $S^1$-spectrum. 

\paragraph{Organization}
We start by reviewing some parts of parametrized $\infty$-category theory in \cref{sec:G_cats_background}.
We hope this short exposition will assist the reader unfamiliar with the theory of $G$-$\infty$-categories.

In \cref{sec:G_manifolds} we construct the $G$-$\infty$-categories of $G$-manifolds and $G$-disks with equivariant tangential structures, and their $G$-symmetric monoidal structure which encodes disjoint unions and topological induction.
These constructions provide a bridge between the geometry of $G$-manifolds and parametrized $\infty$-category theory, and enables the construction of genuine $G$-factorization homology in \cref{sec:G_FA_G_SM_functor}.

Our definition of equivariant tangential structures in \cref{sec:G_framed_mfld} uses an equivariant version of the tangent classifier of \cite{AF} which may be of independent interest, see \cref{sec:G_VB_and_G_tangent_classifier}.
While we focus on framed $G$-manifolds, our definition is flexible enough to consider more general tangential structures such as equivariant orientations, as well as allowing us to restrict our attention to manifolds with a free $G$-action.

We finish \cref{sec:G_manifolds} by studying some aspects of these constructions.
In \cref{sec:G_disks_and_configurations} we study the relation between embedding spaces of $G$-disks and $G$-configuration spaces.
In \cref{sec:DV_and_Gdisk_V_fr} use the work of \cite{Bonventre} to show that the $G$-$\infty$-operad encoding $V$-framed $G$-disk algebras is closely related to the $G$-operad of little disks in a representation $V$.

The technical results and constructions of \cref{sec:G_manifolds} provide a solid foundation for the use of abstract theory of parametrized $\infty$-categories in the following sections. 

In \cref{sec:G_factorization_homology} we define framed $G$-disk algebras and construct $G$-factorization homology, first as a $G$-functor (by $G$-left Kan extension, see \cref{sec:G_factorization_homology_G_functor}) and then as a $G$-symmetric monoidal $G$-functor (\cref{sec:G_FA_G_SM_functor}).

In \cref{sec:G_FH_properties} we study the properties of $G$-factorization homology.
In \cref{sec:G_collar_decomposition} we define $G$-collar decompositions and construct an ``inverse image'' functor.
We use these in \cref{sec:G_tensor_excision}, where we define $G$-$\otimes$-excision for a general $G$-symmetric monoidal functor $\Gmfld \to \ul\C$, and show that $G$-factorization homology satisfies $G$-tensor excision.
In \cref{sec:G_FH_respects_seq_unions} we show that $G$-factorization homology respects sequential unions.

In \cref{sec:axiomatic_characterization} we prove our main result, giving an axiomatic characterization of $G$-factorization homology using equivariant Morse theory. 

In \cref{sec:applications} we describe real topological Hochschild homology using $G$-factorization homology (\cref{sec:Real_THH}), and the relative norm of a genuine $C_n$-ring spectrum as the geometric fixed points of $G$-factorization homology (\cref{sec:Twisted_THH}).

In \cref{MooreOverCat} we show how to model $\infty$-slice categories in the framework of topological categories. 
For the convenience of the reader we recall the definition of $G$-symmetric monoidal categories in \cref{sec:G_SM_cat}.
We collect some Some general statements about mapping spaces in over categories in \cref{sec:maping_spaces}.

\paragraph{Notation.}
In this work we use the quasi-categories as a model $\infty$-categories (with the exception of \cref{G_cat_as_CS_cat}). 
We assume the reader is familiar with the theory of $\infty$-categories, as developed in \cite{HTT} and \cite{HA}.
Explicitly, an $\infty$-category is a simplicial set $\C$ satisfying the left lifting property with respect to inner horns:
for every $0<i<n$, any map $\Lambda^n_i \to \C$ admits an extension to $\Delta^n \to \C$.

All of the manifolds we consider are smooth and $n$-dimensional for a fixed $n\in\N$.
We fix a finite group $G$, and only consider manifolds with actions of subgroups $H<G$.

We frequently construct $\infty$-categories from topological categories by taking their coherent nerve (which is called the topological nerve in \cite[def. 1.1.5.5]{HTT}).
We emphasize that the coherent nerve of a topological category $\C$ is a two step construction.
First, taking the singular nerve of each mapping space, produces a simplicial category $\Sing(\C)$. 
Second, applying the simplicial nerve functor of \cite[def. 1.1.5.5]{HTT} to $\Sing(\C)$ produces an $\infty$-category.
We denote the resulting $\infty$-category by $\NSing(\C)$.

We denote parametrized $\infty$-categories with an underline, for example $\ul\C$. 
In general, if $\ul\C$ is parametrized over an $\infty$-category $S$ we refer to $\ul\C$ as an $S$-$\infty$-category. 
We say that $\ul\C$ is a \emph{$G$-category} (see \cref{def:G_cat}) if it is parametrized over $\OGop$, where $\OG$ is the orbit category of $G$. 
No other notion of $G$-categories is used; a $G$-category is by definition an $\OGop$-$\infty$-category.

\paragraph{Acknowledgements}
This paper is a revised version of my PhD thesis.
I wish to thank my advisors, Emmanuel Farjoun and Yakov Varshavsky, for their continued support throughout this project.
I am grateful to Tomer Schlank for suggesting this project and for valuable insights.

This paper builds on the work of Clark Barwick, Emanuele Dotto, Saul Glasman, Denis Nardin and Jay Shah on parametrized higher category theory.
I would like to thank Clark Barwick for helpful discussions and Jay Shah for generously sharing an early draft of \cite{Parametrized_algebra}.

\section{Background on parametrized $\infty$-category theory} \label{sec:G_cats_background}

In this section we review parametrized $\infty$-category theory of Barwick, Dotto, Glassman, Nardin and Shah, developed in \cite{Expose1,Expose2,Expose4,Nardin_thesis,Parametrized_algebra}.
We recall the notions of $G$-$\infty$-category theory employed below and fix our notation.
We restrict our discussion to the case of $G$-$\infty$-categories, though nothing substantial would change when working over an arbitrary indexing category. 

This section contains no original results, all the results of this section are entirely due to Barwick, Dotto, Glassman, Nardin and Shah.

\subsection{From Elmendorf-McClure's theorem to $G$-categories}
A good starting point to a discussion of $G$-categories is the Elmendorf-McClure theorem, which recasts the equivariant homotopy theory of $G$-spaces as a presheaf category.
Throughout we fix a finite group $G$.
\begin{mydef} \label{def:orbit_cat}
  The orbit category $\OG$ is the full subcategory of $G$-sets supported by transitive $G$-sets.
\end{mydef}
Note that every orbit in $\OG$ is isomorphic to a quotient of $G$ by some subgroup $H<G$. 
This isomorphism depends only on the choice of a base point of the orbit, with $H$ the stabilizer of the chosen basepoint.
We denote the objects of $\OG$ either by $O$ by $G/H$.
Despite the suggestive notation, we try to refrain from a choice of basepoint when possible. 

Define the $\infty$-category of $G$-spaces $\GTop$ as the coherent nerve of the topological category of $G$-CW spaces and $G$-maps. 
\begin{thm}[Elmendorf-McClure, \cite{Elmendorf}]
  There is an equivalence of $\infty$-categories \\
  \[ 
    \GTop \iso \Fun(\OGop, \Ss) ,
  \]
 sending a $G$-space $X$ to its diagram of fixed points, \( G/H \mapsto X^H \).
\end{thm}
Using straightening/unstraightening (\cite[thm. 2.2.1.2]{HTT}) we get a third description of a $G$-space $X$ as the left fibration over $\OGop$ classifying the diagram of fixed points of $X$.

\begin{mydef} \label{def:G_cat}
  A $G$-$\infty$-category is a coCartesian fibration \( \ul\C \fib \OGop \).
\end{mydef}
For sake of readability we refer to $G$-$\infty$-categories simply as $G$-categories.
Other notions of $G$-categories present in the literature are not present in this paper.

A $G$-category $\ul\C$ is classified by a diagram of $\infty$-categories \( \ul\C_\bullet \colon \OGop \to \Cat_\infty \) sending $G/H\in \OGop$ to the fiber $\ul\C_{[G/H]}$ of $\ul\C \fib \OGop$ over $G/H$.
We systematically use the subscript-square bracket notation $\C_{[G/H]}$ for the fiber $\infty$-category in order to avoid confusion with other subscript notations.
As above, straightening/unstraightening (\cite[sec. 3.2]{HTT}) ensures that this is an equivalent description of the $G$-$\infty$-category $\ul\C$.

\begin{rem} \label{G_cat_as_CS_cat}
  Describing $\Cat_\infty$ as a complete Segal object in the $\infty$-category of spaces, we can use the Elmendorf-McClure theorem to get a third equivalent description of a $G$-category as a complete Segal object in $\GTop$.
  This follows from following the Segal conditions and completeness conditions along the equivalences
  \begin{align*}
    \Fun( \Delta^{op}, \GTop) \simeq \Fun( \Delta^{op}, \Fun(\OGop, \Ss)) \simeq  \Fun( \OGop , \Fun( \Delta^{op}, \Ss)) ,
  \end{align*}
  where the first equivalence is induced by the Elmendorf-McClure theorem.
  
  In particular, categories internal to $G$-spaces (and to $G$-sets) are examples of $G$-categories.
    Note that \cite{Guillou_May_Permutative_G_cats} defines $G$-categories as categories internal to $G$-spaces, making them examples of $G$-($\infty$-)categories in the sense of \cite{Expose1}, used here.
  
  While these equivalent descriptions of a $G$-category are good to have in mind, we stick to the definition of a $G$-category as a coCartesian fibration for its explicit nature.
\end{rem}

When we need more general parametrized $\infty$-categories we use the following definition (and the notation of \cite{Expose2}).
\begin{mydef} \label{def:S_cat}
  Let $S$ be an $\infty$-category. 
  An $S$-$\infty$-category is a coCartesian fibration \( \ul\C \fib S \).
  We denote the fiber of $\ul\C$ over $s\in S$ by $\ul\C_{[s]}$.
\end{mydef}
We refer to $S$-$\infty$-categories as $S$-categories.
\begin{rem}
  Most results recalled in this section hold for general $S$-categories. 
  One notable exception is the description of $\ul\Sp^S$, the $S$-stabilization of the $S$-category of $S$-spaces, using spectral Mackey functors.
  Another exception is the uniqueness of $S$-symmetric monoidal structure on $\ul\Sp^S$. 
  However, these results hold under mild conditions on $S$.
  \footnote{Specifically, they hold for $S$ an atomic orbital $\infty$-category. See \cite{GCats_intro} for examples and \cite{Expose1}, \cite{Expose2}, \cite{Nardin_thesis} for the general theory.}
\end{rem}

\paragraph{Handling $H$-categories as $G$-categories}
Occasionally we have to consider $H$-categories for some subgroup $H<G$. 
When doing so we use the slice category $\ul{G/H}:=(\OGop)_{(G/H)/}$, the opposite of the category of $G$-orbits over $G/H$.
The category $\ul{G/H}$ is equivalent to the category of $H$-orbits.
Moreover, the forgetful functor $\ul{G/H} \to \OGop$ is left fibration classified by the representable functor $\Map(-,G/H) \colon \OGop \to \Ss$. 
In particular a $\ul{G/H}$-category $\ul\C \fib \ul{G/H}$ is a $G$-category by postcomposition with the forgetful functor, $\ul{G/H} \fib \OGop$. 
Note that this construction also avoid a choice of basepoint to $G/H$.
When referring to the fibers of $\ul\C \fib \ul{G/H}$ we adopt the notation $\ul\C_{[\varphi]}$ for the fiber over $\varphi \colon  G/K \to G/H$ as an object in the slice category $\ul{G/H}$.

The $G$-category $\ul{G/H}$ has a second role for us, since a $G$-functor \( \ul{G/H} \to \ul\C \) corresponds to an object in the fiber of $\ul\C \fib \OGop$ over $G/H$, as we now explain.
Under straightening/unstraightening the left fibration \( \ul{G/H} \fib \OGop \) corresponds to the representable functor of the orbit $G/H$, given by \( \Hom_{\OG}(-,G/H)  \colon  \OGop \to \Ss\), and therefore by the Yoneda lemma (\cite[lem. 5.1.5.2]{HTT}) corresponds to an object of \( \ul\C_{[G/H]} \)
\footnote{To make this argument precise we need to replace $\ul\C$ with a presheaf of spaces. To achieve that we straighten $\ul\C^\simeq\subseteq \ul\C$, the maximal $G$-subgroupoid of $\ul\C$, given as a left fibration by the full maximal sub-simplicial set supported on the coCartesian edges of $\ul\C$.}.
We denote the $G$-functor corresponding to $x\in \ul\C_{[G/H]}$ by $\sigma_x \colon  \ul{G/H} \to \ul\C$. 
A more explicit construction is given by choosing a section of the trivial fibration \( \Arr^\text{coCart}_{x\to}(\ul\C) \trivfib \ul{G/H} \) of \cite[not. 2.28]{Expose2} and composing with \( ev_1 \colon  \Arr^\text{coCart}_{x\to}(\ul\C) \to \ul\C \).

\paragraph{The $G$-category of $G$-spaces}
The $\infty$-categories $\Fun(\ul{G/H},\Ss)$ assemble as the fibers of a $G$-category $\ulTopG$, the $G$-category of $G$-spaces (\cite[ex. 7.5]{Expose1}).
By the Elmendorf-McClure theorem the fiber over $G/H$ is equivalent to \( \ulTopG_{[G/H]} \= \Fun(\ul{G/H}, \Ss) \simeq \Fun(\mathcal{O}^{op}_H,\Ss) \simeq \mathbf{Top}^H \), the $\infty$-category of $H$-spaces.
By an $H$-space we always mean an $H$-CW space.  
The $G$-category of $G$-spaces is characterized by the following universal property (see \cite[thm. 7.8]{Expose1}).

For any $G$-category $\ul\C$ we have an equivalence of $\infty$-categories \( \Fun_G(\ul\C, \ulTopG) \simeq \Fun(\ul\C, \Ss) \), i.e $\ulTopG$ is the cofree $G$-category co-generated by the $\infty$-category of spaces. 

Taking our cue from the Elmendorf-McClure theorem, we think of a $G$-category as capturing the notion of a $G$-action on an $\infty$-category.
With this intuition in mind one may think of $\ulTopG$ as follows.
Imagine that the $\infty$-category of spaces admits a non-trivial $G$-action, whose $H$-fixed points is the $\infty$-category of $H$-spaces for all $H<G$.
Think of $\ulTopG$ as capturing this imagined $G$-action.

\begin{rem} \label{ulTopG_as_G_spaces_over_orbits}
  In \cref{sec:G_VB_and_G_tangent_classifier} we use the following explicit model for $\ulTopG$. 
  Construct an auxiliary \emph{topological} category $\OGTop$ as follows. 
  An object of $\OGTop$ is $G$-map $X \to O$ where the domain $X$ is a $G$-CW complex and codomain $O\in \OG$ is a $G$-orbit.
  We refer to an object of $\OGTop$ as $\OG$-space, though it should rightfully be called a ``$G$-space over an orbit''. 
  A map of $\OG$-spaces is given by a (strictly) commuting squares of $G$-spaces
  \begin{align} \label{eq:morphism_in_TG}
    \xymatrix{
      X_1 \ar[d] \ar[r] & X_2 \ar[d] \\
      O_1 \ar[r] & O_2 .
    }
  \end{align}
  The mapping spaces of $\OGTop$ are given by 
  \[ \Map_{\OGTop}( X_1 \to O_1, X_2 \to O_2)  = \Map_G(X_1,X_2) \times_{\Map_G(X_1,O_2)} \Map_G(O_1, O_2), \]
  where $\Map_G(X,Y)$ is the space of $G$-maps $X \to Y$ with the compact-open topology.

  We think of an $\OG$-space $X \to G/H$ as representing the $H$-space given by the fiber $X|_{H}$ of $X \to G/H$ over the coset $H$.
  On the other hand, given an $H$-space $X_0$ we can use topological induction to construct a $\OG$-space $ G \times_H X_0$ whose fiber over $H$ is $X_0$.
  Note that the $\OG$-space $X \to O$ \emph{does not} represent the $G$-space $X$ (in fact, choosing an isomorphism $O\=G/H$ for some $H<G$ exhibits the $G$-space $X$ as the topological induction of the $H$-space represented by $X\to G/H$).

  Applying topological nerve construction of \cite[def. 1.1.5.5]{HTT} produces an $\infty$-category $\NSing (\OGTop)$.
  The forgetful  \( \NSing (\OGTop) \to \OG, \, (X\to O) \mapsto O \) is a Cartesian fibration, 
  and a commuting square \eqref{eq:morphism_in_TG} describes a coCartesian edge in $\NSing(\OGTop)$ if it is a pullback square. 
  To see this use \cite[prop. 2.4.1.1 (2)]{HTT} as in the proof of \cref{prop:pullback_is_Cart}.
  The dual coCartesian fibration \( \NSing (\OGTop)^\wedge \to \OGop \), described in \cite{DualizingCarFibs}, is a $G$-category equivalent to $\ulTopG$.
  We can explicitly describe an object of $\ulTopG_{[G/H]}$ in this model as a $G$-map $X \to G/H$, which we interpret as the $H$-space $X|_{eH}$ given by the fiber over the coset $eH$.
  A map in $\ulTopG$ is given by a (strictly) commutative diagram of $G$-spaces 
  \begin{align*}
    \xymatrix{
      X_1 \ar[d] & \ar[l] \pullbackcornerleft X' \ar[d] \ar[r] & Y \ar[d] \\
      O_1 & \ar[l] O_2 \ar[r]^{=} & O_2
    }
  \end{align*}
  in which the left square is a pullback square. 
  It is a coCartesian edge if and only of the $G$-map $X' \to Y$ is a $G$-homotopy equivalence over $O_2$ (see \cite{DualizingCarFibs}).
  Equivalently, if $O_2=G/H$ then the above edge is coCartesian precisely when the map of fibers $X'|_{eH} \to Y|_{eH}$ is an $H$-homotopy equivalence.

  By definition maps in fiber $\ulTopG_{[O]}$ are commutative diagrams as above, with row given by $ O \xfrom{=} O \xto{=} O$.
  Unwinding the definitions we see that $\ulTopG_{[O]}$ is equivalent to $\NSing(\GTop_{/O})$, the coherent nerve of the topological category of $G$-CW-spaces over $O$.
  If $O=G/H$ then restriction to the fiber over $eH$ defines an equivalence of topological categories $\GTop_{/ G/H} \iso \mathbf{Top}^H$ to the topological category of $H$-CW-spaces.

  Finally, we note that $\NSing(\GTop_{/O}) \simeq \NSing(\GTop)_{/O}$ are equivalent $\infty$-categories . 
  We use the Moore over-category of \cref{MooreOverCat} to see this.
  By \cref{Moore_is_overcat} we have 
  \( \NSing (\GTop)_{/O} \simeq \NSing\left( (\GTop)^{Moore}_{/O} \right) \).
  However, since the orbit $O$ is a discrete $G$-space we see that for every $X\in \GTop$ the only Moore paths in $\Map_{\GTop} (X, O)$ are constant, 
  so \( \GTop_{/O} \to (\GTop)^{Moore}_{/O} \) is an equivalence of topological categories. 
  Therefore the fiber $\ulTopG_{[O]}$ is equivalent to the slice category $\NSing( \GTop)_{/O}$. 
  The mapping spaces of $\ulTopG_{[O]} \simeq \NSing( \GTop_{/O})$ will be denoted by $\Map^G_O(X,Y)$.
\end{rem}

\paragraph{The $G$-category of $G$-spectra}
A more interesting example is given by $\ul\Sp^G$, the $G$-category of $G$-spectra, with fiber over $G/H$ is equivalent to $\ul\Sp^G_{[G/H]} \simeq \Sp^H$, the $\infty$-category genuine orthogonal $H$-spectra (see \cite[thm. 2.40]{Nardin_thesis}, with origins in \cite{Guillo_May_Spectral_Mackey_functors}). 
For a construction of $\ul\Sp^G$ as the $G$-stabilization of $\ulTopG$ see \cite[def. 2.35 and thm. 2.36]{Nardin_thesis}.

\subsection{Constructing $G$-categories}
We frequently use the following constructions of $G$-categories.
\begin{construction}
  Given two $S$-categories $\ul\C,\ul\D$, the fiber product $\ul\C \times_S \ul\D$ is an $S$-category, the fiberwise product of $\ul\C$ and $\ul\D$.
  If $\ul\C,\ul\D$ are $G$-categories, we denote the fiberwise product $\ul\C \times_{\OGop} \ul\D$ by $\ul\C \ultimes \ul\D$.
  In particular, we use the fiberwise product to restrict a $G$-category $\ul\C \fib \OGop$ to a $\ul{G/H}$-category $\ul\C \ultimes \ul{G/H} \fib \ul{G/H}$ (``forgetting the $G$-action on $\C$ to get an $H$-action''). 
\end{construction}
\begin{construction}
  Given a $G$-category $\ul\C$ define the \emph{fiberwise arrow category} $\Arr_G(\ul\C)$ as the fiber product $ \OGop \times_{\Fun(\Delta^1,\OGop) } \Fun(\Delta^1, \ul\C)$ (see \cite[not. 4.29]{Expose2}). 
  Note that $\Arr_G(\ul\C)$ is equivalent to the functor $G$-category \( \ulFun_G(\OGop \times \Delta^1, \ul\C) \), where the $G$-category $\OGop \times \Delta^1$ is the constant $G$-category on $\Delta^1$.
    More generally, for any $S$-category $\ul\C \fib S$ define the fiberwise arrow $S$-category $\Arr_S(\ul\C)$ as the fiber product \( S \times_{\Fun(\Delta^1,S)} \Fun(\Delta^1, \ul\C) \).
\end{construction}
\begin{construction} \label{const:G_slice}
    Let $\ul\C$ be a $G$-category and $x\in \ul\C_{[G/H]}$ an object over $G/H$, corresponding to the $G$-functor $ \sigma_x \colon \ul{G/H} \to \ul\C$. 
    Following \cite[not. 4.29]{Expose2}, we define the \emph{parametrized slice-category} $\ul\C_{/\ul{x}} \fib \ul{G/H}$ by pulling back the coCartsian fibration $ev_1 \colon \Arr_G(\ul\C) \fib \ul\C$ along $\sigma_x$,  
    i.e. $\ul\C_{/\ul{x}} := \Arr_G(\ul\C) \times_{\ul\C} \ul{G/H}$.
    We will also consider $\ul\C_{/\ul{x}} \fib \ul{G/H}$ as a $\ul{G/H}$-category. 

    Note that the fiber of $\C_{/\ul{x}}\fib \ul{G/H}$ over $\varphi \colon G/K \to G/H$ is equivalent to the $\infty$-over-category $(\C_{[G/K]})_{/\varphi^* x}$, where $\varphi^* x \in \C_{[G/K]}$ is determined by choosing a coCartesian lift $ x \to \varphi^* x$ of $\varphi$. 
\end{construction}

\begin{construction}
  For $\ul\C \fib S$ an $S$-category, the fiberwise cone $S$-category of $\ul\C$ is defined as the parametrized join $\ul\C \star_S S$ (see \cite[not. 4.2]{Expose2} or \cref{sec:G_SM_cat}).
\end{construction}

\paragraph{Parametrized functors and parametrized functor categories}

\begin{mydef}
  Let $\ul\C , \ul\D $ be $S$-categories, i.e. coCartesian fibrations $\ul\C \fib S, \ul\D \fib S$.
  An $S$-functor is a functor \( \ul\C \to \ul\D \) over $S$ which preserves coCartesian edges.
  Let $\Fun_S(\ul\C,\ul\D)\subseteq \Fun_{/S}(\ul\C, \ul\D)$ be the full subcategory of functors $\ul\C \to \ul\D$ over $S$ which preserve coCartesian edges.
  When $S=\OGop$ we refer to a $\OGop$-functor as a $G$-functor, and denote the $\infty$-category of $G$-functors by $\Fun_G(\ul\C,\ul\D)$.
\end{mydef}
\begin{rem}
  An $S$-functor $\ul\C \to \ul\D$ encodes the data of a coherent natural transformation $\ul\C_\bullet \Rightarrow \ul\D_\bullet$ between the $S$-diagrams $\ul\C_\bullet,\ul\D_\bullet \colon  S \to \Cat_\infty$ classified by the coCartesian fibrations $\ul\C \fib S$ and $\ul\D \fib S$.
\end{rem}
\begin{rem} \label{rem:GH_objects}
  Since the left fibration \( \ul{G/H} \to \OGop \) is corepresentable by construction, we have \( \ul\C_{[G/H]} \simeq \Fun_G (\ul{G/H}, \ul\C) \).
\end{rem}

The $\infty$-category of $G$-categories admits an internal hom, a $G$-category denoted $\ulFun_G(\ul\C,\ul\D)$ see \cite[thm. 9.7]{Expose1} and \cite[def. 9.2]{Expose1} for an explicit construction.
The fiber of $\ulFun_G(\ul\C,\ul\D) \fib \OGop$ over $G/H$ admits the following description. 
Forget the $G$-action on $\ul\C,\ul\D$ to an $H$-action by taking the fiber products \( \ul\C \ultimes \ul{G/H},\, \ul\D \ultimes \ul{G/H} \).
The fiber $\ulFun_G(\ul\C,\ul\D)_{[G/H]}$ is equivalent to the $\infty$-category \( \Fun_{\ul{G/H}}(\ul\C \ultimes \ul{G/H}, \ul\D \ultimes \ul{G/H}) \) of $\ul{G/H}$-functors \( \ul\C \ultimes \ul{G/H} \to  \ul\D \ultimes \ul{G/H} \), (which we think of as modeling ``$H$-equivariant functors from $\ul\C$ to $\ul\D$'').

More generally, for any two $S$-categories $\ul\C\fib S , \ul\D \fib S$ there is an $S$-category of functors $\ulFun_S(\ul\C,\ul\D)$ with fibers \( \ulFun_S(\ul\C,\ul\D)_{[s]} \simeq \Fun_{\ul{s}}(\ul\C \times_S \ul{s} , \ul\D \times_S \ul{s}) \) where $\ul{s}= S_{s/}$.
The $S$-category of functors possesses the universal property of internal hom, from \cite[thm. 9.7]{Expose1}.
\begin{thm}[Barwick-Dotto-Glasman-Nardin-Shah]
  Let $\ul\C,\ul\D,\ul\E$ be $S$-categories. Then there are natural equivalences 
  \begin{align*}
    \ulFun_S(\ul\C,\ulFun_S(\ul\D,\ul\E)) \iso \ulFun_S( \ul\C \times_S \ul\D, \ul\E), \quad 
    \Fun_S(\ul\C,\ulFun_S(\ul\D,\ul\E)) \iso \Fun_S( \ul\C \times_S \ul\D, \ul\E).
  \end{align*}
\end{thm}
Note that if $\ul\C,\ul\D,\ul\E$ are $G$-categories, then the second equivalence follows from the first by restricting to the fiber over the orbit $[G/G]$, the terminal object of $\OG$.

\subsection{Parametrized adjoints, colimits, left Kan extensions}
We follow \cite{Nardin_thesis}, defining parametrized colimits and parametrized left Kan extensions using parametrized adjoints.

\paragraph{Parametrized adjoints}
Let $\ul\C,\ul\D$ be $S$-categories. 
An $S$-adjunction (\cite[def. 8.1]{Expose2}) is a relative adjunction $L \colon \ul\C \adj \ul\D \noloc R$ over $S$ (\cite[def. 7.3.2]{HA}) where both $L$ and $R$ are $S$-functors.
In particular, for each $s\in S$ we have an adjunction $L_{[s]} \colon \ul\C_{[s]} \adj \ul\D_{[s]} \noloc  R_{[s]}$ between the fibers over $s$.
When $S=\OGop$ we will refer to an $\OGop$-adjunction as a $G$-adjunction.

\paragraph{Parametrized colimits}
Let $p \colon \ul{I}\to \ul\C$ be an $S$-functor, which we think of as an $S$-diagram in $\ul\C$.
The $S$-colimit of $p$ is an $S$-object of $\ul\C$, i.e a coCartesian section $S-\colim(p) \colon  S \to \ul\C$ of the structure fibration $\ul\C \fib S$. 
For a general definition of $\colim(p)$ as the $S$-initial $S$-cone under $p$ see \cite[def. 5.2]{Expose2}. 
We define $\ul{I}$-shaped $S$-colimits as the $S$-left adjoint to the ``constant $\ul{I}$-diagram'' $S$-functor, following \cite[def. 2.1]{Expose4}.
This definition is justified by \cite[10.4]{Expose2}, since we only take $S$-colimits in $S$-cocomplete $S$-categories. 

Explicitly, precomposition with the coCartesian fibration $\ul{I} \fib S$ induces an $S$-functor \( \Delta_I \colon  \ul\C \simeq \ulFun_S(\ul{S},\ul\C) \to \ulFun_S(\ul{I},\ul\C) \), where $\ul{S}$ is the terminal $S$-category (given by $id \colon S \to S$).
If $\Delta_{\ul{I}}$ admits an $S$-left adjoint we say that $\ul\C$ admits $\ul{I}$-indexed $S$-colimits, and denote the $S$-left adjoint by \(S-\colim \colon  \ulFun_S(\ul{I},\ul\C) \to \ul\C \). 
Note that for every index $s\in S$ we have an adjunction of $\infty$-categories 
\[ S-\colim \colon  \Fun_{\ul{s}}(\ul{I} \times_S \ul{s},\ul\C \times_S \ul{s}) \adj \Fun_{\ul{s}}(\ul{S} \times_S \ul{s}, \ul\C \times_S \ul{s}) \simeq \ul\C_{[s]}  \noloc  \Delta_{\ul{I}} . \]
Particularly, we will use the following type of $\ul{G/H}$-colimit.
\begin{ex} \label{ex:GmodH_colim}
  Let $\ul\C$ be a $G$-category, $\ul{I} \fib \ul{G/H}$ a $\ul{G/H}$-category and $p \colon  \ul{I} \to \ul\C$ a $G$-functor.
  Since $\ul{G/H} \fib \OGop$ is a left fibration we have \( \Fun_G(\ul{I}, \ul\C) \simeq \Fun_{\ul{G/H}} (\ul{I}, \ul\C \ultimes \ul{G/H}) \), under which $p$ corresponds to a $\ul{G/H}$-functor $p \colon \ul{I} \to \ul\C \ultimes \ul{G/H}$, or in other words \( p\in \Fun_{\ul{G/H}}(\ul{I} , \ul\C \ultimes \ul{G/H}) \). 
  Then $\ul{G/H}-\colim(p) \in \ul\C_{[G/H]}$ is given by applying the left adjoint of
  \[ \ul{G/H}-\colim \colon  \Fun_{\ul{G/H}}(\ul{I} ,\ul\C \ultimes \ul{G/H} ) \adj \Fun_{\ul{G/H}}(\ul{G/H}, \ul\C \ultimes \ul{G/H}) \simeq \ul\C_{[G/H]}  \noloc  \Delta_{\ul{I}} . \]
\end{ex}
We say that an $S$-category $\ul\C$ is $S$-cocomplete if for every $s\in S$ the $\ul{s}$-category $\ul\C \ultimes \ul{s}$ admits $\ul{I}$-indexed $\ul{s}$-colimits for any $\ul{s}$-category $\ul{I}$.

\paragraph{Parametrized left Kan extensions}
We follow \cite[def. 2.12]{Nardin_thesis} and define $S$-left Kan extension using the give a global characterization as a left adjoint. 
For a general definition of pointwise parametrized left Kan extensions see \cite[def. 10.1]{Expose2}, which satisfies the global characterization by \cite[10.4]{Expose2}.
We only use the pointwise definition in the proof of \cref{lem:Operadic_GLan_restrict_to_GLan}, a $G$-categorical statement independent from the rest of the paper.

Let $\iota \colon  \ul\D \to \ul\M$ be an $S$-functor and $\ul\C$ an $S$-category. 
Restriction along $\iota$ induces an $S$-functor \( \iota^* \colon  \ulFun_S(\ul\M,\ul\C) \to \ulFun_S(\ul\D,\ul\C) \). 
The $S$-left Kan extension along $\iota$ is the $S$-left adjoint to $\iota^*$ and denoted by $\phi_!$.

We will use the following propositions from \cite{Expose2}.
\begin{prop}{ \cite[thm. 10.3]{Expose2} }
  Let $A \colon  \ul\D \to \ul\C$ and $\iota \colon  \ul\D \to \ul\M$ be $S$-categories, and suppose that for every $x\in \ul\M$ over $s\in S$ the $\ul{s}$-colimit 
  \begin{align*}
    \ul{s}-\colim \left( \ul\D_{/\ul{x}} \to \ul\D \times_S \ul{s} \xto{A \times_S \ul{s}} \ul\C \times_S \ul{s} \right)
  \end{align*}
  exists.
  Then the $S$-left Kan extension of $A$ along $\iota$ exists (and is essentially unique), and acts on $x\in \ul\D$ by sending it to the $\ul{s}$-colimit above, considered as an object in the fiber $\ul\C_{[s]} $.
\end{prop}
\begin{prop}{ \cite[cor. 10.6]{Expose2} } \label{Shah:fully_faithful_G_Lan}
  Let $\ul\C$ be a $S$-cocomplete $S$-category and $\iota \colon  \ul\D \to \ul\M$ a fully faithful $S$-functor (i.e fiberwise fully faithful, see \cite[def. 1.6]{Expose1}).
  Then the $S$-left Kan extension $\iota_!  \colon  \ul\M \to \ul\C$ exists and is $S$-fully faithful.
\end{prop}
When $S=\OGop$ we refer to $S$-left Kan extensions as $G$-left Kan extensions, which we use to define $G$-factorization homology as a $G$-functor (see \cref{G_FH_functor}).

\paragraph{Parametrized Yoneda embedding}
Another useful tool available to us is the parametrized Yoneda embedding of \cite[sec. 10]{Expose1}, which we use in the construction of the $G$-tangent classifier (see \cref{const:G_tangent_classifier}).
Let $\ul\C$ be a $G$-category, and $\ul\C^{vop}$ the fiberwise opposite $G$-category (with fibers $(\ul\C^{vop})_{[G/H]} \= (\ul\C_{[G/H]})^{op}$, see \cite[def. 3.1]{Expose1}).
According to \cite[def. 10.2]{Expose1} there exists a $G$-functor \(j \colon  \C \to \ulFun_G(\C^{vop}, \ulTopG) \), the parametrized Yoneda embedding, which can be informally described as follows. 
The $G$-functor $j$ takes $x\in \ul\C_{[G/H]}$ the $\ul{G/H}$-functor \( \ul{\Map} (-,x) \colon  \ul\C^{vop} \ultimes \ul{G/H} \to \ulTopG\ultimes \ul{G/H} \) sending an object $y\in ( (\ul\C^{vop}) \ultimes \ul{G/H})_{[\varphi]} \= (\ul\C_{[G/K]})^{op}$ in the fiber over \( \varphi \colon  G/K \to G/H \) to the mapping space $\Map(y, \varphi^* x)$ of the $\infty$-category $\ul\C_{[G/K]}$.

\subsection{$G$-symmetric monoidal structures}

The notion of a $G$-symmetric monoidal structure plays a central role in our presentation of $G$-factorization homology. 
In this subsection we give some intuition for $G$-symmetric monoidal structure, hopefully making it more approachable. 
This subsection is expository in nature, the formal definition of a $G$-symmetric monoidal $G$-category can be found in \cite[sec. 3.1]{Nardin_thesis}, or in \cref{sec:G_SM_cat}.

Informally, the data of a $G$-symmetric monoidal structure on a $G$-category $\ul\C$ is given by collection of symmetric monoidal structures on the fibers $\ul\C_{[G/H]}$, together with symmetric monoidal functors $\ul\C_{[G/K]} \to \ul\C_{[G/H]}$, called norm functors, for each map of orbits $G/K \to G/H$.
We have the following examples in mind.
\begin{itemize}
  \item The coCartesian $G$-symmetric monoidal structure on $\ulTopG$, which is given by disjoint unions in $\ulTopG_{[G/H]}\simeq \mathbf{Top}^H$ and norm functors
    \[ \forall K<H<G: \quad \coprod\limits_{H/K}  \colon  \mathbf{Top}^K \to \mathbf{Top}^H, \quad \coprod\limits_{H/K} X = H \times_K X , \]
    where $H \times_K X$ is the quotient of $G\times X$ by the diagonal action of $K$.
  \item The Cartesian $G$-symmetric monoidal structure on $\ulTopG$, which is given by products of $H$-spaces and norm functors
    \[ \forall K<H<G: \quad \prod\limits_{H/K}  \colon  \mathbf{Top}^K \to \mathbf{Top}^H, \quad \prod\limits_{H/K} X = \Map_K (H,X) , \]
    where $\Map_K (H,X)$ is the space of $K$-equivariant maps $H\to X$ with $K$ acting on $H$ by multiplication from the right.
  \item The $G$-category $\ul\Sp^G$ of $G$-spectra has a $G$-symmetric monoidal which is given by smash products in $\ul\Sp^G_{[G/H]} \simeq \Sp^H$ and the Hill-Hopkins-Ravenel norm functors, informally given by taking $X\in Sp^H$ to the smash product of $|G/H|$ copies of $X$ with induced $G$-action.
    Nardin gave a universal property characterizing this $G$-symmetric monoidal structure by proving that $\ul\Sp^G$ admits an essentially unique $G$-symmetric monoidal structure for which the sphere spectrum is the unit, see \cite[cor. 3.28]{Nardin_thesis}.
\end{itemize}

The data of a $G$-symmetric monoidal structure, along with its coherent compatibility, is encoded by a single coCartesian fibration over the indexing category $\GFin_*$, satisfying certain Segal conditions. 
In what follows, we try to explain how this technical description is related to the intuition presented above.

We regard the symmetric monoidal structure on each fiber and the norm functors on equal footing.
To that end, consider the $G$-symmetric monoidal structure as acting on a $U$-family of objects, where we index our family be a finite $G$-set.
The members of a $U$-family $x_\bullet$ in a $G$-category $\ul\C$ correspond to the orbits of $U$, with $x_W \in \ul\C_{[W]}$ for each orbit $W\in \orb(U)$.
Given a $G$-map $I \colon U\to G/H$, we can use the $G$-symmetric monoidal structure to construct an element $\otimes_I x_\bullet \in \ul\C_{[G/H]}$.
Using the operations $\otimes_I$ we can encapsulate the data $G$-symmetric monoidal structure on $\ul\C$. 

The various operations $\otimes_I$ are subject to certain compatibility conditions, which hold upto coherent homotopy.
In order to encapsulate the compatibility of $\otimes_I$ for various $I$ it is convenient to extend $\otimes_I$ from $I \colon  U \to G/H$ to general $G$-maps of finite $G$-sets $\varphi \colon  U \to V$.
The generalized operation $\otimes_\varphi$ takes a $U$-family to a $V$-family by acting on the fibers of $\varphi$, i.e
\[ \forall W' \in \orb(V): \quad  (\otimes_\varphi x_\bullet)_{W'} = \otimes_{\varphi^{-1}(W')}  \left( x_\bullet|_{\varphi^{-1}(W')} \right) \in \ul\C_{[W']} . \] 
Note that we also need to keep track of restrictions taking a $U$-family $x_\bullet$ to a $U'$-family $x_\bullet|_{U'}$ for each inclusion of $G$-sets $U' \cof U$.

All these operations are encoded by a coCartesian fibration over $\GFin_*$, the $G$-category of finite pointed $G$-sets (see \cref{sec:G_SM_cat}), which we think of as our indexing category. 
Note that the fiber of $\GFin_*$ over $G/H$ is the given by the category of spans of finite $G$-sets $ U \hookleftarrow U' \to V$ over $G/H$, where the wrong way map $ U \hookleftarrow U'$ is an inclusion. 
Restriction to the fiber over $eH$ defines an equivalence 
\( (\GFin_*)_{[G/H]} \iso \Fin^H_* \)
to the category of finite pointed $H$-sets, described here as partly defined $H$-maps given by spans of finite $H$-sets \( \tilde{U} \hookleftarrow \tilde{U}' \to \tilde{V} \), where the wrong way map is an inclusion of finite $H$-sets.

We end this subsection by briefly sketching how to extract the tensor products and norms from a coCartesian fibration $p \colon \ul\C^\otimes \fib \GFin_*$ describing a $G$-symmetric monoidal structure on a $G$-category $\ul\C$.

First we describe the tensor product of two objects $x_1,x_2\in \ul\C_{[G/H]}$.
The $\infty$-category $\ul\C_{[G/H]}$ is given as the fiber of $p$ over $G/H \xto{=} G/H$. 
Let $U= G/H \coprod G/H$ and $I\in \GFin_*$ given by the fold map $I \colon U \to G/H$.
By the Segal conditions we have an equivalence \( \ul\C^\otimes_I \iso \ul\C_{[G/H]} \times \ul\C_{[G/H]} \) from the fiber of $p$ over $I$.
Through this equivalence we identify the ordered pair $(x_1,x_2) \in \ul\C_{[G/H]} \times \ul\C_{[G/H]}$ with an object $x_\bullet \in \ul\C^\otimes_I$ (a $U$-family). 
Choose a $p$-coCartesian lift $x_\bullet \to y$ of the span $ U \xfrom{=} U \xto{I} G/H$ over $G/H$.
The tensor product $x_1 \otimes x_2$ is given by $y\in \ul\C_{[G/H]}$.

Next we describe the norm of an object $x\in \ul\C_{[G/K]}$ along $\varphi \colon  G/K \to G/H$.
As before, the $\infty$-category $\ul\C_{[G/K]}$ is the fiber of $p$ over $G/K \xto{=} G/K$.
Consider the map $\varphi \colon  G/K \to G/H$ as an object of $\GFin_*$.
By the Segal conditions we have an equivalence \( \ul\C^\otimes_\varphi \iso \ul\C_{[G/K]} \) from the fiber of $p$ over $\varphi$.
Through this equivalence we identify $x \in \ul\C_{[G/K]} \times \ul\C_{[G/H]}$ with an object $x_\bullet \in \ul\C^\otimes_\varphi$ (a $\varphi$-family). 
Choose a $p$-coCartesian lift $x_\bullet \to y$ of the span $ G/K \xfrom{=} G/K \xto{\varphi} G/H$ over $G/H$.
The norm $\otimes_\varphi x$ is given by $y\in \ul\C_{[G/H]}$.

\section{$G$-manifolds and $G$-disks} \label{sec:G_manifolds}
Genuine $G$-factorization homology will be constructed in \cref{sec:G_factorization_homology} using parametrized $\infty$-category theory.
Our goal in this section is to construct and study the $G$-$\infty$-categories 
needed there.
Most of this section is devoted to the construction of these $G$-$\infty$-categories and their $G$-symmetric monoidal structures.
These constructions may be of independent interest, as they 
provide a bridge between geometry of manifolds with a finite group action and the theory of parametrized $\infty$-categories. 

In \cref{sec:G_mfld} we construct $\Gmfld$, the $G$-category of $G$-manifolds.
The construction is inspired by the model of $\ulTopG$ described in \cref{ulTopG_as_G_spaces_over_orbits}.
We then turn to study its relation to $G$-vector bundles, and construct an equivariant version of the tangent classifier functor of \cite{AF}.
This $G$-tangent classifier is used in \cref{sec:G_framed_mfld} to construct framed variants of $\Gmfld$.

Next, we turn our attention to $G$-disjoint unions. 
In \cref{sec:G_disj_union} we define a $G$-symmetric monoidal structure on $\Gmfld$ encoding disjoint unions and topological inductions.
The construction is quite explicit, and relies on the unfurling construction Barwick, introduced in \cite{BarwickSMF1}.
In \cref{sec:G_disj_union} we lift $G$-disjoint unions to a $G$-symmetric monoidal structures on the framed variants of $\Gmfld$.
Our main tool will be the $G$-coCartesian structures constructed in \cite{Parametrized_algebra}.

The $G$-symmetric monoidal structure of $G$-disjoint unions will be used in \cref{sec:G_factorization_homology} when defining factorization homology in two ways.
First, the expected interaction of genuine $G$-factorization homology with disjoint unions and topological inductions is expressed by being a $G$-symmetric monoidal functor from $\Gmfld$.
Second, the definition of $G$-disk algebras relies on the definition of the $G$-symmetric monoidal $G$-$\infty$-category of $G$-disks, defined in \cref{sec:G_disk}.

Next, we turn to study our constructions. 
In \cref{sec:G_disk_as_env} we show that $G$-disks are exactly the $G$-manifolds generated from linear representations of subgroups $H<G$ by taking disjoint unions and topological inductions.
In \cref{sec:G_disks_and_configurations} we compare equivariant embeddings of $G$-disks with equivariant configurations spaces.
The results of this comparison will be used in \cref{sec:G_tensor_excision} to show that genuine $G$-factorization homology satisfies $\otimes$-excision.
In \cref{sec:DV_and_Gdisk_V_fr} we define the $G$-$\infty$-operad $\EE_V$ of little $V$-disks, and use the results of \cref{sec:G_disks_and_configurations} to relate $\EE_V$ to $V$-framed $G$-disks.

\subsection{The $G$-category of $G$-manifolds} \label{sec:G_mfld}
The goal of this subsection is to give an explicit model for the $G$-$\infty$-category $\Gmfld$ of $n$-dimensional $G$-manifolds.

Before going into the details of the construction, let us first recall the construction of the $\infty$-category $\mfld^G$ of $G$-manifolds, achieved by a standard procedure.
Let $M_1,M_2$ be smooth $n$-dimensional manifolds equipped with a smooth action of a finite group $G$.
The set $Emb^G(M_1,M_2)$ of smooth $G$-equivariant open embeddings $M_1 \cof M_2$ comes with a natural topology, making the category $\mfld^G$ of $n$-dimensional $G$-manifolds into a topological category.
We consider $\mfld^G$ as an $\infty$-category by taking its coherent nerve (\cite[def. 1.1.5.5]{HTT}). 

We can extend the construction of $\mfld^G$ to construct the $G$-$\infty$-category $\Gmfld$ as follows.
Consider the $\infty$-categories $\mfld^H$ of $n$-dimensional $H$-manifolds and $H$-embeddings for all subgroups $H<G$.
The $\infty$-categories $\mfld^H$ form a diagram of $\infty$-categories, by related by two types of functors:
\begin{enumerate}
  \item 
    First, if $M$ is a $G$-manifold and $H<G$ we can consider $M$ as an $H$-manifold, which defines a functor of topological categories
    \( \mfld^G \to \mfld^H \). 
    Similarly we have \( \mfld^H \to \mfld^K \) for $K<H<G$.
  \item 
    Second, suppose $K,H<G$ are conjugate subgroups, i.e $H = g K g^{-1}$ for some $g\in G$, and $M$ is an $H$-manifold.
    We can consider $M$ as a $K$-manifold by twisting the $H$-action by conjugation,
    defining an isomorphism of topological categories 
\( conj^H_K \colon \mfld^H \to \mfld^K \). 
\end{enumerate}
A standard verification shows that the topological categories $\mfld^H$ define a diagram of topological categories indexed by subgroups $H<G$, with functors indexed contravariantly by $G$-maps $G/K \to G/H$.
Note that this indexing category is equivalent to the orbit category $\OG$ (see \cref{def:orbit_cat}).
Composing with the topological nerve we get a diagram of $\infty$-categories 
\begin{align*}
 \mfld^\bullet \colon \OGop \to \Cat_\infty, \quad G/H \mapsto \NSing(\mfld^H),
\end{align*}
which we can unstraighten to a coCartesian fibration $UnSt( \mfld^\bullet) \fib \OGop$ (see \cite[sec. 3.2]{HTT}). 
The casual reader can use $UnSt(\mfld^\bullet)$ as the definition of the $G$-category of $G$-manifolds, and skip the rest of this subsection. 

The construction of $UnSt(\mfld^\bullet)$ is unsatisfying to us in two respects. 
First, it depends on an implicit choice of an inverse to the inclusion of the full subcategory $ \set{G/H}_{H<G} \subset \OG$ into the category of $G$-orbits (which is equivalent to choosing a basepoint for every transitive $G$-set).
Second, manipulating $UnSt(\mfld^\bullet)$ as a simplicial set is inconvenient, as unstraightening is a right adjoint functor.
Instead of working with $UnSt(\mfld^\bullet)$ we construct an equivalent $G$-$\infty$-category $\Gmfld$ (\cref{def:Gmfld}) which admits a more accessible description as a simplicial set. 
This is the main construction of this subsection.

Let us briefly describe our strategy for constructing $\Gmfld$, inspired by the model of $\ulTopG$ described in \cref{ulTopG_as_G_spaces_over_orbits}.
First we construct a topological category $\OGmfld$ equipped with functor to the orbit category, and show that the topological nerve defines a Cartesian fibration $\NSing(\OGmfld) \to \OG$ of simplicial sets, which classifies a diagram of $\infty$-categories equivalent to $\NSing(\mfld^\bullet)$. 
We then define $\Gmfld$ (\cref{def:Gmfld}) as the coCartesian fibration dual to $\NSing(\OGmfld) \to \OG$, which  classifies the same diagram $\NSing(\mfld^\bullet)$. 
The dual coCartesian fibration admits an explicit construction span categories (see \cite{DualizingCarFibs}) which we use to describe the objects and morphisms of $\Gmfld$
and the coCartesian morphisms of $\Gmfld \fib \OGop$.

\begin{rem}
  In this section we denote objects of the orbit category by $O\in \OGop$, as opposed to $G/H$ elsewhere. 
  This is merely for notational convenience.
\end{rem}

\paragraph{$\OG$-manifolds and their spaces of smooth equivariant embeddings.}
We start by defining $\OG$-manifolds and spaces of smooth equivariant embeddings which will serve as objects and mapping spaces of the topological category $\OGmfld$, see \cref{def:mfldG}.

\begin{mydef}
  An \emph{$\OG$-manifold} $M \to O$ is a smooth $n$-dimensional manifold $M$ with an action of $G$ on $M$ by smooth maps, together with a $G$-map $M\to O$ from the underlying $G$-space of the manifold $M$ to a $G$-orbit $O \in \OG$. 
\end{mydef}
We \emph{always} think of an $\OG$-manifold $M\to G/H$ as encoding a smooth $n$-dimensional manifold with an action of $H$, given by the fiber $M|_{H}$ of the $G$-map $M\to G/H$ over the coset $H$.
Note that a choice of a basepoint $o\in O$ induces an isomorphism $G/H \stackrel{\=}{\to} O, \, gH \mapsto g\cdot o$, where $H<G$ is the stabilizer of $o$. 
We therefore think of an $\OG$-manifold $M\to O$ as encoding the smooth action of $H=Stab(o)$ on the fiber $M|_H$.

\begin{notation}
  Suppose $M,N$ are smooth $n$-dimensional manifolds. 
  Denote by $C^\infty(M,N)$ the space of \emph{smooth} maps $M\to N$ with the \emph{compact-open} topology. 
\end{notation}
\begin{mydef} \label{def:G_mfld_map}
  Let $M_1 \to O_1, M_2 \to O_2$ be $\OG$-manifolds.
  For \( \varphi \colon O_1 \to O_2 \) a map in $\OG$, define $\EmbOG_{\varphi}(M_1,M_2)\subset C^\infty(M_1,M_2)$ as the subspace of smooth maps $f \colon M_1 \to M_2$ such that 
  \begin{enumerate}
    \item $f$ is a $G$-map 
    \item $f$ is over $\varphi$, i.e 
      \begin{align} \label{eq:G_mfld_map}
        \diag{
          M_1 \ar[d] \ar[r]^{f} & M_2 \ar[d] \\
          O_1 \ar[r]^{\varphi} & O_2
        }
      \end{align}
      is a commutative square of $G$-spaces.
    \item the induced map \( M_1 \to O_1 \times_{O_2} M_2 \) is an embedding. 
  \end{enumerate}
  Define the topological space $\EmbOG(M_1,M_2)$ as the coproduct 
  \begin{align}
    \EmbOG(M_1,M_2) := \coprod_{\varphi} \EmbOG_\varphi(M_1,M_2) ,
  \end{align}
  where the coproduct is indexed by the set $\Hom_{\OG}(O_1,O_2)$. 
\end{mydef}

\begin{notation}
  When the orbit map $\varphi$ is an identity $G/H \xto{=} G/H$ we use the notation \( Emb^G_{G/H}(M_1,M_2) \) for the space $\EmbOG_\varphi(M_1,M_2)$ of $G$-equivariant embeddings $M_1 \to M_2$ over $G/H \xto{=} G/H$.
  Restriction to the fiber over $H$ defines a homeomorphism from $Emb^G_{G/H}(M_1,M_2)$ with the space of $H$-equivariant embeddings $M_1|_H \to M_2|_H$ between the fibers over $H$.
\end{notation}

\begin{mydef}
  Let $M_1\to O_1, M_2 \to O_2$ be $\OG$-manifolds.
  A \emph{$G$-isotopy over $\varphi \colon O_1 \to O_2$} is a path in $\EmbOG_{\varphi}(M_1,M_2)$.
  When $M_1 \to O, M_2 \to O$ are over the same orbit we call a path in $\EmbOG_O(M_1,M_2)$ a $G$-isotopy over $O$. 
\end{mydef}
Note that a $G$-isotopy over $G/H$ is equivalent to an $H$-equivariant isotopy between two $H$-equivariant embeddings $M_1|_{H} \to M_2|_{H}$.

\paragraph{The topological category of $\OG$-manifolds.}
We now turn to the definition of the topological category of $\OG$-manifolds.
Note that the pullback of smooth embeddings of $n$-dimensional manifolds is a smooth embedding, therefore we have
\begin{lem}
  Let $M_1 \to O_1,M_2 \to O_2,M_3 \to O_3$ be $\OG$-manifolds. 
  The composition of smooth functions defines a continuous map 
  \[ \EmbOG(M_2,M_3) \times \EmbOG(M_1,M_2) \to \EmbOG(M_1,M_3), \quad (g,f) \mapsto g\circ f . \] 
\end{lem}

\begin{mydef} \label{def:mfldG}
  The \emph{category of $\OG$-manifolds $\OGmfld$} is the topological category whose objects are a $\OG$-manifolds.
  The morphism space from $M_1 \to O_1$ to $M_2 \to O_2$ is given by \( \Map_{\OGmfld}(M_1,M_2) := \EmbOG(M_1,M_2) \).

  Define a forgetful functor \( q \colon \OGmfld \to \OG \) by sending $M\to O$ to the orbit $O$, and the subspace $\EmbOG_\varphi(M_1,M_2)\subset \EmbOG (M_1,M_2)$ to $\varphi \in \Hom_{\OG}(O_1,O_2)$.

  By \cite[ex. 1.1.5.12]{HTT} the topological nerve $\NSing(\OGmfld)$ is an $\infty$-category, and by \cite[ex. 1.1.5.8]{HTT} the topological nerve of $\OG$ can be identified with its ordinary nerve, which we identify with $\OG$ by standard abuse of notation. 

  Applying the topological nerve functor of \cite[1.1.5.5]{HTT} to $q$ produces a functor of $\infty$-categories \( \NSing(q) \colon  \NSing(\OGmfld) \to \OG \). 
\end{mydef}
In particular, an object of the $\infty$-category $\NSing(\OGmfld)$ is an $\OG$-manifold $M\to O$, a map is given by a commutative square \cref{eq:G_mfld_map} satisfying the conditions of \cref{def:G_mfld_map}, and by \cite[thm. 1.1.5.13]{HTT} the mapping spaces of $\NSing(\OGmfld)$ are weakly equivalent to the mapping spaces of $\OGmfld$.

\begin{rem}
  The fiber of $\OGmfld \to \OG$ over an orbit $G/H$ is the topological nerve of the topological category whose objects are $\OG$-manifolds $M \to G/H$ and morphism spaces are $Emb^G_{G/H}(M_1,M_2)$.
  This topological category is equivalent to the category $\mfld^H$ of $H$-manifolds and $H$-equivariant embeddings by restriction to the fibers over $H$.
\end{rem}
\begin{rem}
  We caution the reader not to pass to $\infty$-categories prematurely.
  One can construct the topological category $\OGmfld$ as a subcategory of the topological arrow category $\mfld^G \downarrow \OG$.
  However, the $\infty$-category $\NSing(\OGmfld)$ is \emph{not} a subcategory of the topological nerve $\NSing(\mfld^G \downarrow \OG)$ in the sense of \cite[sec. 1.2.11]{HTT}.
  To see this note that a subcategory of $\NSing(\mfld^G \downarrow \OG)$ is specified by a subcategory of its homotopy category $ho \NSing (\mfld^G \downarrow \OG)$, and therefore given by a choosing connected components of each mapping space of $\mfld^G \downarrow \OG$.
  On the other hand condition (3) of \cref{def:mfldG} is not preserved by $G$-homotopy equivalence, so the subspace 
  \[ Emb^G_{G/H} ( M_1 \to O_1, M_2 \to O_2) \subset \Map_{\mfld^G\downarrow\OG}(M_1\to O_1,M_2\to O_2) \]
  is not given by a set of connected components.  
  The same phenomenon exists in the non-equivariant setting. 
\end{rem}

\paragraph{Equivalences of $\OG$-manifolds.}
Unwinding the definition of equivalence in a nerve of a topological category, we see that a map \( f \colon M_1 \to M_2 \) in $\OGmfld$ is 
an equivalence in $\NSing(\OGmfld)$ if it has a $G$-isotopy inverse: a map \(g \colon M_2 \to M_1\) in $\OGmfld$, 
together with a $G$-isotopy over $id_{q(M_1)}$ from $g\circ f$ to $id_{M_1}$ and a $G$-isotopy over $id_{q(M_2)}$ from $f \circ g$ to $id_{M_2}$.
\begin{mydef}
  We say that a map $f \colon M_1 \to M_2$ of $\OG$-manifolds is a \emph{$G$-isotopy equivalence} if it is an equivalence in the $\infty$-category $\NSing(\OGmfld)$.
\end{mydef}
Note that an equivalence $f$ always lies over an isomorphism of orbits $q(f) \colon O_1 \to O_2$.
Using the homeomorphism between the mapping space $Emb^G_{G/H}(M_1,M_2)$ over an orbit $G/H$ and the space of $H$-equivariant embeddings $M_1|_{H} \to M_2|_{H}$ we see that
a map $f \colon M_1 \to M_2$ over an orbit $G/H$ is an equivalence in $\NSing(\OGmfld)$ if and only if its restriction to the fibers $f|_H  \colon  M_1|_{H} \to M_2|_{H}$ is invertible upto $H$-isotopy.
In particular, $f$ need not induce an equivariant diffeomorphism. 
Nonetheless, its existence is enough to ensure that there exists an equivariant diffeomorphism between underlying manifolds.
We learned the following argument from an answer of Ian Agol on MathOverflow \cite{Agol_isotopic_vs_diffeomorphic}, which we reproduce here (with addition of a $G$-action).
\begin{prop}
  Let $M_1 \to G/H$ and $M_2 \to G/H$ be two $\OG$-manifolds over $G/H$.
  If $f\in Emb^G_{G/H}(M_1,M_2)$ and $g\in Emb^G_{G/H}(M_2,M_1)$ are $G$-isotopy inverses over $G/H$ 
  then there exists a $G$-equivariant diffeomorphism $M_1 \= M_2$ over $G/H$.
\end{prop}
\begin{proof}
  We prove the statement by reduction.
  Since  \( Emb^G_{G/H}(M_1,M_2) \) is homeomorphic to the space of $H$-invariant embeddings between  $M_1|_H \to M_2|_H$ it is enough to consider the case $G=H$.

  Suppose $M,N$ are $n$-dimensional manifolds with smooth actions of $G$, and we are given $G$-equivariant embeddings $f \colon  M \to N, \, g \colon N \to M$.
  Consider the direct limit 
  \[ X = \colim ( M \xto{f} N \xto{g} M   \xto{f} N \xto{g} \cdots ) , \]
  given by the explicit model \( M \times \N \sqcup N \times \N / \sim \) with equivalence relation generated by $(m,k) \simeq (f(m),k)$ and $(n,k) \simeq (g(n),k+1)$.
  Then $X$ is a smooth manifold with an action of $G$, as a sequential union of nested open submanifolds.

  Since $X$ is $G$-diffeomorphic  to \( Y = \colim ( N \xto{g} M   \xto{f} N \xto{g} \cdots ) \) (removing the first term of the sequence does not change the colimit), it is enough to show that $X$ is $G$-diffeomorphic to $M$.

  Note that $X$ is $G$-diffeomorphic to  $\colim ( M \xto{F_1} M \xto{F_1} M \xto{F_1} \cdots ) $ for $F_1 = g \circ f$, and $F_1$ is $G$-isotopic to $id_M$.
  Let $F_t \colon M \to M , t\in [0,1]$ be the $G$-isotopy from $F_0=id_M$ to $F_1= g \circ f$, and define 
  \( X_t = \colim ( M \xto{F_t} M \xto{F_t} \cdots ) \), so that $X_1=X$ and $X_0 = M$.

  Choose a sequence of compact $G$-submanifolds with boundary $K_1 \subset K_2 \subset K_3 \subset \cdots M$ such that $M= \cup_i K_i$ and $F(K_i \times [0,1]) \subset int(K_{i+1})$.
  Such a sequence can be chosen inductively using a $G$-invariant Morse function on $M$ (which exists by \cite[cor. 4.10]{Wasserman}). 
  Define \( Y_t = \colim ( K_1 \xto{F_t} K_2 \xto{F_t} K_3 \xto{F_t} \cdots) \) using the restrictions of the $F_t$ to the subsets $K_i$.
  We claim that $Y_t=X_t$, using the standard model for direct limits.
  Write $X_t =  M \times \N / (x,i) \sim (F_t(x), i+1)$, and note that $Y_t \subseteq X_t$ as the points $(x,i)$ with $x\in K_i$.
  We claim that each point $x\in X_t$ is in $Y_t$. 
  Represent $x$ by $(x,i)\in M \times \N$, then since $M=\cup K_i$ we have $x\in K_j$ for some $j\in \N$.
  If $j\leq i$ then $K_j \subset K_i$, so $x\in K_i$, hence $(x,i)$ represents an point in $Y_t$. 
  Otherwise $(x,i) \sim (F_t^{j-i}(x), j)$ in represents the same point in $X_t$, and since $F_t^{j-i}(K_i) \subset K_j$ we get $F_t^{j-i}(x) \in K_{j}$, so $(F_t^{j-i}(x), j)$ represents an element of $Y_t$.

  We showed that $Y_t = X_t$, so it is enough to prove that $Y_0 \= M$ is $G$-diffeomorphic to $Y_1 \= X$.
  By definition we have $Y_0 = \colim ( K_1 \cof K_2 \cof K_3 \cof \cdots)$ and $Y_1 = \colim ( K_1 \xto{F_1} K_2 \xto{F_1} K_3 \xto{F_1} \cdots )$, hence it is enough to construct compatible $G$-diffeomorphims \( \phi_i \colon  K_i \cof K_i \), i.e satisfying $\phi_{i+1}|_{K_i} = F_1 \circ \phi_i$.
  
  We now inductively construct $G$-equivariant maps $ G^i \colon  K_i \times [0,1] \to K_i$ such that $G_0=Id_{K_i}$, $\forall t\in [0,1]: G_t \colon K_i \to K_i$ is a diffeomorphism and $\forall x\in K_i, t\in [0,1]: \, F_t \circ G^i_t(x) = G^{i+1}_t(x)$, i.e the diagram 
  \begin{align*}
    \xymatrix{
      K_i \times [0,1] \ar@{^(->}[r] \ar[d]^{G^i \times Id} & K_{i+1} \times [0,1] \ar[d]^{G^{i+1}} \\
      K_i \times [0,1] \ar[r]^{F|_{K_i \times [0,1]}} & K_{i+1}
    }
  \end{align*}
  commutes.
  \footnote{The map $G^i$ is an equivariant \emph{diffeotopy} in terminology of \cite{Hirsch} and an \emph{equivariant isotopy starting from the identity} in the terminology of \cite{Bredon_transformation_groups}.}

  We start with setting $G^1_t = Id_{K_1}$.
  Assume that a $G^i$ has been constructed. 
  Consider the isotopy \( K_i \times [0,1] \xto{G^i \times Id} K_i \times [0,1] \xto{F|_{K_i \times [0,1]}} K_{i+1} \).
  Since $K_i \subset K_{i+1}$ is a compact submanifold and $F(K_i) \subset Int(K_{i+1})$ the conditions of the isotopy extension theorem \cite[ch. 8 thm. 1.3]{Hirsch} are satisfied.
  Therefore there exists a diffeotopy $\tilde{G}^{i+1} \colon  K_{i+1} \times [0,1] \to K_{i+1}$ which extends the isotopy \( K_i \times [0,1] \xto{G^i \times Id} K_i \times [0,1] \xto{F|_{K_i \times [0,1]}} K_{i+1} \) and satisfies $\tilde{G}^{i+1}_0 = Id_{K_{i+1}}$, but might not be $G$-equivariant. 
  Since $K_{i+1}$ is compact we can apply \cite[thm 3.1]{Bredon_transformation_groups}, and get a $G$-equivariant diffeotopy $G^{i+1} \colon  K_{i+1} \times [0,1] \to K_{i+1}$ with $G^{i+1}_0 = \tilde{G}^{i+1}_0=Id_{K_{i+1}}$ and which agrees with $\tilde{G}^{i+1}$ on the subset \( \set{ x\in K_{i+1} \, \vert \, \forall g\in G, t\in [0,1]: \, \tilde{G}^{i+1}_t(gx)=g\tilde{G}^{i+1}_t(x) } \).
  In particular, for $x\in K_i$ we have $\tilde{G}^{i+1}_t(x)=F_t G^i_t(x)$, so the $G$-equivariant diffeotopy $G^{i+1}$ agrees with $\tilde{G}^{i+1}$ on $K_i \times [0,1]$.

  Setting $\phi^i=G^i_1$ gives the compatible $G$-diffeomorphisms proving that $Y_0 \= M$ is indeed $G$-diffeomorphic to $Y_1 \= X$.
\end{proof}

\paragraph{Cartesian edges in $\OGmfld$.}
We now identify the Cartesian edges of the forgetful functor \(\NSing(\OGmfld) \to \OG\), as well as the coCartesian edges over isomorphisms.
We start with
\begin{lem}
  The forgetful functor \( \mathbf{N}(q) \colon \NSing(\OGmfld) \to \OG \) is an inner fibration.
\end{lem}
\begin{proof}
  For every pair $M_1\to O_1,M_2\to O_2$ of $\OG$-manifolds, $q$ induces a Kan fibration
  \[ \Map_{\Sing(\OGmfld)}(M_1,M_2) \to \Map_{\Sing(\OG)}(O_1,O_2) = \Hom_{\OG}(O_1,O_2), \] 
  because its a map from a Kan simplicial complex and to discrete simplicial set.
  Therefore by \cite[prop. 2.4.1.10(1)]{HTT} 
  the functor
  $\NSing(q)$ is an inner fibration.
\end{proof}

Note that a map $M\to O$ from an $n$-dimensional manifold to a finite set is always a submersion, so its pullback along any map of finite sets is an $n$-dimensional manifold.
\begin{prop} \label{prop:pullback_is_Cart}
  Suppose that \(\varphi \colon  O_1 \to O_2 \) be a map of orbits, and $M \to O_2$ a $\OG$-manifold. 
  Then the pullback square of topological $G$-spaces
  \begin{align*}
    \diag{
      O_1 \times_{O_2} M \ar[d] \ar[r]^-{f} \pullbackcorner & M \ar[d] \\
      O_1 \ar[r] & O_2
    }
  \end{align*}
  defines a $\mathbf{N}(q)$-Cartesian morphism $f$ in $\OGmfld$. 
  In particular, $\mathbf{N}(q)$ is a Cartesian fibration.
\end{prop}
\begin{proof}
  Checking that $f$ satisfies the conditions of \cref{def:mfldG} is immediate. 
  
  By \cite[prop. 2.4.1.1 (2)]{HTT} the morphism $f$ is $\NSing(q)$-Cartesian if and only if, for every $\OG$-manifold $T \to O$, the square of spaces
  \begin{align*}
    \xymatrix{
      Map_{\Sing(\mfld)}(T,O_1 \times_{O_2} M) \ar@{->>}[d] \ar[r]^-{f_*} & Map_{\Sing(\OGmfld)}(T,M) \ar@{->>}[d] \\
      \Hom_{\OG}(O,O_1) \ar[r]^{q(f)_*} & \Hom_{\OG}(O,O2)
    }
  \end{align*}
  is a homotopy pullback square.
  Since the vertical maps are Kan fibrations, this square is a homotopy pullback if and only if the horizontal map 
  \begin{align*}
    \xymatrix{
      \Map(T,O_1 \times_{O_2} M) \ar@{->>}[dr] \ar[rr] & & \Hom(O,O_1) \times_{\Hom(O,O_2)} \Map(T,M) \ar@{->>}[dl] \\
      & \Hom(O,O_1)
    }
  \end{align*}
  is a homotopy equivalence, or equivalently, if $f_*$ induces an equivalence between the fiber over every $\tau\in\Hom(O,O_1)$.

  Let $\tau:O \to O_1$. Then $f_*$ induces a map of fibers over $\tau$
  \begin{align*}
    & \EmbOG_{\tau}(T,O_1 \times_{O_2} M) \to \set{\tau} \times_{\Hom(O,O_2)} \Map(T,M)=\EmbOG_{q(f)\circ\tau}(T,M), \\
    & \left( \diag{ T \ar[d] \ar[r]^-{g} & O_1 \times_{O_2} M \ar[d] \\ O \ar[r]^{\tau} & O_1 } \right) 
    \mapsto 
    \left( \diag{ T \ar[d] \ar[r]^-{g} & O_1 \times_{O_2} M \ar[d] \ar[r]^-{f} \pullbackcorner & M \ar[d] \\ O \ar[r]^{\tau} & O_1 \ar[r]^{q(f)} & O_2 } \right)
  \end{align*}
  This continuous map is a bijection by the universal property of the pullback.
  We leave to it to the reader to verify it is an open map using the definition of the compact-open topology.
\end{proof}
This gives the following complete description of the Cartesian edges in $\OGmfld$.
\begin{cor} \label{G_iso_PB}
  A morphism \eqref{eq:G_mfld_map} is $\mathbf{N}(q)$-Cartesian if and only if it is equivalent to a pullback, i.e. the morphism 
  \begin{align} \label{eq:G_iso_equiv}
    \diag{ M_1 \ar[d] \ar[r] & O_1 \times_{O_2} M_2 \ar[d] \\ O_1 \ar[r]^{=} & O_1 }
  \end{align}
  is a $G$-isotopy equivalence.
\end{cor}
\begin{proof}
  Factor the morphism \eqref{eq:G_mfld_map} as the composition of \eqref{eq:G_iso_equiv} and a pullback square. 
  Combining \cref{prop:pullback_is_Cart} and \cite[prop. 2.4.1.7]{HTT} we see that the morphism \eqref{eq:G_mfld_map} is $\mathbf{N}(q)$-Cartesian if and only if the map above is $\mathbf{N}(q)$-Cartesian. 
  Since the morphism \eqref{eq:G_iso_equiv} lies over an equivalence it is $\mathbf{N}(q)$-Cartesian if and only if it is an equivalence, by \cite[prop. 2.4.1.5]{HTT}.
\end{proof}

\paragraph{Construction of the $G$-category of $G$-manifolds.}
The construction of the $G$-category $\Gmfld$ now follows easily from the description of the Cartesian fibration $\NSing(q)$ and the explicit construction of \cite{DualizingCarFibs}. 

\begin{mydef} \label{def:Gmfld}
  Let \( p \colon \Gmfld \to \OGop \) be the dual of the Cartesian fibration \( \OGmfld \to \OG \) in the sense of \cite[def. 3.5]{DualizingCarFibs}.
  Explicitly, $\Gmfld$ is the pullback of the effective Burnside category 
  \[ A^{eff}(\OGmfld, \, \OGmfld \times_{\OG} \OG^{\=}, \, \operatorname{q-Cart}(\OGmfld)) \]
  along the equivalence \( \OGop \trivcof A^{eff}(\OG, \OG^{\=}, \OG) \), where $\OG^{\=}$ is the maximal subgroupoid of $\OG$ and $\operatorname{q-Cart}(\OGmfld)\subset \OGmfld$ is the subcategory spanned by all objects and morphisms which are $q$-Cartesian.
\end{mydef}
By \cite[prop. 3.4]{DualizingCarFibs} the map \( p \colon \Gmfld \to \OGop \) is a coCartesian fibration, and we have an explicit description of the objects and morphisms of $\Gmfld$.
The objects of the total $\infty$-category $\Gmfld$ are $\OG$-manifolds $M \to O$.
A morphism in $\Gmfld$ from $M_1 \to O_1$ to $M_2 \to O_2$ is a diagram of the form 
\begin{align} \label{diag:Gmfld_morphism_as_span}
  \xymatrix{
    M_1 \ar[d] & \ar[l] M \ar[d] \ar[r] & M_2 \ar[d] \\
    O_1 & \ar[l] O_2 \ar[r]^{=} & O_2 
  }
\end{align}
where the left square is a coCartesian edge in $\OGmfld$ (in other words, it is equivalent to a pullback square, see \cref{G_iso_PB}).
This arrow is $p$-coCartesian exactly when the right square is a $G$-isotopy equivalence.

Without loss of generality we will represent a morphism in $\Gmfld$ by a span \eqref{diag:Gmfld_morphism_as_span} where the left square is a pullback square. 

\begin{rem} \label{Gmfld_GH_is_mfld_H}
  Let $H<G$ be a subgroup. 
  Topological induction defines a functor 
  \[ G \times_H (-) \colon  \mfld^H \to \Gmfld_{[G/H]}, \quad G \times_H M = \left( (G \times M) / G \to (G \times pt)/H = G/H \right) \]
  where we quotient by the $H$-action \( h \cdot ( g,x)= (gh^{-1}, gx) \). 
  Topological induction if a functor of topological categories, and in fact an equivalence of topological categories
  \( \mfld^H \iso \Gmfld_{[G/H]} \),
  with inverse \( (M \to G/H) \mapsto M|_{eH} \) given by restriction to the fiber over $eH$. 

  Informally, the coCartesian fibration $\Gmfld \to \OGop$ classifies the functor \( \OGop \to \Cat_\infty \) sending $G/H$ to $\mfld^H$.
\end{rem}

\begin{notation}
  We will refer to $\Gmfld$ as \emph{the $G$-category of $G$-manifolds}, to stress its conceptual role and not its technical construction.
  We urge the reader to regards the objects of $\Gmfld$ not as $\OG$-manifolds (which they are), but as a technical means of encoding manifolds with an action of a subgroup of $G$.
  This naming convention is also compatible with \cite[ex. 7.5]{Expose1}, where $\ul{\mathbf{Top}}_T$ is referred to as the $T$-$\infty$-category of $T$-spaces. 
\end{notation}

By construction, we have a simple description of the fiberwise opposite\footnote{The superscript ``vop'' stands for taking ``vertical opposites''.}
category $(\Gmfld)^{vop}$, introduced in \cite[sec. 3]{Expose1}.
It is helpful to keep this description in mind when we use the parametrized Yoneda embedding to construct the equivariant tangent classifier in \cref{const:G_tangent_classifier}.
\begin{prop} \label{Gmfld_vop}
  Applying the opposite $\infty$-category functor $(-)^{op}$ to the Cartesian fibration $\OGmfld \to \OG$ produces a $G$-category $(\OGmfld)^{op} \fib \OGop$ equivalent to  $(\Gmfld)^{vop} \fib \OGop$.
\end{prop}
\begin{proof}
  By \cite[def. 3.1]{Expose1} the opposite $G$-category $(\Gmfld)^{vop} \fib \OGop$ is given by taking the opposite of the dual Cartesian fibration $ (\Gmfld)_{\wedge}  \fib \OG$.
  The result follows, since taking the dual coCartesian fibration is homotopy inverse to taking the dual Cartesian fibration (see \cite[thm. 1.7]{DualizingCarFibs}).
\end{proof}

\subsection{Representations, $G$-vector bundles and the $G$-tangent classifier} \label{sec:G_VB_and_G_tangent_classifier}
In this subsection we study the relation between $G$-vector bundles, $H$-representations of subgroups $H<G$ and the $G$-category of $G$-manifolds, $\Gmfld$, constructed in \cref{sec:G_mfld}.
We do this by identifying $H$-representations with $G$-vector bundles over $G/H$, which in turn span a full $G$-subcategory $\ul{\Rep}^G_n \subset \Gmfld$.
An equivariant version of ``smooth Kister's theorem'' implies that $\ul{\Rep}^G_n$ is in fact a $G$-$\infty$-groupoid, which can be identified with the $G$-space classifying $n$-dimensional $G$-vector bundles, $BO_n(G)$.
We use $\ul{\Rep}^G_n$ to construct an equivariant version of the tangent classifier of \cite[sec 2.1]{AF} (see \cref{const:G_tangent_classifier}), which will be used in \cref{sec:G_framed_mfld} to define equivariant tangential structures on $G$-manifolds. 
It is worth noting that parametrized $\infty$-category theory is essential for \cref{const:G_tangent_classifier}, which relies on the identification of the \emph{$G$-space} $BO_n(G)$ with a full $G$-subcategory of $\Gmfld$.

We start by recalling the standard definition of $G$-vector bundles.
\begin{mydef}[{see \cite[sect. VI.2]{Bredon_transformation_groups}, \cite[ch. I, def. 9.1]{tom_Dieck_transformation_groups}}]
  Let $X$ be a $G$-space. 
  A \emph{$G$-vector bundle} over $X$ is a (real) vector bundle \( p \colon  E \to X \) together with a $G$-action on $E$ by bundle maps (i.e linear action on each fiber) such that $p$ is a $G$-map. 
  We say $p \colon E\to X$ is \emph{smooth} if $E,X$ are (smooth) $G$-manifolds and $p$ is a smooth map.
  Let $G-\vect_{/X}$ denote the category of $G$-vector bundles over $X$.
\end{mydef}
Note that $G$-vector bundles are stable under pullback along $G$-maps, and that a $G$-vector bundle over a point is the same as a $G$-representation.
It is useful to keep in mind the correspondence between representations of subgroups $H<G$ and $G$-vector bundles over the orbit $G/H$:
\begin{prop}\cite[special case of prop. I.9.2]{tom_Dieck_transformation_groups}
  Let $H<G$ be a subgroup.
  Restriction to the fiber over $[eH]$ gives an equivalence 
  \( G-\vect_{/(G/H)} \iso H-\vect_{/pt} \= \Rep^H \)
  from the category of $G$-vector bundles over the orbit $G/H$ to the category of $H$-representations.
  An inverse is given by sending a representation of $H$ on $\R^n$ to its topological induction \( G \times_H \R^n \). 
\end{prop}

The subject of this subsection is the following $G$-subcategory.
\begin{mydef} \label{def:RepGn}
  Let $\ulRep^G_n \subset \Gmfld$ be the full $G$-subcategory spanned by $G$-vector bundles $(E \to G/H)$, i.e $\OG$-manifolds $E\to G/H$ such that $E$ can be endowed with a structure of a $G$-vector bundle over $G/H$. 
\end{mydef}
\begin{rem}
  We will use $G$-vector bundles as a model for ``$G$-disks''.
  Specifically, an embedding of a $G$-disk in an $\OG$-manifold $M\in \Gmfld$ is just a map in $\Gmfld$ with target is $M$ and domain in $\ulRep^G_n$.
  Genuine $G$-factorization homology is defined as a parametrized colimit over finite disjoint unions of $G$-disks in $M$ (see \cref{def:G_FH_as_colimit}).
  In \cref{sec:G_disk} we organize these disjoint unions into a $G$-$\infty$-category $\Gdisk$.
\end{rem}

In order to see the close relation of $\ulRep^G_n$ with representation theory we use the following equivariant version of the ``smooth Kister-Mazur'' theorem (see \cite{Kupers}).
\begin{prop} \label{G_Kister}
  Let $V$ be a finite dimensional real representation of $H<G$. 
  Let \( \Aut_{\Rep^H}(V) \) be the automorphism group of $V$ as an $H$-representation, i.e linear $H$-equivariant isomorphisms.
  Let \(Emb^H_0(V,V) \) denote the subspace of smooth $H$-equivariant embedding fixing the origin, and \( \Aut^H_0(V) \subset Emb^H_0(V,V) \) the subspace of $H$-equivariant diffeomorphisms. 
  Then the inclusions 
  \begin{align*}
    \Aut_{\Rep^H}(V) \cof \Aut^H_0(V) \cof Emb^H_0(V,V)
  \end{align*}
  are homotopy equivalences.  
\end{prop}
\begin{proof}
  The proof of \cite[thm. 2.4]{Kupers} applies verbatim when restricting to subspaces of $H$-equivariant maps after checking that the formulas for \(G^{(1)}_s, G^{(2)}_s \) produce $H$-equivariant homotopies. 
\end{proof}
The central role played by $\ulRep^G_n$ in what follows stems from the following characterization. 
\begin{prop} \label{RepGn_is_G_gpd}
  The $G$-category $\ul{\Rep}^G_n$ is a $G$-$\infty$-groupoid, with fibers $(\ul{\Rep}^H_n)_{[G/H]}$ equivalent to the topological groupoid $\Rep^H_n$ of $n$-dimensional real representations of $H$ and (linear, $H$-equivalent) isomorphisms, 
  where the mapping space \( \operatorname{Iso}_{\Rep^H}(V_0,V_1) \) is endowed with the compact-open topology. 
\end{prop}
\begin{proof}
  In order to show that $\ul{\Rep}^G_n$ is a $G$-$\infty$-groupoid we have to prove that the coCartesian fibration $\ul{\Rep}^G_n \fib \OGop$ is a left fibration. 
  By \cite[prop. 2.4.2.4]{HTT} it is enough to show that the fibers $(\ul{\Rep}^G_n)_{[G/H]}$ are $\infty$-groupoids.
  The equivalence $\Gmfld_{[G/H]} \= \mfld^H$ of \cref{Gmfld_GH_is_mfld_H} takes a $G$-vector bundle \( E \to G/H  \) to an $H$-vector bundle $E|_{eH} \to pt$, i.e. an $n$-dimensional real $H$-representation \( V = ( H \curvearrowright \R^n ) \), so we have to show that for every $V_0,V_1\in \Rep^H_n$ the inclusion \( \Aut^H (V_0,V_1) \subset Emb^H(V_0,V_1) \) 
  is a weak equivalence. 

  Let \( Emb^H_0(V_0,V_1) \subset Emb^H(V_0,V_1) \) denote the subspace of origin fixing maps. 
  Clearly the inclusion \( Emb^H_0(V_0,V_1) \cof Emb^H(V_0,V_1)\) is a homotopy equivalence.
  By \cref{G_Kister} 
  the inclusion  \( \operatorname{Iso}_{\Rep^H}(V_0,V_1) \cof Emb^H_0 (V_0, V_1) \) is a weak equivalence,  so \(  \operatorname{Iso}_{\Rep^H}(V_0,V_1) \trivcof Emb^H(V_0,V_1) \) is a weak equivalence. 
  
  In other words, the functor \( \Rep^H_n \to (\ul\Rep^G_n)_{[G/H]} \) is fully faithful. 
  Since by definition it is essentially surjective it is an equivalence of $\infty$-categories.
  In particular $(\ul\Rep^G_n)_{[G/H]}$ is equivalent to the (coherent nerve of) the topological groupoid $\Rep^H_n$, hence 
  an $\infty$-groupoid.
\end{proof}
By construction of the classifying space of $G$-vector bundles (see \cite{G_smoothing_theory,Waner}) we have the following statement.
\begin{cor} \label{RepGn_is_BOnG}
  The $G$-$\infty$-groupoid $\ulRep^G_n$ corresponds to $BO_n(G)\in \mathbf{Top}^G$, the classifying $G$-space of rank $n$ real $G$-vector bundles.
\end{cor}

We can now construct an equivariant version of the tangent classifier of Ayala-Francis (see \cite[sec. 2.1]{AF}).
\begin{construction}[$G$-tangent classifier] \label{const:G_tangent_classifier}
  Let $j \colon  \Gmfld \to \ulFun_G ( (\Gmfld)^{vop} , \ul{\operatorname{Top}}^G )$ be the parametrized Yoneda embedding $G$-functor of \cite{Expose1} (see \cref{Gmfld_vop} for a description of the fiberwise opposite $(\Gmfld)^{vop}$).
  Define a \emph{$G$-tangent classifier} by the composition of $G$-functors 
  \begin{align*}
    \tau \colon \Gmfld \xto{j} \ulFun_G ( (\Gmfld)^{vop} , \ul{\mathbf{Top}}^G ) \to \ulFun_G ( (\ul\Rep^G_n)^{vop} , \ulTopG ) \simeq \ulTopG_{/\ul{ BO_n(G)}}
  \end{align*}
  where the last equivalence is given by parametrized straightening/unstraightening.
\end{construction}
In order to show that the $G$-tangent classifier sends a $G$-manifold $M$ to the $G$-map classifying its tangent bundle we will use the following description of the $G$-slice category $\ulTopG_{/\ul{B}}$.
\begin{rem} \label{ulTopG_over_B_description}
  A $G$-space $B$ defines a $G$-object $\ul{B} \colon  \OGop \to \ulTopG$ (i.e. a coCartesian section, see \cite[def. 7.1]{Expose1}).
  Using the explicit model of $\ulTopG$ given in \cref{ulTopG_as_G_spaces_over_orbits} we can describe $\ul{B}$ as 
  \[ \ul{B} \colon  \OGop \to \ulTopG, \, [G/H] \mapsto ( B \times G/H \to G/H). \]
  By \cite[lem]{AF} and \cref{ulTopG_as_G_spaces_over_orbits} it follows that the fibers of the parametrized slice category \( \ulTopG_{/\ul{B}} \) are given by 
  \[ 
    \left( \ulTopG_{/\ul{B}} \right)_{[G/H]} \simeq \left( \ulTopG_{[G/H]} \right)_{/  \ul{B}(G/H)} \simeq \left( \GTop_{/G/H} \right)_{/ ( B \times G/H \to G/H )} 
    \iso \GTop_{/B\times G/H}.
  \] 
  In particular an object of $\left( \ulTopG_{/\ul{B}} \right)_{[G/H]}$ is given by a $G$-space over $B \times G/H$, 
  which we consider as an object  $(Y \to G/H) \in \GTop_{/G/H} \simeq \ulTopG_{[G/H]}$, together with a $G$-map $f \colon Y\to B$.
  We write $\bar{f} \colon Y \to B \times G/H$ for the $G$-map corresponding to the pair $(Y \to G/H, Y \xto{f} B)$.

  The mapping spaces of the slice category \( \left( \ulTopG_{[G/H]} \right)_{/\ul{B}(G/H)} \simeq \GTop_{/B \times G/H} \) will be denoted by $\Map^G_{/\ul{B}(G/H)}(X,Y)$.
  An explicit description of these mapping spaces is given by the Moore over category, see \cref{MooreOverCat}.
\end{rem}

\begin{prop} \label{tangent_classifier_classifies_TM}
  Let $(M\to G/H)$ be an $\OG$-manifold, and consider the tangent bundle $TM \to M$ as a $G$-vector bundle.
  Then $\tau_M \in \left( \ulTopG_{/\ul{BO_n(G)}} \right)_{[G/H]}$ is given by $(M \to G/H) \in \ulTopG_{[G/H]}$ together with the $G$-map \( \tau_M \colon  M \to BO_n(G) \) classifying the tangent bundle of $M$.
\end{prop}
\begin{proof}
  Recall that an $\OG$-manifold $M\to G/H$ has an open cover by $G$-embeddings \( E_\alpha \cof M \) over $G/H$, where the patches $(E_\alpha \to G/H)$ are $G$-vector bundles. 
  The mapping space $\Map^G( M, BO_n(G) )$ is the homotopy limit of $\Map^G( E_\alpha, BO_n(G))$, 
  so by the functionality of $\tau$ in $M$ we are reduced to verifying the statement for $ E \to G/H$ a $G$-vector bundle.

  By construction the restriction of $\tau$ to $\ulRep^G_n$ is given by straightening the functor associated to the Yoneda embedding \( \ul\Rep^G_n \cof \ulFun_G( (\ulRep^G_n)^{vop}, \ulTopG ) \). 
  Recalling the construction of the parametrized Yoneda embedding (\cite[sec. 10]{Expose1}) we see that $\tau|_{\ulRep^G_n}$ is associated to the left fibration of the parametrized twisted arrow category \( \widetilde{\mathscr{O}}(\ulRep^G_n / \OGop) \fib (\ulRep^G_n)^{vop} \ultimes \ulRep^G_n \), end $\tau_{E}$ is associated to its pullback
  \begin{align*}
    \xymatrix{
      P_E \ar@{->>}[d] \ar[r] \pullbackcorner & \widetilde{\mathscr{O}}(\ulRep^G_n / \OGop) \ar@{->>}[d] \\
      (\ulRep^G_n)^{vop} \ultimes \ul{G/H} \ar[r]^{id \ultimes \ul{E}} & (\ulRep^G_n)^{vop} \ultimes \ulRep^G_n .
    }
  \end{align*}
  By \cref{RepGn_is_G_gpd} this is a pullback square of $G$-$\infty$-groupoids, and using \cref{RepGn_is_BOnG} we can identify it a homotopy pullback of $G$-spaces given by the top square of the following diagram
  \begin{align*}
    \xymatrix{
      G/H \ar[r]^{\simeq} \ar@/_2pc/[rdd]_{=} & \mathcal{P}_e \ar@{->>}[d] \ar[r] \pullbackcorner & \Map( \Delta^1 , BO_n(G) ) \ar@{->>}[d] \\
      & BO_n(G) \times G/H \ar[r]^{id \times e} \ar[d]^{proj} \pullbackcorner & BO_n(G) \times BO_n(G) \ar[d]^{proj} \\
       & G/H \ar[r]^{e} & BO_n(G) .
    }
  \end{align*}
  Since the bottom square (given by projections to the second coordinate) is a homotopy pullback square it follows that the outer rectangle is a homotopy limit diagram.
  Observe that the composition of the right vertical maps is an equivalence, and therefore the composition of the left vertical maps is an equivalence as well.
  It follows that $\tau_E$ is equivalent to the $G$-map $(e,id) \colon G/H \to BO_n(G) \times G/H)$, where $e \colon G/H \to BO_n(G)$ classifies the $G$-vector bundle $E\to G/H$.

  On the other hand the tangent bundle $TE$ is given by fiber product $TE \= E \times_{G/H} E$ and therefore classified by the composition of the bottom maps in 
  \begin{align*}
    \xymatrix{
      TE \ar[r] \ar[d] \pullbackcorner & E \ar[d] \\
      E \ar[r]^{\simeq} & G/H \ar[r]^-e & BO_n(G),
    }
  \end{align*}
  which is clearly equivalent to \( E \xto{\simeq} G/H \xto{(e,id)} BO_n(G) \times G/H \xto{proj} BO_n(G) \).
\end{proof}

\subsection{The $G$-category of $f$-framed $G$-manifolds} \label{sec:G_framed_mfld}
We now turn to the definition of the $G$-$\infty$-category of $G$-manifolds with additional tangential structure.
Our main interest is in the tangential structure defining $V$-framed $G$-manifolds, for $V$ a $G$-representation. 
However, the definition of equivariant framing on $G$-manifolds supports other interesting tangential structures, including equivariant orientations in the sense of \cite{CMW_equivariant_orientation_theory}, and free $G$-manifolds (an example not a priori associated with tangential structures). 

The specific type of $G$-tangential structure, such as equivariant framing or equivariant orientation, is specified by a $G$-space $B$ and a $G$-map $f \colon  B\to BO_n(G)$, in the following manner. 
An $f$-framing on a $G$-manifold $M$ is given by a $G$-map $M \to B$ such that the composition $M \to B \xto{f} BO_n(G)$ classifies the tangent bundle of $M$.
Similarly, if $H<G$ is a subgroup and $M$ is an $H$-manifold, we say that $M$ is $f$-framed its tangent bundle is classified by the $H$-map $M \to B \xto{f} BO_n(G)$.

The $\infty$-categories of $f$-framed $H$-manifold for $H<G$ can be arranged into an $\OGop$-diagram, encoded by a $G$-$\infty$-category $\ulmfld^{G,f-fr}$.
We start by giving a precise definition of $\ulmfld^{G,f-fr}$ and the $G$-functor $\ulmfld^{G,f-fr} \to \Gmfld$ that forgets the tangential structure.
\begin{mydef} \label{def:Gmfld_framed}
  Let $B\in \mathbf{Top}^G$ be a $G$-space and $f \colon B \to BO_n(G)$ be a $G$-map.
  Define the $G$-categories of $f$-framed $G$-manifolds as the pullback
  \begin{align*}
    \xymatrix{
      \ulmfld^{G,f-fr} \pullbackcorner \ar[d] \ar[r] & \ulTopG_{/\ul{B}} \ar[d]^{f_*} \\ 
      \Gmfld \ar[r]^-\tau & \ulTopG_{/\ul{BO_n(G)}} .
    }
  \end{align*}
\end{mydef}

\begin{rem} \label{rem:f_framed_mapping_spaces}
  Unwinding the definition, an object of \( (\ulmfld^{G, f-fr})_{[G/H]} \) is given by $(M\to G/H)\in \Gmfld_{[G/H]}$, a $G$-map \( f_M \colon  M  \) and a $G$-homotopy between $f \circ f_M$ exhibiting 
  \begin{align*}
    \xymatrix{
      & B \ar[d]^{f } \\
      M \ar[r]^-{\tau_M} \ar[ru]^{f_M} & BO_n(G) 
    }
  \end{align*}
  as homotopy coherent diagram of $G$-spaces.

  The mapping spaces of $\ulmfld^{G, f-fr}_{[G/H]}$ are given by homotopy pullbacks \\
  \begin{tikzcd}
    Emb^{G,f-fr}_{G/H}(M,N) \pbcorner \arrow[d] \arrow[r] & \Map^G_{/ B \times G/H }  (M \xto{\bar{f_M}} B \times G/H ,N \xto{\bar{f_N}} B \times G/H) \arrow{d}{(f\times G/H)_*} \\
    Emb^G_{G/H}(M,N) \arrow{r}{\tau} &  \Map^G_{/ BO_n(G) \times G/H } (M \xto{\bar{\tau_M}} BO_n(G) \times G/H, N \xto{\bar{\tau_N}} BO_n(G) \times G/H) . 
  \end{tikzcd}
\end{rem}

We finish this subsection with some examples of equivariant tangential structures on $G$-manifolds. 
We are primarily interested in equivariantly framed $G$-manifolds, which is our first example.
\begin{ex}[$V$-framed $G$-manifolds] \label{ex:V_fr_G_mfld}
  Let $B=pt$. 
  A $G$-map $f \colon pt \to BO_n(G)$ factors through the space of $G$-fixed points $(BO_n(G))^G = \coprod_V B \Aut_{\Rep^G_n} (V)$, so choosing $f$ is equivalent to choosing a connected component, i.e a real $n$-dimensional $G$-representation $V$. 
  A $V$-framing of an $H$-manifold $M$ is therefore a homotopy lift 
  \begin{align*}
    \xymatrix{
      & pt \ar[d]^{V} \\
      M \ar[r]^-{\tau_M} \ar[ru] & \BAut(V) ,
    }
  \end{align*}
  which under \cref{tangent_classifier_classifies_TM} and restriction to fibers over the coset $eH$ is equivalent to a choice of trivialization $ TM \= M \times V$ as an $H$-vector bundle. 
\end{ex}
\begin{ex}[$G$-manifolds with no tangential structure]
  Apply \cref{def:Gmfld_framed} for the $G$-space $B=BO_n(G)$ and $id \colon BO_n(G) \to BO_n(G)$ constructs $\ulmfld^{G,id-fr} \= \Gmfld$.
\end{ex}
\begin{ex}[$G$-orientated $G$-manifolds]
  Orientations of $G$-vector bundles were studied by Costenoble, May and Waner in \cite{CMW_equivariant_orientation_theory}\footnote{see \cite[def. 2.8]{CMW_equivariant_orientation_theory} for a precise definition},
  and used in \cite{CostenobleWaner} to prove equivariant versions of Poincar\'e duality.

  Let us recall the relevant results from \cite{CMW_equivariant_orientation_theory}.
  First, there exists a universal oriented $G$-$n$-plane bundle, given by a $G$-map $ EO_n(G,S) \to BO_n(G,S)$, see \cite[thm. 22.4]{CMW_equivariant_orientation_theory}.
  Second, there is a $G$-map $f \colon  BO_n(G,S) \to BO_n(G)$ representing the forgetful functor from oriented $n$-plane bundles to $G$-$n$-plane bundles.
  Therefore an orientation on a $G$-vector bundle is given by a $G$-homotopy lift of its classifying map along the $G$-map $f$.
  
  Applying \cref{def:Gmfld_framed} to $B=BO_n(G,S)$ and $f \colon  BO_n(G,S) \to BO_n(G)$ we get a $G$-$\infty$-category $\ulmfld^{G,or}$ of oriented $G$-manifolds. 
\end{ex}
\begin{rem}
  The notion of an oriented $G$-manifold seems not to agree with the notion of oriented global orbifold (see, for example, \cite[p. 34]{OrbifoldsStringyTop}).
\end{rem}
Finally, we can use equivariant tangential structures to restrict the class of $G$-manifolds we consider, an idea introduced in \cite[rem. 1.1.9]{AFT_stratified_local_structures}.
\begin{ex}
  Applying \cref{def:Gmfld_framed} with $B=BO_n(G) \times EG$ and a $G$-map given by the projection $pr \colon  BO_n(G) \times EG \to BO_n(G)$ produces a $G$-$\infty$-category $\ulmfld^{G,pr-fr}$.
  In this example the forgetful $G$-functor $\ulmfld^{G,pr-fr}\to \Gmfld$ is fully faithful, and exhibits $\ulmfld^{G,pr-fr}$ as the full $G$-subcategory of $\Gmfld$ spanned by $\OG$-manifolds $M\to O$ where $M$ is a free $G$-manifold. 
  We now give a quick sketch the argument.
  
  We consider a manifold $M$ with an action of $G$, describing an object $(M\to G/G) \in \Gmfld_{[G/G]}$ (the argument for an $\OG$-manifold $M \to G/H$ is similar). 
  A homotopy lift of $\tau_M \colon  M \to BO_n(G)$ along the projection the same as a $G$-map $M \to EG$.
  A $G$-map $M\to EG$ exists if and only if the action of $G$ on $M$ is free, in which case the space of $G$-maps $\Map_G(M, EG)$ is contractible.
  This is easily seen by using the Elmendorf-McClure theorem; the presheaves that represents $M$ and $EG$ send
  \begin{align*}
    M,EG \colon  \OGop \to \Ss, \quad M \colon  G/H \mapsto M^H, \, EG \colon  G/H \mapsto \begin{cases} pt, \quad & H=e , \\ \emptyset \quad & H\neq e, \end{cases} 
  \end{align*}
  and a map $M^H \to \emptyset$ exists if and only if $M^H$ is empty. 
  It follows that a map of $\OG$-presheaves $M\to EG$ exists if and only if the action of $G$ on $M$ is free.
  Finally, if $G$ acts freely on $M$ then $\Map_{\Fun(\OGop,\Ss)}(M,EG) \simeq \Map(M, pt) \simeq pt$.
\end{ex}

\subsection{$G$-disjoint union of $G$-manifolds} \label{sec:G_disj_union}
The goal of this section is to endow the $G$-$\infty$-category $\Gmfld$ with a $G$-symmetric monoidal structure associated to disjoint unions.

Recall that $n$-dimensional manifolds with $G$-action and $G$-equivariant embedding can be organized into a topological category $\mfld^G$.
Despite the fact that $\mfld^G$ does not have coproducts\footnote{ Note that $Emb^G(M_1 \sqcup M_2, M) \not \simeq Emb^G(M_1,M) \times Emb^G(M_2,M)$.}, 
we can still endow $\mfld^G$ with a symmetric monoidal structure by taking disjoint unions.
Therefore the $\infty$-category $\NSing(\mfld^G)$ admits a symmetric monoidal structure $\NSing^\otimes(\mfld^G) \fib \Fin_*$, given by applying the \emph{operadic nerve} construction of \cite[def. 2.1.1.23]{HA}.

Similarly, disjoint unions endow the $\infty$-category $\NSing(\mfld^H)$ with a symmetric monoidal structure, making the restriction and conjugation functors symmetric monoidal.
We can therefore enhance $\mfld^\bullet$ from a diagram of $\infty$-categories to a diagram of symmetric monoidal $\infty$-categories $\NSing^\otimes(\mfld^\bullet)$. 
However, this construction does not encode the operation of topological induction and its coherent compatibility with the symmetric monoidal structure and restriction and conjugation of the action.
The main point of this subsection is that all of the structure we are interested in can be encoded as a $G$-symmetric monoidal structure on the $G$-category of $G$-manifolds (see \cref{def:Gmfld}).
It would be preferable to define this $G$-symmetric monoidal structure by an appropriate variant of the operadic nerve construction, however we are unaware of such construction. 
We therefore define the $G$-symmetric monoidal structure by explicitly constructing a coCartesian fibration $\GmfldD \fib \GFin_*$ (see \cref{def:GmfldD_as_coCart_fib}).

Our construction can be briefly described as follows. 
The category $\GFin_*$ is constructed as a category of spans in the category of finite $G$-sets over an orbit, $\OGFin$, (see \cref{lem:S_adeq_triple}), so it is natural to construct $\GmfldD$ as category of spans of an auxiliary $\infty$-category $\mfldGD$, defined over $\OGFin$.
In \cref{def:mfldGD_as_top_cat} we construct $\mfldGD$ as a topological category over $\OGFin$.
We want to apply Barwick's unfurling construction, see \cite{BarwickSMF1}, to the functor $\NSing(\mfldGD)\to \OGFin$, in order to produce a coCartesian fibration $\GmfldD \fib \GFin_*$ between the respected $\infty$-categories of spans. 
There is a simple criterion, described in \cite{BarwickSMF1}, that ensures that the unfurled functor is a coCartesian fibration: 
\begin{enumerate}
  \item Egressive arrows in $\OGFin$, serving as the ``wrong way arrows'' in the span category $\GFin_*$, have Cartesian lifts (verified in \cref{lem:Cart_PB}).
  \item Ingressive arrows in $\OGFin$, serving as the ``right way arrows'' in the span category $\GFin_*$, have coCartesian lifts (verified in \cref{lem:coCart_OrbitIso}).
  \item The pullback squares appearing in the definition of composition in the span category $\GFin_*$ satisfies a ``Beck-Chevalley condition'' (verified in \cref{prop:adeq_fib}).
\end{enumerate}
The resulting ``unfurled'' $\infty$-category $\GmfldD$ (see \cref{def:GmfldD_as_coCart_fib}) admits an explicit description as an $\infty$-category of spans.
In particular we have a description of the objects, morphisms and coCartesian morphisms of $\GmfldD$.
Using the explicit description of $\GmfldD \fib \GFin_*$ we show that it satisfies the $G$-Segal conditions and that its underlying $G$-$\infty$-category is $\Gmfld$ (\cref{def:GmfldD}).

\subsubsection*{Construction of the auxiliary category $\mfldGD$}
In this subsection we define a topological category $\mfldGD$ with a functor to the category $\OGFin$ of finite $G$-sets over orbits. 
The topological category $\mfldGD$ serves as input to the unfurling construction (\cite[sec. 11]{BarwickSMF1}), producing a coCartesian fibration \(\GmfldD \fib \GFin_* \) that defines the $G$-symmetric monoidal structure of $G$-disjoint union on $\Gmfld$ (see \cref{def:GmfldD_as_coCart_fib}).

We start with a definition of the category $\OGFin$, which serves as the base category of the unfurling construction. 
\begin{mydef} 
  The category $\OGFin$ is the pullback 
  \( \OGFin:=\Fun(\Delta^1,Fin^G) \times_{\Fun(\set{1},Fin^G)} \OG .  \)
  The category $\OGFin$ is a full subcategory of the arrow category $\Fun(\Delta^1,Fin^G)$, whose objects are arrows $U \to O$ in $Fin^G$ such that $O\in \OG$. 
  A morphism in $\OGFin$ is a \emph{summand-inclusion}  (\cite[def. 4.12]{Expose4}) if it factors as an inclusion over orbit-identity followed by a pullback square
  \begin{align} \label{eq:summand_inclusion}
    \diag{
      U_1 \ar[d] \ar@{^(->}[r] & \varphi^* U_2 \ar[d] \ar[r] \pullbackcorner & U_2 \ar[d] \\
      O_1 \ar[r]^{=} & O_1 \ar[r]^{\varphi} & O_2 
    } .
  \end{align}
  Note that we the inclusion of $G$-sets \( U_1 \cof \varphi^* U_2 \) exhibits $\varphi^* U_2$ as the coproduct of the $G$-sets $U_1$ and $U' = \varphi^* U_2 \setminus U_1$.
  We can therefore identify $\varphi^* U_2 \= U_1 \coprod U'$.

  Let \( \OGFin^{\dagger} \subset \OGFin \) be the subcategory consisting of all objects while morphisms are summand-inclusions.
\end{mydef}
It is straightforward to see that the $G$-category $\GFin_*$ of \cite[def. 4.12]{Expose4} can be defined by the following unfurling construction.
\begin{lem} \label{lem:S_adeq_triple}
  The triple $(\OGFin,\OGFin\times_{\OG} \OG^{\=}, \OGFin^\dagger)$ is an adequate triple in the sense of \cite[def. 5.2]{BarwickSMF1},
  and its effective Burnside category fits into a pullback square
  \begin{align} \label{eq:GFin_into_Aeff}
    \diag{ 
      \GFin_* \ar[d] \ar@{^(->}[r] \pullbackcorner & A^{eff}(\OGFin,\OGFin\times_{\OG} \OG^{\=}, \OGFin^\dagger) \ar[d] \\
      \OG \ar@{^(->}[r] & A^{eff}(\OG,\OG^{\=}, \OG)
    }
  \end{align}
\end{lem}

We now define a topological category of ``parametrized $\OG$-manifolds'' over $\OGFin$.
\begin{mydef}
  An \emph{$\OGFin$-manifold} $M \to U \to O$ is 
  \begin{enumerate}
    \item a smooth $n$-dimensional manifold $M$ with an action of $G$ on $M$ by smooth maps, 
    \item together with a $G$-map $M\to U$ from the underlying $G$-space of the manifold $M$ to a $G$-finite set $U\in Fin^G$,
    \item and an arrow $U \to O$ in $\Fin^G$ such that $O\in \OG$. 
  \end{enumerate}
  An morphism of $\OGFin$-manifolds is given by a commuting square of $G$-spaces
  \begin{align*}
    \xymatrix{
      M_1 \ar[d] \ar[r]^{f} & M_2 \ar[d] \\
      U_1 \ar[d] \ar[r]^{\overline{\varphi}} & U_2 \ar[d] \\
      O_1 \ar[r]^{\varphi} & O_2 ,
    }
  \end{align*}
  such that the induced map $M_1 \to O_1 \times_{O_2} M_2$ is an embedding.  
\end{mydef}

\begin{mydef} \label{def:EmbG}
  Let $M_1 \to U_1 \to O_1,\, M_2\to U_2 \to O_2$ be $\OGFin$-manifolds and $ \varphi \colon   I_1 \to I_2$ a morphism in $\OGFin$ given by 
  \( \vcenter{
    \xymatrix{
      U_1 \ar[d] \ar[r]^{\overline{\varphi}} & U_2 \ar[d] \\
      O_1 \ar[r]^{\varphi} & O_2 
    }
  } \). 
  Define $Emb^{\OGFin}_{\varphi}(M_1,M_2)\subset C^\infty(M_1,M_2)$ as the subspace of smooth maps $f \colon M_1 \to M_2$ such that $(f,\overline{\varphi},\varphi)$ is a morphism of $\OGFin$-manifolds from $M_1 \to U_1 \to O_1$ to $M_2 \to U_2 \to O_2$.
\end{mydef}

\begin{mydef} \label{def:mfldGD_as_top_cat}
  The \emph{Category of $\OGFin$-manifolds $\mfldGD$} is the topological category whose objects are $\OGFin$-manifolds. 
  The morphism space from $M_1$ to $M_2$ is given by 
  \begin{align*}
    \Map_{\mfldGD}(M_1,M_2) :=  \coprod_{\varphi} Emb^{\OGFin}_\varphi(M_1,M_2),
  \end{align*}
  where the coproduct is indexed by  $\Hom_{\OGFin}(U_1 \to O_1,U_2 \to O_2)$.

  Define a forgetful functor \( p \colon \mfldGD \to \OGFin \) by sending $M\to U\to O$ to $U\to O$, and the subspace $Emb^{\OGFin}_\varphi(M_1,M_2)\subset Map(M_1,M_2)$ to $\varphi \in \Hom_{\OGFin}(U_1 \to O_1,U_2 \to O_2)$.
\end{mydef}
From here on we will abuse notation, writing $\mfldGD$ for both the topological category $\mfldGD$, its incarnation as a fibrant simplicial category $\Sing(\mfldGD)$ and its incarnation as an $\infty$-category $\NSing(\mfldGD)$, distinguishing between these incarnations by context. 

\begin{rem}
  Note that an equivalence $f \colon M \to N$ in $\mfldGD$ is always an embedding of smooth manifolds, since it lies over an isomorphism of orbits. 
  Moreover, it is $G$-isotopic to an identity-of-manifolds over the isomorphism $p(f)$.
  On the other hand, if $f$ is $G$-isotopic to an identity-of-manifolds over an isomorphism of finite $G$-sets then $f$ is an equivalence in $\mfldGD$, so we have a complete characterization of equivalences in $\mfldGD$.
\end{rem}

\subsubsection*{Some Cartesian and coCartesian edges of $\mfldGD\to \OGFin$}
We characterize $p$-Cartesian edges of $\mfldGD$ over summand-inclusions and $p$-coCartesian edges over isomorphisms of orbits. 
We summarize the results of this subsection as follows.
\begin{prop} \label{prop:Cart_coCart_edges}
  A morphism $f$ of $\mfldGD$ over $\OGFin^\dagger$ is $p$-Cartesian if and only if it is equivalent to a pullback over a summand-inclusion.
  A morphism $g$ of $\mfldGD$ over $\OGFin\times_{\OG} \OG^{\=}$ is $p$-coCartesian if and only if it is $G$-isotopic to an identity-of-manifolds over an orbit-isomorphism.
\end{prop}
The characterization of $p$-Cartesian edges is given in \cref{cor:Cart_SI_equiv}, and the characterization of $p$-coCartesian edges is given in \cref{cor:coCart_overIso_equiv}.

\begin{rem}
  By \cite[prop. 2.4.1.10(1)]{HTT} the map $\mfldGD \to \OGFin$ is an inner fibration.
\end{rem}
\begin{lem} \label{lem:Cart_PB}
  Let \(\varphi\in\Hom_{\OGFin}(U_1\to O_1, U_2\to O_2)\) be a morphism in $\OGFin$ given by a pullback square \( \diag{ U_1 \ar[d] \ar[r] \pullbackcorner & U_2 \ar[d] \\ O_1 \ar[r] & O_2} \),
  and \(N \to U_2 \to O_2 \) a $\OGFin$-manifold over its target.
  Then the pullback 
  \begin{align*}
    \diag{
      M \ar[d] \ar[r]^{f} \pullbackcorner & N \ar[d] \\ 
      U_1 \ar[d] \ar[r] \pullbackcorner & U_2 \ar[d] \\ 
      O_1 \ar[r] & O_2
    }
  \end{align*}
  defines a $p$-Cartesian morphism $f$ in $\mfldGD$ lifting $\varphi$.
\end{lem}
\begin{proof}
  According to \cite[prop. 2.4.1.10(2)]{HTT} we have to show that for every $\OGFin$-manifold \(T \to U \to O \) the commutative square
  \begin{align*}
    \diag{
      \Map(T,M) \ar@{->>}[d] \ar[r]^{f_*} & \Map(T,N) \ar@{->>}[d] \\
      \Hom_{\OGFin}(p(T),p(M)) \ar[r]^{p(f)_*} & \Hom_{\OGFin}(p(T),p(N) 
    }
  \end{align*}
  is a homotopy pullback.
  Since the vertical maps are Kan fibrations, this square is a homotopy pullback if and only if $f_*$ induces an equivalence between the fibers over every vertex of the base $\Hom_{\OGFin}(p(T),p(M))$.

  Let $\tau\in\Hom_{\OGFin}(p(T),p(M))$.
  The functor $f_*$ induces a map of the fibers over $\tau$ 
  \[ (f_*)|_{\tau} \colon Emb^{\OGFin}_{\tau}(T,M) \to \set{\tau} \times_{\Hom_{\OGFin}(p(T),p(M)} \Map(T,N) .\]
  Unwinding the definition of the mapping space in $\mfldGD$, we have 
  \begin{align*}
    \set{\tau} \times_{\Hom_{\OGFin}(p(T),p(M)} \Map(T,N) & = \set{\tau} \times_{\Hom_{\OGFin}(p(T),p(M)} \left( \coprod_{\varphi} Emb^{\OGFin}_\varphi(T,N) \right) \\
      & = Emb^{\OGFin}_{p(f)\circ\tau} (T,N) ,
  \end{align*}
  where the last equality holds since pullback along a fixed map preserve coproducts. 

  Suppose that the $\OGFin$-manifold $T$ is given by \( T \to U \to O \) and $\tau \colon p(T)\to p(N)$ is given by the square
  \( \diag{ U \ar[d] \ar[r] & U_1 \ar[d] \\ O \ar[r] & O_1} \).
  Then the map \( (f_*)|_{\tau} \colon Emb^{\OGFin}_{\tau}(T,M) \to Emb^{\OGFin}_{p(f)\circ\tau} (T,N) \) sends $h \colon  T \to M$ to $f\circ h$: 
  \begin{align*}
    (f_*)|_{\tau} \colon 
    h=\left( \diag{ 
      T \ar[d] \ar[r]^{h} & M \ar[d] \\
      U \ar[d] \ar[r] & U_1 \ar[d] \\
      O \ar[r] & O_1
    } \right) 
    \mapsto
    \left( \diag{ 
      T \ar[d] \ar[r]^{h} & M \ar[d] \ar[r]^{f} \pullbackcorner & N \ar[d] \\
      U \ar[d] \ar[r] & U_1 \ar[d] \ar[r] \pullbackcorner & U_2 \ar[d] \\
      O \ar[r] & O_1 \ar[r] & O_2
    } \right).
  \end{align*}
  The universal property of the pullback $M = N \times_{U_2} U_1$ shows that $(f_*)|_{\tau}$ is a continuous bijection: 
  injectivity follows from uniqueness of maps to the pullback. 
  Surjectivity: suppose \( g\in Emb^{\OGFin}_{p(f)\circ \tau} (T,N) \), by existence of a map to the pullback we have a candidate map \(h \colon T\to M \) over $\tau$ such that \( g=f\circ g\). 
  We have to show that $h\in Emb^{\OGFin}_{\tau}$.
  Clearly $h$ is a smooth $G$-map, so we only have to verify condition (3) of \cref{def:EmbG}: $h$ induces an embedding \( T \to O\times_{O_1} M \).
  To see that observe that $g$ induces an embedding \(T \cof O \times_{O_2} N \) which factors as the map induced by $h$ followed by the isomorphism \(O \times_{O_1} M = O \times_{O_1} (O_1 \times_{O_2} N) \= O \times_{O_2} N \).

  We leave it as an exercise to the reader to verify that $(f_*)|_{\tau}$ is an open map, and therefore a homeomorphism.
\end{proof}

Note that every $G$-map \( M \to U_1 \coprod U_2 \) from a manifold with $G$-action to a coproduct of $G$-sets factors as coproduct of $G$-maps \( M = M_1 \coprod M_2 \to U_1 \coprod U_2 \). 

\begin{lem} \label{lem:Cart_Inclusion}
  Let $\varphi$ be an inclusion of finite $G$-sets over $id_{O}$ in ${\OGFin}$, given by the diagram 
  \( \diag{ U_1 \ar[d] \ar@{^(->}[r] & U_1 \coprod U_2 \ar[d] \\ O_1 \ar[r]^{=} & O_1 } \),
  and \( M_1 \coprod M_2 \to U_1 \coprod U_2 \to O_2 \) a $\OGFin$-manifold over its target.
  Then the pullback 
  \begin{align*}
    \diag{
      M_1 \ar[d] \ar@{^(->}[r]^-{i} \pullbackcorner & M_1 \coprod M_2 \ar[d] \\ 
      U_1 \ar[d] \ar@{^(->}[r] & U_1 \coprod U_2 \ar[d] \\ 
      O_1 \ar[r]^{=} & O_1
    }
  \end{align*}
  defines a $p$-Cartesian morphism $i$ in $\mfldGD$ lifting $\varphi$.
\end{lem}

\begin{proof}
  As in \cref{lem:Cart_PB}, we have to show that for every $\OGFin$-manifold \( T \to U \to O \) and every $\tau \colon p(T) \to p(M_1)$ the map $i_*$ induces equivalence of the fibers
  \[ Emb^{\OGFin}_{\tau}(T,M_1) \to Emb^{\OGFin}_{p(i)\circ\tau}(T,M_1\coprod M_2). \]
  As above, we use the universal property of the pullback to show this map is a bijection, and leave it to the reader to verify it is an open map.

  The only part which is different is the verification of condition (3) of \cref{def:EmbG}: $g$ induces an embedding  $T \cof O\times_{O_1} (M_1 \coprod M_2)$, which factors as the composition of the map induced by $h$, an inclusion and an isomorphism
  \[ T \to O\times_{O_1} M_1 \cof O\times_{O_1} M_1 \coprod O\times_{O_1} M_2 \= O\times_{O_1} (M_1 \coprod M_2) . \]
  Since the composition is an embedding, the map $T \to O\times_{O_1} M_1$ induced by $h$ is an embedding. 
\end{proof}
Together, the lemmas above show the existence of $p$-Cartesian lifts over summand-inclusions and characterizes them.
\begin{cor} \label{cor:Cart_SI}
  Let \(\varphi\in\Hom_{\OGFin}(U_1\to O_1, U_2\to O_2)\) be a morphism in ${\OGFin}^\dagger$  
  and \(N \to U_2 \to O_2 \) an $\OGFin$-manifold over its target.
  Then the pullback 
  \begin{align*}
    \diag{
      M \ar[d] \ar[r]^{f} \pullbackcorner & N \ar[d] \\ 
      U_1 \ar[d] \ar[r] & U_2 \ar[d] \\ 
      O_1 \ar[r] & O_2
    }
  \end{align*}
  defines a $p$-Cartesian morphism $f$ in $\mfldGD$ lifting $\varphi$.
\end{cor}
\begin{proof}
  Factor the summand-inclusion $\varphi$ as in \eqref{eq:summand_inclusion}, apply \cref{lem:Cart_PB} and \cref{lem:Cart_Inclusion}.
\end{proof}
By \cite[prop. 2.4.1.7 and 2.4.1.5]{HTT}, we have
\begin{cor} \label{cor:Cart_SI_equiv}
  A morphism $f$ of $\mfldGD$ over ${\OGFin}^\dagger$ is $p$-Cartesian if and only if it is equivalent to a pullback over a summand-inclusion, i.e the left map in the factorization 
  \begin{align*}
    f= \left( \diag{
      M_1 \ar[d] \ar@{^(->}[r] & M_2 \times_{U_2} U_1 \ar[d] \ar[r] \pullbackcorner  & M_2 \ar[d] \\
      U_1 \ar[d] \ar[r]^{=} & U_1  \ar[d] \ar[r] & U_2 \ar[d] \\
      O_1 \ar[r]^{=} & O_1  \ar[r] & O_2
    } \right)
  \end{align*}
  is an equivalence in $\mfldGD$ (a $G$-isotopy equivalence over $U_1$). 
\end{cor}
Next, we construct $p$-coCartsian lifts over isomorphism of orbits.
\begin{lem} \label{lem:coCart_OrbitIso}
  Let \( \varphi = \left( \diag{ U_1 \ar[d] \ar[r] & U_2 \ar[d] \\ O_1 \ar[r]^{\=} & O_2 } \right) \) be a morphism of ${\OGFin}\times_{\OG} \OG^{\=}$ and $M \to U_1 \to O_1$ an $\OGFin$-manifold.
  Then 
  \( 
    f= \left( \diag{
          M \ar[d] \ar[r]^{=} & M \ar[d] \\
          U_1 \ar[d] \ar[r] & U_2 \ar[d] \\ 
          O_1 \ar[r]^{\=} & O_2
        } \right) 
  \)
  is a $p$-coCartesian lift of $\varphi$. 
\end{lem}
\begin{proof}
  By the dual version of \cite[prop. 2.4.1.10(2)]{HTT} we have to show that for every $\OGFin$-manifold $T \to U \to O$ the square
  \begin{align*}
    \diag{
      Map(M \to U_2 \to O_1, T \to U \to O) \ar@{->>}[d] \ar[r]^{f^*} & Map( M \to U_1 \to O_1, T \to U \to O ) \ar@{->>}[d] \\
      \Hom_{\OGFin}( U_2 \to O_1, U \to O) \ar[r]^{p(f)^*} & \Hom_{\OGFin}(U_1 \to O_1,  U \to O) 
    }
  \end{align*}
  is a homotopy pullback square.
  Since the vertical maps are Kan fibrations, this square is a homotopy pullback if and only if $f^*$ induces an equivalence between the fibers.
  Next, note that the map $f^*$ is induced by composition with $id_M$, and the fibers over  
  $\tau \in \Hom_{\OGFin} \left( U_2 \to O_1, U \to O \right)$ and $ \tau \circ p(f) \in \Hom_{\OGFin} \left( U_1 \to O_1, U \to O \right)$ are both subspaces of the space of smooth maps $C^\infty(M,T)$:
  \[
    Emb^{\OGFin}_{\tau}(M,T) \subset C^\infty(M,T) , \quad Emb^{\OGFin}_{\tau\circ p(f)}(M,T) \subset C^\infty(M,T) .
  \] 
  We finish the proof by observing that these subspaces are equal: conditions (1),(3) of \cref{def:EmbG} coincide, while the equivalence of condition (2) follows from the commutativity of the square \( \diag{ M \ar[d] \ar[r]^{=} & M \ar[d] \\ U_1 \ar[r] & U_2 } \), the top square of $f$.
\end{proof}
We therefore have a characterisation of $p$-coCartesian edges over orbit isomorphisms.
\begin{cor} \label{cor:coCart_overIso_equiv}
  A morphism $f$ of $\mfldGD$ over an orbit-isomorphism is $p$-Cartesian if and only if it is equivalent to an identity-of-manifolds, 
  i.e. the right map in the factorization 
  \begin{align*}
    f= \left( \diag{
      M_1 \ar[d] \ar[r]^{=} & M_1 \ar[d] \ar@{^(->}[r] & M_2 \ar[d] \\
      U_1 \ar[d] \ar[r] & U_2 \ar[d] \ar[r]^{=} & U_2 \ar[d] \\ 
      O_1 \ar[r]^{\=} & O_2 \ar[r]^{=} & O_2
    } \right) 
  \end{align*}
  is an equivalence in $\mfldGD$ (a $G$-isotopy equivalence over $U_2$). 
\end{cor}

\subsubsection*{Construction of the $G$-symmetric monoidal category $\GmfldD$}
We now turn to the goal of this subsection, the construction of a $G$-symmetric monoidal structure on the $G$-category of $G$-manifolds.
In \cref{def:GmfldD_as_coCart_fib} we use the unfurling construction of \cite[sect. 11]{BarwickSMF1} to define a coCartesian fibration \( \GmfldD \fib \GFin_* \),
and in \cref{def:GmfldD} we verify the Segal conditions, showing that it defines a $G$-symmetric monoidal structure on $\Gmfld$.

We first make sure that the conditions for applying Barwick's unfurling construction hold.
Since Cartesian lifts of egressive morphisms and coCartesian lifts of ingressive morphisms were constructed in \cref{prop:Cart_coCart_edges} it remains to verify the Beck-Chevalley conditions.
\begin{prop} \label{prop:adeq_fib}
  The inner fibration $\mfldGD \to \OGFin$ is adequate over the triple $(\OGFin,\OGFin\times_{\OG} \OG^{\=}, \OGFin^\dagger)$ (\cite[def. 10.3]{BarwickSMF1}).
\end{prop}
\begin{proof}
  Conditions \cite[cond. (10.3.1),(10.3.2)]{BarwickSMF1} follow from \cref{prop:Cart_coCart_edges}.
  To verify condition \cite[cond. (10.3.3)]{BarwickSMF1} construct the natural map \( i_! \circ q^* (\tilde{N}) \to q'^* \circ j_! (\tilde{N}) \) by choosing appropriate $p$-Cartesian and $p$-coCartesian lifts, and show that map is the universal map between two models of the same pullback, hence a diffeomorphism over an identity map.

  Let 
  \( \diag{ 
    *+[r]{s \,\,} \ar@{->>}[d]_{q} \ar@{>->}[r]^{i} \pullbackcorner & s' \ar@{->>}[d]_{q'} \\
    *+[r]{t \,\,} \ar@{>->}[r]^{j} & t'
  } \) 
  be an ambigressive pullback square in $\OGFin$, whose objects and morphisms are given by
  \begin{align*}
  & s  = \left( \diag{ \tilde{U} \ar[d] \\ \tilde{O}_1 } \right),\quad 
    s' = \left( \diag{ U \ar[d] \\ O_1 } \right),\quad 
    t  = \left( \diag{ \tilde{V} \ar[d] \\ \tilde{O}_2 } \right),\quad 
    t' = \left( \diag{ V \ar[d] \\ O_2 } \right), \quad \quad \quad
    i  = \left(\diag{ \tilde{U} \ar[d] \ar[r] & U \ar[d] \\ \tilde{O}_1 \ar[r]^{\=} & O_1 }\right), \\
  & j  = \left(\diag{ \tilde{V} \ar[d] \ar[r] & V \ar[d] \\ \tilde{O}_2 \ar[r]^{\=} & O_2 }\right), \quad
    q  = \left(\diag{ \tilde{U} \ar[d] \ar[r] & \tilde{V} \ar[d] \\ \tilde{O}_1 \ar[r] & \tilde{O}_2 }\right), \quad
    q' = \left(\diag{ U \ar[d] \ar[r] & V \ar[d] \\ O_1 \ar[r] & O_2 }\right)
  \end{align*}
  And \( \tilde{N} = ( \tilde{N} \to \tilde{V} \to \tilde{O}_2 ) \) an object in the fiber of $p$ over $t$.
  We compute \( i_! \circ q^* (\tilde{N}) ,\, q'^* \circ j_! (\tilde{N}) \) and the  map \( i_! \circ q^* (\tilde{N}) \to q'^* \circ j_! (\tilde{N}) \) (natural in $\tilde{N}$) by choosing appropriate $p$-Cartesian and $p$-coCartesian lifts. 
  
  Let \( \tilde{M} := \tilde{N} \times_{\tilde{N}} \tilde{U} \). 
  Since $q$ is a summand-inclusion by \cref{cor:Cart_SI} the map
  \begin{align*}
    \diag{
      \tilde{M} \ar[d] \ar[r] \pullbackcorner & \tilde{N} \ar[d] \\
      \tilde{U} \ar[d] \ar[r] & \tilde{V} \ar[d] \\
      \tilde{O}_1 \ar[r] & \tilde{O}_2 
    }
  \end{align*}
  is $p$-Cartesian over $q$ , so \( q^*(\tilde{N}) := (\tilde{M} \to \tilde{U} \to \tilde{O}_1 ) \). 

  Since $i$ is over an isomorphism of orbits, by \cref{lem:coCart_OrbitIso} the map
  \begin{align*}
    \diag{
      \tilde{M} \ar[d] \ar[r]^{=} & \tilde{M} \ar[d] \\
      \tilde{U} \ar[d] \ar[r] & U \ar[d] \\
      \tilde{O}_1 \ar[r]^{\=} & O_1 
    }
  \end{align*}
  is $p$-coCartesian over $i$, so \( i_! \circ q^*(\tilde{N}) := (\tilde{M} \to U \to O_1 ) \).

  Since $j$ is over an isomorphism of orbits, by \cref{lem:coCart_OrbitIso} the map
  \begin{align*}
    \diag{
      \tilde{N} \ar[d] \ar[r]^{=} & \tilde{N} \ar[d] \\
      \tilde{V} \ar[d] \ar[r] & V \ar[d] \\
      \tilde{O}_2 \ar[r]^{\=} & O_2
    }
  \end{align*}
  is $p$-coCartesian over $j$, so \( j_!(\tilde{N}) := (\tilde{N} \to V \to O_2 ) \).

  Let \( M:= \tilde{N} \times_{V} U \).
  Since $q'$ is a summand-inclusion by \cref{cor:Cart_SI} the map
  \begin{align*}
    \diag{
      M \ar[d] \ar[r] \pullbackcorner & \tilde{N} \ar[d] \\
      U \ar[d] \ar[r] & V \ar[d] \\
      O_1 \ar[r] & O_2 
    }
  \end{align*}
  is $p$-Cartesian over $q'$, so \( q'^* \circ j_! (\tilde{N}) := (M \to U \to O_1) \). 
  
  Next, we choose a map \( \xi \colon  q^*(\tilde{N}) \to q'^*\circ j_!(\tilde{N}) \) over $i$ by composing the lifts of $q$ and $j$ above and using the universal property of the pullback $M$
  \begin{align*}
    \diag{
      \tilde{M} \ar[d] \ar[r] \pullbackcorner & \tilde{N} \ar[d] \ar[r]^{=} & \tilde{N} \ar[d] \\
      \tilde{U} \ar[d] \ar[r] & \tilde{V} \ar[d] \ar[r] & V \ar[d] \\
      \tilde{O}_1 \ar[r] & \tilde{O}_2 \ar[r]^{\=} & O_2
    } \quad \Rightarrow \quad 
    \diag{
      \tilde{M} \ar[d] \ar@{-->}[r]^{\exists !\xi} & M \ar[d] \ar[r] \pullbackcorner & \tilde{N} \ar[d] \\
      \tilde{U} \ar[d] \ar[r] &  U\ar[d] \ar[r] & V \ar[d] \\
      \tilde{O}_1 \ar[r]^{\=} & O_1 \ar[r] & O_2
    }.
  \end{align*}
  The map $\xi$ induces the natural map \( \overline{\xi} \colon  i_! \circ q^*(\tilde{N}) \to q'^*\circ j_!(\tilde{N}) \) over $id_{s'}$ by  
  \begin{align*}
    \diag{
      \tilde{M} \ar[d] \ar[r]^{=} & \tilde{M} \ar[d] \ar@{-->}[r]^{\exists !\xi} & M \ar[d] \\
      \tilde{U} \ar[d] \ar[r] & U \ar[d] \ar[r]^{=} & U \ar[d] \\
      \tilde{O}_1 \ar[r]^{\=} & O_1 \ar[r]^{=} & O_1
    }.
  \end{align*}
  In order to verify \cite[cond. (10.3.3)]{BarwickSMF1} we have to show that $\overline{\xi}$ is an equivalence in the fiber over $s'$.
  We show that $\xi$ is a diffeomorphism.
  Consider the diagram
  \begin{align*}
    \diag{
      \tilde{M} \ar[d] \ar[r] \pullbackcorner & \tilde{N} \ar[d]  \\
      \tilde{U} \ar[d] \ar[r] \pullbackcorner & \tilde{V} \ar[d] \\
      U \ar[r] & V
    }
  \end{align*}
  the top square is a pullback square by definition of $\tilde{M}$, and the bottom square is a pullback square by assumption.
  Therefore the outer rectangle is a pullback square. 
  By the universal property of \( M = \tilde{N} \times_{U} V \) the induced map $\xi$ is a diffeomorphism, as claimed.

  This ends the proof of \cref{prop:adeq_fib}.
\end{proof}
We can now define $\GmfldD$ by applying the unfurling construction to $\mfldGD \to \OGFin$.
\begin{mydef}
  Define a subcategory \( (\mfldGD)^\dagger \subset \mfldGD \) with the same objects as $\mfldGD$, and with morphisms the $p$-Cartesian edges over summand-inclusions (i.e over edges over $\OGFin^\dagger$).
  Define a subcategory \( (\mfldGD)_\dagger\subset \mfldGD \) by
  \[ (\mfldGD)_\dagger := \mfldGD \times_{\OGFin} (\OGFin \times_{\OG} \OG^{\=}) \= \mfldGD \times_{\OG} \OG^{\=} . \]
\end{mydef}
\begin{construction}
  By \cref{lem:S_adeq_triple}, \cref{prop:adeq_fib} and \cite[prop. 11.2]{BarwickSMF1} the triple \\
  \( (\mfldGD, (\mfldGD)_\dagger, (\mfldGD)^\dagger ) \)
  is adequate.
  This condition ensures we can form the $\infty$-category of spans $A^{eff}(\mfldGD,(\mfldGD)_\dagger,(\mfldGD)^\dagger)$.
  Applying the effective Burnside construction to $p \colon \mfldGD \to \OGFin$ we get a functor 
  \[
    \begin{tikzcd}
      A^{eff}(\mfldGD,(\mfldGD)_\dagger,(\mfldGD)^\dagger) \ar{d}{\Upsilon(p)} \\
      A^{eff}(\OGFin,\OGFin\times_{\OG} \OG^{\=}, \OGFin^\dagger) ,
    \end{tikzcd}
  \]
  called the unfurling of $p$ in \cite[def. 11.3]{BarwickSMF1}. 
\end{construction}
\begin{lem}
  The functor $\Upsilon(p)$ is a coCartesian fibration. 
\end{lem}
\begin{proof}
  The functor $\Upsilon(p)$ is an inner fibration by \cite[lem. 11.4]{BarwickSMF1}, and a coCartesian fibration by \cite[lem. 11.5]{BarwickSMF1}
  and \cref{prop:Cart_coCart_edges}.
\end{proof}

\begin{mydef} \label{def:GmfldD_as_coCart_fib}
  Define a coCartesian fibration \( \GmfldD \fib \GFin_* \) by pulling $\Upsilon(p)$ along the inclusion \( \GFin_* \cof A^{eff}(\OGFin, \OGFin\times_{\OG} \OG^{\=}, \OGFin^\dagger) \) of \eqref{eq:GFin_into_Aeff}.
\end{mydef}
\begin{rem} \label{GmfldD_description}
  Unwinding the definition of the effective Burnside category, we see that the objects of $\GmfldD$ are $\OGFin$-manifolds, and a morphism $f \colon M_1\to M_2$ is represented by a span 
  \begin{align*}
    f=\left(\diag{ 
      M_1 \ar[d] & \ar[l] M  \ar[d] \ar[r] & M_2\ar[d] \\ 
      U_1 \ar[d] & \ar[l] U \ar[d] \ar[r] & U_2 \ar[d] \\ 
      O_1 & \ar[l] O_2 \ar[r]^{=} & O_2 
    } \right),
  \end{align*}
  where the 'backwards arrow' is equivalent to a pullback over a summand-inclusion. 
  The morphism $f$ is coCartsian exactly when the 'forward arrow' is equivalent to an identity-of-manifolds
  (see \cref{prop:Cart_coCart_edges} and \cite[lem. 11.5]{BarwickSMF1}).
\end{rem}

\begin{prop} \label{def:GmfldD}
  The coCartesian fibration \( \GmfldD \fib \GFin_* \) of \cref{def:GmfldD_as_coCart_fib} is $G$-symmetric monoidal category whose underlying $G$-category is isomorphic to the $G$-category $\Gmfld$ of \cref{def:Gmfld}.
  We call this $G$-symmetric monoidal structure \emph{$G$-disjoint union} of $G$-manifolds.
\end{prop}
\begin{proof}
  By \cref{def:underlying_G_cat} the underlying $G$-category of \( \GmfldD \) has objects $\OGFin$-manifolds of the form \( (M \to O \xto{=} O) \) and maps represented by spans of the form
  \begin{align*}
    \left( \diag{ M_1 \ar[d] & \ar[l] M \ar[d] \ar[r] & M_2 \ar[d] \\ O_1 \ar[d]^{=} & \ar[l] O_2 \ar[d]^{=} \ar[r]^{=} & O_2 \ar[d]^{=} \\ O_1 & \ar[l] O_2 \ar[r]^{=} & O_2 } \right) 
  \end{align*}
  with left square equivalent to a pullback.
  This $G$-category is isomorphic to $\Gmfld$ by the forgetful functor \( (M\to O \xto{=} O) \mapsto ( M \to O) \). 

  By \cref{G_SM_cat_Fiberwise_char} it is enough to show that for every \( I=(U\to O) \in \GFin_* \) the induced functor \( \prod \rho^W_* \colon  \GmfldD_I \to \prod_{W\in \orb(U)} \Gmfld_{[W]} \) is an equivalence of $\infty$-categories,
  where $\rho^W_*$ is induced by the fibration \( \GmfldD \fib \GFin_* \) and the inert edge 
  \begin{align*}
   \rho^W = \left( \vcenter{ \xymatrix{
     U \ar[d] & \ar@{_(->}[l] W \ar[r]^{=}  \ar[d]^{=} & W \ar[d]^{=} \\
     O & \ar[l] W \ar[r]^{=} & W .
   }} \right), \quad
   \rho^W \in \GFin_* .
  \end{align*}
  Let \( ( M \to U \to O ) \in \Gmfld_I \) be an $\OGFin$-manifold.
  The decomposition $U=\coprod_{W \in \orb(U)} W$ into orbits induces a decomposition of $M$ into a disjoint union \( M = \sqcup_{W \in \orb(U)} M_W \).
  The action of $\rho^W_*$ on \( (M \to U \to O ) \) is specified by a choice of coCartesian lift over $\rho^W$. 
  By the above description of coCartesian edges we see that 
  \begin{align*}
    \left( \vcenter{\xymatrix{
      M \ar[d] & \ar@{_(->}[l] M_W \ar[d] \ar[r]^{=} & M_W \ar[d] \\
      U \ar[d] & \ar@{_(->}[l] W \ar[r]^{=}  \ar[d]^{=} & W \ar[d]^{=} \\
      O & \ar[l] W \ar[r]^{=} & W .
    } } \right)
  \end{align*}
  is such a coCartesian edge, therefore the functor \( \prod \rho^W_* \) is given by 
  \begin{align*}
   \prod \rho^W_* \colon  \GmfldD_I & \to \prod_{W\in \orb(U)} \Gmfld_{[W]}, \\
   \prod \rho^W_* \colon  
    \left( \diag{ M \ar[d] \\ U  \ar[d] \\ O } \right) =
    \left( \diag{ \bigsqcup_{W\in\orb(U)} M_W \ar[d] \\ \coprod_{W\in\orb(U)} W  \ar[d] \\ O } \right) 
    & \mapsto \left( \diag{ M_W \ar[d] \\ W  \ar[d]^{=} \\ W } \right)_{W\in \orb(U)} 
  \end{align*}
  which is an equivalence by inspection.
\end{proof}

\subsection{$G$-disjoint union of $f$-framed $G$-manifolds}
In this subsection we 
lift $G$-disjoint union of $G$-manifolds to a $G$-symmetric monoidal structure on $\ulmfld^{G,f-fr}$.
Recall that $\ulmfld^{G,f-fr}$ was defined as the pullback of $G$-$\infty$-categories (see \cref{def:Gmfld_framed}).
We will show that the $G$-symmetric monoidal structure of $\Gmfld$ lifts to $\ulmfld^{G,f-fr}$ by exhibiting the pullback square of \cref{def:Gmfld_framed} as underlying a pullback square of $G$-symmetric monoidal $G$-$\infty$-categories and $G$-symmetric monoidal functors.

In addition to $G$-disjoint unions of $G$-manifolds we will use the $G$-coCartesian structure, constructed in \cite{Parametrized_algebra} and given by $G$-coproducts.
In general the $G$-coCartesian structure on a $G$-category $\ul\C$ is given by a $G$-$\infty$-operad $\ul\C^\amalg$.
However, we will only use this construction for $\ul\C$ with finite $G$-coproducts, in which case $\ul\C^\amalg$ is a $G$-symmetric monoidal $G$-$\infty$-category.

We show show that the $G$-functors in the pulback square of \cref{def:Gmfld_framed} extend to $G$-symmetric monoidal functors in two steps.
By a formal argument these $G$-functors extend to \emph{lax} $G$-symmetric monoidal functors.
It then remains to verify that these lax $G$-symmetric monoidal functors are in fact $G$-symmetric monoidal.

The following claim allows us to extend $G$-functors to $\ul\C$ from certain $G$-$\infty$-operads to lax $G$-symmetric monoidal functors. 
\begin{lem} \label{coCart_GSM_extension_of_functors}
  Let $\ul\C$ be a $G$-category and $\ul{O}^\otimes$ a unital $G$-$\infty$-operad.
   Restriction to the underlying $G$-category induces an equivalence 
   \begin{align*}
     Alg_G(\ul{O},\ul\C) \to \Fun_G (\ul{O},\ul\C)
   \end{align*}
   between the $\infty$-category of morphisms of $G$-$\infty$-operads from $\ul{O}^\otimes$ to $\ul\C^\amalg$ and the $\infty$-category of $G$-functors between the underlying $G$-categories. 
\end{lem}

Let $B\in \mathbf{Top}^G$ be a $G$-space and $f \colon B\to BO_n(G)$ be a $G$-map.
Endow the parametrized slice $G$-categories $\ulTopG_{/\ul{B}},\,\ulTopG_{/\ul{BO_n(G)}}$ with the $G$-coCartesian $G$-symmetric monoidal structure.
By \cref{coCart_GSM_extension_of_functors} the $G$-functors 
\begin{align*}
  f_* \colon \ulTopG_{/\ul{B}} \to \ulTopG_{/\ul{BO_n(G)}}, \quad 
  \tau  \colon  \Gmfld \to \ulTopG_{/\ul{BO_n(G)}}
\end{align*}
admit an essentially unique lift to lax $G$-symmetric monoidal functors 
\begin{align*}
  f_* \colon \ulTopG_{/\ul{B}} \to (\ulTopG_{/\ul{BO_n(G)}})^\amalg, \quad 
  \tau  \colon  \GmfldD \to (\ulTopG_{/\ul{BO_n(G)}})^\amalg
\end{align*}

The following description of the $G$-coCartesian structure  $(\ulTopG_{/\ul{B}})^\amalg$ is useful when verifying that the lax $G$-symmetric monoidal functors $\tau, f_*$ constructed above are in fact $G$-symmetric monoidal.
\begin{rem} \label{G_coprod_in_ulTopG_over_B}
  Let $I=(U \to G/H)\in \GFin_*$. 
  Then a $U$-family \( x_\bullet  \colon  \ul{U} \to \ulTopG \) can be described by a $G$-map $X \to U$.
  Moreover, under this description the parametrized coproduct $\coprod_I x_\bullet  \colon \ul{G/H} \to \ulTopG$ is given by the $G$-map $X \to U \to G/H$.

  To see this first construct the left fibration associated to $x_\bullet$, and then notice it is a map of $G$-$\infty$-groupoids and therefore can identified with a map of $G$-spaces $X \to U$.
  One should think of the family $X \to U$ as assigning to each $W \in \orb(U)$ the $G$-map $(X|_W \to W) \in \ulTopG_{[W]}$, where we use the explicit model of \cref{ulTopG_as_G_spaces_over_orbits}.
  In order to see that $\coprod_I x_\bullet$ is given by $(X \to U \to G/H) \in \ulTopG_{[G/H]}$ recall that $\coprod_I$ is given by $G$-left Kan extension along $\ul{U} \to \ul{G/H}$, which by \cite[prop. 10.9]{Expose2} is given by (unparametrized) left Kan extension along $\ul{U} \to \ul{G/H}$.
  Applying straightening/unstraightening, we see that $\coprod_I$ is let adjoint to pulling back along $\ul{U} \to \ul{G/H}$, and therefore given by post-composition with $\ul{U} \to \ul{G/H}$.

  Let $B$ be a $G$-space. 
  Combining 
  \cref{ulTopG_over_B_description} with the description of $U$-families in $\ulTopG$ above, we get the following description of  $G$-coproducts in $\ulTopG_{/\ul{B}}$.
  A $U$-family \( x_\bullet \colon  \ul{U} \to \ulTopG_{/\ul{B}} \) is given by a $G$-map $X \to U$ together with a collection of $G$-maps $\set{X|_W \to B}$ indexed by $W\in \orb(U)$.
  Equivalently, \( x_\bullet \colon  \ul{U} \to \ulTopG_{/\ul{B}} \) is given by a pair of $G$-maps $(X \to U,\, X \to B)$.
  The $G$-coproduct $\coprod_I x_\bullet \in \left( \ulTopG_{/\ul{B}} \right)_{[G/H]} $ is given by $(X\to U \to G/H) \in \ulTopG_{[G/H]}$ together with the $G$-map $X \to B$.
\end{rem}

\begin{lem}
  The functor \( \tau \colon  \GmfldD \to (\ulTopG_{/\ul{BO_n(G)}})^\amalg \) is a $G$-symmetric monoidal functor.
\end{lem}
\begin{proof}
  By \cref{tangent_classifier_classifies_TM} and the Segal conditions we have a concrete description of $\tau$. 
  Namely, if $I=(U\to G/H) \in \GFin_*$ and $(M \to U \to G/H) \in \GmfldD_{I}$ is a $\OGFin$-manifold then $\tau ( M \to U \to G/H) \in ( \ulTopG_{/\ul{BO_n(G)}})^\amalg$ 
  is given by $(M\to U \to G/H) \in \ulTopG_I$ together with the $G$-map \( M \to BO_n(G) \) classifying $TM\to M$. 
  Therefore the $G$-coproduct $\coprod_I \tau( M \to U \to G/H)$ is given by 
  $(M \to U \to G/H) \in \ulTopG_{[G/H]}$ together with the $G$-map $M \to BO_n(G)$ classifying $TM \to M$.

  On the other hand, by \cref{GmfldD_description} the $G$-disjoint union $ \sqcup_I M \in \Gmfld_{[G/H]}$ is the $\OG$-manifold given by the composition $M \to U \to G/H$, therefore \( \tau ( \sqcup_I M) \) is given by the $\OG$-manifold $(M \to U \to G/H) \in \Gmfld_{[G/H]}$ together with the $G$-map $ M \to BO_n(G)$ classifying $TM \to M$.
\end{proof}

\begin{prop}
  The $G$-functor \(  f_*  \colon \ulTopG_{/\ul{B}} \to \ulTopG_{/\ul{BO_n(G)}}\) extends to a $G$-symmetric monoidal functor \(  f_*  \colon  (\ulTopG_{/\ul{B}})^\amalg \to (\ulTopG_{/\ul{BO_n(G)}})^\amalg \) .
\end{prop}
\begin{proof}
  This is an immediate consequence of the description of $G$-coproducts in $\ulTopG_{/\ul{B}}$ and $\ulTopG_{/\ul{BO_n(G)}}$:
  for $I=(U\to G/H)$ the diagram
   \begin{align*}
     \xymatrix{
        \left( \ulTopG_{/\ul{B}} \right)^\amalg_I \ar[r]^-{f_*} \ar[d]^{\sqcup_I} & \left( \ulTopG_{/\ul{BO_n(G)}} \right)^\amalg_I \ar[d]^{\coprod_I} \\
        \left( \ulTopG_{/\ul{B}} \right)_{[G/H]} \ar[r]^-{f_*} & \left( \ulTopG_{/\ul{BO_n(G)}} \right)_{[G/H]} 
     } 
   \end{align*}
   is commutativity, since \cref{G_coprod_in_ulTopG_over_B} implies it is given by
   \begin{align*}
     \xymatrix{
       (X \to U , X \to B) \ar@{|->}[d]^{\sqcup_I} \ar@{|->}[r]^-{f_*} & (X \to U , X \to B \xto{f} BO_n(G) ) \ar@{|->}[d]^{\coprod_I} \\
       (X \to U \to G/H, X \to B) \ar@{|->}[r]^-{f_*} & (X \to U \to G/H, X \to B \xto{f} BO_n(G) ) .
     }
   \end{align*}
\end{proof}

It follows that given a $G$-map $f \colon B \to BO_n(G)$ over $G/H$ we can endow $\ulmfld^{G,f-fr}$ with a $G$-symmetric monoidal structure. 
\begin{cor} \label{GmfldD_framed}
  The $G$-symmetric monoidal structure of $G$-disjoint union on $\Gmfld$ lifts to a $G$-symmetric monoidal structure on $\ulmfld^{G,f-fr}$, given by the pullback
  \begin{align*}
    \xymatrix{
      \ulmfld^{G,f-fr,\sqcup} \pullbackcorner \ar[d] \ar[r] &  (\ul{\mathbf{Top}}^{G}_{/\ul{B}})^\amalg \ar[d]^{f_*} \\ 
      \GmfldD \ar[r]^-\tau &  (\ul{\mathbf{Top}}^{G}_{/\ul{BO_n(G)}})^\amalg . 
    }
  \end{align*}
\end{cor}
\begin{proof}
  The $\infty$-category $\Cat^{G,\otimes}_\infty$ of $G$-symmetric monoidal categories admits limits, and the forgetful $G$-functor \( \Cat^{G,\otimes}_{\infty} \to \Cat_\infty^G \) sending a $G$-symmetric monoidal category $\ul\C^\otimes \fib \GFin_*$ to its underlying $G$-category $\ul\C = \ul\C^\otimes \times_{\GFin_*} \OGop$ preserves limits.
\end{proof}
\begin{rem} \label{GmfldDfr_description}
  Informally, we can describe an object of $\ulmfld^{G,f-fr,\sqcup}$ over $(U \to G/H) \in \GFin_*$ as an $\OGFin$-manifold $(M \to U \to G/H)$ together an $f$-framing $f_M \colon M \to B \times G/H$ .
\end{rem}

\begin{mydef} \label{def:ulRepG_framed_operad}
  Let $\ulRep^{G,f-fr,\sqcup}_n \subset \ulmfld^{G,f-fr,\sqcup}$ be the full $G$-subcategory of $\ulmfld^{G,f-fr,\sqcup}$ given by the pullback
  \begin{align*}
    \xymatrix{
      \ulRep^{G,f-fr,\sqcup}_n \pullbackcorner \ar[d] \ar@{^(->}[r] & \ulmfld^{G,f-fr,\sqcup} \ar[d]  \\ 
      \ulRep^{G,\sqcup}_n \ar@{^(->}[r] & \GmfldD . 
    }
  \end{align*}
  It follows that $\ulRep^{G,f-fr,\sqcup}_n \subset \ulmfld^{G,f-fr,\sqcup}$ is the full subcategory of $f$-framed $\OGFin$-manifolds $(E \to U \to G/H)$ where $E \to U$ is a $G$-vector bundle.
  Note that $\ulRep^{G,f-fr,\sqcup}_n \fib \GFin_*$ is a $G$-$\infty$-operad.
\end{mydef}

\subsection{The $G$-category of $G$-disks and the definition of $G$-disk algebras} \label{sec:G_disk}

Our next goal is to define the $G$-symmetric monoidal $G$-$\infty$-category of $G$-disks $\GdiskD$, and its framed variants $\ul\disk^{G,f-fr, \sqcup}$.
These $G$-$\infty$-categories are the point of contact between equivariant algebra and equivariant geometry.

On the one hand,  we use $\GdiskD$ to define $G$-disk algebras, which serve as coefficients for genuine equivariant factorization homology.
In a nutshell, the algebraic structure of a $G$-disk algebra is indexed by equivariant embeddings of $G$-disks. 

On the other hand, $G$-disks capture the local geometry of $G$-manifolds: $G$-disks are designed to be the $G$-tubular neighbourhoods of a configuration of orbits in a $G$-manifold.
We will therefore define $G$-disks as a full $G$-subcategory $\Gdisk \subset \Gmfld$ of the $G$-$\infty$-category of $G$-manifolds.

After defining $\Gdisk$ we show that $G$-disjoint unions endow it with $G$-symmetric monoidal structure (see \cref{def:GdiskD_as_coCart_fib} and \cref{G_disj_union_of_G_disks}).
Finally, we construct $\ul\disk^{G,f-fr}$, a framed version of the $G$-$\infty$-category of $G$-disks (see \cref{def:Gdisk_framed}) and define $f$-framed $G$-disk algebras (see \cref{def:G_disk_alg}).

\begin{mydef}[$G$-disks] \label{def:Gdisk}
  A $G$-disk is a $G$-vector bundle \( E \to O \) rank $n$, where $O\in \OG$ is an orbit.
  Clearly a $G$-disk is an $\OG$-manifold.

  Let \( \Gdisk \subset \Gmfld, \, \disk^G \subset \mfld^G \) be the full subcategories spanned by $\OG$-manifolds equivalent to a composition  \( E \to U \to O \) of $G$-vector bundle $E \to U$ of rank $n$ over a finite $G$-set \( (U \to O) \in \GFin \).
\end{mydef}

\begin{rem}
  The $G$-subcategory \( \Gdisk \subset \Gmfld \) is the full $G$-subcategory generated from $G$-disks by finite $G$-disjoint unions. 
  We think of \( (E \to U \to O) \in \Gmfld \) as a $G$-disjoint union of $G$-disks:
  the decomposition \( U = \sqcup_{W\in\orb(U)} W \) into orbits decomposes $E$ into a disjoint union of $G$-vector bundles \( E_W \to W \), and each composition \( E_W \to W \to O \) exhibits $E_W\to O$ as the topological induction of $ E_W \to W $ along $W \to O$. 

  In fact, $\Gdisk$ is the free $G$-category generated from $H$-representations for $H<G$, considered as $G$-vector bundles over $G/H$, by disjoint unions and topological induction (see \cref{Gdisk_is_Env_G_RepG} below).
\end{rem}

We first verify that $\Gdisk$ is a $G$-$\infty$-category.
\begin{prop} \label{Gdisk_a_G_subcat}
  The subcategory \( \Gdisk \subset \Gmfld \) is a $G$-subcategory stable under equivalences. 
\end{prop}
\begin{proof}
  By \cite[lem. 4.5]{Expose1} it is enough to show that for any coCartesian edge \( x\to y \) in $\Gmfld$ if $x\in \Gdisk$ then $y\in \Gdisk$. 
  Recall that an edge
  \begin{align*}
  \left(\diag{
    M_1 \ar[d] & \ar[l] M \ar[d] \ar[r] & M_2 \ar[d] \\
    O_1 & \ar[l]_{\varphi} O_2 \ar[r]^{=} & O_2
  }\right)
  \end{align*}
  in $\Gmfld$ is coCartesian if and only if the left square is equivalent to a pullback square and the right square is a $G$-isotopy equivalence.
  Let \( (M_1 \to O_1) \in \Gdisk \), then by definition it is equivalent to \( E \to U \to O_1\) for $U$ a finite $G$-set and $E\to U$ a $G$-vector bundle.
  Pulling back along \( \varphi \) shows that \( M \to O_2 \) is equivalent to \( \varphi^* E \to \varphi^* U \to O \), a $G$-vector bundle over a finite $G$-set. 
  Since $M_2 \to O_2$ is equivalent to $M\to O_2$ it follows that \( (M_2 \to O ) \in \Gdisk \). 
\end{proof}

\begin{rem}
  The coCartesian fibration \( \Gdisk \fib \OGop \) is dual to the Cartesian fibration \( \disk^G \to \OG \). 
\end{rem}

\paragraph{$G$-disjoint union of $G$-disks}
We now show (\cref{G_disj_union_of_G_disks}) that $G$-disjoint union of $G$-manifolds (see \cref{def:GmfldD}) induces a $G$-symmetric monoidal structure on $\Gdisk$.
\begin{mydef} \label{def:GdiskD_as_coCart_fib}
  Define \( \GdiskD \subset \GmfldD \) to be the full subcategory spanned by the $\OGFin$-manifolds $ M\to U \to O$ 
  equivalent to $ E \to U' \to U \to O$ where $E\to U'$ is a $G$-vector bundle over a finite $G$-set $U'$.
\end{mydef}
\begin{rem} \label{rem:disk_bun_connected_components}
  Note that if $ M\to U \to O$ is equivalent to $ E \to U' \to U \to O$ where $E\to U'$ is a $G$-vector bundle over a finite $G$-set $U'$, then $U'=\pi_0(E) \= \pi_0(M)$ is the set of connected components of $M$ with the induced action.
\end{rem}
\begin{lem} \label{GdiskD_is_a_G_subcat}
  The subcategory \( \GdiskD \subset \GmfldD \) is a $G$-subcategory stable under equivalences. 
\end{lem}
\begin{proof}
  The proof follows from the characterization of coCartesian edges of \( \GmfldD \fib \GFin_* \) as spans of $\OGFin$-manifolds where the 'backwards arrow' is equivalent to a pullback over a summand-inclusion and the 'forwards arrow' is equivalent to an identity-of-manifolds, following the outline of \cref{Gdisk_a_G_subcat}.
\end{proof}
\begin{cor} \label{G_disj_union_of_G_disks}
  The operation of $G$-disjoint union on $\Gmfld$ induces a $G$-symmetric monoidal structure on the $G$-subcategory $\Gdisk$.
\end{cor}
\begin{proof}
  By \cref{GdiskD_is_a_G_subcat} it is enough to show that to show that the underlying $G$-category of $\GdiskD \fib \GFin_*$ is equivalent to $\Gdisk$.
  Indeed, pulling back along the $G$-functor
  \[ \sigma_{<G/G>}  \colon  \OGop \to \GFin_* ,\quad O \mapsto (O  \xto{=} O) \] 
  we see that the underlying category \( \GdiskD_{<G/G>} \) has objects $\OGFin$-manifolds equivalent to $(E \to U' \to O \xto{=} O) $ for $E \to U'$ a $G$-vector bundle over a $G$-finite set. 
  Therefore the full $G$-subcategory $\GdiskD_{<G/G>}\subset \GmfldD_{<G/G>}$ corresponds to \( \Gdisk \subset \Gmfld \) under the identification 
  \[ \GmfldD_{<G/G>} \simeq \Gmfld , \quad (M \to O \xto{=} O) \mapsto (M\to O). \qedhere \]
\end{proof}

\paragraph{Framed $G$-disks} 
We now define $f$-framed $G$-disks by restricting the underlying $\OG$-manifolds of $f$-framed $\OG$-manifolds to $G$-disks.
\begin{mydef} \label{def:Gdisk_framed}
  Let $B\in \mathbf{Top}^G$ be a $G$-space and $f \colon B \to BO_n(G)$ be a $G$-map.
  Define the $G$-categories of $f$-framed $G$-disks as the pullback on the left.
  \begin{align*}
    \xymatrix{
      \ul{\disk}^{G,f-fr} \pullbackcorner \ar[d] \ar@{^(->}[r] & \ulmfld^{G,f-fr} \ar[d] \\ 
      \Gdisk \ar@{^(->}[r] & \Gmfld  ,
    } \quad
    \xymatrix{
      \ul\disk^{G,f-fr,\sqcup} \pullbackcorner \ar[d] \ar@{^(->}[r] & \ulmfld^{G,f-fr,\sqcup} \ar[d] \\ 
      \GdiskD \ar@{^(->}[r] & \GmfldD . 
    }
  \end{align*}
  The $G$-symmetric monoidal structure of $G$-disjoint union on $\Gdisk$ lifts to a $G$-symmetric monoidal structure on $\ul\disk^{G,f-fr}$, given by the right pullback above.
\end{mydef}

\paragraph{$G$-disk algebras}
We define $G$-disk algebras using $G$-symmetric monoidal functors.
\begin{notation}
  Let \( p \colon \ul\C^\otimes \fib \GFin_*, \, q \colon \D^\otimes \fib \GFin_* \) be two $G$-symmetric monoidal categories. 
  A $G$-symmetric monoidal functor from $\ul\C$ to $\ul\D$ is a functor of $\infty$-categories \( f \colon \ul\C^\otimes \to \D^\otimes \) over $\GFin_*$ that takes $p$-coCartesian edges to $q$-coCartesian edges.
  Denote the $\infty$-category of $G$-symmetric monoidal functors from $\ul\C$ to $\ul\D$ by
  \(
    \Fun^\otimes_{G}(\ul\C,\ul\D) := \Fun_{\GFin_*}(\ul\C^\otimes, \ul\D^\otimes) .
  \)
\end{notation}

\begin{mydef} \label{def:G_disk_alg}
  Let $\ul\C^\otimes \fib \GFin_*$ be a $G$-symmetric monoidal category.
  A \emph{$G$-disk algebra} with values in $\ul\C$ is a $G$-symmetric monoidal functor
  $A \colon \GdiskD \to \ul\C^\otimes$ (see \cref{def:GdiskD_as_coCart_fib}).
  Denote  the $\infty$-category of $G$-disk algebras in $\ul\C$ by \( \Fun^\otimes_G( \Gdisk, \ul\C) \).

  Let  $f \colon B \to BO_n(G)$ a $G$-map, as in \cref{def:Gmfld_framed}.
  An \emph{$f$-framed $G$-disk algebra} with values in $\ul\C$ is a $G$-symmetric monoidal functor
  $A \colon \ul\disk^{G,f-fr,\sqcup} \to \ul\C^\otimes$ (see \cref{GmfldD_framed}). 
  Denote  the $\infty$-category of $G$-disk algebras in $\ul\C$ by \( \Fun^\otimes_G( \ul\disk^{G,f-fr}, \ul\C) \).
\end{mydef}
We will use $G$-disk algebras as coefficients in the definition of $G$-factorization homology in \cref{sec:G_factorization_homology}.

\begin{ex}
  Let $V: pt \to BO_n(G)$ be the $G$-map corresponding to a real $n$-dimensional $G$-representation $V$ (see \cref{ex:V_fr_G_mfld}), and $\ul\disk^{G,V-fr,\sqcup}$ be the $G$-symmetric monoidal category of $V$-framed $G$-disks.
  A $V$-framed $G$-disk algebra is a $G$-symmetric monoidal functor $\ul\disk^{G,V-fr,\sqcup} \to \ul\C^\otimes$.
  In \cref{EV_algs_vs_V_disk_algs} we will see that $V$-framed $G$-disk algebras are equivalent to $\EE_V$-algebras.
\end{ex}

\subsection{$G$-disks as a $G$-symmetric monoidal envelope} \label{sec:G_disk_as_env}
There is a close relationship between $\Gdisk$ and the $G$-$\infty$-category $\ulRep^G_n$ of \cref{def:RepGn}. 
To state it we first define a $G$-$\infty$-operad $\ul{\Rep}^{G,\sqcup}_n$ whose underlying $G$-$\infty$-category is $\ulRep^G_n$ (see \cref{def:ulRepG_framed_operad}), 
and then show that $\Gdisk$ is the $G$-symmetric monoidal envelope of $\ul{\Rep}^{G,\sqcup}_n$.

See \cite{Parametrized_algebra} for the construction and universal property of the $G$-symmetric monoidal envelope.
\begin{mydef}
  Let $\ulRep^{G,\sqcup}_n \subset \GmfldD$ be the full $G$-subcategory on the objects of $\ulRep^G_n$ (using the Segal conditions on the fibers of $\GmfldD$).
  Note that $\ulRep^{G,\sqcup}_n \fib \GFin_*$ is a $G$-$\infty$-operad.
  Equivalently, $\ulRep^{G,\sqcup}_n$ is the full subcategory on $\OGFin$-manifolds $E \to U \to O$ where $E \to U$ is a $G$-vector bundle.
\end{mydef}

\begin{lem} \label{Gdisk_is_Env_G_RepG}
  The $G$-symmetric monoidal $G$-category of $G$-disks, $\GdiskD$, is equivalent to $Env_G(\ulRep^{G,\sqcup}_n)$, the $G$-symmetric monoidal envelope of $\ulRep^{G,\sqcup}_n$.
\end{lem}
\begin{proof}
  Recall that $Env_G(\ulRep^{G,\sqcup}_n)$ is given by the fiber product $\ulRep^{G,\sqcup}_n \times_{\GFin_*} \Arr_G^{act}(\GFin_*)$, where $\Arr_G^{act}(\GFin_*)\subset \Arr_G(\GFin_*)$ is the full subcategory of fiberwise active arrows.
  Unwinding the definition, we identify the objects of $Env_G(\ulRep^{G,\sqcup}_n)$ with 
  \begin{align*}
    \xymatrix{
      E \ar[d] \\ 
      U_1 \ar[d] & \ar[l]_{=} U_1 \ar[d] \ar[r] & U_2 \ar[d] \\
      O & \ar[l]_{=} O  \ar[r]^{=} & O ,
    }
  \end{align*}
  where $E \to U_1$ is a $G$-vector bundle.

  The inclusion $\ulRep^{G,\sqcup}_n \cof \GmfldD$ is a morphism of $G$-$\infty$-operads, so by the universal property of the enveloping $G$-symmetric monoidal $G$-category induces a $G$-symmetric monoidal $G$-functor \( Env_G(\ulRep^{G,\sqcup}_n) \to \GmfldD \), taking an object 
  \begin{align*}
    \xymatrix{
      E \ar[d] \\ 
      U_1 \ar[d] & \ar[l]_{=} U_1 \ar[d] \ar[r] & U_2 \ar[d] \\
      O & \ar[l]_{=} O  \ar[r]^{=} & O  
    }
  \end{align*}
  to the $\OGFin$-manifold $E \to U_1 \to U_2 \to O$.
  Therefore the essential image of $Env_G(\ulRep^{G,\sqcup}_n) \to \GmfldD$ is $\GdiskD$.

  We have to show that the $G$-functor $Env_G(\ulRep^{G,\sqcup}_n) \to \GdiskD$ is a fully faithful (i.e. that it is fiberwise fully faithful). 
  However, the mapping spaces of \( Env_G( \ulRep^{G,\sqcup}_n ) = \ulRep^{G,\sqcup}_n \times_{\GFin_*} \Arr_G^{act}(\GFin_*)\) are given by homotopy pullbacks of the mapping spaces of $\ulRep^{G,\sqcup}_n$ and $\Arr_G^{act}(\GFin_*)$ over $\GFin_*$.
  This follows from the definition of the mapping spaces of $\GmfldD$ after decomposing the mapping spaces of $\GdiskD$ using the Segal conditions.
\end{proof}

It follows that $G$-disk algebras (see \cref{def:G_disk_alg}) are equivalent to algebras over the $G$-$\infty$-operad $\ulRep^{G,\sqcup}_n$.
\begin{cor} \label{G_disk_algs_as_G_operad_algs}
  Let $\ul\C^\otimes$ be a $G$-symmetric monoidal category. 
  The $\infty$-category $\Fun_G^\otimes(\Gdisk, \ul\C)$ of $G$-symmetric monoidal functors $A \colon  \GdiskD \to \ul\C^\otimes$ is equivalent to the $\infty$-category $Alg_G(\ulRep^{G},\ul\C)$ of morphisms of $G$-$\infty$-operads $\ulRep^{G,\sqcup}_n \to \ul\C^\otimes$, i.e algebras of the $G$-$\infty$-operad $\ulRep^{G,\sqcup}_n$ in $\ul\C$.
\end{cor}

A similar result holds for $f$-framed $G$-disks, for $B$ a $G$-space and $f \colon B\to BO_n(G)$ a $G$-map as in \cref{def:Gdisk_framed}.  
\begin{prop} \label{Gdisk_is_Env_G_RepG_framed}
  The $G$-symmetric monoidal category  $\ul\disk^{G,f-fr,\sqcup}$ is equivalent to the $G$-symmetric monoidal envelope of $\ulRep^{G,f-fr,\sqcup}$.
\end{prop}

\subsection{Embedding spaces of $G$-disks and equivariant configuration spaces} \label{sec:G_disks_and_configurations}
We compare the mapping spaces of $f$-framed $\OG$-manifolds with equivariant configuration spaces.
\begin{notation}
  Let \( (M\to G/H) \in \Gmfld \) be an $\OG$-manifold over $G/H$ and $(U\to G/H)\in \GFin $ a finite $G$-set over $G/H$. 
  Denote by \( \Conf^G_{G/H}(U;M) \subset \Map^G_{G/H} (U,M) \) the space of injective $G$-equivariant functions \( U \to M \) over $G/H$ with compact-open topology.
\end{notation}

\begin{rem}
  The space $\Conf^G_{G/H}(U;M)$ can be identified with the space of configurations of disjoint orbits in the $H$-manifold $M|_{eH}$ (the fiber of $M\to G/H$ over the base coset $eH$), where the orbits of the configurations are indexed by the orbits of $U|_{eH}$, with stabilizers specified by  $Stab(W), \, W\in \orb(U|_{eH})$.
\end{rem}
In order to compare equivariant embedding spaces of $G$-disks in $M$ with equivariant configuration spaces we first study the equivariant embedding space of a single $G$-disk. 
\begin{mydef}
  Let $E \to U$ be a $G$-vector bundle over a finite $G$-set, and choose a $G$-equivariant metric on $E$. 
  For $t>0$ define \( B_t(E) \subset E, \, B_t(E) = \set{ v\in E \, \big\vert \, \lVert v \rVert <t } \), so $ B_t(E) \to U$ is the ``open ball of radius $t$'' subbundle. 
  Define $Germ(E,M)= \colim_n Emb^G_{G/H}(B_{\frac{1}{2^n}}(E),M)$. 
\end{mydef}

\begin{lem}
  For $s<t$ the restriction map \( Emb^G_{G/H} (B_t(E), M) \to Emb^G_{G/H} (B_s(E),M) \) is a homotopy equivalence.
\end{lem}
\begin{proof}
  By radial dilation we see that the inclusion \( B_s(E) \cof B_t(E) \) is $G$-isotopic over $G/H$ to a $G$-equivariant homeomorphism.
\end{proof}
\begin{cor} \label{GVB_Emb_is_Germ}
  The restriction map \( Emb^G_{G/H}(E,M) \to Germ(E,M) \) is a homotopy equivalence.
\end{cor}

Let $(E \to U \to G/H)\in \Gdisk$ be a finite $G$-disjoint union of $G$-disks, i.e. \( E \to U \) a $G$-vector bundle,  $U=\pi_0 E$, and \( (M \to G/H) \in \Gmfld \) an $\OG$-manifold.
Precomposition with the zero section inclusion \( U \to E \) defines a fibration
\begin{align}  \label{center_map}
  c \colon  Emb^G_{G/H}(E,M) \fib \Conf^G_{G/H}(U;M) ,
\end{align}
which we think of as sending a configuration of $G$-disks in $M$ to the configuration of points which are in the centers these $G$-disks.

Similarly, for $t>0$ we have fibrations \( Emb^G_{G/H}(B_t(E), M) \fib \Conf^G_{G/H}(U;M) \), whose colimit forms a fibration \(c \colon  Germ(E,M) \fib \Conf^G_{G/H}(U;M) \). 

The following corollary is used in the proof of the axiomatic properties of $G$-factorization homology (see the proofs of \cref{ev0_cofinal} and \cref{lem:disks_over_Mi}).
\begin{cor} \label{Emb_Conf_PB}
  Let $(E \to U \to G/H)\in \Gdisk$ be a finite $G$-disjoint union of $G$-disks, i.e. \( E \to U \) a $G$-vector bundle,  $U=\pi_0 E$.
  Let \( (M \to G/H) \in \Gmfld \) be an $\OG$-manifold and \( N \subset M \) an open $G$-submanifold.
  Then 
  \begin{align*}
    \diag{
      Emb^G_{G/H}(E,N) \ar@{->>}[d]^{c} \ar[r] & Emb^G_{G/H}(E,M) \ar@{->>}[d]^{c} \\ 
      \Conf^G_{G/H}(U;N) \ar[r] & \Conf^G_{G/H}(U;M)
    }
  \end{align*}
  is a homotopy Cartesian square of spaces, where the vertical maps are given by \eqref{center_map}.
\end{cor} 
\begin{proof}
  By \cref{GVB_Emb_is_Germ} the left horizontal maps in the diagram \\
  \begin{equation} \label{eq:Emb_Germ_Conf_along_embedding}
    \begin{tikzcd}
      Emb^G_{G/H}(E,N) \arrow{d} \arrow{r}{\sim} & Germ(E,N) \arrow{d} \arrow[r, two heads, "c"] & \Conf^G_{G/H}(U;N) \arrow[d] \\ 
      Emb^G_{G/H}(E,M) \arrow{r}{\sim} & Germ(E,M) \arrow[r, two heads, "c"] & \Conf^G_{G/H}(U;M)
    \end{tikzcd}
  \end{equation}
  are homotopy equivalences, so we have to show the right square is a homotopy pullback square.
  Since the right horizontal arrows are fibrations, it is enough to show that the right square is a pullback square.
  
  Let $x_\bullet \in \Conf^G_{G/H}(U;N)$, given by an injective $G$-map $x_\bullet \colon  U \to N$.
  For $t\in \R$ denote the fiber of $Emb^G_{G/H}(B_t(E),N) \fib Conf_U(N)$ by $Emb^G_{G/H}(B_t(E),N)_{x_\bullet}$.
  We have a map of fibrations \\
  \begin{tikzcd}
    Emb^G_{G/H}(B_t(E),N)_{x_\bullet} \arrow[d] \arrow[r] & Emb^G_{G/H}(B_t(E),N) \arrow[d] \arrow[r, two heads] & \Conf^G_{G/H}(U;N) \arrow[d] \\
    Emb^G_{G/H}(B_t(E),M)_{x_\bullet} \arrow[r] & Emb^G_{G/H}(B_t(E),M) \arrow[r, two heads] & \Conf^G_{G/H}(U;M) .
  \end{tikzcd} \\
  For small enough $t>0$ the left vertical arrow is an isomorphism, hence the right square is a pullback square. 
  Since filtered colimits commute with pullbacks in $\mathbf{Top}$, we see that the right square of diagram \eqref{eq:Emb_Germ_Conf_along_embedding} is indeed a pullback square. 
\end{proof}

Our goal for the rest of this subsection is to study the framed version of the map \eqref{center_map},
and show that its $V$-framed variant is an equivalence (\cref{ex:framed_GVB_emb_is_conf}).
This fact will be used in \cref{sec:DV_and_Gdisk_V_fr} to compare the $G$-$\infty$-operad $\ulRep^{G,V-fr, \sqcup}_n$ (\cref{def:ulRepG_framed_operad}) with the classical $G$-operad of little disks in $V$.

We begin by showing that the decomposition of the configuration of $G$-disks $E$ into orbits of $G$-disks induces a decomposition on its space of $G$-embeddings into $M$.
\begin{prop} \label{Emb_Conf_product_hPB}
  Let $(M \to G/H) \in \Gmfld_{[G/H]}$ be an $\OG$-manifold over $G/H$ and $(E \to U \to G/H)\in \Gdisk_{[G/H]}$. 
  For $W\in \orb(U)$ let $E_W\in \Gdisk_{[G/H]}$ denote $(E|_W \to W \to G/H)$, the restriction of the vector bundle $E\to U$ to the orbit $W \subseteq U$.
  Then the commutative square of spaces
  \begin{align*}
    \xymatrix{
      Emb^G_{G/H}(E,M) \pullbackcorner \ar[d]^c \ar[r] & \prod\limits_W Emb^G_{G/H}(E_W,M) \ar[d]^c \\
      \Conf^G_{G/H}(U;M) \ar[r] & \prod\limits_{W} \Conf^G_{G/H}(W;M)
    }
  \end{align*}
  is a homotopy pullback square, where the products are indexed by $W\in \orb(U)$, 
  and the vertical maps are given by \eqref{center_map}.
\end{prop}
\begin{proof}
  By \cref{GVB_Emb_is_Germ} the left horizontal maps in the diagram \\
  \begin{equation*} 
    \begin{tikzcd}
      Emb^G_{G/H}(E,M) \arrow{d} \arrow{r}{\sim} & Germ(E,M) \arrow{d} \arrow[r, two heads, "c"] & \Conf^G_{G/H}(U;M) \arrow[d] \\ 
      \prod\limits_{W} Emb^G_{G/H}(E_W,M) \arrow{r}{\sim} & \prod\limits_{W} Germ(E_W,M) \arrow[r, two heads, "c"] & \prod\limits_{W} \Conf^G_{G/H}(W;M)
    \end{tikzcd}
  \end{equation*}
  are equivalences, so it is enough to show that right square is a homotopy pullback square.
  The right horizontal maps are Kan fibrations, therefore it is enough to show this square is a pullback square. 
  This is clear, since the induced map on the fibers of the horizontal maps is a homeomorphism.
\end{proof}
\begin{prop} \label{Emb_Conf_product_hPB_framed}
  Let $M\in \ulmfld^{G,f-fr}_{[G/H]}$ be an $f$-framed $\OG$-manifold over $G/H$, given by a $\OG$-manifold $M\to G/H$ together with an $f$-framing $f_M \colon  M \to B $ lifting $\tau_M \colon  M \to BO_n(G)$.
  Let $E\in \ul\disk^{G,f-fr}_{[G/H]}$, given by $(E \to U \to G/H)\in \Gdisk_{[G/H]}$ and $f$-framing $f_E \colon  E \to B $.
  For $W\in \orb(U)$ denote $E_W\in \ul\disk^{G,f-fr}_{[G/H]}$ denote the restricted $G$-vector bundle $(E|_W \to W \to G/H)$, with the restricted framing $f_W \colon  E|_W \subset E \xto{f_E} B \times G/H$.

  Then the commutative square of spaces
  \begin{align*}
    \xymatrix{
      Emb^{G,f-fr}_{G/H}(E,M) \pullbackcorner \ar[d]^c \ar[r] & \prod\limits_W Emb^{G,f-fr}_{G/H}(E_W,M) \ar[d]^c \\
      \Conf^G_{G/H}(U;M) \ar[r] & \prod\limits_{W} \Conf^G_{G/H}(W;M)
    }
  \end{align*}
  is a homotopy pullback square, where the products are indexed by $W\in \orb(U)$, 
  and the vertical maps are given by precomposition with the zero section.
\end{prop}
\begin{proof}
  Recall the notation of \cref{ulTopG_over_B_description}, 
  \begin{align*}
    & \ul{B}(G/H) ,\, \ul{BO_n(G)}(G/H)  \in \GTop_{/G/H} , \\ 
    & \ul{B}(G/H) = ( B \times G/H \to G/H ) , \quad \ul{BO_n(G)}(G/H) = ( BO_n(G) \times G/H \to G/H ) .
  \end{align*}
  Consider the commutative diagram
  \[
    \begin{tikzcd}[row sep=1.5em, column sep = -3em]
      Emb^{G,f-fr}_{G/H}(E,M)  \ar[dd] \ar[rr] \ar[dr] & & \Map^G_{/ \ul{B}(G/H) }  ( E ,M ) \ar[dd] \ar[dr] \\
      & \prod\limits_W Emb^{G,f-fr}_{G/H}(E_W,M) \ar[rr, crossing over] & & \prod\limits_W \Map^G_{/ \ul{B}(G/H)} ( E_W ,M ) \ar[dd] \\
      Emb^G_{G/H}(E,M) \ar[rr] \ar[dr] & &  \Map^G_{/ \ul{BO_n(G)}(G/H) }  ( E ,M )  \ar[dr] \\
      & \prod\limits_W Emb^G_{G/H}(E_W,M) \ar[from=uu, crossing over] \ar[rr] & & \prod\limits_W \Map^G_{/ \ul{BO_n(G)}(G/H) } (E_W , M ) ,
    \end{tikzcd}
  \]
  where
  \begin{align*}
    \Map^G_{/ \ul{B}(G/H) }  ( E ,M ) &= \Map^G_{/ \ul{B}(G/H) }  ( E \xto{\overline{f_E}} B \times G/H ,M \xto{\overline{f_M}} B \times G/H ) , \\
    \Map^G_{/ \ul{BO_n(G)}(G/H) }  ( E ,M ) &= \Map^G_{/ \ul{BO_n(G)}(G/H) }  ( E \xto{\overline{f_E}} BO_n(G) \times G/H, M \xto{\overline{f_M}} BO_n(G) \times G/H ) .
  \end{align*}
  The forward and backward faces of the cube are homotopy pullback squares by definition (see \cref{rem:f_framed_mapping_spaces}).
  The diagonal morphisms on the right are equivalences, since 
  \begin{align*}
    & \left( E \xto{\overline{f_E}} B \times G/H \right)  = \coprod_W \left( E_W \xto{\overline{f_W}} B \times G/H \right) , \\
    & \left( E \xto{\overline{\tau_E}} BO_n(G) \times G/H \right) = \coprod_W \left( E_W \xto{\overline{\tau_W}} BO_n(G) \times G/H \right) ,
  \end{align*}
  and in particular the right face is a homotopy pullback square.
  By \cite[lem. 4.4.2.1]{HTT} the left face is a homotopy pullback square.

  Note that the left face is above diagram the same as the top square of the diagram
  \begin{align*}
    \xymatrix{
      Emb^{G,f-fr}_{G/H}(E,M) \pullbackcorner \ar[d] \ar[r] & \prod\limits_W Emb^{G,f-fr}_{G/H}(E_W,M) \ar[d] \\
      Emb^G_{G/H}(E,M) \pullbackcorner \ar[d]^c \ar[r] & \prod\limits_W Emb^G_{G/H}(E_W,M) \ar[d]^c \\
      \Conf^G_{G/H}(U;M) \ar[r] & \prod\limits_{W} \Conf^G_{G/H}(W;M)
    }
  \end{align*}
  By \cref{Emb_Conf_product_hPB} the bottom square is a homotopy pullback square, hence by \cite[lem. 4.4.2.1]{HTT} so is the outer rectangle. 
\end{proof}
 
\paragraph{The endomorphism space of a single framed $G$-disk} 
We identify the endomorphism space of a single framed $G$-disk as a loop space.
Let $E\xto{\pi} G/H$ be a $G$-vector bundle.
Note that as an object of $\ulTopG_{[G/H]}$ it is equivalent to the terminal object $(G/H \xto{=} G/H)$, so the $G$-tangent classifier \( \tau_E  \colon  E \to BO_n(G) \) is given by a choice of connected component of $(BO_n(G))^H \simeq \coprod_V \BAut_{\Rep^H} (V)$, i.e an $H$-representation $V$ of dimension $n$.
In particular, we have an isomorphism $ E \= V \times_H G$ of $G$-vector bundles over $G/H$.

An $f$-framing on $E$ is given by a $G$-map \( e \colon  E \to B  \) lifting $V \colon  E \to BO_n(G) $ up to $G$-homotopy.
Using the equivalence $\Map^G(E, B ) \simeq \Map^G(G/H, B) \simeq \Map^H(pt, B) \simeq B^H$ we can consider $e$ as a point in $B^H$.

\begin{prop} \label{framed_rep_emb_is_loopspace}
  Let $E\to G/H$ be a $G$-vector bundle with $f$-framing $e \colon  E \to B$. 
  Then the endomorphism space of $(E\to G/H) \in \mfld^{G,f-fr}_{[G/H]}$ is weakly equivalent to the loop space of $B^H$ with base point $e$, 
\( Emb^{G,f-fr}_{G/H}(E, E) \simeq  \Omega_e B^H. \)
\end{prop}
\begin{proof}
  The endomorphism space of $E$ is given by the homotopy pullback 
  \begin{equation} \label{PB_disk_framed_aut}
    \begin{tikzcd}[column sep = 0.5em]
      Emb^{G,f-fr}_{G/H}(E,E) \pbcorner \arrow[d] \arrow[r] & \Map^G_{/ \ul{B}(G/H) } (E \xto{\overline{e}} B \times G/H ,E \xto{\overline{e}} B \times G/H) \arrow{d}{(f\times G/H)_*} \\
      Emb^G_{G/H}(E,E) \arrow{r}{\tau} &  \Map^G_{/ \ul{BO_n(G)}(G/H) } (E \xto{\overline{V}} BO_n(G) \times G/H ,E \xto{\overline{V}} BO_n(G) \times G/H) . 
    \end{tikzcd}
  \end{equation}
  We prove our claim by identifying the mapping spaces on the right column with loop spaces and showing that the horizontal maps are equivalences.

  Since \( \Map^G_{/ \ul{B}(G/H) } (E \xto{\overline{e}} B \times G/H ,E \xto{\overline{e}} B \times G/H) \) is a mapping space in the slice category 
  $\left( \ulTopG_{[G/H]} \right)_{/\ul{B}(G/H)}$ 
  it is equivalent to the homotopy pullback
  \[
    \begin{tikzcd}[column sep = 0.5em]
      \Map^G_{/ \ul{B}(G/H) } (E \xto{\overline{e}} B \times G/H ,E \xto{\overline{e}} B \times G/H) \pbcorner \arrow[d] \arrow[r] & \Map^G_{G/H} (E \to G/H, E\to G/H) \arrow{d}{\overline{e}_*} \\
      \ast \arrow{r}{\overline{e}} & \Map^G_{G/H}(E\to G/H,B \times G/H \to G/H) .
    \end{tikzcd}
  \]
  Since $(E \to G/H) \in \ulTopG_{[G/H]}$ is terminal we have
  \begin{align*}
    \Map^G_{G/H} (E \to G/H, E \to G/H) & \simeq \Map^G_{G/H}( G/H \xto{=} G/H, G/H \xto{=} G/H) = \ast, \\
    \Map^G_{G/H}( E \to G/H, B \times G/H \to G/H) & \simeq \Map^G_{G/H}(G/H \xto{=} G/H, B \times G/H \to G/H) \\
    & \simeq \Map^G( G/H, B) \= B^H ,
  \end{align*}
  hence \( \Map^G_{/ \ul{B}(G/H) } (E \xto{\overline{e}} B \times G/H ,E \xto{\overline{e}} B \times G/H)  \simeq \Omega_e B^H \).

  Replacing $B$ with $BO_n(G)$, the same calculation shows 
  \begin{align*}
    \Map^G_{/ \ul{BO_n(G)}(G/H) } (E \xto{\overline{V}} BO_n(G) \times G/H ,E \xto{\overline{V}} BO_n(G) \times G/H) 
     & \simeq \Omega_V \left( \coprod_{\rho \colon  H \curvearrowright \R^n} \BAut_{\Rep^H}(\rho) \right) \\
     & = \Omega \BAut_{\Rep^H}(V) . 
  \end{align*}
  Identify \( Emb^G_{G/H}(E,E) \= Emb^G_{G/H}(V\times G/H,V \times G/H) \= Emb^H(V, V) \) in the homotopy pullback square \eqref{PB_disk_framed_aut} we get a homotopy pullback square 
  \[
    \begin{tikzcd}
        & & Emb^{G,f-fr}_{G/H}(E,E) \pbcorner \arrow[d] \arrow[r] & \Omega_e B^H \arrow{d} \\
        \Aut_{\Rep^H}(V) \arrow[r, hook] & Emb^H_0(V,V) \arrow[r, hook] & Emb^H(V,V) \arrow{r}{\tau} &  \Omega \BAut_{\Rep^H}(V)  , 
    \end{tikzcd}
  \]
  where $Emb^H_0(V,V)$ is the subspace of $H$-equivariant self embeddings $V \cof V$ that fix the origin.

  By \cref{G_Kister} the bottom left map is a weak equivalence, and the middle bottom arrow is clearly a homotopy equivalence.
  Since the composition of the bottom maps is the known equivalence \( \Aut_{\Rep^H}(V) \to \Omega \BAut_{\Rep^H} (V) \), we conclude that $\tau$ is a weak equivalence, and therefore the top map of the homotopy pullback square is a weak equivalence as well. 
\end{proof}

Finally, we return to the $V$-framed variant of the map \eqref{center_map}.
\begin{ex} \label{ex:framed_GVB_emb_is_conf}
  Consider $V$-framed manifolds for $V$ be a real $n$-dimensional $G$-representation (\cref{ex:V_fr_G_mfld}), and let
  $E\in \ul\disk^{G,V-fr}_{[G/H]}$ be given by $(E \to U \to G/H)\in \Gdisk_{[G/H]}$, with $V$-framing inducing a trivialization \( E \= U \times V \) of the $G$-vector bundle $E \to U$.
  Consider the mapping space from $E$ to $(V\times G/H \to G/H) \in \ul\disk^{G,V-fr}_{[G/H]}$.
  For every orbit $W\in \orb(U)$ we have $E_W \= W \times V$ as $G$-vector bundles over $W$. 
  Therefore the homotopy fiber of $c \colon  Emb^{G,V-fr}_{G/H} (E_W , V \times G/H) \to \Conf^G_{G/H}(W;V \times G/H)$ is equivalent to the loop space of a point (see \cref{framed_rep_emb_is_loopspace}),
  hence contractible.
  It follows that the map
  \( \prod\limits_W Emb^{G,V-fr}_{G/H}(E_W,V \times G/H) \to \prod\limits_{W} \Conf^G_{G/H}(W;V \times G/H) \)
  is an equivalence, since its homotopy fibers are contractible.
  By \cref{Emb_Conf_product_hPB_framed} precomposition with the zero section of $E \to U$ induces a homotopy equivalence 
  \begin{align} \label{V_framed_emb_is_conf}
    c \colon Emb^{G,V-fr}_{G/H}(E,V \times G/H) \iso \Conf^G_{G/H}(U;V \times G/H) .
  \end{align}
  More generally, for any $V$-framed $\OG$-manifold $M \in \ulmfld^{G,V-fr}_{[G/H]}$ we have
  \begin{align} 
    c \colon Emb^{G,V-fr}_{G/H}(E,M) \iso \Conf^G_{G/H}(U;M) ,
  \end{align}
  since by \cref{Emb_Conf_product_hPB_framed} the homotopy fibers are contractible.
\end{ex}
We will use \cref{ex:framed_GVB_emb_is_conf} in \cref{sec:DV_and_Gdisk_V_fr}.

\subsection{Comparison of the equivariant little disks $G$-operad and the $G$-$\infty$-operad of $V$-framed representations} \label{sec:DV_and_Gdisk_V_fr}
Let $V$ be a real $n$-dimensional representation of $G$, 
and $\ulRep^{G,V-fr,\sqcup}$ the $G$-$\infty$-operad of \cref{def:ulRepG_framed_operad}.
In this subsection we define the $G$-$\infty$-operad $\EE_V$ of little $G$-disks (\cref{def:EV_G_operad}) using the genuine operadic nerve construction of Bonventre,
and show that it is equivalent to $\ulRep^{G,V-fr,\sqcup}$ (\cref{DV_is_ulRepG_framed_operad}), hence $\EE_V$-algebras are equivalent to $V$-framed $G$-disk algebras. 

We first review the relevant details of Bonventre's construction. 
This construction is best understood in the light of \cite[thm. III]{Bonventre_Pereira_genuine_G_operads} which gives a (right) Quillen equivalence 
\[ i_*  \colon  sOp^G \to sOp_G \]
between the $G$-graph model structure on simplicial $G$-operads (where weak equivalences is detected on graph-subgroup fixed points) and the projective model structure on genuine $G$-operads.

\begin{construction}[The genuine equivariant category of operators, see {\cite[def. 4.1]{Bonventre}}]
  Let $\mathcal{P}\in sOp_G$ be a genuine $G$-operad. 
  Define a simplicial category $\mathcal{P}^\otimes$ as follows.
  The objects of $\mathcal{P}^\otimes$ are objects of $\GFin_*$, i.e. $G$-maps $U\to G/H$ from a finite $G$-set to a $G$-orbit.
  The simplicial space of maps $\mathcal{P}^\otimes(U_1 \to G/H , U_2 \to G/K)$ is given by 
  \begin{align*}
    \Map_{\mathcal{P}^\otimes} \left( \substack{ U_1 \\ \downarrow \\ G/H}, \substack{ U_2 \\ \downarrow \\ G/K} \right) 
    = \coprod\limits_{\varphi} \prod\limits_{W\in \orb(U_2)} \mathcal{P} \left( \substack{ f^{-1}(W) \\ \downarrow \\ W} \right) ,
  \end{align*}
  where the coproduct is indexed by \( \varphi \in \Map_{\GFin_*}  \left( \substack{ U_1 \\ \downarrow \\ G/H}, \substack{ U_2 \\ \downarrow \\ G/K} \right) \).
  Composition in $\mathcal{P}^\otimes$ is defined using coproducts of the composition maps of the genuine $G$-operad $\mathcal{P}$.
\end{construction}
\begin{thm}[{\cite[thm. 4.10]{Bonventre}}]
  Let $\mathcal{P}\in sOp_G$ be a genuine $G$-operad, and $N^\otimes(\mathcal{P})$ the coherent nerve of the $\mathcal{P}^\otimes$.
  If $\mathcal{P}\in sOp_G$ is locally fibrant, then $N^\otimes(\mathcal{P})$ is a $G$-$\infty$-operad.
\end{thm}
We call $N^\otimes(\mathcal{O})$ as the \emph{genuine operadic nerve} of $\mathcal{O}$.
\begin{cor}[{\cite[cor 6.3]{Bonventre}}]
  Let $\mathcal{O} \in sOp^G$ be a graph-fibrant simplicial $G$-operad with a single color. 
  Then $i_* \mathcal{O} \in sOp_G $ is locally fibrant, and thus there exists a $G$-$\infty$-operad $N^\otimes(\mathcal{O})$ associated to $\mathcal{O}$.
\end{cor}
In particular, the genuine coherent nerve construction associates a $G$-$\infty$-operad to the equivariant little disk operad.
\begin{ex}[{\cite[ex. 6.5]{Bonventre}}] \label{ex:coherent_nerve_of_D_V}
  Let $V$ be a real orthogonal $n$-dimensional $G$-representation, and $D(V)$ the open unit disk of $V$.
  For $H<G$ and $U$ a finite $H$-set let $Emb^{Aff,H}( U \times D(V), D(V) )$ denote the space of $H$-equivariant affine embeddings $U\times D(V) \cof D(V)$.
  Let $\mathcal{D}_V$ be the \emph{little $V$-disks operad} (see e.g \cite[def. 1.1]{Guillou_May_Permutative_G_cats} or \cite[def. 3.11(ii)]{Blumberg_Hill}).
  Applying the functor $\Sing$ to the spaces $(\mathcal{D}_V)_n$ we get a locally fibrant simplicial $G$-operad, hence an associated $G$-$\infty$-operad $N^\otimes(\mathcal{D}_V)$.

  The mapping spaces of $N^\otimes(\mathcal{D}_V)$ are given by  
  \begin{align*}
    \Map_{ N^\otimes(\mathcal{D}_V) } \left( \substack{ U_1 \\ \downarrow \\ G/H } , \substack { U_2 \\ \downarrow \\ G/K}  \right) 
    = \coprod\limits_{\varphi} \prod\limits_{Gx \in \orb(U_2)} Emb^{Aff,Stab(x)} (f^{-1}(x) \times D(V), D(V) ) .
  \end{align*}
\end{ex}
\begin{mydef} \label{def:EV_G_operad}
  Fix a real orthogonal $G$-representation $V$, and let $\mathcal{D}_V$ denote the $G$-operad of little $V$-disks.
  Let $\EE_V^\otimes$ denote the genuine operadic nerve $N^\otimes(\mathcal{D}_V)$ of \cite[ex. 6.5]{Bonventre}.
\end{mydef}
Before defining $\EE_V$-algebras we recall the definition of a $G$-$\infty$-operad map.
\begin{notation}
Let $\P^\otimes \to \GFin_*$, $\Q^\otimes \to \GFin_*$ be $G$-$\infty$-operads (see \cite[def. 3.1]{Nardin_thesis}).
A map of $G$-$\infty$-operads from $\P^\otimes$ to $\Q^\otimes$ is a map of simplicial sets $f \colon \P^\otimes \to \Q^\otimes$ such that 
\begin{enumerate}
  \item The diagram
    \[
      \begin{tikzcd}
        \P^\otimes \ar[dr] \ar[rr, "f"] & & \Q^\otimes \ar[dl] \\
        & \GFin_*
      \end{tikzcd}
    \]
    commutes. 
  \item The functor $f$ carries coCartesian edges over inert morphisms to coCartesian edges. 
\end{enumerate}
\end{notation}
\begin{mydef} \label{def:EV_algs}
  Let $\ul\C^\otimes \fib \GFin_*$ be a $G$-symmetric monoidal category.
  An \emph{$\EE_V$-algebra} in $\ul\C$ is a map of $G$-$\infty$-operads $A \colon \EE_V \to \ul\C^\otimes$.
  Let $Alg_{\EE_V}(\ul\C) \subseteq \Fun_{/\GFin_*} (\EE_V^\otimes, \ul\C^\otimes)$ denote the full subcategory spanned by $\EE_V$-algebras.
\end{mydef}

\paragraph{Comparison with $\ulRep^{G,V-fr,\sqcup}$.}
We can now easily compare 
the $G$-$\infty$-operads $\EE_V$ of \cref{def:EV_G_operad} and  $\ulRep^{G,V-fr,\sqcup}$ of \cref{def:ulRepG_framed_operad}.

We start with some observations.
Fix a $V$-framed $G$-diffeomorphism $D(V) \= V $.
Let $H<G$ and $U'$ a finite $H$-set and $U = G \times_H U'$ its topological induction. 
Then the topological induction of $U' \times D(V)$ from $H$ to $G$ is $ U \times D(V) $.
Note that the induced map $U \times D(V) \to G/H$ is $G$-vector bundle equivalent to $U \times V \to U$ by our chosen diffeomorphism, and hence $V$-framed.
Let \( Emb^{Aff,G}_{G/H} (U \times D(V) , G/H \times D(V) ) \) denote the space of affine $G$-embeddings over $G/H$.
Note that restriction to the fiber over $eH$ defines a homeomorphism 
\[ Emb^{Aff,H}(U' \times D(V) , D(V) ) \= Emb^{Aff,G}_{G/H} (U \times D(V) , G/H \times D(V) )  . \]
On the other hand, affine $G$-embeddings $U \times D(V) \cof G/H \times D(V)$ over $G/H$ are clearly $V$-framed (using the chosen $G$-diffeomorphism \(D(V) \= V \).
Therefore we have a map 
\begin{align*}
  Emb^{Aff,G}_{G/H} \left( U \times D(V) , G/H \times D(V) \right) \cof Emb^{G,V-fr}_{G/H} \left( U \times D(V) , G/H \times D(V) \right) .
\end{align*}

Construct a functor \( F \colon  \EE_V  \to \ulRep^{G,V-fr,\sqcup} \) over $\GFin_*$ as follows.
For every finite $G$-set $U$ 
define $F(U \to G/H) = (U \times D(V) \to U \to G/H) $.
Define $F$ on mapping spaces by the embeddings
\begin{align*}
  \Map_{ \EE_V } \left( \substack{ U_1 \\ \downarrow \\ G/H } , \substack { U_2 \\ \downarrow \\ G/K}  \right) 
  & = \coprod\limits_{\varphi} \prod\limits_{Gx \in \orb(U_2)} Emb^{Aff,Stab(x)} (f^{-1}(x) \times D(V), D(V) ) \\
   & \cof \coprod\limits_{\varphi} \prod\limits_{W\in \orb(U_2)} Emb^{G,V-fr}_W (f^{-1}(W) \times D(V), G/H \times D(V) ) \\
   & = \Map_{\ulRep^{G,V-fr,\sqcup}}  ( U_1 \times D(V), U_2 \times D(V) ) .
\end{align*}

\begin{prop} \label{DV_is_ulRepG_framed_operad}
  The functor \( F \colon  \EE_V  \to \ulRep^{G,V-fr,\sqcup} \) is an equivalence of $G$-$\infty$-operads.
\end{prop}
\begin{proof}
  By construction $F$ is a functor over $\GFin_*$, therefore it is enough to show that $F$ is an equivalence of $\infty$-categories. 
  Clearly $F$ is essentially surjective, since any $V$-framed $G$-vector bundle over a finite $G$-set $U$ is equivariant to $U \times V \= U \times D(V)$.
  We therefore have to show that $F$ is fully faithful. 

  By the Segal conditions it is enough to show that $F$ induces an equivalence of spaces 
  \begin{align*}
    F \colon  \Map_{ \EE_V } \left( \substack{ U \\ \downarrow \\ G/H } , \substack { G/H \\ \downarrow \\ G/H}  \right) 
      \to \Map_{\ulRep^{G,V-fr,\sqcup}} \left( \substack{ U \times D(V) \\ \downarrow \\ U \\ \downarrow \\ G/H } , \substack { G/H \times D(V) \\ \downarrow \\ G/H \\ \downarrow \\ G/H}  \right)  
  \end{align*}
  on the mapping spaces over \( \varphi \in \Map_{\GFin_*}  \left( \substack{ U \\ \downarrow \\ G/H}, \substack{ G/H \\ \downarrow \\ G/H} \right) \).
  By \cref{ex:coherent_nerve_of_D_V} we have 
  \begin{align*}
    \Map_{ \EE_V } \left( \substack{ U \\ \downarrow \\ G/H } , \substack { G/H \\ \downarrow \\ G/H}  \right) \simeq \prod\limits_{W\in \orb(G/H)} Emb^{Aff,Stab(W)} (f^{-1}(W) \times D(V), D(V) ) ,
  \end{align*}
  and since $\ulRep^{G,V-fr,\sqcup} \subset \ulmfld^{G,V-fr,\sqcup}$ is a full $G$-subcategory we have 
  \begin{align*}
      \Map_{\ulRep^{G,V-fr,\sqcup}} \left( \substack{ U \times D(V) \\ \downarrow \\ U \\ \downarrow \\ G/H } , \substack { G/H \times D(V) \\ \downarrow \\ G/H \\ \downarrow \\ G/H}  \right)  
      & = Emb^{G,V-fr}_{G/H} ( U \times D(V),  G/H \times D(V) )  \\
      & \= Emb^{G,V-fr}_{G/H} ( U \times V,  G/H \times V ).
  \end{align*}

  Consider the commutative diagram
  \begin{align*}
    \xymatrix{
      \prod\limits_{W\in\orb(U)} Emb^{Aff,Stab(W)} (f^{-1}(W) \times D(V), D(V) ) \ar[d]^{c} \ar[r]^-F & Emb^{G,V-fr}_{G/H}( U \times V, G/H  \times V ) \ar[d]^{c} \\
      \prod\limits_{W\in\orb(U)} \mathbf{Inj}^{Stab(W)} (f^{-1}(W) , V ) \ar[r] & \Conf^G_{G/H} (U;V \times G/H) 
    }
  \end{align*}
  where the vertical map is given by taking the centers of disks, and the right vertical map is given by precomposition with the zero section.
  We wish to prove that the top horizontal map is an equivalence.
  The left vertical map is known to be an equivalence (see \cite[prop. 4.19]{Blumberg_Hill} and \cite[lem 1.2]{Guillou_May_Permutative_G_cats}).
  The right vertical map is an equivalence by \cref{ex:framed_GVB_emb_is_conf}, and the bottom horizontal map is a homeomorphism by inspection.
\end{proof}
We immediately see that $\EE_V$-algebras are equivalent to $V$-framed $G$-disk algebras.
\begin{cor} \label{EV_algs_vs_V_disk_algs}
  There is an equivalence of $\infty$-categories $Alg_{\EE_V}(\ul\C) \simeq \Fun^\otimes_G( \ul\disk^{G,V-fr} , \ul\C)$.
\end{cor}
\begin{proof}
  Precomposition with the equivalence $F: \EE_V \iso \ulRep^{G,V-fr,\sqcup}$ of \cref{DV_is_ulRepG_framed_operad} induces an equivalence 
  \[
    Alg_{\EE_V}(\ul\C) \iso Alg_{\ulRep^{G,V-fr}}(\ul\C).
  \]
  By \cref{Gdisk_is_Env_G_RepG_framed} the $G$-symmetric monoidal envelope of $\ulRep^{G,V-fr,\sqcup}$ is equivalent to $\ul\disk^{G,V-fr,\sqcup}$,
  so by its universal property we have 
  \[
     Alg_{\ulRep^{G,V-fr}}(\ul\C) \simeq \Fun^\otimes_G( \ul\disk^{G,V-fr}, \ul\C).
  \]
\end{proof}

\section{Genuine $G$-factorization homology} \label{sec:G_factorization_homology}
In this section we use the $G$-categories $\ulmfld^{G,f-fr}$ and $\ul\disk^{G,f-fr}$ to define genuine equivariant factorization homology. 
We define $G$-factorization homology, first as a parametrized colimit (\cref{def:G_FH_as_colimit}), then as a $G$-functor (\cref{G_FH_functor}) and finally as a $G$-symmetric monoidal functor (\cref{def:G_FH_as_G_SM_functor}).

\subsection{The definition of $G$-disk algebras and $G$-factorization homology as a $G$-functor} \label{sec:G_factorization_homology_G_functor}
In this subsection we define equivariant factorization homology (see \cref{G_FH_functor}).
This is an smooth equivariant version of the factorization homology of \cite{AF} and of topological chiral homology of \cite[7.5.2]{HA}.

In order to define genuine $G$-factorization homology using parametrized $\infty$-colimits we first recall the definition of a parametrized over-category from \cite{Expose2}.
The parametrized over-category plays the role of an indexing category in the $G$-colimit defining factorization homology below (see \cref{def:G_FH_as_colimit}), and more generally in the $G$-colimit formula for $G$-left Kan extensions (see \cite[thm. 10.3]{Expose2}).

Let $\ul\C$ be a $G$-category and $x\in \ul\C_{[G/H]}$ an object over $G/H$, classified by the $G$-functor $\sigma_x \colon \ul{G/H} \to \ul\C$.
Define the parametrized over-category $\ul\C_{/\ul{x}} \fib \ul{G/H}$  (see \cite[not. 4.29]{Expose2}) as the fiber product $\Arr_G(\ul\C) \times_{\ul\C} \ul{G/H}$, considered as a $\ul{G/H}$-category by pulling back the coCartesian fibration $ev_1 \colon \Arr_G(\ul\C) \to \ul\C$ along $\sigma_x \colon  \ul{G/H} \to \ul\C$. 
Note that the fiber of $\ul\C_{/\ul{x}}\fib \ul{G/H}$ over $\varphi \colon G/K \to G/H$ is equivalent to the $\infty$-over-category $(\ul\C_{[G/K]})_{/\varphi^* x}$, where $\varphi^* x \in \ul\C_{[G/K]}$ is determined by choosing a coCartesian lift $ x \to \varphi^* x$ of $\varphi$. 

If $\ul\C' \subseteq \ul\C$ is a full $G$-subcategory we abuse notation and write $\ul\C'_{/\ul{x}}$ for the restricted $G$-over-category, given by the fiber product \( \ul\C' \times_{\ul\C} \ul\C_{/\ul{x}} \).

We now return to the definition of genuine $G$-factorization homology.

Let $A\in\Fun^\otimes_{G}(\ul\disk^{G,f-fr},\ul\C)$ be an $f$-framed $G$-disk algebra with values in $\ul\C$, and \( M \in \ulmfld^{G,f-fr}_{[G/H]} \) an $f$-framed $\OG$-manifold. 
In the following definition we use the parametrized over-category $\ul\disk^{G,f-fr}_{/\ul{M}}$ associated to $M\in \ulmfld^{G,f-fr}_{[G/H]}$ and $\ul\disk^{G,f-fr} \subset \ulmfld^{G,f-fr}$.

\begin{rem}
  Note that $\Gdisk_{/\ul{M}} \to \ul{G/H}$ is the coCartesian fibration dual to the Cartesian fibration \( (\disk^G)_{/M} \to (\OG)_{/[G/H]} \) (see \cite[prop 2.4.3.1]{HTT}, compare \cite[prop. 4.31]{Expose2}), and therefore can be modeled by the topological Moore over category (see \cref{MooreOverCat}). 
\end{rem}

Construct a $G$-functor over $\ul{G/H}$ by composing
\begin{align}
  \diag{
    \ul\disk^{G,f-fr}_{/\ul{M}} \ar@{->>}[dr] \ar[r] & \ul\disk^{G,f-fr} \ul\times \ul{G/H} \ar@{->>}[d] \ar[rr]^{A \ul\times id} & & \ul\C \ul\times \ul{G/H} \ar@{->>}[dll] \\
      & \ul{G/H} .
  }
  \label{G_disk_diag}
\end{align}
Consider the functor \eqref{G_disk_diag} as an $\ul{G/H}$-diagram in the $\ul{G/H}$-category \( \ul\C \ul\times \ul{G/H} \). 
Note that the $\ul{G/H}$-colimit of the above diagram is a coCartesian section of \( \ul\C \ul\times \ul{G/H} \fib \ul{G/H} \), or equivalently a $G$-functor $\ul{G/H} \to \ul\C$ and that a $G$-functor $\ul{G/H} \to \ul\C$ represents an object of $\ul\C$ over $[G/H]$. 
\begin{mydef} \label{def:G_FH_as_colimit}
  Let $M\in\ulmfld^{G,f-fr}_{[G/H]}$ be an $f$-framed $\OG$-manifold, and $A$ an $f$-framed $G$-disk algebra. 
  Define the \emph{$G$-factorization homology} of $M$ with coefficients in $A$ by the parametrized colimit
  \begin{align} \label{G_FH_colimit_formula}
    \int_M A \in \ul\C,\quad \int_M A := \ul{G/H}-\colim \left( \ul\disk^{G,f-fr}_{/\ul{M}} \to \ul\disk^{G,f-fr} \ul\times \ul{G/H} \xto{A\ul\times id} \ul\C \ul\times \ul{G/H} \right).
  \end{align}
\end{mydef}

In what follows, assume that $\ul\C$ is a $G$-cocomplete $G$-category (i.e $\ul\C$ has all $\ul{G/H}$-colimits for every $H<G$, see \cite[def. 5.12]{Expose2}), so that all the parametrized colimits of \cref{G_FH_functor} exist.
Next we show that 
the assignment \( M \mapsto \int_M A \) extends to a $G$-functor \( \int_- A  \colon  \ulmfld^{G,f-fr} \to \ul\C \), and that the $G$-functors $\int_- A$ are in turn functorial in $A$ (\cref{G_FH_functor}).
\begin{construction}
  Let \( \iota \colon  \ul\disk^{G,f-fr} \hookrightarrow \ulmfld^{G,f-fr} \) denote the inclusion of the full $G$-subcategory of finite $G$-disjoint unions of $G$-disks and $\ul\C$ be a cocomplete $G$-symmetric monoidal category. 
  The inclusion $G$-functor $\iota$ induces a restriction $G$-functor \( \iota^*  \colon  \ulFun_{G} (\ulmfld^{G,f-fr}, \ul\C) \to \ulFun_{G}(\ul\disk^{G,f-fr}, \ul\C) \).
  By \cite[cor. 10.6]{Expose2} (\cref{Shah:fully_faithful_G_Lan}) the restriction $G$-functor has a fully faithful left $G$-adjoint
  \begin{align*}
    \iota_! \colon \ulFun_{G}(\ul\disk^{G,f-fr}, \ul\C) \adj \ulFun_{G}(\ulmfld^{G,f-fr},\ul\C) \noloc \iota^*  .
  \end{align*}
  In particular, define $\iota_!$ to be the fully faithful left adjoint of  
  \begin{align} \label{iota_adj}
    \iota_! \colon \Fun_{G}(\ul\disk^{G,f-fr}, \ul\C) \adj \Fun_{G}(\ulmfld^{G,f-fr},\ul\C) \noloc \iota^*  ,
  \end{align}
  the adjunction of $\infty$-categories between the fibers over the terminal orbit $[G/G]$
\end{construction}

\begin{prop} \label{G_FH_functor}
  Let $\ul\C$ be a cocomplete $G$-symmetric monoidal category. 
  Then the functor 
  \begin{align*}
    \Fun^\otimes_G ( \ul\disk^{G,f-fr},\ul\C) \to \Fun_G (\ul\disk^{G,f-fr}, \ul\C) \xto{\iota_!} \Fun_G (\ulmfld^{G,f-fr}, \ul\C), \\
    (A \colon \ul\disk^{G,f-fr}D \to \ul\C^\otimes) \mapsto ( A \colon \ul\disk^{G,f-fr} \to \ul\C ) \mapsto ( \iota_! A  \colon  \ulmfld^{G,f-fr} \to \ul\C ) 
  \end{align*}
  sends a $G$-disk algebra $A$ to a $G$-functor
  \begin{align*}
     \iota_! A  \colon  \ulmfld^{G,f-fr} \to \ul\C ,\quad,  M \mapsto (\iota_! A)(M) = \int_M A.
  \end{align*}
\end{prop}
\begin{proof}
  By \cite[thm. 10.4]{Expose2} for every $G$-disk algebra $A\in \Fun_G (\ulmfld^{G,f-fr},\ul\C)$ the left $G$-adjoint $\iota_!(A) \colon \ulmfld^{G,f-fr} \to \ul\C) $ is given by left $G$-Kan extension of $A$ along $\iota$. 
  By \cite[thm 10.3]{Expose2} applying the $\iota_!(A)$ to $M \in \ulmfld^{G,f-fr}_{[G/H]}$ is given by the $\ul{G/H}$-colimit 
  \begin{align*}
    (\iota_! A)(M) = \ul{G/H}-\colim \left( \ul\disk^{G,f-fr}_{/\ul{M}} \to \ul\disk^{G,f-fr} \ul\times \ul{H/G} \xto{A\ul\times id} \ul\C \ul\times \ul{G/H} \right) = \int_M A .
  \end{align*}
\end{proof}

\subsection{Extensding $G$-factorization homology to a $G$-symmetric monoidal functor } \label{sec:G_FA_G_SM_functor}
In this subsection we prove (\cref{FH_G_SM}) that $G$-factorization homology with values in a presentable $G$-symmetric monoidal category extends to a $G$-symmetric monoidal functor (see \cref{def:G_FH_as_G_SM_functor}). 

\begin{mydef} \label{def:G_SM_presentable_cat}
  Let \( \ul\C^\otimes \fib \GFin_* \) be a $G$-symmetric monoidal category. 
  We say $\ul\C^\otimes$ is \emph{a presentable $G$-symmetric monoidal category} if the underlying $G$-category is presentable and for every active map $\alpha  \colon  I \to J$ in $\GFin_*$ the $G$-functor \( \otimes_\alpha  \colon  \ul\C^\otimes_{<I>} \to \ul\C^\otimes_{<J>} \) is distributive (\cite[sec. 3.3]{Nardin_thesis}).
\end{mydef}

\begin{prop} \label{FH_G_SM}
  Let \( \ul\C^\otimes \fib \GFin_* \) be a presentable $G$-symmetric monoidal category.

  Then the adjunction \cref{iota_adj} lifts to an adjunction 
  \begin{align}
    \label{operadic_Lan_extends_Lan}
    \diag{
      (\iota^\otimes)_! \colon \Fun^\otimes_G (\ul\disk^{G,f-fr},\ul\C)  \ar[d] \ar@/_/[r] &   \Fun^\otimes_G (\ulmfld^{G,f-fr}, \ul\C) \noloc (\iota^\otimes)^\ast \ar@/_/[l] \ar[d]  \\
      \iota_!  \colon \Fun_{G} (\ul\disk^{G,f-fr}, \ul\C) \ar@/_/[r] & \Fun_{G} (\ulmfld^{G,f-fr}, \ul\C ) \noloc \iota^\ast \ar@/_/[l] 
    }
  \end{align}
  where \( \iota^\otimes \colon \ul\disk^{G,f-fr,\sqcup} \to \ulmfld^{G,f-fr.\sqcup} \) is the inclusion of the subcategory of $f$-framed indexed disks (see \cref{GmfldD_framed}).
\end{prop}
Note that since $\iota_!$ is fully faithful the Segal conditions imply that \( (\iota^\otimes)_! \colon \Fun^\otimes_G (\ul\disk^{G,f-fr},\ul\C) \to \Fun^\otimes_G (\ulmfld^{G,f-fr}, \ul\C) \) is fully faithful.

\begin{mydef} \label{def:G_FH_as_G_SM_functor}
  For \( \ul\C^\otimes \fib \GFin_*, \, A \colon \ul\disk^{G,f-fr,\sqcup} \to \ul\C^\otimes \) as in \cref{FH_G_SM}, denote the $G$-symmetric monoidal functor $(\iota^\otimes)_!$ by 
  \begin{align*} 
    \diag{
      \int_- A \colon \ulmfld^{G,f-fr.\sqcup} \ar[dr] \ar[rr] & & \ul\C^\otimes \ar[dl] \\
        & \GFin_* .
    }
  \end{align*}
  Commutativity of the diagram \eqref{operadic_Lan_extends_Lan} shows that $\int_- A$ extends the $G$-functor \( \iota_! A \colon \ulmfld^{G,f-fr} \to \ul\C \) of \cref{iota_adj} which sends an $\OG$-manifold $(M\to G/H)$ to its $G$-factorization homology \( \iota_! A (M) = \int_M A \) (\cref{G_FH_functor}) to a $G$-symmetric monoidal functor.
  We call the $G$-symmetric monoidal functor $\int_- A \colon \ulmfld^{G,f-fr,\sqcup} \to \ul\C^\otimes$ the \emph{$G$-factorization homology functor with coefficients in $A$}.
\end{mydef}

In the remainder of this subsection we prove \cref{FH_G_SM}.
The proof has two parts, the first is a general $G$-categorical lemma, \cref{G_Lan_extends_to_G_SM}, giving conditions ensuring that a $G$-left Kan extension lifts to a $G$-symmetric monoidal functor, and the second is a verification of these conditions. 

We start with by recalling the notion of a $G$-lax monoidal functor and stating a useful proposition from \cite{Parametrized_algebra}.
Let $\ul\D^\otimes, \, \ul\C^\otimes$ be $G$-symmetric monoidal categories.
Recall that a lax $G$-symmetric monoidal $G$-functor $F$ from $\ul\D$ to $\ul\C$ is a functor \( F \colon \ul\D^\otimes \to \ul\C^\otimes \) over $\GFin_*$ which preserves inert edges (i.e. coCartesian edges over inert morphisms). 
Let \( Alg(\ul\D,\ul\C) \subset \Fun_{/\GFin_*} ( \ul\D^\otimes , \ul\C^\otimes) \) be the full subcategory of functors over $\GFin_*$ which are lax $G$-symmetric monoidal.
\begin{prop} \label{lem:Operadic_GLan_restrict_to_GLan}
  Let $\ul\C^\otimes \fib \GFin_*$ be a presentable $G$-symmetric monoidal category, let \( \ul\M^\otimes \fib \GFin_* \) be a small $G$-symmetric monoidal category and \( \iota^\otimes  \colon \ul\D^\otimes \cof \ul\M^\otimes \) an inclusion of a full $G$-symmetric monoidal subcategory. 
  Denote by $\iota \colon \ul\D \to \ul\M$ the induced $G$-functor on the underlying categories.

  Then the restriction along $ \iota^\otimes$ has a left adjoint $(\iota^\otimes)_! \colon Alg(\ul\D,\ul\C) \to Alg(\ul\M,\ul\C)$.
  Moreover, the adjunction \( (\iota^\otimes)_! \colon Alg(\ul\D,\ul\C) \adj Alg(\ul\M,\ul\C) \noloc (\iota^\otimes)^* \) restricts to the adjunction \( \iota_! \colon \Fun(\ul\D,\ul\C) \adj \Fun(\ul\M,\ul\C) \noloc \iota^* \), 
  where $ \iota_! \colon \Fun(\ul\D,\ul\C) \to \Fun(\ul\M,\ul\C)$ is left adjoint to the restriction along $\iota$. 

  In particular we have a commuting square of $\infty$-categories
  \begin{align*}
    \diag{
       Alg (\ul\D,\ul\C)  \ar[d] \ar[r]^{(\iota^\otimes)_! } &  Alg (\ul\M, \ul\C) \ar[d] \\ 
       \Fun_{G} (\ul\D, \ul\C) \ar[r]^{\iota_!} & \Fun_{G} (\ul\M, \ul\C ) . 
    }
  \end{align*}
\end{prop}

We will prove \cref{FH_G_SM} by applying the following $G$-categorical lemma (a $G$-categorical version of \cite[lem. 2.16]{AyalaFrancisTanaka}).
\begin{lem} \label{G_Lan_extends_to_G_SM}
  Let $\ul\C^\otimes \fib \GFin_*$ be a presentable $G$-symmetric monoidal category, let \( \ul\D^\otimes, \ul\M^\otimes \fib \GFin_* \) be small $G$-symmetric monoidal categories and \( \iota^\otimes \colon \ul\D^\otimes \cof \ul\M^\otimes \) 
  be an inclusion of a full $G$-symmetric monoidal subcategory. 
  Denote by $\iota \colon \ul\D \to \ul\M$ the induced $G$-functor on the underlying categories.
  
  If for every active morphism $\psi \colon I \to J$ in a fiber $(\GFin_*)_{[G/H]}$ and every coCartesian lift $x \to y$ of $\psi$ to $\ul\M^\otimes$ the $\ul{G/H}$-functor 
  \( \otimes_\psi \colon  (\ul\D^\otimes_{<I>})_{/\ul{x}} \to (\ul\D^\otimes_{<J>})_{/\ul{y}} \)
  is $\ul{G/H}$-cofinal then the diagram
  \begin{align*}
    \diag{
      \Fun^\otimes_G (\ul\D,\ul\C)  \ar[d] \ar[r]^{(\iota^\otimes)_!} &   \Fun^\otimes_G (\ul\M, \ul\C)  \ar[d]  \\
      \Fun_{G} (\ul\D,\ul\C) \ar[r]^{\iota_! } & \Fun_{G} (\ul\M, \ul\C )  
    }
  \end{align*}
  commutes, where 
  $(\iota^\otimes)_!$ and  $\iota_!$ the left adjoins to the restrictions along $\iota^\otimes$ and $\iota$, respectively.
\end{lem}
\begin{proof}
  Applying \cref{lem:Operadic_GLan_restrict_to_GLan} we have:
  \begin{align*}
    \diag{
      (\iota^\otimes)_!  \colon  Alg (\ul\D,\ul\C)  \ar[d] \ar@/_/[r] &  Alg (\ul\M, \ul\C) \noloc (\iota^\otimes)^\ast \ar@/_/[l] \ar[d]  \\
      \iota_! \colon \Fun_{G} (\ul\D, \ul\C) \ar@/_/[r] & \Fun_{G} (\ul\M, \ul\C ) \noloc \iota^\ast \ar@/_/[l] 
    }
  \end{align*}
  We need to show that the adjunction \( (\iota^\otimes)_! \colon Alg (\ul\D,\ul\C) \adj Alg(\ul\M,\ul\C) \noloc (\iota^\otimes)^* \) restricts to an adjunction between the full subcategories 
  \[ \Fun^\otimes_{G} ( \ul\D,\ul\C) \subset Alg(\ul\D, \ul\C) ,\quad  \Fun^\otimes_{G} ( \ul\M,\ul\C) \subset Alg(\ul\M, \ul\C) . \]
  Clearly precomposition with the $G$-symmetric monoidal functor $\iota^\otimes \colon \ul\D^\otimes \to \ul\M^\otimes$ takes $G$-symmetric monoidal functors to $G$-symmetric monoidal functors, so the right adjoint restricts to a functor 
  \[
    (\iota^\otimes)^* \colon \Fun^\otimes_{G}(\ul\M,\ul\C) \to \Fun^\otimes_{G} (\ul\D,\ul\C) . 
  \]

  Let $F^\otimes \colon  \ul\D^\otimes \to \ul\C^\otimes$ be a $G$-symmetric monoidal functor, with $F \colon  C\to D$ the induced $G$-functor on the underlying categories. 
  Applying the left adjoint $(\iota^\otimes)_!$ to $F^\otimes$ we get a lax $G$-symmetric monoidal functor \( (\iota^\otimes)_! F^\otimes  \colon  \ul\M^\otimes \to \ul\C^\otimes \), in other words $(\iota^\otimes)_! F^\otimes$ preserves coCartesian edges over inert morphisms. 
  We have to show that \( (\iota^\otimes)_! F^\otimes \) preserves all coCartesian edges. 
  Using the inert-fiberwise active factorization system on $M^\otimes$ (which exists on any $G$-$\infty$-operad, see \cite{Parametrized_algebra}), we are reduced to showing that \( (\iota^\otimes)_! F^\otimes \) preserves fiberwise active coCartesian edges.
  By the Segal conditions it is enough to show  \( (\iota^\otimes)_! F^\otimes \) preserves arrows over maps $I \to J$ in $\GFin_*$ with $J=(G/H \xto{=} G/H)$. 

  Before showing that $(\iota^\otimes)_! F^\otimes$ preserve these coCartesian edges, let us first recall how the functor $(\iota^\otimes)_! F^\otimes$ acts on morphisms.

  By definition $\iota_! \colon  \ul\M\to \ul\C$ is a left $G$-Kan extension.
  Using the construction of \cite[def. 10.1]{Expose2} we have a $G$-functor 
  \[
    (\ul\D \times_{\ul\M} \Arr_G(\ul\M) ) \star_{\ul\M} \ul\M \to \ul\C 
  \]
  which is an $\ul\M$-parametrized $G$-colimit diagram, where 
  \[
    \Arr_G(\ul\M) = \OGop \times_{\Fun(\ul\Delta^1,\OGop) } \Fun(\ul\Delta^1, \ul\M) \simeq \ulFun_G(\OGop \times \ul\Delta^1, \ul\M) 
  \]
  is the fiberwise arrow category (see \cite[not. 4.29]{Expose2}). 
  Note that by definition the restriction to the first coordinate  
  $ \ul\D \times_{\ul\M} \Arr_G(\ul\M) \to \left( \ul\D \times_{\ul\M} \Arr_G(\ul\M) \right) \star_{\ul\M} \ul\M \to \ul\M$ factors as \( \ul\D \times_{\ul\M} \Arr_G(\ul\M) \xto{\pi_{\ul\D}} \ul\D \xto{F} \ul\C \).
  and the restriction to the second coordinate is the left $G$-Kan extension functor $\iota_! F$, i.e $ \iota_! F \colon  \ul\M \to (\ul\D \times_{\ul\M} \Arr_G(\ul\M) ) \star_{\ul\M} \ul\M \to \ul\C $.

  Let  $x\in \ul\M^\otimes$ be an object over $I= (U \to G/H)$ and $\psi \colon I \to J$ be an active morphism in the fiber $(\GFin_*)_{[G/H]}$ with target $J=( G/H \xto{=} G/H )$, given by the span
  \begin{align*}
    \psi = \left( \vcenter{
      \xymatrix{
       U \ar[d]^f & \ar[l]_{=} U \ar[d]^f \ar[r]^{f} & G/H \ar[d]^{=} \\
       G/H & \ar[l]_{=} G/H \ar[r]^{=} & G/H 
     }
   } \right) .
  \end{align*}

  Denote the $G$-functor classified by $x$ by $x_\bullet \colon  \ul{U} \to \ul\M$  (see \cref{rem:object_as_family}).
  Pulling back the coCartesian fibration \( (\ul\D \times_{\ul\M} \Arr_G(\ul\M) ) \star_{\ul\M} \ul\M  \fib \ul\M \) along $x_\bullet$ we get a $\ul{U}$-parametrized $G$-colimit diagram 
  \( (\ul\D \times_{\ul\M} \Arr_G(\ul\M) \times_{\ul\M} \ul{U} ) \star_{\U} \U \to (\ul\D \times_{\ul\M} \Arr_G(\ul\M) ) \star_{\ul\M} \ul\M \to  \ul\C \) (implicitly using \cite[lem. 4.4]{Expose2}),
  and therefore a $\ul{U}$-colimit diagram 
  \[ 
    \overline{p} \colon  (\ul\D \times_{\ul\M} \Arr_G(\ul\M) \times_{\ul\M} \ul{U} ) \star_{\U} \U \to \ul\C \ultimes \ul{U} .
  \]
  Denote the $\ul{U}$-category indexing the colimit diagram above by \( \ul\D_{/\ul{x_\bullet}} := \ul\D \times_{\ul\M} \Arr_G(\ul\M) \times_{\ul\M} \ul{U}\).
  Note that by definition the restriction of \( \overline{p}\) to $ \ul\D_{/\ul{x_\bullet}} $ factors as the $\ul{U}$-functor 
  \[
    \ul\D_{/\ul{x_\bullet}} \to \ul\D \ultimes \ul{U} \xto{F\ultimes \ul{U}} \ul\C \ultimes \ul{U} 
  \]
  and the restriction to $\ul{U}$ is the $\ul{U}$-functor 
  \[
    \iota_! F(x_\bullet) \colon  \ul{U} \xto{x_\bullet \ultimes \ul{U} } \ul\M \ultimes \ul{U} \xto{ \iota_! F \ultimes \ul{U}} \ul\C \ultimes \ul{U} .
  \]
  
  Since $\ul\C^\otimes$ is a presentable $G$-symmetric monoidal category the tensor product functor 
  \[ 
    \otimes_\psi \colon  \prod_I \ul\C \ultimes \ul{U} \to \ul\C \ultimes \ul{G/H}
  \]
  of \cref{def:param_tensor_prod} is a distributive $\ul{G/H}$-functor (see \cite[def. 3.15]{Nardin_thesis}).
  Therefore the $\ul{U}$-colimit diagram $\overline{p}$ induces a $\ul{G/H}$-colimit diagram 
  \begin{align} \label{diag:GLAN_of_tensor_prod}
    \left( \prod\limits_I \ul\D_{/\ul{x_\bullet}} \right) \star_{\ul{G/H}} \ul{G/H} \to 
    \prod\limits_I \left( \ul\D_{/\ul{x_\bullet}} \star_{\ul{U}} \ul{U} \right) \xto{\prod\limits_\psi \overline{p} } 
    \prod\limits_I \left( \ul\C \ultimes \ul{U} \right) \xto{\otimes_\psi} 
    \ul\C \ultimes \ul{G/H} 
  \end{align}
  exhibiting the $\ul{G/H}$-object 
  \[
    \otimes_\psi \left(\prod_I \iota_! F ( x_\bullet) \right)  \colon  \ul{G/H} \xto{\simeq} \prod_I \ul{U}  \xto{\prod_I \iota_!F (x_\bullet)} \prod_I \ul\C \ultimes \ul{U} \xto{\otimes_\psi} \ul\C \ultimes \ul{G/H} 
  \]
  as the $\ul{G/H}$-colimit of 
  \[
    p \colon \prod\limits_I \ul\D_{/\ul{x_\bullet}} \to \prod\limits_I \ul\D \ultimes \ul{U} \to \prod\limits_I \ul\C \ultimes \ul{U}  \xto{\otimes_\psi} \ul\C \ultimes \ul{G/H} .
  \]

  First, note that we can express the $\ul{G/H}$-colimit \(\otimes_\psi \left(\prod_I \iota_! F ( x_\bullet) \right) \) of \eqref{diag:GLAN_of_tensor_prod} in simpler terms.
  Since $(\iota^\otimes)_! F^\otimes \colon  D^\otimes \to \ul\C^\otimes$ is a lax $G$-symmetric monoidal functor we have 
  \begin{align*}
    \diag{
      \ul{G/H} \ar[d]^{\simeq} \ar[rr]^{x} \ar@/^20pt/[rrrr]^{(\iota^\otimes)_! F^\otimes (x)} & & \ul\M^\otimes_{<I>} \ar[d]^{\simeq} \ar[rr]^{ ( (\iota^\otimes)_! F^\otimes )_{<I>} } & & \ul\C^\otimes_{<I>} \ar[d]^{\simeq} \\
      \prod_I \ul{U} \ar[rr]^-{ \prod_I x_\bullet \ultimes \ul{U}} \ar@/_10pt/[rrrr]_{\prod_I \iota_! F (x_\bullet) } & & \prod_I \ul\M \ultimes \ul{U} \ar[rr]^-{\prod_I \iota_! F \ultimes \ul{U}} & & \prod_I C \ultimes \ul{U} ,
    }
  \end{align*}
therefore \(\otimes_\psi \left(\prod_I \iota_! F ( x_\bullet) \right) \simeq \otimes_\psi \left( (\iota^\otimes)_! F^\otimes (x) \right) \).

  On the other hand, we can also express the $\ul{G/H}$-diagram $p$ in simpler terms. 
  To see this observe the commutative diagram
  \begin{align} \label{diag:tensor_first}
    \diag{
      (\ul\D^\otimes_{<I>})_{/\ul{x}} \ar[d]^{\simeq} \ar[r] & \ul\D^\otimes_{<I>} \ar[d]^{\simeq} \ar[r]^{F^\otimes_{<I>}} & \ul\C^\otimes_{<I>} \ar[d]^{\simeq} \\
      \prod\limits_I \ul\D_{/\ul{x_\bullet}} \ar[d]^{\otimes_\psi} \ar[r] & \prod\limits_I \ul\D \ultimes \U  \ar[d]^{\otimes_\psi} \ar[r]^{\prod_I F \ultimes \ul{U}} & \prod\limits_I \ul\C \ultimes \U  \ar[d]^{\otimes_\psi} \\
      \ul\D_{/\ul{\otimes_\psi x}} \= \left( \ul\D \ultimes \ul{G/H} \right)_{/\ul{\otimes_\psi x}} \ar[r] & \ul\D \ultimes \ul{G/H}  \ar[r]^{ F \ultimes \ul{G/H}} & \ul\C \ultimes \ul{G/H}
    }
  \end{align}
  where the left vertical column is induced by taking the $\ul{G/H}$-limit of the rows of the following diagram of $\ul{G/H}$-categories:
  \begin{align*}
    \diag{
      \ul\D^\otimes_{<I>} \ar[d]^{\simeq} \ar[r] & \ul\M^\otimes_{<I>} \ar[d]^{\simeq} & \ar[l] \Arr_{\ul{G/H}} ( \ul\M^\otimes_{<I>} )  \ar[d]^{\simeq} \ar[r] & \ul\M^\otimes_{<I>} \ar[d]^{\simeq} & \ar[l]_{x} \ul{G/H}  \ar[d]^{\simeq} \\
      \prod_I \ul\D \ultimes \ul{U} \ar[r] \ar[d]^{\otimes_\psi} & \prod_I \ul\M \ultimes \ul{U} \ar[d]^{\otimes_\psi} & \ar[l] \Arr_{\ul{G/H}} (\prod_I \ul\M \ultimes \ul{U}) \ar[d]^{\otimes_\psi} \ar[r] & \prod_I \ul\M \ultimes \ul{U} \ar[d]^{\otimes_\psi} & \ar[l]_-{\prod_I x_\bullet} \prod_I \ul{U} \ar[d]^{\simeq} \\
      \ul\D \ultimes \ul{G/H} \ar[r] & \ul\M \ultimes \ul{G/H} & \ar[l] \Arr_{\ul{G/H}} ( \ul\M \ultimes \ul{G/H} ) \ar[r] & \ul\M \ultimes \ul{G/H} & \ar[l]_-{\otimes_\psi x} \ul{G/H}.
    }
  \end{align*}
  Note that $p$ is the composition of the middle row of diagram \eqref{diag:tensor_first} followed by the lower right vertical $\ul{G/H}$-functor $\otimes_\psi$. 
  Therefore $p$ is equivalent to the composition of the left vertical column of diagram \eqref{diag:tensor_first} followed by the bottom row: 
  \begin{align*}
    (\ul\D^\otimes_{<I>})_{/\ul{x}} \xto{\otimes_\psi} \ul\D_{/\ul{\otimes_\psi x}} \to \ul\D \ultimes \ul{G/H} \xto{ F \ultimes \ul{G/H} }  \ul\C \ultimes \ul{G/H}.
  \end{align*}

  Finally, by the assumption of the lemma the $\ul{G/H}$-functor \( \otimes_\psi \colon   (\ul\D^\otimes_{<I>})_{/\ul{x}} \to \left( \ul\D \ultimes \ul{G/H} \right)_{/\ul{\otimes_\psi x}} \) is $\ul{G/H}$-cofinal, therefore
  \begin{align*}
    \otimes_\psi \left( (\iota^\otimes)_! F^\otimes (x) \right) \simeq
    \ul{G/H}-\colim \left( (\ul\D^\otimes_{<I>})_{/\ul{x}} \xto{\otimes_\psi} \ul\D_{/\ul{\otimes_\psi x}} \to \ul\D \ultimes \ul{G/H} \xto{ F \ultimes \ul{G/H} }  \ul\C \ultimes \ul{G/H} \right) \\
    \iso \ul{G/H}-\colim \left( \ul\D_{/\ul{\otimes_\psi x}} \to \ul\D \ultimes \ul{G/H} \xto{ F \ultimes \ul{G/H} }  \ul\C \ultimes \ul{G/H}  \right) \simeq \iota_! F (\otimes_\psi x),
  \end{align*}
  so we have a coCartesian edge \( e \colon  (\iota^\otimes)_! F^\otimes (x) \to \iota_! F (\otimes_\psi x) \) in $\ul\C^\otimes$ over $\psi$.

  We can now show that \( (\iota^\otimes)_! F^\otimes  \colon  \ul\M^\otimes \to \ul\C^\otimes \) preserves coCartesian edges over $\psi \colon I \to J$ as above. 
  Let $e' \colon x\to y$ be a coCartesian edge in $\ul\M^\otimes$ over $\psi$.
  By definition of $\otimes_\psi$ this coCartesian edge factors as \( x \to \otimes_\psi x \iso y \) over \( I \xto{\psi} J \xto{=} J \) (See \cite[rem. 2.4.1.4 and prop. 2.4.1.5]{HTT}). 
  Applying $(\iota^\otimes)_! F^\otimes$ we get
  \( (\iota^\otimes)_! F^\otimes (e') \colon  (\iota^\otimes)_! F^\otimes (x) \to  (\iota^\otimes)_! F^\otimes(y) = \iota_! F(y) \), 
  and we need to show $ (\iota^\otimes)_! F^\otimes (e')$ is a coCartesian lift of $\psi$.
  However, we already have a coCartesian lift of $\psi$, the edge $e$ we constructed above. 
  Therefore $ (\iota^\otimes)_! F^\otimes (e')$ factors through $e$ as \( (\iota^\otimes)_! F^\otimes (x) \xto{e} \iota_! F (\otimes_\psi x) \to \iota_! F(y) \). 
  Note that the morphism \( \iota_! F (\otimes_\psi x) \to \iota_! F(y) \) is induced from \( \otimes_\psi x \iso y \), and therefore an equivalence. 
  Hence $ (\iota^\otimes)_! F^\otimes (e')$ is coCartesian as a composition of a coCartesian edge and an equivalence. 

  This ends the proof of \cref{G_Lan_extends_to_G_SM}.
\end{proof}

We can now prove \cref{FH_G_SM} by verifying the cofinality conditions of \cref{G_Lan_extends_to_G_SM}. 
In fact, we prove the cofinality of the maps in \cref{G_Lan_extends_to_G_SM} by showing that they are equivalences. 

We rely on the following result to reduce our calculations to the non-framed case $B=BO_n(G)$.
\begin{prop} \label{Gmfld_framed_over_M}
  Let $f \colon B\to BO_n(G)$ be a $G$-map as in \cref{def:Gmfld_framed}, and $M \in \ulmfld^{G,f-fr}_{[G/H]}$ an $f$-framed $\OG$-manifold over $G/H$.
  Then the $\ul{G/H}$-functor 
  \begin{align*}
    ( \ulmfld^{G,f-fr} )_{/\ul{M}} \to \Gmfld_{/\ul{M}}
  \end{align*}
  is an equivalence of $\ul{G/H}$-categories. 

  In particular, every $G$-submanifold $N\subseteq M$ has an essentially unique lift to $N\in \ulmfld^{G,f-fr}_{/\ul{M}}$ (informally, $M$ induces an $f$-framing of $N$).
\end{prop}
\begin{proof}
  We show that for every $\varphi\in \ul{G/H}, \, \varphi \colon G/K \to G/H$ the induced functor on the fibers over $\varphi$,
  \(    \left( (\ulmfld^{G,f-fr})_{/\ul{M}} \right)_{[\varphi]} \to \left( \Gmfld_{/\ul{M}} \right)_{[\varphi]} \),
  is an equivalence.
  By construction the fibers of the parametrized over category are equivalent to the over categories
  \begin{align*}
    \left( (\ulmfld^{G,f-fr})_{/\ul{M}} \right)_{[\varphi]} \simeq \left( \ulmfld^{G,f-fr}_{[G/K]} \right)_{/\varphi^* M} , \quad 
    \left( (\ulmfld^{G})_{/\ul{M}} \right)_{[\varphi]} \simeq \left( \ulmfld^{G}_{[G/K]} \right)_{/\varphi^* M} .
  \end{align*}
  By \cref{def:Gmfld_framed} the fiber $\ulmfld^{G,f-fr}_{[G/K]}$ is given by the pullback of $\infty$-categories
  \begin{align*}
    \xymatrix{
      \ulmfld^{G,f-fr}_{[G/K]} \ar[d] \ar[r] \pullbackcorner & (\ulTopG_{[G/K]})_{/ \ul{B}(G/K)} \ar[d] \\
      \ulmfld^{G}_{[G/K]} \ar[r] & (\ulTopG_{[G/K]})_{/ \ul{BO_n(G)}(G/K)} ,
    }
  \end{align*}
  where 
  \[
    \ul{B}(G/K) = (B \times G/K \to G/K), \quad \ul{BO_n(G)}(G/K) = (BO_n(G) \times G/K \to G/K), 
  \]
  see \cref{ulTopG_as_G_spaces_over_orbits}.
  We can simplify the pullback square above using the equivalences of \cref{ulTopG_over_B_description}:
  \[
    (\ulTopG_{[G/K]})_{/ \ul{B}(G/K)} \iso \GTop_{/B \times G/K}, \quad
    (\ulTopG_{[G/K]})_{/ \ul{BO_n(G)}(G/K)} \iso \GTop_{/BO_n(G) \times G/K}.
  \]
  We can now express the slice category $(\ul\mfld^{G,f-fr}_{[G/K]})_{/\varphi^* M}$ as a pullback of slice categories
  \begin{align*}
    \xymatrix{
      (\ulmfld^{G,f-fr}_{[G/K]})_{/\varphi^* M} \ar[d] \ar[r] \pullbackcorner & \left( \GTop_{/B \times G/K} \right)_{/(\varphi^* M \to B \times G/K)} \ar[d] \\
      (\ulmfld^{G}_{[G/K]})_{/\varphi^* M} \ar[r] & \left( \GTop_{/BO_n(G) \times G/K} \right)_{/(\varphi^* M \to BO_n(G) \times G/K)} \ar[d] \\
      & \GTop_{/\varphi^* M} .
    }
  \end{align*}
  By \cite[lem. 2.5]{AF} both the bottom right vertical arrow and the composition of the right vertical arrows are equivalences of $\infty$-categories.
  By the two-out-of-three property we see that the top vertical arrow is an equivalence of $\infty$-categories, and therefore the left vertical arrow is also an equivalence, as claimed.
\end{proof}

\begin{proof}[Proof of \cref{FH_G_SM}]
  By the Segal conditions and \cref{Gmfld_framed_over_M} the $\ul{G/H}$-functors 
  \begin{align*}
    (\ul\disk^{G,f-fr,\sqcup}_{<I_i>})_{/\ul{M_i}}  \to (\GdiskD_{<I_i>})_{/\ul{M_i}}   , \quad i=1,2
  \end{align*}
  are equivalences of $\ul{G/H}$-categories, therefore it is enough to prove the non-framed case. 

  Let $\psi \colon I \to J$ be an active morphism in the fiber $(\GFin_*)_{[G/H]}$. 
  Without loss of generality, $\psi$ is represented by the span
  \begin{align*}
    \psi=\left(\diag{ 
      U_1 \ar[d] & \ar[l]_{=} U_1 \ar[d] \ar[r] & U_2 \ar[d] \\ 
      G/H & \ar[l]_{=} G/H \ar[r]^{=} & G/H 
    } \right).
  \end{align*}
  By \cref{GmfldD_description} a coCartesian lift $f \colon M_1 \to M_2$ of $\psi$ is represented by a span
  \begin{align*}
    f=\left(\diag{ 
      M_1 \ar[d] & \ar[l]_{=} M_1  \ar[d] \ar[r]^{\sim} & M_2\ar[d] \\ 
      U_1 \ar[d] & \ar[l]_{=} U_1 \ar[d] \ar[r] & U_2 \ar[d] \\ 
      G/H & \ar[l]_{=} G/H \ar[r]^{=} & G/H 
    } \right).
  \end{align*}
  By \cref{G_Lan_extends_to_G_SM} it is enough to show that the $\ul{G/H}$-functor \(  (\GdiskD_{<I>})_{/\ul{M_1}} \to (\GdiskD_{<J>})_{/\ul{M_2}} \) is $\ul{G/H}$-cofinal.
  We prove that it is in fact an equivalence, by showing that it is induces fiberwise equivalences. 
  
  Consider the induced functor \(  \left( (\GdiskD_{<I>})_{/\ul{M_1}} \right)_{[\varphi]}  \to \left( (\GdiskD_{<J>})_{/\ul{M_2}} \right)_{[\varphi]} \) between the fibers over 
  \( \varphi \in \ul{G/H}, \, \varphi \colon  G/K \to G/H \).

  We now inspect each fiber.
  By definition, we have 
  \begin{align*}
    (\GdiskD_{<I>})_{/\ul{M_i}} :=  \GdiskD_{<I>} \times_{\GmfldD_{<I>}} (\GmfldD_{<I>})_{/\ul{M_i}}, \quad i=1,2 ,
  \end{align*}
  and therefore the fibers over $\varphi$ are given by
  \begin{align*}
    \left( (\GdiskD_{<I>})_{/\ul{M_i}} \right)_{[\varphi]} = \left( \GdiskD_{<I>}  \right)_{[\varphi]} \times_{\left(\GmfldD_{<I>} \right)_{[\varphi]} } \left((\GmfldD_{<I>})_{/\ul{M_i}} \right)_{[\varphi]} , \quad i=1,2 .
  \end{align*}
  By the definition of parametrized slice category \cite[not. 4.29]{Expose2} we have 
  \begin{align*}
    \left((\GmfldD_{<I>})_{/\ul{M_i}} \right)_{[\varphi]} \= \left((\GmfldD_{<I>})_{[\varphi]} \right)_{/ \varphi^* M_i}, \quad i=1,2,
  \end{align*}
  where $\varphi^* M_i,\, i=1,2$ is the pullback of $M_i \to U_i \to G/H$ along $\varphi \colon  G/K \to G/H$.
  
  Next, note that \( (\GmfldD_{<I_i>})_{[\varphi]} \= (\GmfldD)_{\varphi^* I_i} \) is the fiber of \( \GmfldD \fib \GFin_* \) over \(\varphi^* I_i = ( U_i \times_{G/H} G/K  \to G/H ) \in \GFin_* \).

  However, using the definition of the coCartesian fibration $\GmfldD \fib \GFin_*$ (\cref{def:GmfldD_as_coCart_fib}) and the definition of the unfurling construction (see \cite[prop. 11.6]{BarwickSMF1} and the description of the fibers following it) we see that 
  \( (\GmfldD_{<I_i>})_{[\varphi]} \) 
  is equivalent to $(\mfldGD)_{\varphi^* I_i}$, 
  (the coherent nerve of) the full topological subcategory of $\OGFin$-manifolds, $\mfldGD$ (\cref{def:mfldGD_as_top_cat}) spanned by $\OGFin$-manifolds over $\varphi^* I_i$. 
  It follows that the $\infty$-category $\left((\GmfldD_{<I_i>})_{/\ul{M_i}} \right)_{[\varphi]}$ is equivalent to the slice category $\left( (\mfldGD)_{\varphi^* I_i} \right)_{/\varphi^* M_i}$, modeled by the coherent nerve of the Moore over category $\left( (\mfldGD)_{\varphi^* I_i} \right)^{\Moore}_{/\varphi^* M_i}$. 
  Therefore, the fiber \( \left( (\GdiskD_{<I_i>})_{/\ul{M_i}} \right)_{[\varphi]} \) is equivalent to the full subcategory of $\left( (\mfldGD)_{\varphi^* I_i} \right)^{\Moore}_{/\varphi^* M_i}$ spanned by objects represented by morphisms \( (E \to U' \to U_i \times_{G/H} G/K \to G/H ) \to ( \varphi^* M_i \to U_i \times_{G/H} G/K \to G/H ) \) over $\phi^* I_i$, where $E \to U'$ is a $G$-vector bundle. 
  Recall that \( U' = \pi_0(E) \) (\cref{rem:disk_bun_connected_components}).
  Unwinding the definition of morphisms in $\mfldGD$ 
  over $\varphi^* I_i = U_i \times_{G/H} G/K$, we see that such morphisms are represented by commutative diagrams
  \begin{align*}
    \diag{ 
      E \ar[d] & \ar[l]_{\sim} E'  \ar[d] \ar@{^(->}[r] & M_i \ar[dd] \\ 
      \pi_0(E) \ar[d] & \ar[l]_{=} \pi_0(E)  \ar[d]  &  \\ 
      \varphi^* I_i   \ar[d] & \ar[l]_{=} \varphi^* I_i \ar[d] \ar[r]^{=} & \varphi^* I_i \ar[d] \\ 
      G/K & \ar[l]_{=} G/K \ar[r]^{=} & G/K ,
    } 
  \end{align*}
  or equivalently, by a $G$-equivariant embedding $E \cof M_i$ over $\varphi^* I_i$. 

  With this concrete description of the fibers at hand, the induced functor between the fibers \(  \left( (\GdiskD_{<I_1>})_{/\ul{M_1}} \right)_{[\varphi]}  \to \left( (\GdiskD_{<I_2>})_{/\ul{M_2}} \right)_{[\varphi]} \)  is given by composition with 
  \begin{align*}
    \varphi^* f = \left(\diag{ 
      \varphi^* M_1 \ar[d] & \ar[l]_{=} \varphi^* M_1  \ar[d] \ar[r]^{\sim} & \varphi^* M_2 \ar[d] \\ 
      \varphi^* I_1 \ar[d] & \ar[l]_{=} \varphi^* I_1 \ar[d] \ar[r] & \varphi^* I_2 \ar[d] \\ 
      G/K & \ar[l]_{=} G/K \ar[r]^{=} & G/K
    } \right).
  \end{align*}
  By inspection the induced functor  
  \begin{align*}
    (\varphi^* f) \circ - \colon  \left( (\mfldGD)_{\varphi^* I_1} \right)^{\Moore}_{/\varphi^* M_1} & \to \left( (\mfldGD)_{\varphi^* I_2} \right)^{\Moore}_{/\varphi^* M_2} \\
    \left( \diag{
      E \ar[d] \ar@{^(->}[r] & \varphi^* M_1 \ar[dd] \\
      \pi_0(E) \ar[dr]  \\
        & \varphi^* I_1 \ar[d] \\
        & G/K
    } \right) & \mapsto 
     \left(\diag{ 
       E \ar[d] \ar@{^(->}[r] & \varphi^* M_1 \ar[dd] & \ar[l]_{=} \varphi^* M_1  \ar[dd] \ar[r]^{\sim} & \varphi^* M_2 \ar[dd] \\ 
      \pi_0(E) \ar[dr]  \\
      & \varphi^* I_1 \ar[d] & \ar[l]_{=} \varphi^* I_1 \ar[d] \ar[r] & \varphi^* I_2 \ar[d] \\ 
      & G/K & \ar[l]_{=} G/K \ar[r]^{=} & G/K
    } \right) 
  \end{align*}
  is an equivalence of topological categories.
\end{proof}

\section{Properties of $G$-factorization homology} \label{sec:G_FH_properties}
In this subsection we prove two properties of $G$-factorization homology: it satisfies $G$-$\otimes$-excision (\cref{G_tensor_excision}) and respects $G$-sequential colimits (\cref{G_FH_respects_seq_unions}).

\subsection{Collar decomposition of $G$-manifolds} \label{sec:G_collar_decomposition}
We define $G$-collar decompositions of $G$-manifolds and construct inverse image functors (\cref{inverse_image_functor}).
In the next subsection we use these constructions to define $G$-$\otimes$-excision and prove that $G$-factorization homology satisfies $G$-$\otimes$-excision (\cref{G_tensor_excision}).

We begin with an equivariant version of collar-gluing, see \cite[def. 3.13]{AF}.
The same definition is given in \cite[def. 4.20]{Weelinck}.
\begin{mydef}
  Let $M\in \mfld^G$ be an $n$-dimensional  $G$-manifold. 
  A \emph{$G$-collar decomposition} of $M$ is a smooth $G$-invariant function \(f \colon  M \to [-1,1] \) to the closed interval for which the restriction  \( f|_{(-1,1)} \colon M|_{(-1,1)} \to  (-1,1) \) is a manifold fiber bundle, with a choice of trivialization \( M|_{(-1,1)} \= M_0 \times (-1,1) \).
  Here $M_{(-1,1)} = f^{-1}(-1,1), \, M_0 = f^{-1}(0)$.
  For such a decomposition, denote \( M_+ := f^{-1} (-1,1], \, M_- := f^{-1}[-1,1)\). 

  A $G$-collar decomposition of an $f$-framed $\OG$-manifold $M\in \ulmfld^{G,f-fr}_{[G/H]}$ is a $G$-collar decomposition of the underlying $G$-manifold $M$.
\end{mydef}

\begin{rem}
  Note that $M|_{(-1,1)}$ is a tubular neighborhood of the codimension one $G$-submanifold $M_0$, and that $M_0$ splits $M$ into two $G$-manifolds, i.e. there exists a continuous $G$-invariant function \( M \setminus M_0 \to [-1,1] \setminus \set{0} \to  \set{-1,1} \) to the set with two elements.
  On the other hand, a $G$-submanifold $M_0 \subset M$ of codimension one that splits $M$ into two $G$-manifolds has an equivariant tubular neighbourhood $T\subset M$ equivalent to the total space $\nu(M_0)$ of the normal bundle of $M_0$ (compatible with the $M_0 \subset M$ and the zero section \( M_0 \to \nu(M)\) ). 
  By assumption, the normal bundle of $M_0$ is a trivial vector bundle of rank 1 with trivial $G$ action.
  A choice of $G$-diffeomorphisms $T \= \nu(M) \= M_0 \times \R \= M_0 \times (-1,1)$ (compatible with $M_0$) determines a $G$-collar decomposition of $M$. 
\end{rem}
\begin{rem}
  A $G$-collar decomposition $f \colon M \to [-1,1]$ defines a decomposition of $M$ into a union of open $G$-submanifolds \( M = M_- \cup M_+ \) with a chosen isomorphism \( M_- \cap M_+ \= M_0 \times (-1,1) \).
  The purpose of the above definition is to specify these decompositions among all decompositions $M=U \cup V$ of $M$ as a union of two open $G$-submanifolds. 
  We will see that $G$-equivariant homology is compatible with $G$-collar decompositions (\cref{def:G_tensor_excision} and \cref{G_tensor_excision}).
  This should be compared with Bredon homology, which is compatible with all decompositions $M=U \cup V$ into two equivariant open subsets (the equivariant Mayer-Vietoris property).
\end{rem}

Next we construct an ``inverse image'' functor \(f^{-1} \colon  \mfld^{\bnd,or}_{/[-1,1]} \to (\ulmfld^{G,f-fr}_{[G/H]} )_{/M} \) from the $\infty$-category of $1$-dimensional oriented manifolds with boundary over the interval $[-1,1]$ (see \cite{AF}). 
By \cref{Gmfld_framed_over_M} we have an equivalence of $\infty$-categories 
\begin{align} \label{Gmfld_framed_fiber_over_M}
  (\ulmfld^{G,f-fr}_{[G/H]} )_{/M} \iso (\Gmfld_{[G/H]} )_{/(M\to G/H)} ,
\end{align}
so it is enough to construct \(f^{-1} \colon  \mfld^{\bnd,or}_{/[-1,1]} \to (\Gmfld_{[G/H]} )_{/(M\to G/H)} \) for $(M \to G/H) \in \Gmfld_{[G/H]}$ the underlying $\OG$-manifold of $M \in \ulmfld^{G,f-fr}_{[G/H]}$.

Note that both the domain and the codomain of the functor $f^{-1}$ can be described using coherent nerve of the Moore over categories (see \cref{MooreOverCat}),
since the $\infty$-categories 
\[ 
  \mfld^{bnd, or}, \quad \Gmfld_{[G/H]} \= \mfld^G_{[G/H]} 
\]
are coherent nerves of topological categories.
We construct the inverse image functor as the coherent nerve of a functor of topological categories between the Moore over categories. 
\begin{construction} \label{inverse_image_functor}
  Let $(M\to G/H)$ be an $\OG$-manifold with a collar decomposition \( f \colon M \to [-1,1] \). 
  Define a topological functor \( (\mfld^{\bnd,or})^{\Moore}_{/[-1,1]} \to (\Gmfld_{[G/H]})^{\Moore}_{/(M\to G/H)} \) between the Moore over categories: 
  \begin{enumerate}
    \item Send an object of $(\mfld^{\bnd,or})^{\Moore}_{/[-1,1]}$ given by oriented embedding \( \varphi \colon  V \hookrightarrow [-1,1] \) to its inverse image, \( f^{-1}(\varphi)  \colon  f^{-1} V \hookrightarrow M \) given by the pullback of $\varphi$ along $f$.
      Since the function $f$ is $G$-invariant the embedding $f^{-1}(\varphi)$ is $G$-equivariant. 
      The composition 
      \begin{align*}
        \diag{ 
          f^{-1} V \ar@{^(->}[rr]^-{f^{-1}(\varphi)} & & M \ar[r] & G/H 
        }
      \end{align*}
      makes $f^{-1}V$ a $G$-manifold over $G/H$, hence \( f^{-1} (\varphi) \) is a point in the topological space 
      \( Emb^G_{G/H}( f^{-1} V, M) \), i.e an object of the Moore over category $(\Gmfld_{[G/H]})^{\Moore}_{/(M\to G/H)}$ . 
    \item Let \( \varphi \colon  V \hookrightarrow [-1,1]\) and \( \varphi' \colon  V' \hookrightarrow [-1,1]\) be two objects of the Moore over category \((\mfld^{\bnd,or})^{\Moore}_{/[-1,1]}\).
      Let $(h,(r,\gamma))$ be a point in $\Map_{(\mfld^{\bnd,or})^{\Moore}_{/[-1,1]}} (\varphi,\varphi')$, where $h \colon V \hookrightarrow V'$ is an oriented embedding and 
      $(r,\gamma) \in [0,\infty) \times \left( Emb^{\bnd,or}(V,[-1,1]) \right)^{[0,\infty)}$ is a Moore path from $\varphi$ to $\varphi'\circ h$.
      Define a continuous function
      \begin{align*}
        & f^{-1} \colon \Map_{(\mfld^{\bnd,or})^{\Moore}_{/[-1,1]}} (\varphi,\varphi') \to \Map_{(\Gmfld_{[G/H]})^{\Moore}_{/(M\to G/H)}} \left( f^{-1}(\im\varphi) \subset M, f^{-1}(\im\varphi') \subset M \right) , \\ 
        & f^{-1}(h,(r,\gamma)) := ( f^{-1}(h) ,(r,\alpha)), \quad 
        \diag{
          f^{-1} V \pullbackcorner \ar[d] \ar@{^(->}[r]^-{f^{-1}(h)} & f^{-1} \pullbackcorner \ar[d] \ar@{^(->}[r]^-{f^{-1}(\varphi')} & M \ar[d]^{f} \\
          V \ar@{^(->}[r]^{\varphi} & V \ar@{^(->}[r]^{h} & [-1,1]
        }
      \end{align*}
      where $f^{-1}(h)$ is given by the pullback 
      \begin{align*}
        \diag{
          f^{-1} V \pullbackcorner \ar[d] \ar@{^(->}[r]^-{f^{-1}(h)} & f^{-1} \pullbackcorner \ar[d] \ar@{^(->}[r]^-{f^{-1}(\varphi')} & M \ar[d]^{f} \\
          V \ar@{^(->}[r]^{\varphi} & V \ar@{^(->}[r]^{h} & [-1,1]
        }
      \end{align*}
      and \( \alpha \colon [0,\infty) \to   Emb^G(f^{-1}(V), M) \) is the Moore path of length $r$ defined as follows.
      If \( x\in M|_{(-1,1)} \= M_0 \times (-1,1) \) corresponds to \( (y,s) \in M_0 \times (-1,1) \) define 
      \[ \alpha_t (x)= (y,\gamma_t \circ \varphi^{-1}(s)) \in M_0 \times (-1,1) \= M|_{(-1,1)} , \] 
      otherwise (i.e. $f(x)= \pm 1$) define \( \alpha_t(x)=x \). 
      Verification that $(r,\alpha)$ is a smooth $G$-equivariant isotopy depending continuously on $\gamma$ is left to the reader.
  \end{enumerate}
  Clearly $f^{-1}$ preserve disjoint unions. 
\end{construction}

\begin{rem}
  More generally, one can try to define an inverse image functor along a general smooth invariant map $M \to N$ to a oriented manifold with boundary $N$. 
  However, not every map $f$ will do.
  First, in order to define the isotopy lift $\alpha$ assume that the restrictions of $f$ to  \( f^{-1}(N\setminus\bnd N) \) and $f^{-1}(\bnd N)$ are smooth fiber bundles, and use $G$-equivariant parallel transport between the fibers.
  The connections on the fiber bundles need to be compatible in order for $\alpha$ to be continuous and smooth.
  However, such parallel transport defines functions which are only continuous in the $C^1$-topology on \( Emb^G(f^{-1} V, M) \), since they depend on the time derivative of the isotopy $\gamma$. 
  Nevertheless, if the connections chosen are flat then parallel transport depends only on the end points, and therefore defines a continuous function relative to the compact-open topology. 
  All these conditions can be can be captured together by assuming that $f \colon M\to N$ is a $G$-invariant flat complete Riemannian submersion. 
  This condition implies that the restrictions to $N\setminus\bnd N$ and $\bnd N$ are flat fiber bundles, with compatibly chosen flat $G$-equivariant Ehresmann connections (i.e. a constructible fiber bundle relative to the boundary stratification). 
\end{rem}

\subsection{$G$-$\otimes$-excision} \label{sec:G_tensor_excision}
We define an equivariant version of $\otimes$-excision of \cite[def. 3.15]{AF} (see \cref{def:G_tensor_excision}), and prove it is satisfied by $G$-factorization homology (\cref{G_tensor_excision}).

Given a $G$-symmetric monoidal functor \( F \colon \ulmfld^{G,f-fr} \to \ul\C \) and a $G$-collar decomposition of 
an $f$-framed $\OG$-manifold 
$M\in \ulmfld^{G,f-fr}_{[G/H]}$ 
we construct a comparison map \( F(M_-) \otimes_{F(M_0 \times (-1,1))} F(M_+) \to F(M) \) in $\ul\C_{[G/H]}$.
This construction depends on the ``inverse image'' functor of \cref{inverse_image_functor}.

\begin{construction}
  Let \( F \colon \ulmfld^{G,f-fr} \to \ul\C \) be a $G$-symmetric monoidal functor.
  Let $M\in \ulmfld^{G,f-fr}_{[G/H]}$ with underlying $\OG$-manifold $(M \to G/H)\in \Gmfld_{[G/H]}$, and $f \colon M\to [-1,1]$ a $G$-collar decomposition.
  Consider the $\disk^{\bnd,or}_{/[-1,1]}$-shaped diagram in $\ul\C_{[G/H]}$ given by the functor 
  \begin{align} \label{tensor_diag}
    \disk^{\bnd,or}_{/[-1,1]} \to \mfld^{\bnd}_{/[-1,1]} \xto{f^{-1}} (\Gmfld_{[G/H]})_{/(M\to G/H)} \simeq (\ulmfld^{G,f-fr}_{[G/H]})_{/M}  \xto{F} (\ul\C_{[G/H]})_{/F(M)}
  \end{align}
  where the first functor is the embedding of disks in manifolds followed by the functor forgetting orientation (see \cite[def. 2.18]{AF}),  
  the second functor is the inverse image functor defined in \cref{inverse_image_functor},
  followed by the equivalence of \cref{Gmfld_framed_fiber_over_M},
  and the third functor is induced by the action of $F$ on the over categories.
  By \cite[lem. 3.11]{AF} there is a cofinal map \( \Delta^{op} \to \disk^{\bnd,or} \), 
  therefore the colimit of \cref{tensor_diag} in $\ul\C_{[G/H]})_{/F(M)}$ is given by 
  \begin{align*}
    \left( \vcenter{\xymatrix{ \colim\left( \SimplicialDiagram{F(M_-) \otimes F(M_0\times (-1,1)) \otimes F(M_+)}{F(M_-) \otimes F(M_+)} \right) \ar[d] \\  F(M) }} \right) \in (\ul\C_{[G/H]})_{/F(M)} ,
  \end{align*}
  known as the two sided bar construction.
  Assume that $\ul\C_{[G/H]}$ admits sifted colimits and that the tensor product functor of $\ul\C_{[G/H]}$ 
  preserves sifted colimits separately in each variable (i.e the coCartesian fibration $\ul\C_{[G/H]}^\otimes \to \Fin_*$ is compatible with sifted colimits in the sense of \cite[def. 3.1.1.18]{HA}). 
  Then the relative tensor product \( F(M_-) \otimes_{F(M_0 \times (-1,1))} F(M_+) \) can be identified with the colimit of this two sided bar construction  (see \cite[thm. 4.4.2.8]{HA}).
  Hence we identify the colimit of the diagram \cref{tensor_diag} with
  \begin{align} \label{tensor_excision_map}
    \left( \vcenter{\xymatrix{  F(M_-) \otimes_{F(M_0 \times (-1,1))} F(M_+) \ar[d] \\  F(M) }} \right) \in (\ul\C_{[G/H]})_{/F(M)} .
  \end{align}
\end{construction}

\begin{mydef} \label{def:G_tensor_excision}
  A $G$-symmetric monoidal functor \( F \colon \ulmfld^{G,f-fr} \to \ul\C \) \emph{satisfies $G$-$\otimes$-excision} if 
  for every  $M \in \ulmfld^{G,f-fr}$ 
  with underlying $\OG$-manifold $(M\to G/H)$ 
  together with a $G$-collar decomposition $f \colon M\to [-1,1]$ 
  the morphism \eqref{tensor_excision_map} 
  is an equivalence in $\ul\C_{[G/H]}$.
\end{mydef}

The main result of this subsection is
\begin{prop} \label{G_tensor_excision}
  Let \( A \colon  \ul\disk^{G,f-fr,\sqcup} \to \ul\C^\otimes \) be an $f$-framed $G$-disk algebra. 
  Then the $G$-factorization homology functor $\int A \colon  \ulmfld^{G,f-fr,\sqcup} \to \ul\C^\otimes$ of \cref{def:G_FH_as_G_SM_functor} satisfies $G$-$\otimes$-excision.
\end{prop}
\begin{rem}
  We view the proof of \cref{G_tensor_excision} as an instance of ``equivariant pushforward''.
  We conjecture that the pushforward paradigm of \cite[sec. 3.4]{AF} and \cite[sec. 2.5]{AyalaFrancisTanaka} has an equivariant generalization to a smooth constructible $G$-fiber bundle between equivariantly-framed $\OG$-manifolds with boundary. 
  However, the definition of equivariantly framed $\OG$-manifolds with boundary is beyond the scope of this work. 
  
  Instead, we are able to prove \cref{G_tensor_excision} without these definitions because the action of $G$ on the oriented manifold $[-1,1]$ is trivial.

  We could have followed a slightly more general approach, considering a $G$-constructible bundle \( M \to N \) where $N$ has boundary and trivial $G$-action.
  To do this, note that since $G$ acts trivially on $N$, any $G$-embedding $V \cof N$ must have a trivial $G$-action as well, so the slice category of $G$-disks over $N$ is a constant $G$-diagram.
  This allows us to harness the definition of (nonequivariant) framing given in \cite{AF} to construct a replacement for the expected ``$G$-slice category of $f$-framed $G$-embeddings $V\cof N$'' needed to preform pushforward.
  We chose not to prove this generalization since we do not currently need it, and we believe it would further obfuscate the proof.
\end{rem}

In order to prove \cref{G_tensor_excision} we need the following auxiliary construction.
\newcommand{\Xf}{\ul{\mathbf{X}}_f}
\begin{construction}
  Let $M \to G/H$ be an $\OG$-manifold and $f \colon M \to [-1,1]$  be a $G$-collar decomposition of $M$. 
  Define a $\ul{G/H}$-category $\Xf \to \ul{G/H}$ and $\ul{G/H}$-functors 
  \begin{align*}
    ev_0  \colon  \Xf \to \ul\disk^{G,f-fr}_{/\ul{M}} , \quad 
    ev_1  \colon  \Xf \to \ul{G/H}\times \disk^{\bnd,or}_{/[-1,1]} 
  \end{align*}
  by the taking the limit of the following diagram of $\ul{G/H}$-categories. 
  \[
    \begin{tikzcd}[column sep = -1em]
      \ul\disk^{G,f-fr}_{/\ul{M}} \ar[dr] & & \ar[dl] \ulFun_{\ul{G/H}} (\ul{G/H}\times \Delta^1, \ulmfld^{G,f-fr}_{/\ul{M}}) \ar[dr] & & \ar[dl,"f^{-1}"] \ul{G/H}\times \disk^{\bnd,or}_{/[-1,1]} \\
      & \ulmfld^{G,f-fr}_{/\ul{M}}  & & \ulmfld^{G,f-fr}_{/\ul{M}} 
    \end{tikzcd}
  \]
  where $f^{-1}$ is the inverse image functor of \cref{inverse_image_functor}.
\end{construction}

\begin{rem} \label{ev0_explained}
  Using \cref{Gmfld_framed_over_M} and unwinding the definitions shows that the fiber of $\Xf \to \ul{G/H}$ over \( \varphi \colon G/K \to G/H \) is given by the limit of 
  \[
    \begin{tikzcd}[column sep = -1em]
      (\disk^G_{[G/K]})_{/\varphi^* M} \ar[dr] & & \ar[dl] \Fun(\Delta^1, (\disk^G_{[G/K]})_{/\varphi^* M} ) \ar[dr] & & \ar[dl,"f^{-1}"] \disk^{\bnd,or}_{/[-1,1]} \\
      & (\mfld^G_{[G/K]})_{/\varphi^* M}  & & (\mfld^G_{[G/K]})_{/\varphi^* M} ,
    \end{tikzcd}
  \]
  where the $\infty$-over categories \( (\disk^G_{[G/H]})_{/ \varphi^* M}, \, (\mfld^G_{[G/H]})_{/ \varphi^* M} \) can be modeled as the coherent nerve of the Moore over category (see \cref{MooreOverCat}). 
  Explicitly, an object of $(\Xf)_{[\varphi]}$ is given by \( (g \colon V \cof [-1,1], h \colon  E \cof \varphi^* M, h' \colon  E \cof f^{-1} V, \gamma) \) where 
  \begin{itemize}
    \item $V$ is a finite disjoint union of 1-dimensional oriented disks with boundary, i.e oriented open intervals equivalent to $\mathbb{R}$ and oriented half open intervals equivalent to \( [0,1) \) or \( (0,1] \),
    \item $g$ is an orientation preserving embedding of $V$ into the closed interval $[-1,1]$,
    \item $E \to U \to G/K$  is a finite $G$-disjoint union of $G$-disks (i.e $E\to U$ a $G$-vector bundle over a finite $G$-set), 
    \item $h$ is a $G$-equivariant embedding over $G/K$ of $E$ into the pullback of $M\to G/H$ along $\varphi$, 
    \item $h'$ is a $G$-equivariant embedding over $G/K$ of $E$ into the preimage $f^{-1} V$
    \item $\gamma$ is a Moore path in \( Emb^G_{[G/K]}(E,\varphi^*M) \) from \( h\) to  \(  E \xto{h'} f^{-1} V \xto{f^{-1}(g)} \varphi^* M \).
  \end{itemize}
  The functor $ev_0$ sends the object \( (g \colon V \cof [-1,1], h \colon  E \cof \varphi^* M, h' \colon  E \cof f^{-1} V, \gamma) \) described above to \( (h \colon E \cof \varphi^* M) \in (\mfld^G_{[G/K]})_{/\varphi^* M} \), while the functor $ev_1$ sends it to \( (g \colon V \cof [-1,1]) \in \disk^{\bnd,or}_{/[-1,1]} \) .

  By \cite[prop. 2.4.7.12]{HTT} it follows that for every $\varphi\in \ul{G/H}$ the functor 
  \[
    (ev_0)_{[\varphi]} \colon  (\Xf)_{[\varphi]} \to (\disk^G_{[G/H]})_{/\varphi^* M} 
  \]
  is a Cartesian fibration (and therefore that $ev_1$ is a $\ul{G/H}$-Cartesian fibration, see \cite[def. 7.1]{Expose2}).
\end{rem}

The following lemma is the main ingredient in the proof of \cref{G_tensor_excision}.
\begin{lem} \label{ev0_cofinal}
  The $\ul{G/H}$-functor $ev_0 \colon  \Xf \to \ul\disk^{G,f-fr}_{/\ul{M}}$ is $\ul{G/H}$-cofinal.
\end{lem}

The following proof is an adaptation of \cite[thm. 5.5.3.6]{HA}, \cite[lem. 3.21]{AF} and \cite[lem. 2.27]{AyalaFrancisTanaka} to the equivariant setting. 
\begin{proof}[Proof of \cref{ev0_cofinal}.] 
  By \cref{Gmfld_framed_over_M} we have to prove that $ev_0 \colon  \Xf \to \Gdisk_{/\ul{M}}$ is $\ul{G/H}$-cofinal.
  By \cite[thm. 6.7, def. 6.8]{Expose2} the $\ul{G/H}$ functor $ev_0$ is $\ul{G/H}$-cofinal if and only if for each $(\varphi \colon  G/K \to G/H) \in \ul{G/H}$ the functor 
  \( (ev_0)_{[\varphi]}  \colon  (\Xf)_{[\varphi]} \to (\Gdisk_{/\ul{M}})_{[\varphi]} \) is cofinal. 
  
  By replacing $f \colon M \to [-1,1]$ with $\varphi^* M \to M \xto{f} [-1,1]$ we reduce to $\varphi=(G/H \xto{=} G/H)\in \ul{G/H}$: it is enough to prove that $(ev_0)_{[G/H]}$ is cofinal.

  By \cref{ev0_explained} the functor $(ev_0)_{[G/H]}$ is a Cartesian fibration, therefore by \cite[prop. 4.1.3.2]{HTT} it is enough to show that for each $(E \cof M)\in (\disk^G_{[G/H]})_{/M}$ the fiber \( (ev_0)^{-1}( E \cof M) \) is weakly contractible. 

  Note that the category $(ev_0)^{-1}(E\cof M)$ has a functor to $\disk^{\bnd,or}_{/[-1,1]}$ by construction:
  \[
    \begin{tikzcd}
      (ev_0)^{-1}(E \cof M) \ar[d] \ar[r] \pbcorner & \disk^{\bnd,or}_{/[-1,1]} \= \set{E \cof M} \times \disk^{\bnd,or}_{/[-1,1]} \ar[d] \\
      (\Xf)_{[G/H]} \ar[r,"{( (ev_0)_{[G/H]} , (ev_1)_{[G/H]} )}"] \ar[d] \pbcorner & (\disk^G_{[G/H]})_{/M} \times \disk^{\bnd,or}_{/[-1,1]} \ar[d, "{(\iota,f^{-1})}"] \\
      Fun(\Delta^1, (\mfld^G_{[G/H]})_{/M}) \ar[r] & (\mfld^G_{[G/H]})_{/M} \times (\mfld^G_{[G/H]})_{/M}
    \end{tikzcd}
  \]
  The top horizontal functor \( (ev_0)^{-1}(E\cof M) \to \disk^{\bnd,or}_{/[-1,1]} \) is pullback of a left fibration, since it can be written as
  \[
    \begin{tikzcd}
      (ev_0)^{-1}(E \cof M) \ar[r,->>] \ar[d] \pbcorner & \disk^{\bnd,or}_{/[-1,1]} \ar[d,"{f^{-1}}"] \\
      \left( (\mfld^G_{[G/H]})_{/M} \right)_{(E \cof M)/} \ar[r,->>] \ar[d] \pbcorner & (\mfld^G_{[G/H]})_{/M}  \ar[d,"{(\set{E \cof M},id)}"] \\
      Fun(\Delta^1, (\mfld^G_{[G/H]})_{/M}) \ar[r]  & (\mfld^G_{[G/H]})_{/M} \times (\mfld^G_{[G/H]})_{/M} ,
    \end{tikzcd}
  \]
  where the middle horizontal arrow is a left fibration by \cite[cor. 2.1.2.2]{HTT}.

  The left fibration \( (ev_0)^{-1}(E \cof M) \to \disk^{\bnd,or}_{/[-1,1]} \) classifies the functor
  \begin{align*}
    (\mfld^G_{[G/H]})_{/M} \to \Ss, \quad (V\cof [-1,1]) \mapsto \Map_{(\mfld^G_{[G/H]})_{/M}}  (E \cof M, f^{-1} V \cof M) ,
  \end{align*}
  and by \cite[3.3.4.5]{HTT} we have to show that the colimit
  \begin{align*}
   \colim_{(V\cof[-1,1])\in \disk^{\bnd,or}_{[-1,1]}} \Map_{(\mfld^G_{[G/H]})_{/M}}  (E \cof M, f^{-1} V \cof M) 
  \end{align*}
  is weakly contractible. 
  
  Let \( \disk^{\bnd,or}([-1,1]) \) denote the ordinary category with the same objects as $\disk^{\bnd,or}_{/[-1,1]}$ and sets of morphisms given by forgetting the topology of the mapping spaces of $\disk^{\bnd,or}_{/[-1,1]}$ (see \cite[def. 2.8]{AF}).
  Note that the category \( \disk^{\bnd,or}([-1,1]) \) is equivalent to the partial ordered set of open subsets $V \subsetneq [-1,1]$ for $V$ a finite disjoint union of intervals in $[-1,1]$, possibly containing the edge points $-1,1$, after excluding the whole interval $[-1,1] \subseteq [-1,1]$. 

  By \cite[prop. 2.19]{AF} the functor \( \disk^{\bnd,or}([-1,1]) \to \disk^{\bnd,or}_{/[-1,1]} \) is cofinal, hence it is enough to show that the homotopy colimit
  \begin{align*}
    \hocolim_{(V\subsetneq [-1,1])\in \disk^{\bnd,or}([-1,1])} \Map_{(\mfld^G_{[G/H]})_{/M}}  (E \cof M, f^{-1} V \cof M) 
  \end{align*}
  is contractible. 

  Using \cref{overcat_map_is_h_fiber} we see that the space $\Map_{(\mfld^G_{[G/H]})_{/M}} (E \cof M, f^{-1} V \cof M)$ is the homotopy fiber of \( Emb^G_{[G/H]}(E, f^{-1} V) \to Emb^G_{[G/H]}(E,M) \), hence by \cite{HTT} it is enough to show that the map 
  \begin{align*}
    \hocolim_{(V\subsetneq [-1,1])\in \disk^{\bnd,or}([-1,1])} Emb^G_{[G/H]}(E, f^{-1} V) \to Emb^G_{[G/H]}(E,M) 
  \end{align*}
  is an equivalence. 

  By \cite[thm. 6.1.0.6]{HTT} colimits in $\Ss$ are universal, therefore by \cref{Emb_Conf_PB} it is enough to prove that the map 
  \begin{align*}
    \hocolim_{(V\subsetneq [-1,1])\in \disk^{\bnd,or}([-1,1])} \ul\Conf^G_{G/H}(U;f^{-1} V) \to \ul\Conf^G_{G/H}(U;M) 
  \end{align*}
  is an equivalence.
  Since \( \set{\ul\Conf^G_{G/H}(U;f^{-1} V)}_{(V\subsetneq[-1,1])\in \disk^{\bnd,or}([-1,1])} \) is a complete open cover of $\ul\Conf^G_{G/H}(U;M)$ it follows from \cite[cor. 1.6]{DuggerIsaksen} that the above map is an equivalence. 
\end{proof}

We will also need a simple cofinality lemma. 
Assume we have a coCartesian fibration \(p \colon \ul\D \fib \ul\C \) between $S$-categories $\ul\C,\ul\D$ (i.e. an $S$-coCartesian fibration, see \cite[rem. 7.3]{Expose2}), and an $S$-object $x \colon S \to \ul\C$.
Let \( p^{-1}(x) := S \times_{\ul\C} \ul\D \) be the pullback of $p$ along $x$.
Since \( p^{-1}(x) \fib S \) is a coCartesian fibration we can considered $p^{-1}(x)$ as an $S$-category, which we denote by $\ul{p^{-1}(x)}$.
\begin{lem} \label{fiber_to_slice_cofinal}
  Let $\ul\C,\ul\D$ be $S$-categories and \(p \colon \ul\D \fib \ul\C \) be a coCartesian fibration, and $x \colon S \to \ul\C$ an $S$-object of $\ul\C$.
  Then the $S$-functor \( \ul{p^{-1}(x)} \to \ul\D_{/\ul{x}} \) is $S$-cofinal.
\end{lem}
\begin{proof}
  By \cite[thm. 6.7]{Expose2} we have to show that for each $s\in S$ the functor \( \ul{p^{-1}(x)}_{[s]}  \to (\ul\D_{\ul{x}})_{[s]} \) between the fibers of \( \ul{p^{-1}(x)} \to \ul\D_{/\ul{x}} \) over $s$ is cofinal. 
  Since $p_{[s]} \colon  \ul\D_{[s]} \fib \ul\C_{[s]}$ is a coCartesian fibration it follows that \( (p_{[s]}^{-1}(x(s)) \to (\ul\D_{[s]})_{/x(s)} \) is cofinal. 
  The result now follows from the equivalence \( p^{-1}_{[s]} \= (p_{[s]})^{-1}(x(s)) \). 
\end{proof}

With \cref{ev0_cofinal} at hand we turn to the proof of \cref{G_tensor_excision}.
The proof follows the outline of the proof of \cite[prop. 3.23]{AF} (``pushforward'').
\begin{proof}[Proof of \cref{G_tensor_excision}.]
  By \cref{G_FH_colimit_formula} and \cref{ev0_cofinal} we have 
  \begin{align*}
    \int_M A & := \ul{G/H}-\colim ( \ul\disk^{G,f-fr}_{/\ul{M}} \to \ul\disk^{G,f-fr} \ultimes \ul{G/H} \xto{A \ultimes \ul{G/H}} \ul\C \ultimes \ul{G/H} ) \\
    & = \ul{G/H}-\colim ( \Xf \xto{ev_0} \ul\disk^{G,f-fr}_{/\ul{M}} \to \ul\disk^{G,f-fr} \ultimes \ul{G/H} \xto{A \ultimes \ul{G/H}} \ul\C \ultimes \ul{G/H} ) .
  \end{align*}
  Using the characterization of parametrized Kan extensions as parametrized left adjoints (see \cite[thm. 10.4]{Expose2}, and also \cite[def. 2.10 and def. 2.12]{Expose4}) 
  we can express the above $\ul{G/H}$-colimit as a left $\ul{G/H}$-Kan extension of $L \colon \Xf \to \ul\C \ultimes \ul{G/H}$ along the structure map $\Xf \to \ul{G/H}$, where $L$ is the $\ul{G/H}$-functor given by the composition 
  \begin{align} \label{Xf_diag}
    L  \colon  \Xf \xto{ev_0} \ul\disk^{G,f-fr}_{/\ul{M}} \to \ul\disk^{G,f-fr} \ultimes \ul{G/H} \xto{A \ultimes \ul{G/H}} \ul\C \ultimes \ul{G/H} .
  \end{align}
  Equivalently the $\ul{G/H}$-colimit over $\Xf$ is given by the left $\ul{G/H}$-adjoint to restriction along the structure map \( \Xf \to \ul{G/H} \),
  \begin{align*}
    \ul{G/H}-\colim  \colon  
      \ulFun_{\ul{G/H}}( \Xf, \ul\C \ultimes \ul{G/H} )
      \adj
      \ulFun_{\ul{G/H}}( \ul{G/H} ,  \ul\C \ultimes \ul{G/H} ) \simeq \ul\C \ultimes \ul{G/H} .
  \end{align*}
  By construction $ev_1 \colon  \Xf \to \ul{G/H} \times \disk^{\bnd,or}_{/[-1,1]}$ is a  $\ul{G/H}$ category, therefore the structure map \( \Xf \to \ul{G/H} \) factors as 
  \( \Xf \xto{ev_1} \ul{G/H} \times \disk^{\bnd,or}_{/[-1,1]} \to \ul{G/H} \). 
  We can now extend $L$ along $\Xf \to \ul{G/H}$ in two steps, again using \cite[thm. 10.4]{Expose2}, as the composition of left $\ul{G/H}$-adjoints
  \[
    \begin{tikzcd}
      (ev_1)_! \colon \ulFun_{\ul{G/H}}( \Xf, \ul\C \ultimes \ul{G/H} ) \ar[r, shift right] &
      \ulFun_{\ul{G/H}}( \ul{G/H} \times \disk^{\bnd,or}_{/[-1,1]} ,  \ul\C \ultimes \ul{G/H} ) \noloc (ev_1)^*  \ar[l, shift right], 
    \end{tikzcd} 
  \]
  \[
    \begin{tikzcd}
      \ul{G/H}-\colim \colon \ulFun_{\ul{G/H}}( \ul{G/H} \times \disk^{\bnd,or}_{/[-1,1]} ,  \ul\C \ultimes \ul{G/H} )  \ar[r, shift right] &
      \ulFun_{\ul{G/H}}( \ul{G/H} ,  \ul\C \ultimes \ul{G/H} ) \simeq \ul\C \ultimes \ul{G/H}   \ar[l, shift right],
    \end{tikzcd}
  \]
  where $(ev_1)_!$ is the left $\ul{G/H}$-Kan extension of \eqref{Xf_diag} along $ev_1$.
  In particular restricting to fibers over $(G/H \xfrom{=} G/H) \in \ul{G/H}$ we get composition of unparametrized left adjoints (see \cite[prop. 7.3.2.6]{HA} and \cite[def. 8.1]{Expose2}):
  \[
    \begin{tikzcd}
      (ev_1)_! \colon \Fun_{\ul{G/H}}( \Xf, \ul\C \ultimes \ul{G/H} ) \ar[r, shift right] &
      \Fun_{\ul{G/H}}( \ul{G/H} \times \disk^{\bnd,or}_{/[-1,1]} ,  \ul\C \ultimes \ul{G/H} ) \noloc (ev_1)^* \ar[l, shift right] ,
    \end{tikzcd} 
  \]
  \[
    \begin{tikzcd}
      \ul{G/H}-\colim \colon \Fun_{\ul{G/H}}( \ul{G/H} \times \disk^{\bnd,or}_{/[-1,1]} ,  \ul\C \ultimes \ul{G/H} ) \ar[r, shift right] &
      \Fun_{\ul{G/H}}( \ul{G/H} ,  \ul\C \ultimes \ul{G/H} ) \simeq \ul\C_{[G/H]} \ar[l, shift right] . 
    \end{tikzcd}
  \]
  Applying both left adjoints to the $\ul{G/H}$-functor  $L \colon \Xf \to \ul\C \ultimes \ul{G/H}$ of \eqref{Xf_diag} produces the $G$-factorization homology \( \int_M A = \ul{G/H}-\colim(L \colon \Xf \to \ul\C \ultimes \ul{G/H})  \).
  Let \(L' := (ev_1)_! (L) \in \Fun_{\ul{G/H}}(\ul{G/H}\times \disk^{\bnd,or}_{/[-1,1]}, \ul\C \ultimes \ul{G/H}) \) be the left $\ul{G/H}$-Kan extension of $L$ along $ev_1$.
  Then the $\ul{G/H}$-colimit of $L'$ is 
  \begin{align*}
    \ul{G/H}-\colim(L') = \ul{G/H}-\colim \left( (ev_1)_!(L) \right) \simeq \ul{G/H}-\colim(L) = \int_M A .
  \end{align*}
  
  Next, note that the $\ul{G/H}$-colimit over the constant diagram \( \ul{G/H}\times \disk^{\bnd,or}_{/[-1,1]} \) is equivalent to an unparametrized colimit over \(\disk^{\bnd,or}_{/[-1,1]} \).
  To see this, use the equivalence 
  \begin{align*}
    \Fun_{\ul{G/H}}(\ul{G/H}\times \disk^{\bnd,or}_{/[-1,1]}, \ul\C \ultimes \ul{G/H}) & \iso \Fun ( \disk^{\bnd,or}_{/[-1,1]}, C_{[G/H]} ) , \\
    L' & \mapsto L'|_{ \set{G/H \xfrom{=} G/H} \times \disk^{\bnd,or}_{/[-1,1]} } 
  \end{align*}
  and the global definition of a colimit as a left adjoint
  \begin{align*}
    \colim  \colon   \Fun ( \disk^{\bnd,or}_{/[-1,1]}, \ul\C_{[G/H]} ) \adj \ul\C_{[G/H]}.
  \end{align*}
  Therefore, we have  
  \begin{align*}
    \int_M A \simeq \ul{G/H}-\colim(L') \simeq \colim \left( L'\vert_{ \set{G/H \xfrom{=} G/H} \times \disk^{\bnd,or}_{/[-1,1]} } \right) ,
  \end{align*}
  which we write as 
  \[
    \int_M A \simeq \colim\limits_{ (V \cof [-1,1]) \in \disk^{\bnd,or}_{/[-1,1]} }  L'(G/H \xto{=} G/H,V \cof [-1,1] ) .
  \]
  
  Out next goal is to calculate $L'(y)$ for \(y =(G/H \xto{=} G/H,V \cof [-1,1] ) \) where  \( (V \cof [-1,1]) \in \disk^{\bnd,or}_{/[-1,1]} \) is an oriented embedding.
  We claim that \( L'(y) \simeq \int_{f^{-1} V} A \in \ul\C_{[G/H]} \) is the $G$-factorization homology of 
  \( f^{-1} (V)  \in \ul\mfld^{G,f-fr}_{[G/H]}\). 
  After asserting our claim we use the cofinal map \( \Delta^{op} \to \disk^{\bnd,or}_{/[-1,1]} \) of \cite[lem. 3.11]{AF} to deduce $\int_M A$ is equivalent to the colimit of the simplicial diagram \( \Delta^{op} \to \set{G/H \xfrom{=} G/H} \times \disk^{\bnd,or}_{/[-1,1]} \xto{L'} \ul\C_{[G/H]} \).
  Since the functor $L'(V\cof [-1,1]) \simeq \int_{f^{-1} V} A$ takes disjoint unions over $[-1,1]$ to tensor product in $\ul\C_{[G/H]}$ (see \cref{FH_G_SM}), and using the equivalence of oriented open and half open intervals over $[-1,1]$ (as objects of $\disk^{\bnd,or}_{/[-1,1]}$) we see that $\int_M A \in \ul\C_{[G/H]}$ is equivalent to the realization of the two sided bar construction 
  \begin{align*}
    \SimplicialDiagram{ \int_{M_-} A \otimes \int_{M_0 \times (-1,1)} A \otimes \int_{M_-} A }{ \int_{M_-} A \otimes \int_{M_-} A }, \\
    [n] \mapsto \left(\int_{M_-} A\right) \otimes \left(\int_{M_0 \times (-1,1)} A\right)^{\otimes n}  \otimes \left(\int_{M_-} A \right).
  \end{align*}
  By \cite[4.4.2.8-11]{HA} we see that $G$-factorization homology of $M$ is equivalent to the relative tensor product  \( \int_M A \simeq  (\int_{M_-} A) \otimes_{(\int_{M_0 \times (-1,1)} A)} (\int_{M_-} A )\).
  Therefore it is enough to prove our claim that \( L'(y) \simeq \int_{f^{-1} V} A \in \ul\C_{[G/H]} \).

  Since $L'$ is the left $\ul{G/H}$-Kan extension of $L \colon \Xf \to \ul\C \ultimes \ul{G/H}$ along $ev_1$, it is given by the following $\ul{G/H}$-colimit
  \begin{align*}
    L'(y) & = \ul{G/H}-\colim \left( (\Xf)_{/\ul{y}} \to \Xf \xto{L} \ul\C \ultimes \ul{G/H} \right) \\
    & = \ul{G/H}-\colim \left( (\Xf)_{/\ul{y}} \to \Xf \xto{ev_0} \ul\disk^{G,f-fr}_{/\ul{M}} \to \ul\disk^{G,f-fr} \ultimes \ul{G/H} \xto{A \ultimes \ul{G/H}} \ul\C \ultimes \ul{G/H} \right) \\
  \end{align*}
  Next we replace the $\ul{G/H}$-category $\Xf$ indexing the above colimit by a $\ul{G/H}$-category which is more closely related to $G$-disks in $f^{-1} V$. 
  Note that \( ev_1 \colon  \Xf \to \ul{G/H} \times \disk^{\bnd,or} \) is a coCartesian fibration, and let 
  \(
  \ul{(ev_1)^{-1}(y)} 
  \)
  denote the pullback of $\Xf$ along the $\ul{G/H}$-functor 
  \[
    \ul{G/H} \to \ul{G/H} \times \disk^{\bnd,or}, \quad ( G/K \to G/H ) \mapsto ( G/K \to G/H , V \cof [-1,1])
  \] 
  corresponding to $y=(G/H \xto{=} G/H, V \cof [-1,1]) \in \ul{G/H} \times \disk^{\bnd,or}$.
  By \cref{fiber_to_slice_cofinal} the $\ul{G/H}$-functor \( \ul{(ev_1)^{-1}(y)} \to (\Xf)_{/\ul{y}} \) is $\ul{G/H}$-cofinal, 
  hence $L'(y)$ is the $\ul{G/H}$-colimit of the $\ul{G/H}$-diagram
  \[
    \ul{(ev_1)^{-1}(y)} \to (\Xf)_{/\ul{y}} \to \Xf \xto{ev_0} \ul\disk^{G,f-fr}_{/\ul{M}} \to \ul\disk^{G,f-fr} \ultimes \ul{G/H} \xto{A \ultimes \ul{G/H}} \ul\C \ultimes \ul{G/H}  .
  \]
  Since $ev_1 \colon \Xf \to \ul{G/H} \times \disk^{\bnd,or}$ factors through $ (\ul\disk^{G,f-fr}_{/\ul{M}}) \times_{\ul{G/H}} \left( \ul{G/H} \times \disk^{\bnd,or}_{/[-1,1]} \right)$ we can express $\ul{(ev_1)^{-1}(y)}$ as the iterative pullback
  \[
    \begin{tikzcd}
      \ul{(ev_1)^{-1}(y)} \pbcorner \ar[d] \ar[r] & 
        \Xf \ar[d,"{(ev_0,ev_1)}"] \\
        \ul\disk^{G,f-fr}_{/\ul{M}} \times_{\ul{G/H}} \ul{G/H} \pbcorner \ar[d] \ar [r,"id\ultimes y"]  &
        (\ul\disk^{G,f-fr}_{/\ul{M}}) \times_{\ul{G/H}} \left( \ul{G/H} \times \disk^{\bnd,or}_{/[-1,1]} \right) \ar[d]  \\
      \ul{G/H} \ar[r]^-{y} &
        \ul{G/H} \times \disk^{\bnd,or}_{/[-1,1]} .
    \end{tikzcd}
  \]
  On the other hand we can express $\Xf$ is the pullback
  \[
    \begin{tikzcd}
        \Xf \pbcorner \ar[r] \ar[d,"{(ev_0,ev_1)}"] & 
        \ulFun_{\ul{G/H}} ( \ul{G/H} \times \Delta^1 , \ulmfld^{G,f-fr}_{/\ul{M}} ) \ar[d] \\
        (\ul\disk^{G,f-fr}_{/\ul{M}}) \times_{\ul{G/H}} \left( \ul{G/H} \times \disk^{\bnd,or}_{/[-1,1]} \right) \ar[r,"{\iota \times f^{-1}}"] &
        \ulmfld^{G,f-fr}_{/\ul{M}} \times_{\ul{G/H}} \ulmfld^{G,f-fr}_{/\ul{M}} .
    \end{tikzcd}
  \]
  Notice that the composition 
  \[
    \begin{tikzcd}
        \ul\disk^{G,f-fr}_{/\ul{M}} \times_{\ul{G/H}} \ul{G/H}  \ar [d,"id\ultimes y"]  \\
        (\ul\disk^{G,f-fr}_{/\ul{M}}) \times_{\ul{G/H}} \left( \ul{G/H} \times \disk^{\bnd,or}_{/[-1,1]} \right) \ar[d,"{\iota \times f^{-1}}"] \\
        \ulmfld^{G,f-fr}_{/\ul{M}} \times_{\ul{G/H}} \ulmfld^{G,f-fr}_{/\ul{M}} 
    \end{tikzcd}
  \]
  is equivalent to 
  \begin{align*}
    \ul\disk^{G,f-fr}_{/\ul{M}} \times_{G/H} \ul{G/H} 
    \xto{(\iota,f^{-1}(y)) } \ulmfld^{G,f-fr}_{/\ul{M}} \ultimes \ulmfld^{G,f-fr}_{/\ul{M}} ,
  \end{align*}
  and therefore that 
  \[ 
    \ul{(ev_1)^{-1}(y)} \= \left(\ul\disk^{G,f-fr}_{/\ul{M}} \right)_{/\ul{(f^{-1} V \cof M)}} \simeq \ul\disk^{G,f-fr}_{/\ul{f^{-1} V}} 
  \]
  (compare \cite[lem. 2.1]{AF}).
  Finally, since the diagram
  \begin{align*}
    \diag{
      \ul{(ev_1)^{-1}(y)} \ar[d]^{\simeq} \ar[r] & (\Xf)_{/\ul{y}} \ar[r] &  \Xf \ar[d] \\
      \ul\disk^{G,f-fr}_{/\ul{f^{-1} V}} \ar[r] & \ul\disk^{G,f-fr}_{/\ul{M}} \ar[r] & \ul\disk^{G,f-fr} \ultimes \ul{G/H} \ar[r]^-{A\ultimes\ul{G/H}} & \ul\C \ultimes \ul{G/H}
    }
  \end{align*}
  commutes, we get 
  \begin{align*}
    L'(y) 
    & \simeq \ul{G/H}-\colim \left( \ul\disk^{G,f-fr}_{/\ul{f^{-1} V}} \to \ul\disk^{G,f-fr} \ultimes \ul{G/H} \xto{A \ultimes \ul{G/H}} \ul\C \ultimes \ul{G/H} \right).
  \end{align*}
  Therefore, by the definition of left $\ul{G/H}$-Kan extension we see that indeed \( L'(y) \simeq \int_{f^{-1} V} A \).
\end{proof}

\subsection{$G$-sequential unions} \label{sec:G_FH_respects_seq_unions}
\begin{mydef} \label{def:G_seq_unions}
  Let $M$ be a $G$-manifold. 
  A \emph{$G$-sequential union} of $M$ is a sequence of open $G$-submanifolds \( M_1 \subset M_2 \subset \cdots \subset M\) with $M = \cup_{i=1}^{\infty} M_i$. 
  A $G$-sequential union of an $f$-framed $\OG$-manifold $M \in \ulmfld^{G,f-fr}_{[G/H]}$ is a $G$-sequential union of its underlying $G$-manifold.
\end{mydef}
If \( F \colon \ulmfld^{G,f-fr} \to \ul\C \) is a $G$-symmetric monoidal functor and $M=\cup_{i=1}^{\infty} M_i$ is a $G$-sequential union of 
$M\in \ulmfld^{G,f-fr}_{[G/H]}$,
then we have a comparison morphism 
\( \colim F(M_i) \to F(M) \) 
in $\ul\C_{[G/H]}$.
\begin{mydef} \label{def:F_respects_G_seq_unions}
  We say that $G$-symmetric monoidal functor \( F \colon \ulmfld^{G,f-fr} \to \ul\C \) \emph{respects $G$-sequential unions} if for every $G$-sequential union $M=\cup_{i=1}^{\infty} M_i$ the comparison morphism
  \[
    \colim F(M_i) \to F(M) 
  \]
  is an equivalence in $\ul\C_{[G/H]}$.
\end{mydef}

\begin{prop} \label{G_FH_respects_seq_unions}
  Let $\ul\C^\otimes \fib \GFin_*$ be a $G$-symmetric monoidal $G$-category and $A$ be an $f$-framed $G$-disk algebra with values in $\ul\C$. 
  Then $G$-factorization homology \( \int_- A  \colon  \ulmfld^{G,f-fr} \to \ul\C \) of \cref{def:G_FH_as_G_SM_functor} respects $G$-sequential unions. 
\end{prop}

The proof of \cref{G_FH_respects_seq_unions} relies on the following lemma. 
\begin{lem} \label{lem:disks_over_Mi}
  Let $M\in\ulmfld^{G,f-fr}_{[G/H]}$ be an $f$-framed $\OG$-manifold over $G/H$, and \( M = \cup_{i=1}^{\infty} M_i \) a $G$-sequential union of $M$. 
  Then the $\ul{G/H}$-functor
  \( \colim \ul\disk^{G,f-fr}_{/\ul{M_i}} \iso \ul\disk^{G,f-fr}_{/\ul{M}} \)
  is an equivalence of $\ul{G/H}$-categories.
\end{lem}
\begin{proof}
  By \cref{Gmfld_framed_over_M} it is enough to prove that the $\ul{G/H}$-functor
  \( \colim \Gdisk_{/\ul{M_i}} \to \Gdisk_{/\ul{M}} \)
  is a fiberwise equivalence. 
  Without loss of generality we show that the functor between the fiber over $(G/H \xfrom{=} G/H) \in \ul{G/H}$ is an equivalence.
  Since colimits of $G$-categories are computed fiberwise, we have to show that \( \colim_i (\Gdisk_{[G/H]})_{/M_i} \to (\Gdisk_{[G/H]})_{/M} \) is an equivalence of $\infty$-categories.
  
  In order to show that this functor is fully faithful we first show that $Emb^G_{G/H}(E,M)$ is equivalent to the homotopy colimit $\hocolim_i Emb^G_{G/H}(E,M_i)$.
  Let $(E \to U \to G/H)\in \Gdisk$ be a finite $G$-disjoint union of $G$-disks, i.e. \( E \to U \) a $G$-vector bundle,  $U=\pi_0 E$.
  By \cref{Emb_Conf_PB} the square
  \begin{align*}
    \diag{
      Emb^G_{G/H}(E,M_i) \ar[d] \ar[r] & Emb^G_{G/H}(E,M) \ar[d] \\ 
      \ul\Conf^G_{G/H}(U;M_i) \ar[r] & \ul\Conf^G_{G/H}(U;M)
    }
  \end{align*}
  is a homotopy pullback square for each $i\in \N$.
  Since filtered homotopy colimits preserves homotopy pullbacks, the square
  \begin{align*}
    \diag{
      \hocolim_i Emb^G_{G/H}(E,M_i) \ar[d] \ar[r] & Emb^G_{G/H}(E,M) \ar[d] \\ 
      \hocolim_i \ul\Conf^G_{G/H}(U;M_i) \ar[r] & \ul\Conf^G_{G/H}(U;M)
    }
  \end{align*}
  is also a homotopy pullback square. 
  However, \( \set{ \ul\Conf^G_{G/H}(U;M_i) }_i\in \N \) is a complete open cover of $\ul\Conf^G_{G/H}(U;M)$, so by \cite[cor. 1.6]{DuggerIsaksen} the bottom map is a weak equivalence. 
  Therefore the map \( \hocolim_i Emb^G_{G/H}(E,M_i) \iso Emb^G_{G/H}(E,M) \) is a weak equivalence. 
  
  We now show that \( \colim_i (\Gdisk_{[G/H]})_{/M_i} \to (\Gdisk_{[G/H]})_{/M} \) is fully faithful.

  Let \( (E' \to U' \to G/H) , \, (E'' \to U'' \to G/H) \in \Gdisk_{ [G/H] } \) and \( f' \colon  E' \cof M_{i'} \, f'' \colon E'' \cof M_{i''}) \) be two $G$-embeddings over $G/H$, representing two objects in \( \colim_i (\Gdisk_{[G/H]})_{/M_i} \) .
  For $i$ greater then $i'$ and $i''$ they represent objects of the same slice category 
  \[ (f'_i \colon  E' \xto{f'} M_{i'} \subseteq M_i ), \, ( f''_i \colon  E'' \xto{f''} M_{i''} \subseteq M_i) \in (\Gdisk_{[G/H]})_{/M_i} ,\]
  with mapping space \( \Map_{(\Gdisk_{[G/H]})_{/M_i}} ( f'_i  \colon  E' \cof M_i, f''_i \colon  E'' \cof M_i) \)  given by the homotopy fiber of \( (f''_i)_*  \colon  Emb^G_{G/H} (E' , E'') \to Emb^G_{G/H}(E',M_i) \) over $f'_i \in Emb^G_{G/H}(E',M_i)$.
  Homotopy fibers are preserved by filtered homotopy colimits, so the homotopy fiber of the map 
  \[ Emb^G_{G/H} (E' , E'') \to \hocolim_i Emb^G_{G/H}(E',M_i) \simeq Emb^G_{G/H}(E',M) \]
  induced by post composition with $f'' \colon  E'' \cof  M_i \subset M$ over $f' \colon  E' \cof M_i\subseteq M$ is equivalent to \( \hocolim_i \Map_{(\Gdisk_{[G/H]})_{/M_i}} ( f'_i  \colon  E' \cof M_i, f''_i \colon  E'' \cof M_i) \).
  On the other hand, this homotopy fiber is equivalent to the mapping space of the slice category $ (\Gdisk_{[G/H]})_{/M}$, hence 
  \[ 
    \colim_i \Map_{(\Gdisk_{[G/H]})_{/M_i}} ( f'_i \colon  E' \cof M_i, f''_i  \colon E'' \cof M_i) 
  \]
  is homotopy equivalent to 
  \[
    \Map_{(\Gdisk_{[G/H]})_{/M}} ( f' \colon  E' \cof M, f''  \colon E'' \cof M) , 
  \]
  so the functor \( \colim_i (\Gdisk_{[G/H]})_{/M_i} \to (\Gdisk_{[G/H]})_{/M} \) is fully faithful.

  It remains to show that  \( \colim_i (\Gdisk_{[G/H]})_{/M_i} \to (\Gdisk_{[G/H]})_{/M} \) is essentially surjective. 
  Let \( (E\to U \to G/H) \in \Gdisk_{[G/H]}, \, (f \colon E\cof M) \in (\Gdisk_{[G/H]})_{/M}\) for $E\to U$ a $G$-vector bundle. 
  Choose $t>0$ small enough so that the restriction of $f$ to the open ball of radius $t$ bundle, \( B_t(E) \cof E \xto{f} M \), factors through some $M_i\subseteq M$.
  By radial dilation we see that the inclusion \( (B_t(E) \to G/H) \to (E\to G/H) \) is an equivalence in $\Gdisk_{[G/H]}$.
  Postcomposition with $f \colon E \cof M$ induces an equivalence \( (f \colon E \cof M ) \simeq (B_t(E) \cof E \xto{f} M) \) of objects in the slice category $(\Gdisk_{[G/H]})_{/M} $.
  On the other hand, since $(B_t(E) \cof E \xto{f} M)$ factors through $M_i$ this object is clearly in the image of the functor \( \colim_i (\Gdisk_{[G/H]})_{/M_i} \to (\Gdisk_{[G/H]})_{/M} \), showing the functor is indeed essentially surjective.
\end{proof}

We now show that $G$-factorization homotopy respects sequential colimits.
\begin{proof}[Proof of \cref{G_FH_respects_seq_unions}]
  Let $M \in \ulmfld^{G,f-fr}_{[G/H]}$ be an $f$-framed $\OG$-manifold and \( M = \cup_{i=1}^{\infty} M_i \) a $G$-sequential union of $M$. 
  The assembly map \( \colim_i \int_{M_i} A \to \int_M A \) factors as a sequence of equivalences
  \begin{align*}
    \colim_i \int_{M_i} A & = \colim_i \left( \ul{G/H}-\colim \left( \ul\disk^{G,f-fr}_{/\ul{M_i}} \to \ul\disk^{G,f-fr} \ultimes \ul{G/H} \xto{A \ultimes id} \ul\C \ultimes \ul{G/H} \right) \right) \\
    & \simeq \ul{G/H}-\colim \left( \colim_i \left( \ul\disk^{G,f-fr}_{/\ul{M_i}} \to \ul\disk^{G,f-fr} \ultimes \ul{G/H} \xto{A \ultimes id} \ul\C \ultimes \ul{G/H} \right) \right) \\
    & \iso \ul{G/H}-\colim \left( \ul\disk^{G,f-fr}_{/\ul{M}} \to \ul\disk^{G,f-fr} \ultimes \ul{G/H} \xto{A \ultimes id} \ul\C \ultimes \ul{G/H} \right) = \int_M A,
  \end{align*}
  where the second equivalence is induced by the equivalence \( \colim \ul\disk^{G,f-fr}_{/\ul{M_i}} \iso \ul\disk^{G,f-fr}_{/\ul{M}} \) of \cref{lem:disks_over_Mi}.
\end{proof}

\section{Axiomatic characterization of $G$-factorization homology theories} \label{sec:axiomatic_characterization}
In this subsection we give an axiomatic characterization of $G$-factorization homology theories with values in a presentable $G$-symmetric monoidal $G$-category (\cref{def:G_SM_presentable_cat}), as $G$-symmetric monoidal functors that satisfy $G$-$\otimes$-excision (\cref{def:G_tensor_excision}) and respects $G$-sequential unions (\cref{def:F_respects_G_seq_unions}).

\begin{mydef} 
  Let \( \ul\C^\otimes \fib \GFin_* \) be a $G$-symmetric monoidal category and $B \to BO_n(G)$ a $G$-map, as in \cref{def:Gmfld_framed}. 
  An \emph{equivariant homology theory of $G$-manifolds} is a $G$-symmetric monoidal functor \( F \colon \ulmfld^{G,f-fr,\sqcup} \to \ul\C^\otimes \) which satisfies $G$-$\otimes$-excision and respects $G$-sequential unions.
  We denote the full subcategory of equivariant homology theories by $\mathcal{H}(\ulmfld^{G,f-fr},\ul\C) \subset \Fun^\otimes_{G} (\ulmfld^{G,f-fr}, \ul\C)$.
\end{mydef}

The main result in this subsection is the following characterization of $G$-factorization homology.
\begin{thm} \label{thm:G_FH_axiomatic_characterization}
  Let \( \ul\C^\otimes \fib \GFin_* \) be a presentable $G$-symmetric monoidal category. 
  Then the full subcategory $\mathcal{H}(\ulmfld^{G,f-fr},\ul\C) \subset \Fun^\otimes_{G} ( \ulmfld^{G,f-fr}, \ul\C)$ is spanned by objects for which the counit map of the adjunction 
  \[
    \diag{
      (\iota^\otimes)_! \colon \Fun^\otimes_{G} (\ul\disk^{G,f-fr},\ul\C)  \ar@/_/[r] &   \Fun^\otimes_{G} (\ulmfld^{G,f-fr}, \ul\C) \noloc (\iota^\otimes)^\ast \ar@/_/[l]  \\
    }
  \]
  of \eqref{operadic_Lan_extends_Lan} is an equivalence.
  In particular, the adjunction restricts to an equivalence 
  \begin{align*}
    (\iota^\otimes)_! \colon \Fun^\otimes_{G} (\ul\disk^{G,f-fr}, \ul\C) \iso \mathcal{H}(\ulmfld^{G,f-fr},\ul\C), \quad
     A \mapsto \int_- A
  \end{align*}
  sending an $f$-framed $G$-disk algebra $A$ to $G$-factorization homology with coefficients in $A$. 
\end{thm}

\begin{proof}
  Let $A$ be a $G$-disk algebra.
  By \cref{G_tensor_excision} and \cref{G_FH_respects_seq_unions} the functor 
  \[
    (\iota^\otimes)_!  \colon \Fun^\otimes_{G} (\ul\disk^{G,f-fr}, \ul\C) \to \Fun^\otimes_{G} (\ulmfld^{G,f-fr}, \ul\C) 
  \] 
  factors though the full $G$-subcategory \(\mathcal{H}(\ulmfld^{G,f-fr}, \ul\C) \subset \Fun^\otimes_{G} (\ulmfld^{G,f-fr}, \ul\C) \).
  
  On the other hand, let \( F \in \mathcal{H}(\ulmfld^{G,f-fr},\ul\C) \) be an equivariant homology theory of $G$-manifolds. 
  Denote by $A \colon \ul\disk^{G,f-fr,\sqcup} \to \ul\C^\otimes$ the restriction of $F$ along $\iota^\otimes$.
  We have to show that the counit \( \int_- A \to F \) is an equivalence.
  Since $F,\, \int_- A$ are $G$-symmetric monoidal functors it is enough to show that for every $f$-framed $\OG$-manifold 
  \(M \in \ulmfld^{G,f-fr} \)
  the counit map \( \int_M A \to F(M) \) is an equivalence in $\ul\C$.
  We proceed by induction.

  For $k=0,1,\ldots, n$ let  $\mathcal{F}_{\leq k} \subseteq \ulmfld^{G,f-fr}$ be the full $G$-subcategory of $f$-framed $\OG$-manifolds whose underlying $\OG$-manifold is of the form \( (M \times_{G/H} D \to G/H) \)
  where $G/H \in \OG$ is a $G$-orbit, $M \to G/H$ is a \emph{$k$-dimensional} $\OG$-manifold
  and $(D\to G/H)$ is a finite $G$-disjoint union of $(n-k)$-dimensional $G$-disks, i.e. 
  equivalent to $D \to U\to G/H$ where $U$ is a finite $G$-set and $D \to U$ is a $G$-vector bundle of rank $n-k$ (and therefore $U=\pi_0(D)$).
  
  We now prove that the counit map is an equivalence on objects of $\mathcal{F}_{\leq k}$ by induction on $k$. 

  For $k=0$ the underlying $\OG$-manifold of $M \in \mathcal{F}_{\leq 0}$ is simply a finite $G$-disjoint union of $G$-disks, $(D\to G/H)\in\Gdisk$, therefore 
  \( M \in \ul\disk^{G,f-fr} \) and
  \(  \int_M A \simeq A(M) = F(M) \), since $\iota_!$ is fully faithful $A$ is the restriction of $F$ along $\iota$. 

  For $k\geq 1$, let $N\in \mathcal{F}_{\leq k}$ with underlying $\OG$-manifold \( (M \times_{G/H} D \to G/H) \). 
  We show that the counit map \( \int_N A \to F(N) \) is an equivalence using equivariant Morse theory.
  In what follows we only consider $G$-submanifolds of $M \times_{G/H} D \to G/H$, which by \cref{Gmfld_framed_over_M} have an essentially unique $f$-framing induced from the inclusion into $N$.
  Therefore we omit the identification of such $G$-submanifolds with their $f$-framed lift to $\ulmfld^{G,f-fr}$.

  Choose a $G$-equivariant Morse function \( f \colon  M \to \R \) with $f^{-1}(-\infty,r]$ a compact $G$-submanifold for every $r\in \R$ (see \cite[thm. 4.10]{Wasserman}).
  Choose an increasing sequence of regular values \( r_0 < r_1 < r_2 <  \cdots \) such that $f^{-1}(-\infty,r_0)=\emptyset$, the interval \( (r_i,r_{i+1}) \) contains a single critical value and \( r_i \to \infty \).

  Let $M_i := f^{-1}( (-\infty,r_i) ),$ then $M= \bigcup\limits_{i=0}^\infty M_i$ and therefore 
  \( M \times_{G/H} D \= \left( \bigcup\limits_{i=0}^\infty M_i \right) \times_{G/H} D \= \bigcup\limits_{i=0}^\infty \left(  M_i \times_{G/H} D \right) \) is a $G$-sequential union of $(M\times_{G/H} D \to G/H)$ (\cref{def:G_seq_unions}).
  Since both $F\in\mathcal{H}(\ulmfld^{G,f-fr},\ul\C)$ and $\int_- A$ respect $G$-sequential unions (\cref{def:F_respects_G_seq_unions} and \cref{G_FH_respects_seq_unions}) we have 
  \( F(M \times_{G/H} D) \simeq \colim F(M_i \times_{G/H} D),\, \int_{M \times_{G/H} D} A \simeq \colim \int_{M_i \times_{G/H} D} A \).
  Therefore it is enough to prove that the counit map \( \int_{M_i \times_{G/H} D} A \to F( M_i \times_{G/H} D) \) is an equivalence, which we prove by induction on $i$. 
  
  Let $\overline{M_i} := f^{-1}(-\infty,r_i]$.
  Since $\overline{M_i}$ is compact $\overline{M_{i+1}} \setminus M_i$ has only a finite number of critical orbits, \( x_j \colon  W_j \cof M,\, j=1,\ldots,s \).
  Note that the tangent bundle \( T_{x_j} M \to W_j \) over the critical orbit $x_j$ is a $G$-vector bundle which decomposes as a direct sum of two $G$-bundles \( T_{x_j} \= P_j \oplus E_j \) on which the Hessian is negative definite (called the index $E_j$) and positive definite (called the co-index $P_j$). 

  By \cite[thm. 4.6]{Wasserman} $\overline{M_{i+1}}$ is equivariantly diffeomorphic to $\overline{M_i}$ with $s$ handle-bundles $N_1,\ldots,N_s$ disjointly attached, where the handle-bundle \( N_j := \cDbun(P_j) \times_{W_j} \cDbun(E_j) \) is the fiberwise product of the closed unit disk bundles \( \cDbun(P_j) \to W_j, \, \cDbun(E_j) \to W_j\), attached to $M_i$ along \( \cDbun(P_j) \times_{W_j} \Sbun(E_j) \) where \( \Sbun(E_j) \to W_j \) is the unit sphere bundle of the negative definite $G$-subbundle (the index). 

  Since the handle-bundles are attached disjointly and $F,\int_- A$ are $G$-symmetric monoidal we can reduce to the case of a single handle-bundle by attaching one handle-bundle at a time. 
  Therefore we assume that there is a single critical orbit $x \colon  W \cof M$ in $\overline{M_{i+1}} \setminus M_i$ with \( T_x M \= P\oplus E \), and 
  \( \overline{M_{i+1}} \= \overline{M_i} \bigcup_{ \cDbun(P)\times_{W} \Sbun(E) } \left( \cDbun(P) \times_{W} \cDbun(E) \right) \).

  Let \( \Abun(E) \to W \) denote the unit annulus bundle of $E$, i.e. the open unit disk bundle minus the zero section.
  Note that $\Abun(E)$ is a $G$-tubular neighbourhood of $\Sbun(E)$, therefore 
  \( \overline{M_{i+1}} \= \overline{M_i} \bigcup_{ \cDbun(P)\times_{W} \Abun(E) } \left( \cDbun(P) \times_{W} \Abun(E) \right) \)
  a union of $k$-dimensional $G$-manifolds with boundary along a $k$-dimensional manifold with boundary.

  Discarding boundary points we see that the $M_{i+1}$ is equivariantly diffeomorphic to the union of $M_i$ with the $G$-manifold \( \Dbun(P) \times_{W} \Dbun(E) \) along the $G$-manifold \( \Dbun(P) \times_{W} \Abun(E) \).
  After taking fibered product with the fibration map \( D \to G/H \) we have 
  \begin{align} \label{morse_handle}
    M_{i+1} \times_{G/H} D \= ( M_i \times_{G/H} D ) \bigcup_{ \left( (\Dbun(P) \times_{W} \Abun(E)) \times_{G/H} D \right)} \left( (\Dbun(P) \times_{W} \Dbun(E)) \times_{G/H} D \right) .
  \end{align}
  This decomposition has the following properties:
  \begin{enumerate}
  \item The decomposition of \cref{morse_handle} is in fact a $G$-collar decomposition.
    Intuitively, the codimension one $G$-submanifold \( (\Dbun(P) \times_W \Sbun(E) ) \times_{G/H} D \) splits the handle bundle of \cref{morse_handle} to two $G$-submanifolds, $M_i$ and the handle bundle.
    Explicitly, construct a $G$-collar decomposition by defining a $G$-invariant smooth $G$-invariant function $M_{i+1} \to [-1,1]$ for which the restriction to the open interval $(-1,1)$ is a manifold bundle as follows.
    Compose the $G$-diffeomorphism of \cref{morse_handle} with the restriction of the Morse function $f \colon M\to \R$ to the handle-bundle of \cref{morse_handle}, followed by a smooth function \( \Psi \colon  \R \to [-1,1] \) such that
    \begin{enumerate}
      \item it sends the closed interval \( (-\infty,a+\epsilon] \) to $-1$ for some small $\epsilon>0$. 
      \item it sends \( [c-\epsilon,\infty) \) to $1$ for $c$ the unique critical value of $f$ in the interval $[a,b]$.
      \item it has a positive derivative in the open interval \( (a+\epsilon, c-\epsilon) \).
    \end{enumerate}
    Note that the fibers of $M_{i+1}\to [-1,1]$ over $(-1,1)$ are $\left( (\Dbun(P) \times_{W} \Sbun(E,r)) \times_{G/H} D \right)$ where $\Sbun(E,r)$ is the radius-$r$-sphere bundle, for various radii $r$. 
    \item The induced handle-bundle \( \left( (\Dbun(P) \times_{W} \Dbun(E)) \times_{G/H} D \to G/H \right) \in \Gdisk \) is a finite $G$-disjoint union of $G$-disks,
    since the open unit disk bundle of a $G$-vector bundle is equivalent to the entire vector bundle.
  \end{enumerate}

  We now distinguish between two cases, according to the rank of the bundle $E\to W$.
  \begin{enumerate}
    \item If the critical orbit $x$ has zero index, i.e. the Hessian is positive definite on $T_x M$, then $E \to W$ is a rank zero $G$-vector bundle, and its unit annulus \( \Abun(E) =\emptyset \) is empty. 
      In this case the $G$-collar decomposition of \cref{morse_handle}
      is a disjoint union 
      \begin{align*}
        M_{i+1} \times_{G/H} D \= ( M_i \times_{G/H} D ) \sqcup \left( (\Dbun(P) \times_{W} \Dbun(E)) \times_{G/H} D \right) .
      \end{align*}
      Since $F,\, \int_- A$ are $G$-symmetric monoidal functors we have
      \begin{align*}
        \int_{M_{i+1}} A \simeq \left( \int_{M_i} A \right) \otimes \left( \int_{  \left( (\Dbun(P) \times_{W} \Dbun(E)) \times_{G/H} D \right)} A \right), \\
        F(M_{i+1}) \simeq F(M_i) \otimes F\left( (\Dbun(P) \times_{W} \Dbun(E)) \times_{G/H} D \right) 
      \end{align*}
      where \( \left( (\Dbun(P) \times_{W} \Dbun(E)) \times_{G/H} D \right) \simeq \left( (P \times_{W} E) \times_{G/H} D \right) \in \Gdisk \) is a finite $G$-disjoint union of $G$-disks. 
      Therefore \( \int_{M_{i+1}} A \iso F(M_{i+1}) \) by induction on $i$. 

    \item Otherwise the critical orbit $x$ has positive index, i.e. $\rank(E)>0$. 
      In this case, \( \Abun(E) \= \Sbun(E) \times (-1,1) \) where $G$ acts trivially on the open interval $(-1,1)$, since the Morse function $f$ is $G$-invariant.
      It follows that 
      \begin{align*}
        (\Dbun(P) \times_{W} \Abun(E)) \times_{G/H} D  \=  \Abun(E) \times_{W} (P \times_{G/H} D) \=  \Sbun(E) \times_W ( (-1,1) \times P \times_{G/H} D) ,
      \end{align*}
      hence \( (\Sbun(E)\to W \to G/H) \) is a $G$-manifold of dimension 
      \[ 
        \dim \Sbun(E) = \rank(E)-1 \leq \dim M -1 = k-1, 
      \]
      so we have 
      \( (\Dbun(P) \times_{W} \Abun(E)) \times_{G/H} D  \in \mathcal{F}_{k-1} \).
      It follows by induction on $k$ that the counit map \( \int_{(\Dbun(P) \times_{W} \Abun(E)) \times_{G/H} D} A \iso  F((\Dbun(P) \times_{W} \Abun(E)) \times_{G/H} D) \) is an equivalence. 

      The $G$-functor $\int_- A$ satisfies $G$-$\otimes$-excision by \cref{G_tensor_excision} and $F$  satisfies $G$-$\otimes$-excision by assumption, therefore applying  $F, \int_- A$ to the $G$-collar decomposition of \cref{morse_handle} we get
      \begin{align*}
        & F(M_i \times_{G/H} D) \otimes_{F\left( (\Dbun(P) \times_{W} \Abun(E)) \times_{G/H} D \right)} F\left( (\Dbun(P) \times_{W} \Dbun(E)) \times_{G/H} D \right) \iso F(M_{i+1}) , \\
        & \left( \int_{(M_i \times_{G/H} D)} A \right) \otimes_{\left( \int_{ (\Dbun(P) \times_{W} \Abun(E)) \times_{G/H} D } A \right)} \left( \int_{(\Dbun(P) \times_{W} \Dbun(E)) \times_{G/H} D } \right) A \iso \int_{M_{i+1}} A 
      \end{align*}
      and by induction on $i$ the map \( \int_{M_{i+1}} A \iso F(M_{i+1}) \) is an equivalence. \qedhere
  \end{enumerate}
\end{proof}

\section{Equivariant versions of Hochschild homology} \label{sec:applications}
As an application of the $G$-$\otimes$-excision property (\cref{G_tensor_excision}) we describe two variants of topological Hochschild homology using $G$-factorization homology.

\subsection{Real topological Hochschild homology as $G$-factorization homology} \label{sec:Real_THH}
Let $C_2$ denote the cyclic group of order two and let $\sigma$ be its one dimensional sign representation.

\paragraph{The structure of an $\EE_\sigma$-algebra in $\ul\Sp^{C_2}$.}
Let us first describe the algebraic structure of an $\EE_\sigma$-algebra $A$ in $\ul\Sp^{C_2}$.
We will use this description in the proof of \cref{prop:Real_THH_as_GFH}.
\begin{itemize}
  \item 
    By \cref{EV_algs_vs_V_disk_algs} we have an equivalence
    \[
      Alg_{\EE_\sigma}(\ul\Sp^{C_2}) \simeq \Fun^\otimes_G ( \ul\disk^{C_2,\sigma-fr}, \ul\Sp^{C_2}) ,
    \]
    so $A$ corresponds to a $C_2$-symmetric monoidal functor $A \colon \ul\disk^{C_2,\sigma-fr} \to \ul\Sp^{C_2}$.
    In particular the $G$-symmetric monoidal functor $A$ restricts to symmetric monoidal functors 
    \begin{equation} \label{eq:E_sigma_functors}
      A_{[C_2/C_2]} \colon \ul\disk^{C_2,\sigma-fr}_{[C_2/C_2]} \to \Sp_{C_2}, \quad A_{[C_2/C_2]} \colon \ul\disk^{C_2,\sigma-fr}_{[C_2/e]} \to \Sp.
    \end{equation}
  \item
    By abuse of notation, we use $A$ to denote the ``underlying'' genuine $C_2$-spectrum,
    \[
        A_{[C_2/C_2]} (\R^\sigma) \in \Sp^{C_2} ,
    \]
    where \( \R^\sigma\in \ul\disk^{C_2,\sigma-fr} \) is the one dimensional sign representation of $C_2$, considered as a $\sigma$-framed $C_2$-manifold.
  \item
    Unwinding the definitions we see that $\ul\disk^{C_2,\sigma-fr}_{[C_2/e]}$ is equivalent to the $\infty$-category $\disk^{fr}_1$ of \cite[rem 2.10]{AF}.
    Since $A \colon \ul\disk^{C_2,\sigma} \to \ul\Sp^{C_2}$ is a $G$-functor it is compatible with the forgetful functors 
    \begin{equation*} 
      Res^{C_2}_e  \colon \ul\disk^{C_2,\sigma}_{[C_2/C_2]} \to \ul\disk^{C_2,\sigma}_{[C_2/e]} \simeq \disk^{fr}_1 , \quad Res^{C_2}_e \colon \Sp_{C_2} \to \Sp,
    \end{equation*}
    therefore $A_{[C_2/e]}(\R^1) = A_{[C_2/e]}(Res^{C_2}_e \R^\sigma) \simeq Res^{C_2}_e A_{[C_2/C_2]}(\R^\sigma) = Res^{C_2}_e A$.
  \item 
    Observe that $Res^{C_2}_e A$ is endowed with a structure of an $\EE_1$-sing spectrum.
    To see this, recall that $\R^1\in \disk^{fr}_1$ is an $\EE_1$-algebra in $\disk^{fr}_1$, which induces an equivalence between  the symmetric monoidal envelope of $\EE_1$ and $\disk^{fr}_1$ (see \cite[prop. 2.12]{AyalaFrancisTanaka}).

  \item
    Let 
    \[
      \sqcup_{C_2} \R^1 \in \ul\disk^{C_2,\sigma-fr}_{[C_2/C_2]}, \quad \sqcup_{C_2} \R^1 = C_2 \times \R^1
    \]
    denote the topological induction of $\R^1 \in \disk^{fr}_1$.
    The compatibility of the $G$-symmetric monoidal functor  $A \colon \ul\disk^{C_2,\sigma-fr} \to \ul\Sp^{C_2}$ with with topological induction and the Hopkins-Hill-Ravenel norm, 
    \[ 
       \sqcup_{C_2}  \colon \disk^{fr}_1 \simeq \ul\disk^{C_2,\sigma}_{[C_2/C_2]} \to \ul\disk^{C_2,\sigma}_{[C_2/e]} , \quad 
       N^{C_2}_e  \colon \Sp \to \Sp_{C_2},
    \]
    implies that $A_{[C_2/C_2]}(\sqcup_{C_2} \R^1) \simeq N^{C_2}_e A $.
  \item 
    Note that $N^{C_2}_e A$ is an $\EE_1$-algebra in $\Sp_{C_2}$, since $N^{C_2}_e  \colon \Sp \to \Sp_{C_2}$ is a symmetric monoidal functor and $Res^{C_2}_e A$ is an $\EE_1$-ring spectrum.
  \item
    The ``underlying'' $C_2$-spectrum $A$ has the structure of a module over $N^{C_2} A$.
    To see this structure, note that an equivariant oriented embeddings 
    \[
       \left( \sqcup_{C_2} \R^1 \right) \sqcup \R^\sigma  \cof \R^\sigma
     \]
    induces a map 
    \[
      N^{C_2}_e A  \otimes A \to A.
    \]
\end{itemize}

\begin{prop} \label{prop:Real_THH_as_GFH}
  For $A$ an $\EE_\sigma$-algebra
  in $\ul\Sp^{C_2}$ 
  there is an equivalence of genuine $C_2$-spectra
  \begin{align*}
    \int_{S^1} A \simeq A \otimes_{N_e^{C_2} A } A .
  \end{align*}
  where $C_2$ acts on $S^1$ by reflection.
\end{prop}
\begin{proof}
  Consider the $C_2$-collar gluing $S_1 = \R^\sigma \cup_{\sqcup_{C_2} \R^1} \R^\sigma$ into two hemispheres, where each hemisphere is reflected onto itself by the action of $C_2$.
  Note that the intersection $\sqcup_{C_2} \R^1$ consists of two segments interchanged by the action of $C_2$.
  Applying \cref{G_tensor_excision} we get an equivalence of genuine $\C_2$-spectra 
  \begin{align*}
    \int_{S^1} A \simeq \left( \int_{ \R^\sigma } A \right) \otimes_{ \left( \int_{ \sqcup_{C_2} \R^1 } A \right)  } \left( \int_{ \R^\sigma } A \right) 
      \simeq A \otimes_{N_e^{C_2} A } A .
  \end{align*}
\end{proof}
\begin{rem} \label{rem:THR_as_GFH}
  The tensor product $A \otimes_{N_e^{C_2} A } A$ appearing in \cref{prop:Real_THH_as_GFH} is equivalent to the derived smash product $A \wedge^{\mathbf{L}}_{N_e^{C_2} A } A$ of left and right $N_e^{C_2}$-modules.
  Dotto, Moi, Patchkoria and Reeh (\cite{Real_THH}) show that for $A$ a flat ring spectrum with anti-involution 
  there is a stable equivalence of genuine $C_2$-spectra
  \begin{align*}
    THR(A) \simeq A \wedge^L_{N_e^{C_2} A} A ,
  \end{align*}
  where $THR(A)$ is the B\"okstedt model for real topological Hochschild homology.

  By \cite[def. 2.1]{Real_THH} we can interpret a ring spectrum with anti-involution as an algebra over an operad $Ass^\sigma$ in $C_2$-sets. 
  Direct inspection shows $Ass^\sigma$ is equivalent to $G$-operad  $\mathcal{D}_\sigma$ of the little $\sigma$-disks
  \footnote{This also follows from a direct analysis of the mapping spaces of $\ulRep^{C_2,\sigma-fr,\sqcup}$, which are homotopically discrete.},
  whose genuine operadic nerve is $\EE_\sigma$.
  Regarding a flat ring spectrum with anti-involution $A$ as an $\EE_\sigma$-algebra in $\ul\Sp^{C_2}$,
  we can reinterpret \cref{prop:Real_THH_as_GFH} as an equivalence 
  \[ 
    \int_{S^1} A \simeq THR(A)
  \]
  of genuine $C_2$-spectra.
\end{rem}

\subsection{Twisted Topological Hochschild Homology of genuine $C_n$-ring spectra} \label{sec:Twisted_THH}

We start with a general lemma relating trivially framed $G$-disk algebras to $\EE_n$-algebras.
Let $G$ be a finite group acting trivially on $\R^n$, and $\ulmfld^{G,\R^n-fr}$ the $G$-category of trivially framed $G$-manifolds.
\begin{lem}
  Let $\ul\C^\otimes \fib \GFin_*$ be a $G$-symmetric monoidal $\infty$-category.
  The $\infty$-category $ \Fun_G^\otimes(\ul\disk^{G,\R^n-fr}, \ul\C)$ of trivially framed $G$-disk algebras in $\C$ is equivalent to the $\infty$-category $Alg_{\EE_n}(\ul\C_{[G/G]})$ of $\EE_n$-algebras in the fiber $\ul\C_{[G/G]}$. 
\end{lem}

\paragraph{The structure of a trivially framed $C_n$-disk algebra}
Let $C_n$ the cyclic group of order $n$ and $\ul\C = \ul\Sp^{C_n}$, the $C_n$-$\infty$-category of genuine $C_n$-spectra.
We will use the following an explicit description of the trivially framed $C_n$-disk algebra corresponding to $A$.
The $C_n$-functor \( A \colon \ul\disk^{C_n,\R^n-fr,\sqcup} \to \ul\Sp^{C_n} \) sends 
\begin{align*}
  \forall H < C_n: \quad A_{[C_n/H]} \colon \sqcup_{C_n/H} \R^1 \mapsto N_H^{C_n}(A) \in \Sp_{H} ,
\end{align*}
where $N_H^{C_n}(A)$ denotes the Hill-Hopkins-Ravenel norm applied to the restriction of the genuine $C_n$-spectrum $A \in \Sp_{C_n}$ to $\Sp_H$.
In particular,  \( A \colon \ul\disk^{C_n,\R^n-fr,\sqcup} \to \ul\Sp^{C_n} \) sends $\R^n$ with trivial $C_n$-action to $A\in \Sp_{C_n}$ and the topological induction $\sqcup_{C_n} \R^1 = C_n \times \R^1 \in$ to $N_e^{C_n}(A)\in \Sp_{C_n}$.

We will need some notation for our next statement. 
Let $A$ be an $\EE_n$-ring spectrum in $\Sp_{C_n}$.
Define an $A-A^{op}$-bimodule structure on $A\in \Sp_{C_n}$ with ``twisted'' left multiplication, given by first acting on the scalar by the generator $\tau\in C_n$:
\[
  A \otimes A^\tau \otimes A \to A^\tau , \quad x \otimes a \otimes y \mapsto \tau x \cdot a \cdot y.
\]
We denote this ``twisted'' $A-A$-bimodule by $A^\tau$.
Let $THH(A;A^\tau)$ denote the topological Hochschild homology of $A$ with coefficients in $A^\tau$.

\begin{prop} \label{prop:twisted_THH_via_GFH}
  Let $A$ be an $\EE_1$-ring spectrum in $\Sp_{C_n}$,
  and $C_n \curvearrowright S^1$ be the standard action. 
  Then there exists an equivalence of spectra 
  \begin{align*}
    \left( \int_{S^1} A \right)^{\Phi C_n} \simeq  THH(A;A^\tau) . 
  \end{align*}
  In particular,  $THH(A;A^\tau)$ admits a natural circle action.
\end{prop}
\begin{proof}
  Consider $S^1$ as the $n$-fold covering space $p \colon S^1 \to S^1$, with the standard $C_n$-action given by deck transformations.
  Let \( S^1 = U \cup_{U \cap V} V \) be the standard collar decomposition of the base $S^1$ by hemispheres. 
  Construct a $C_n$-collar decomposition  \( S^1 = p^{-1}(U) \cup_{p^{-1}(U\cap V)} p^{-1}(V) \) of the covering space by taking preimages.
  Observe that the pieces of this $C_n$-collar decomposition are given by topological induction,
  \begin{align*}
    & p^{-1}(U) = \sqcup_{C_n} U \= \sqcup_{C_n} \R^1, \quad p^{-1}(V) = \sqcup_{C_n} V \= \sqcup_{C_n} \R^1, \\
    & p^{-1}(U\cap V) = \sqcup_{C_n} (U \cap V) \= \sqcup_{C_n} (\R^1 \sqcup \R^1 ) = ( \sqcup_{C_n} \R^1 ) \sqcup ( \sqcup_{C_n} \R^1 ) .
  \end{align*}
  Therefore by $C_n$-$\otimes$-excision
  \begin{align*}
    \int_{S^1} A & \simeq \left( \int_{p^{-1}(U)} A \right) \bigotimes\limits_{ \int_{p^{-1}(U\cap V)} A } \left( \int_{p^{-1}(V)} A \right)
    \simeq  \left( \int_{\sqcup_{C_n} \R^1 } A \right) \bigotimes\limits_{ \int_{ ( \sqcup_{C_n} \R^1 ) \sqcup ( \sqcup_{C_n} \R^1 )} A } \left( \int_{\sqcup_{C_n} \R^1 } A \right) \\
    & \simeq  (N_e^{C_n} A) \bigotimes\limits_{ (N_e^{C_n} A) \otimes (N_e^{C_n} A)^{op} } (N_e^{C_n} A)^\tau .
  \end{align*}
  
  Let us pause and explain the superscript decorations in the last term.
  The $\left(\int_{p^{-1}(U \cap V)} A \right)$-module structure of $\int_{p^{-1}(U)} A$ is induced by the inclusion $p^{-1} (U \cap V) \cof p^{-1}(U)$.
  When we identify $p^{-1} (U \cap V) \= ( \sqcup_{C_n} \R^1 ) \sqcup ( \sqcup_{C_n} \R^1 )$ the module structure on $\int_{p^{-1}(U)} A$ is naturally identified with an $(N_e^{C_n} A)-(N_e^{C_n} A)$-bimodule structure, or equivalently a right $(N_e^{C_n} A) \otimes (N_e^{C_n} A)^{op}$-module structure. 
  Similarly, $\int_{p^{-1}(V)} A$ is naturally a left $N_e^{C_n}(A)-N_e^{C_n}(A)^{op}$-module. 
  However the left module structure is induced by an embedding $ \sqcup_{C_n} \R^1 \cof p^{-1} V$ which defers from the standard embedding (the topological induction of $\R^1 \cof V$) by a deck transformation. 
  Therefore the left multiplication is ``twisted'',  i.e. given by first acting on the scalar by the generator $\tau\in C_n$.
  In order to remember this twist in the module structure of the right hand side we add the superscript $\tau$.

  Next we take geometric fixed points of $\int_{S^1} A$. 
  Since the geometric fixed points functor $(-)^{\Phi C_n} \colon  \Sp_{C_n} \to \Sp$ is symmetric monoidal and preserve homotopy colimits 
  \begin{align*}
    \left( \int_{S^1} A \right)^{\Phi C_n} \simeq  (N_e^{C_n} A)^{\Phi C_n} \bigotimes\limits_{ (N_e^{C_n} A)^{\Phi C_n} \otimes ((N_e^{C_n} A)^{op})^{\Phi C_n} } ( (N_e^{C_n} A)^\tau )^{\Phi C_n} \simeq A \bigotimes_{A \otimes A^{op} } A^\tau .
  \end{align*}
  The right hand side is equivalent to the topological Hochschild homology $THH(A;A^\tau)$ of $A\in \Sp$ with coefficients in the $A-A$-bimodule $A^\tau$.

  Finally, we describe the natural circle action on $\left( \int_{S^1} A \right)^{\Phi C_n}$.
  Note that the automorphism space of  $S^1 \in \ulmfld^{C_n,1-fr}_{[C_n/C_n]}$ acts on $S^1$, so by functoriality it induces a natural action on $\left( \int_{S^1} A \right)^{\Phi C_n}$.
  The endomorphism space of $S^1 \in \ulmfld^{C_n,1-fr}_{[C_n/C_n]}$ is the space of $C_n$-equivariant oriented embeddings $Emb^{C_n}(S^1,S^1)$.
  In particular the endomorphism space $S^1 \in \ulmfld^{C_n,1-fr}_{[C_n/C_n]}$ includes rotations of $S^1$, therefore the circle group acts on $\int_{S^1} A$ by rotations, and by functoriality on $\left( \int_{S^1} A \right)^{\Phi C_n}$.
\end{proof}
\begin{rem}
  The inclusion of the circle group into $Emb^{C_n}(S^1,S^1)$ is in fact a deformation retract.
\end{rem}
\begin{rem}
  This theorem can be seen as an instance of a more general principle: factorization homology with local coefficients on a manifold $M$ can be constructed as the fixed points of $G$-factorization homology on a cover of $M$.
\end{rem}

\paragraph{Relation to the relative norm construction} 
The spectrum $THH(A;A^\tau)$ and its circle action have been used to define the relative norm in  \cite{TC_via_the_norm}. 
In order to give a precise statement we recall the notation of \cite{TC_via_the_norm}.

Fix $U$ a complete universe of the circle group (in the sense of orthogonal spectra), and define a complete $C_n$-universe $\tilde{U}=\iota^*_{C_n}U$.
Let $R$ be an associative ring orthogonal $C_n$-spectrum indexed on the universe $\tilde{U}$.
Let $I^{\R^\infty}_{\tilde{U}}, I^U_{\R^\infty}$ denote the ``change of universe'' functors.
The relative norm \( N^{S^1}_{C_n} R \) of \cite[def. 8.2]{TC_via_the_norm} is the genuine $S^1$-spectra defined as 
\[
  I^U_{\R^\infty} \left| N^{cyc,C_n}_{\wedge}( I^{\R^\infty}_{\tilde{U}} R ) \right| ,
\]
where $N^{cyc,C_n}_{\wedge}(-)$ is the ``twisted cyclic bar construction'' of \cite[def. 8.1]{TC_via_the_norm}.

Note that the geometric realization $| N^{cyc,C_n}_{\wedge} (I^{\R^\infty}_{\tilde{U}} R) |$ is equivalent to $THH(R;R^\sigma)$, computed using the standard bar resolution.
By \cref{prop:twisted_THH_via_GFH} there exists an equivalence of spectra 
\begin{align*}
  \left( \int_{S^1} R \right)^{\Phi C_n} 
  \simeq \left| N^{cyc,C_n}_{\wedge} \left( I^{\R^\infty}_{\tilde{U}} R \right) \right| ,
\end{align*}
where one the left hand side we consider $R$ as an $\EE_1$-algebra in $\Sp_{C_n}$.

Moreover, by inspection the above equivalence respects the circle action, hence after applying the change of universe functor $I^{U}_{\R^{\infty}}$ we get an equivalence of genuine $S^1$-spectra
\[
  N^{S^1}_{C_n} R \simeq I^U_{\R^\infty} \left( \left( \int_{S^1} R \right)^{\Phi C_n} \right) .
\]

\begin{appendices}
  \appendix
  
\section{The Moore over category} \label{MooreOverCat} 
Let $\C$ be a topological category and $x\in\C$ an object.
Denote by \( N(\C) \in \Cat_\infty \) the coherent nerve of $\C$, and by \( N(\C)_{/x} \in \Cat_\infty \) the over category.
Note that $N(\C)_{/x}$ is not equivalent to the coherent nerve of \( \C_{/x} \), the topological over category:
both have the same objects, but a point in \( \Map_{\C_{/x}} ( y_0 \xto{f_0} x, y_1 \xto{f_1} x ) \) is an given by a map \( h \in \Map_\C(y_0,y_1) \) satisfying \( f_0 = f_1 \circ h \), 
while a point in \( \Map_{N(\C)_{/x}} ( y_0 \xto{f_0} x, y_1 \xto{f_1} x ) \) is given by a map \( h \in \Map_\C(y_0,y_1) \) together with a \emph{path} in $ \Map_\C(y_0, x) $ from $f_0$ to $f_1 \circ h$.
Nevertheless, it is useful to have a topological category whose coherent nerve is equivalent to \( N(\C)_{/x} \). 
Of course, this could be achieved by applying homotopy coherent realization to $N(\C)_{/x}$, but unwinding the construction one sees that an explicit description of topological category involves a lot of simplicial combinatorics. 
In what follows we construct a topological category \( \C^\Moore_{/x} \) whose coherent nerve is equivalent to \( N(\C)_{/x} \), which avoids simplicial combinatorics. 

An obvious candidate for the mapping space $\Map_{\C^\Moore_{/x}} ( y_0 \xto{f_0} x, y_1 \xto{f_1} x )$ is the space  of maps $h \colon y_0 \to y_1$ in $\C$ together with a path from $f_0$ to $ f_1 \circ f$ in \( \Map_\C (y_0,x) \), formally given by the fiber product \( \Map_{\C} (y_0,y_1) \times_{\Map_\C (y_0,x)} P(\Map_\C(y_0, x)) \). 
However, one runs into trouble when trying to define composition functions which are strictly associative, since the composition action uses concatenation of paths.
The problem of defining strictly associative concatenation of paths has a classical solution, namely replacing the space of paths with the homotopy equivalent space of Moore paths.
Defining the mapping space $\Map_{\C^\Moore_{/x}} ( y_0 \xto{f_0} x, y_1 \xto{f_1} x )$ using Moore paths leads to a simple construction of a topological category $\C^\Moore_{/x}$, the Moore over category (\cref{def:Moore_over_cat}), whose coherent nerve is equivalent to \( N(\C)_{/x} \) (\cref{Moore_is_overcat}). 

We first recall the definition of the Moore path space and concatenation of Moore paths.
Let $X$ be a topological space. 
The Moore path space of $X$ is the subspace
\begin{align*}
  M(X) \subset [0,\infty) \times X^{[0,\infty)}, \quad M(X) = \set{ (r,\gamma) | \text{ the restriction } \gamma|_{[r,\infty)} \text{ is a constant function} },
\end{align*}
where $X^{[0,\infty)}$ is the space of functions $[0,\infty) \to X$ endowed with the compact-open topology. 
The ``starting point'' and ``finishing point'' fibrations \( \alpha, \omega  \colon  M(X) \fib X \) are the given by \( \alpha(r,\gamma) = \gamma(0), \,\omega(r,\gamma) = \gamma(r) \).
Moreover, the ``ends points'' map \( (\alpha,\omega) \colon  M(X) \fib X \times X \) is also a Serre fibration.
Concatenation of Moore paths is defined by
\begin{align*}
  \ast \colon  M(X) \times_X M(X) \to M(X) , \quad  (r_0,\gamma_0) \ast (r_1,\gamma_1) =
  \left( r_0+r_1, t \mapsto 
  \begin{cases}
    \gamma_0(t)  & t \leq r_0 \\
    \gamma_1(t-r_0) & t \geq r_0
  \end{cases} \right) .
\end{align*}
It is straightforward to verify that concatenation of paths is associative, i.e. 
\[ \left( (r_0,\gamma_0) \ast (r_1,\gamma_1) \right) \ast (r_2,\gamma_2) = (r_0,\gamma_0) \ast \left( (r_1,\gamma_1) \ast (r_2,\gamma_2) \right) . \] 
For $x\in X$  a point, the ``constant instant Moore path'' \( (0, t\mapsto x) \in M(X) \) is a neutral element for concatenation.

With the definition of Moore paths at hand, we can define the Moore path category.
\begin{mydef} \label{def:Moore_over_cat}
  Let $\C$ be a topological category and $x\in \C$ an object. 
  Define a topological category \( \C^\Moore_{/x} \) with objects arrows \( f \colon  y \to x\), i.e pairs $(y,f)$ where \(y\in \C, \, f\in \Map_\C(y,x) \), and morphism spaces $\Map_{\C^\Moore_{/x}} ( y_0 \xto{f_0} x, y_1 \xto{f_1} x ) $ given by the fiber products
  \begin{align*}
    \set{f_0} \times_{\Map_\C(y_0,x)} & M(\Map_\C(y_0,x)) \times_{\Map_\C(y_0,x), (f_1 \circ_\C (-))} \Map_\C(y_0,y_1) \\
    & = \set{ ((r,\gamma),h) \, | \, \gamma(0)=f_0, \, \gamma(r) = f_1 \circ_\C h } .
  \end{align*}
  Define composition in $\C^\Moore_{/x}$ by
  \begin{align*}
    \circ \colon  & \Map_{\C^\Moore_{/x}} ( y_0 \xto{f_0} x, y_1 \xto{f_1} x ) \times \Map_{\C^\Moore_{/x}} ( y_0 \xto{f_1} x, y_1 \xto{f_2} x )  \to \Map_{\C^\Moore_{/x}} ( y_0 \xto{f_0} x, y_1 \xto{f_2} x ) \\
    & \left( ((r,\gamma),h) , ((r',\gamma'),h') \right)  \mapsto \left( (r,\gamma) \ast (r', \gamma' \circ_\C h), h' \circ_\C h \right)
  \end{align*}
  and identity of $f \colon y\to x$ by \( ((0, t \mapsto f),id_y ) \in \Map_{\C^\Moore_{/x}} ( y \xto{f} x, y \xto{f} x) \), using the constant instant Moore path at $f$. 
  We call $\C^\Moore_{/x}$ the \emph{Moore over category} of $\C$ over $x$.
\end{mydef}

\begin{observation} \label{overcat_map_is_h_fiber}
  The mapping space $\Map_{\C^\Moore_{/x}}$ is the homotopy fiber of  
  \[ f_1 \circ (-)  \colon  \Map_\C(y_0,y_1) \to \Map_\C (y_0,x) . \]
\end{observation}

\begin{rem}
  If the mapping spaces \( \Map_\C(y,x) \) of $\C$ has a smooth structure one can replace the Moore spaces of continuous Moore paths by spaces of piecewise smooth Moore paths, without changing the $\infty$-category represented by $N(\C^\Moore_{/x})$. 
\end{rem}

\begin{lem}
  The coherent nerve of the Moore category \(\C^\Moore_{/x} \) has a terminal object  \( (x \xto{=} x) \in \C^\Moore_{/x} \). 
\end{lem}
\begin{proof}
  For every object \( (y\xto{f} x) \in \C^\Moore_{/x} \) the mapping space \( \Map_{\C^\Moore_{/x}} (y \xto{f} x, x \xto{=} x ) \) is the space of Moore paths in $\Map_\C (y,x)$ starting at ${f}$, a contractible space. 
\end{proof}

Define a functor of topological categories $U \colon  \C^\Moore_{/x} \to \C$ sending \( U \colon  (y \xto{f} x) \mapsto y \) and 
\begin{align*}
  U \colon  \Map_{\C_{/x}} ( y_0 \xto{f_0} x , y_1 \xto{f_1} x) \to \Map_\C(y_0,y_1) ,\quad U \colon ((r,\gamma),h) \mapsto h
\end{align*}
on mapping spaces by projection.

\begin{lem}
  The induced map of coherent nerves $N(U) \colon N(\C^\Moore_{/x}) \to N(\C)$ is a right fibration.
\end{lem}
\begin{proof}
  First we observe that $N(U)$ is an inner fibration.
  For each pair of objects $(y_0 \xto{f_0} x), (y_1\xto{f_1} x) \in \C^\Moore_{/x}$ the map $ U \colon \Map_{\C_{/x}} ( y_0 \xto{f_0} x , y_1 \xto{f_1} x) \to \Map_\C(y_0,y_1)$ is a pullback of the ``end points'' fibration $(\alpha,\omega)$ along $\set{f_0} \times (f_1 \circ_\C (-))$, and therefore a fibration. 
  By \cite[prop. 2.4.1.10 (1)]{HTT} it follows that $N(U)$ is an inner fibration. 

  By \cite[prop. 2.4.2.4]{HTT} we need to show that every morphism \( ((r,\gamma),h)  \colon  (y_0 \xto{f_0} x) \to (y_1 \xto{f_1} x) \) in $\C^\Moore_{/x}$ is $U$-Cartesian.
  By \cite[prop. 2.4.1.10 (2)]{HTT} we have to show that for every $(y\xto{f} x)$ in $\C^\Moore_{/c}$ the diagram 
  \begin{align*}
    \diag{
      \Map_{\C^\Moore_{/x}}(f,f_0) \ar@{->>}[d]^{U} \ar[rr]^{((r,\gamma),h) \circ -} & & \Map_{\C^\Moore_{/x}}(f,f_1) \ar@{->>}[d]^{U} \\
      \Map_\C(y,y_0) \ar[rr]^{h\circ_\C -} & & \Map_\C (y,y_1)
    }
  \end{align*}
  is homotopy Cartesian.
  We show that the induced map between the fibers is a homotopy equivalence. 
  For every point  \( (h' \colon y \to y_0) \in \Map_\C (y,y_0) \), the fiber over $h'$ is the space 
  of Moore paths in $\Map_\C(y,x)$ starting at $f$ and ending at $f_0\circ h$,
  the fiber over $h\circ h'$ is the space 
  of Moore paths in $\Map_\C(y,x)$ starting at $f$ and ending at $f_1\circ h' \circ h$,
  and the map between the fibers is given by concatenation with the Moore path \( (r,\gamma \circ h') \) starting at $f_0 \circ h' $ and ending at $f_1 \circ h \circ h'$, a homotopy equivalence. 
\end{proof}
\begin{cor} \label{Moore_is_overcat}
  Let $\C$ be a topological category and $x\in \C$ an object.
  The coherent nerve $N(\C^\Moore_{/x})$ is equivalent to the $\infty$-over category \( N(\C)_{/x} \).
\end{cor}
\begin{proof}
  The right fibration \( N(U)  \colon  N(\C^\Moore_{/x}) \to N(\C) \) takes the terminal object \( (x \xto{=} x) \in N(\C^\Moore_{/x}) \) to $x\in \C$.
  By \cite[prop. 4.4.4.5]{HTT} the right fibrations \( N(U)  \colon  N(\C^\Moore_{/x}) \to N(\C) \) and \( N(\C)_{/x} \fib N(\C) \) are equivalent fibrant objects of the contravariant model structure on \( sSet_{/N(\C)} \) (both right fibrations classify the representable functor \( \Map(-,x) \colon  \C^{op} \to \mathcal{S} \) ), and claim the follows.
\end{proof}

\section{The definition of a $G$-Symmetric Monoidal category} \label{sec:G_SM_cat}
This appendix contains no original results or definitions.
The notion of $G$-symmetric monoidal $\infty$-category, developed by Barwick, Dotto, Glasman, Nardin and Shah, is central to our treatment of $G$-factorization homology. 
For the convenience of the reader we include the definition here (see \cref{def:G_SM_cat}), which is equivalent to the definition given in \cite{Nardin_thesis}.

\paragraph{Parametrized join} First we recall the parametrized version of the join construction.
\begin{mydef}
  Let $S$ be an $\infty$-category. 
  Restricting along \( S \times \partial \Delta^1 \to S\times \Delta^1 \) defines a functor \( sSet_{/ S\times \Delta^1} \to sSet_{/ S\times \partial \Delta^1} \= sSet_{/ S \times\set{0}} \times sSet_{/ S \times\set{1}} \) which carries coCartesian fibrations over $S\times \Delta^1$ to coCartesian fibrations over $S\times \partial\Delta^1 = S \times \set{0} \coprod S \times \set{1}$.
  This functor has a right adjoint which is called the $S$-parametrized join and denoted by 
  \begin{align*}
    sSet_{/ S\times \Delta^1} \leftrightarrows sSet_{/ S\times \set{0} } \times sSet_{/ S\times \set{1} } \noloc \star_S .
  \end{align*}
\end{mydef}
By \cite[prop. 4.3]{Expose2}, if $\ul\C \fib S, \ul\D \fib S$ are coCartesian fibrations (i.e $S$-categories), then $\ul\C \star_S \ul\D \fib S$ is a coCartesian fibration. 

It follows from \cite[thm. 4.16]{Expose2} that the parametrized join carries coCartesian fibrations over $S\times \partial\Delta^1$ to inner fibrations over $S\times\Delta^1$ with coCartesian lifts over $S\times\partial \Delta^1$.

The parametrized join $X \star_S Y \to S \times \Delta^1$ of $ X \to S, \, Y \to S$ can be informally described as follows (see \cite[lem. 4.4]{Expose2}):
its restriction to $S\times\set{0}$ is $X\to S$, its restriction to $S\times\set{1}$ is $Y\to S$, and for each $s\in S$ its restriction to $\set{s}\times\Delta^1$ is the join $X_{[s]} \star Y_{[s]} $, where $X_{[s]},Y_{[s]}$ are the fibers of $X\to S,\, Y\to S$ over $s\in S$. 

Fact: for the case $Y=S$ one gets a coCartesian fibration \( X \star_S S \fib S \). 

\paragraph{Finite pointed $G$-sets} We denote by $\GFin_*$ the $G$-category of finite pointed $G$-sets of \cite[def. 4.12]{Expose4}.
An object $I \in \GFin_*$ over the orbit $G/H$ is a $G$-equivariant map \( I= (U \to G/H) \) from a finite $G$-set $U$.
A morphism in $\GFin_*$ over \(\varphi \colon G/K \to G/H\) is a span of the form
\begin{align*}
 \xymatrix{
   U \ar[d] & \ar[l] U' \ar[d] \ar[r] & V \ar[d] \\
   G/H & \ar[l]_{\psi} G/K \ar[r]^{=} & G/K
 }
\end{align*}
where the left square is a summand inclusion, i.e it induces an inclusion of $U'$ into the pullback \( \psi^* U = G/K \times_{G/H} U\).
The span above is a coCartesian edge if the left square is Cartesian and the map \( U' \to V \) is an isomorphism of finite $G$-sets (\cite[lem. 4.9, def. 4.12]{Expose4}).
We call the span above \emph{inert} if \( U' \to V \) is an isomorphism.

\begin{notation}
  Let \( G/K \in \OGop \) be an orbit.
  Denote by \( I_+(G/K) = (G/K \xto{=} G/K ) \in \GFin_* \) the finite pointed set given by the identity map of $G/K$.
\end{notation}
\begin{mydef} \label{corr_G_functor}
  Let \( I \in \GFin_*, \,  I = ( U \to G/H ) \) be a finite pointed $G$-set over $G/H$. 
  Recall that the left fibration 
  \[ \ul{G/H} = (\OGop)_{[G/H]/} \to \OGop, \quad (G/H \leftarrow G/K) \mapsto G/K \] classifies the representable functor \( \Hom(-,G/H) \colon \OGop \to Set \) (see \cite[ex. 2.4]{Expose1}).
  By Yoneda's lemma the set $(\GFin_*)_{[G/H]}$ of finite $G$-sets over $G/H$ is in bijection with the set of natural transformations \( Nat\left( \Hom(-,G/H) , (\GFin_*)_{(-)} \right) \),
  which in turn is in bijection with the set of $G$-functors \( \ul{G/H} \to \GFin_* \). 
  Define \( \sigma_{<I>}  \colon \ul{G/H} \to \GFin_* \) as the $G$-functor corresponding to $I$ under the bijection above. 
  Explicitly, $\sigma_{<I>}$ acts on objects by \( \sigma_{<I>} \colon ( G/H \xfrom{\varphi} G/K) \mapsto (\varphi^* U \to G/K)  \) .
\end{mydef}

\paragraph{The underlying $G$-categories of the $G$-diagram classified by $\ul\C^\otimes \fib \GFin_*$}
By straightening/unstraightening for $G$-categories (\cite[prop. 8.3]{Expose1}) the coCartesian fibration \( \ul\C^\otimes \fib \GFin_* \) corresponds to a $G$-functor \( \GFin_* \to \ul{\Cat}_{\infty, G} \), which we can interpret as a $\GFin_*$-shaped $G$-diagram in $\ul{\Cat}_{\infty,G}$.
The functor \( \GFin_* \to \ul{\Cat}_{\infty,G} \) assigns to each $I\in (\GFin_*)_{[G/H]}$ an object of \( (\ul{\Cat}_{\infty,G})_{[G/H]} = \Fun({\ul{G/H}},\Cat_\infty) \simeq (\Cat_\infty)_{/ \ul{G/H}}^{coCart} \) (see \cite[ex. 7.5]{Expose1}), i.e a coCartesian fibration over $\ul{G/H}$
\footnote{
  We think of a coCartesian fibration over $\ul{G/H}$ as representing an $H$-category, since the category $\ul{G/H}= (\OGop)_{/[G/H]}$ is equivalent to $\mathcal{O}_H^{op}$. 
}, 
which can be constructed as follows.
\begin{mydef} \label{def:underlying_G_cat}
  Let $\ul\C^\otimes \fib \GFin_*$ be a coCartesian fibration, and $I\in \GFin_*$ a $G$-set over $G/H$. 
  Define a coCartesian fibration \( \ul\C^\otimes_{<I>} \fib \ul{G/H}\) by pulling back $\ul\C^\otimes$ along $\sigma_{<I>}$,
  \begin{align*}
    \xymatrix{ 
      \ul\C^\otimes_{<I>} \ar@{->>}[d] \ar[r] \pullbackcorner & \ul\C^\otimes \ar@{->>}[d] \\
      \ul{G/H} \ar[r]^{\sigma_{<I>}} & \GFin_* .
    }
  \end{align*}
  In particular, for $I_+(G/G)=(G/G \xto{=} G/G)$, the terminal object of $\GFin_*$, denote by \( \ul\C := \ul\C^\otimes_{<I_+(G/G)>} \) the \emph{underlying $G$-category of $\ul\C^\otimes$}. 
\end{mydef}

\paragraph{An inert diagram in $\GFin_*$} 
Let $I=(U\to G/H)$ be a finite pointed $G$-set over $G/H$, as before. 
Applying the parametrized join construction for \( S=\ul{G/H} \) and the left fibrations \( \ul{G/H} \xto{=} \ul{G/H}, \,  \ul{U} \fib \ul{G/H} \)
\footnote{
  Since $\mathbf{Fin}^G$ is a category, a map of finite $G$-sets $U \to V$ induces a $G$-functor $\ul{U} \to \ul{V}$. 
  By comparison, a map \( f \colon x \to y \) in an $\infty$-category $C$ induces a span 
  \( \ul{x}=C_{/x} \xfrom{\sim} C_{/f} \to C_{/y} =\ul{y} \) where both arrows are left fibrations and the left arrow is an equivalence of $\infty$-categories. 
}
we get a coCartesian fibration \( \ul{U} \star_{\ul{G/H}} \ul{G/H} \fib \ul{G/H} \times \Delta^1 \), which we can consider as a $G$-category by composing with 
the coCartesian fibration \( \ul{G/H} \times \Delta^1 \fib \ul{G/H} \) and the left fibration \( \ul{G/H} \fib \OGop \).

For each $I\in \GFin_*$ we construct a $G$-functor \( \Phi_{<I>} \colon \ul{U} \star_{\ul{G/H}} \ul{G/H}  \to \GFin_* \) (a $G$-diagram in $\GFin_*$):
\begin{mydef}
  Let $I=(U\to G/H)$ be a finite pointed $G$-set over $G/H$. 
  We define a $G$-functor  
  \begin{align*}
    \xymatrix{
      \Phi_{<I>} \colon \ul{U} \star_{\ul{G/H}} \ul{G/H} \ar@{->>}[d] \ar[rr] & & \GFin_* \ar@{->>}[d] \\
      \ul{G/H} \times \Delta^1 \ar@{->>}[r] & \ul{G/H} \ar@{->>}[r] & \OGop 
    }
  \end{align*}
  by specifying its restrictions to \( \ul{U} \fib \ul{G/H}\times\set{0}\) and \(\ul{G/H}\fib \ul{G/H} \times \set{1} \), together with its action on morphisms over \( (id, 0 \to 1) \in \ul{G/H} \times \Delta^1 \):
  \begin{enumerate}
    \item The $G$-functor \( \ul{U} \to \GFin_* \) is the composition \( \diag{ \ul{U} \ar@{->>}[r] & \ul{G/H} \ar[r]^{\sigma_{<I>}} & \GFin_* } \), where the first map is the left fibration induced by \( U \to G/H \).
    \item The $G$-functor  \( \ul{G/H} \to \GFin_* \) is the composition \( \diag{ \ul{G/H} \ar@{->>}[r] & \ul{G/G} \ar[rr]^-{\sigma_{<I+(G/G)>}} & & \GFin_* } \), where the first map is the structure map \( \ul{G/H} \fib \OGop=\ul{G/G} \) and the second map is the $G$-functor corresponding to \( I_+(G/G)\) 
      (in fact, the composition is just $\sigma_{<I_+(G/H)>}$). 
    \item Let $( G/H \xfrom{\psi} G/K) \in \ul{G/H}$, then the fiber of \( \ul{U} \star_{\ul{G/H}} \ul{G/H} \to \ul{G/H} \times \Delta^1 \) over \( ( \set{\psi} ,0\to 1) \) is 
      \( \left( \ul{U} \star_{\ul{G/H}} \ul{G/H} \right)_{\psi} = \ul{U}_{\psi} \star \set{\psi} \),
      a co-cone diagram on the finite set of maps \( \varphi \colon G/H \to U \) such that \( \diag{ G/K \ar[r]^{\varphi} \ar[dr]_{\psi} & U \ar[d] \\ & G/H } \) commutes.
      Therefore, morphisms of \( \ul{U} \star_{\ul{G/H}} \ul{G/H} \fib \ul{G/H} \times \Delta^1\) over \( (id_{\psi},0\to 1)\in \ul{G/H} \times \Delta^1 \) are in bijection to $\varphi \colon U\to G/H$ making the above diagram commute.
      Let \( \bar{\varphi}  \colon  G/K \to \psi^* U \) be the unique map given by
      \begin{align*}
        \vcenter{\xymatrix{
          G/K \ar@/_/[ddr]_{=} \ar@{-->}[dr]^{\exists ! \bar{\varphi}} \ar@/^/[drr]^{\varphi} \\
           & \psi^* U \ar[r] \ar[d] \pullbackcorner & U \ar[d] \\
           & G/K \ar[r]^{\psi} & G/H \\
         }}
      \end{align*}

      The functor $\Phi_{<I>}$ sends the morphism over \( (id_\psi,0\to 1) \)  corresponding to $\varphi \colon U \to G/H$  to the span of finite pointed $G$-sets
      \begin{align*}
       \xymatrix{
         \psi^* U \ar[d] & \ar[l]_{\bar{\varphi}} G/K \ar[r]^{=}  \ar[d]^{=} & G/K \ar[d]^{=} \\
         G/K & \ar[l]_{=} G/K \ar[r]^{=} & G/K .
       }
      \end{align*}
      Using the fact that $\OG$ is atomic (i.e orbits have no non-trivial retracts) one can check that the left square is a summand inclusion.
  \end{enumerate}
  Steps 1 and 2 define $\Phi_{<I>}$ on every morphism over \( \ul{G/H} \times (0\to 1) \), since every such morphism uniquely decomposes as a morphism in $\ul{U}$ followed by a morphism over \( (\set{\psi},0\to 1) \) for some $\psi\in \ul{G/H}$.
  Verifying that $\Phi_{<I>}$ is well defined is a straightforward calculation,
  using the fact that every morphism of $\ul{U} \star_{\ul{G/H}} \ul{G/H}$ can be uniquely decomposed as a morphism over \( (\set{\psi'},0\to 1) \) followed by a morphism in  $\ul{G/H}$.
\end{mydef}

\paragraph{Construction of Segal maps and definition of a $G$-symmetric monoidal $\infty$-category} 
For any coCartesian fibration over $\GFin_*$ we construct 'Segal maps':
\begin{mydef}
  Let  \( \ul\C^\otimes \fib \GFin_* \) be a coCartesian fibration and $I=(U\to G/H)$ a finite pointed $G$-set over $G/H$.
  Construct a $G$-functor over $\ul{G/H}$ by the following steps:
  \begin{enumerate}
    \item Pulling $\ul\C^\otimes$ along $\Phi_{<I>}$ produces a coCartesian fibration \( (\Phi_{<I>})^* \ul\C^\otimes \fib \ul{U} \star_{\ul{G/H}} \ul{G/H} \),
      which we can consider as a coCartesian fibration over \( \ul{G/H} \times \Delta^1 \) by the composition 
      \begin{align} \label{Segal_classifier}
        (\Phi_{<I>})^* \ul\C^\otimes \fib \ul{U} \star_{\ul{G/H}} \ul{G/H} \fib \ul{G/H} \times \Delta^1 .
      \end{align}
    \item The restriction of the coCartesian fibration \eqref{Segal_classifier} to \( \ul{G/H} \times \set{0} \) is given by 
      \begin{align*}
        \ul{U} \times_{\ul{G/H}} \ul\C^\otimes_{<I>} \fib \ul{U} \fib \ul{G/H} \times \set{0} , 
      \end{align*}
      as it is the pullback of  \( \ul\C^\otimes \fib \GFin_* \) along  \( \diag{ \ul{U} \ar@{->>}[r] & \ul{G/H} \ar[r]^{\sigma_{<I>}} & \GFin_* } \).
    \item The restriction of the coCartesian fibration \eqref{Segal_classifier} to \( \ul{G/H} \times \set{1} \) is given by 
      \begin{align*}
        \ul{G/H} \ultimes \ul\C  \fib \ul{G/H} \xto{=} \ul{G/H} \times \set{1} , 
      \end{align*}
      as it is the pullback of  \( \ul\C^\otimes \fib \GFin_* \) along \( \diag{ \ul{G/H} \ar@{->>}[r] & \ul{G/G} \ar[rr]^-{\sigma_{<I+(G/G)>}} & & \GFin_* } \).
    \item Therefore, the coCartesian fibration \eqref{Segal_classifier} classifies a $G$-functor over $\ul{G/H}$
      \begin{align*}
        \xymatrix@=2ex{
          \ul{U} \times_{\ul{G/H}} \ul\C^\otimes_{<I>} \ar@{->>}[dr] \ar[rr] & & \ul{G/H} \ultimes \ul\C  \ar@{->>}[dl] \\
          & \ul{G/H},
        }
      \end{align*}
      which by \cite[thm. 9.7]{Expose1} is equivalent to a $G$-functor over $\ul{G/H}$
      \begin{align} \label{bare_Segal_map}
        \xymatrix@=2ex{
           \phi_{<I>} \colon  \ul\C^\otimes_{<I>} \ar@{->>}[dr] \ar[rr] & & \ulFun_{\ul{G/H}} ( \ul{U}, \ul{G/H} \ultimes \ul\C ) \ar@{->>}[dl] \\
          & \ul{G/H}.
        }
      \end{align}
      We call \eqref{bare_Segal_map} the Segal map of $I$. 
  \end{enumerate}
\end{mydef}
We can now give the definition of a $G$-symmetric monoidal $G$-category. 
\begin{mydef} \label{def:G_SM_cat}
  A \emph{$G$-symmetric monoidal $G$-category} is a coCartesian fibration \( \ul\C^\otimes \fib \GFin_* \) such that for every finite pointed $G$-set \( I = (U \to G/H) \) the Segal map \( \phi_{<I>} \) of \cref{bare_Segal_map} is an equivalence of $\ul{G/H}$-categories. 
\end{mydef}

\begin{rem} \label{rem:object_as_family}
  Let $\ul\C^\otimes \to \GFin_*$ be a $G$-symmetric monoidal $G$-category.
  The Segal conditions imply that an object $x\in \ul\C^\otimes$ over $I=(U\to G/H) \in \GFin_*$ classifies a $G$-functor \( x_\bullet \colon \ul{U} \to \ul\C \).
  To see this, first note that by Yoneda's lemma $x$ defines a $\ul{G/H}$ object \( \sigma_x  \colon  \ul{G/H} \to \ul\C^\otimes \).
  Since $x\in \ul\C^\otimes$ is over $I\in\GFin_*$ the composition \( \ul{G/H} \xto{\sigma_x} \ul\C^\otimes \to \GFin_* \) is equivalent to $\sigma_{<I>} \colon  \ul{G/H} \to \GFin_*$, so $\sigma_x$ factors as \( \sigma_x \colon  \ul{G/H} \to \ul\C^\otimes_{<I>}  \to \ul\C^\otimes \).
  Therefore we can regard $\sigma_x$ as a $\ul{G/H}$-object of \( \ul\C^\otimes_{<I>} \). 
  Using the Segal conditions we identify $\sigma_x$ with a $\ul{G/H}$-object of \(\ulFun_{\ul{G/H}} ( \ul{U}, \ul\C \ultimes \ul{G/H} ) \).
  Finally we use the equivalence  
  \begin{align*}
    \Fun_{\ul{G/H}} ( \ul{G/H}, \ulFun_{\ul{G/H}} ( \ul{U}, \ul\C \ultimes \ul{G/H} ) \simeq \Fun_{\ul{G/H}} ( \ul{U}, \ul\C \ultimes \ul{G/H} ) \simeq \Fun_G(\ul{U},\ul\C) . 
  \end{align*}
  to identify \( \sigma_x \colon  \ul{G/H} \to \ulFun_{\ul{G/H}} ( \ul{U}, \ul\C \ultimes \ul{G/H} ) \) with a $G$-functor \(  x_\bullet \colon \ul{U} \to \ul\C\). 
\end{rem}

\begin{rem} \label{rem:Segal_map_product}
  The codomain of the above Segal map is equivalent to a parametrized product:
  The ``internal hom'' $\ul{G/H}$-functor \( \ulFun_{\ul{G/H}}(\ul{U},-)  \colon  \Cat_\infty^{\ul{G/H}} \to \Cat_\infty^{\ul{G/H}} \) is right adjoint to the composition 
  \begin{align*}
    \xymatrix@=1ex{
      \Cat_\infty^{\ul{G/H}} \ar[rr] & & \Cat_\infty^{\ul{U}} \ar[rr] & &   \Cat_\infty^{\ul{G/H}}, \\
      (\ul\D \fib \ul{G/H}) \ar@{|->}[rr] & & (\ul\D \times_{\ul{G/H}} \ul{U} \fib \ul{U} ) \ar@{|->}[rr] & & (X \times_{\ul{G/H}} \ul{U} \fib \ul{U} \to \ul{G/H} ) .
    }
  \end{align*}
  Therefore, it decomposes as the composition of the right adjoints: 
  \begin{align*}
    \xymatrix@=1ex{
    }
  \end{align*}
  Under this equivalence, the Segal map of \( I=(U\to G/H) \) is given by 
  \begin{align} \label{Segal_map}
    \xymatrix@=2ex{
       \phi_{<I>} \colon  \ul\C^\otimes_{<I>} \ar@{->>}[dr] \ar[rr] & & \prod\limits_I \ul\C \ultimes \ul{U}  \ar@{->>}[dl] \\
      & \ul{G/H}.
    }
  \end{align}
  In particular, we can identify an object $x\in \ul\C^\otimes$ over $I$ with a $\ul{G/H}$-object of \( \prod\limits_I \ul\C \ultimes \ul{U} \) as follows.
  Since $\GFin_* \fib \OGop, \, I \mapsto [G/H]$ the object $x$ belongs to the fiber \( \ul\C^\otimes_{[G/H]} \), and by Yoneda's lemma is classified by a $G$-functor \( \sigma_x \colon  \ul{G/H} \to \ul\C^\otimes \). 
  Since $x\in \ul\C^\otimes$ is over $I\in \GFin_*$, the $G$-functor \( \ul{G/H} \to \ul\C^\otimes \fib \GFin_* \) classifies $I\in \GFin_*$, and is therefore equivalent to $\sigma_{<I>}$.
  Therefore it induces a $\ul{G/H}$-functor \( \ul{G/H} \to \ul\C^\otimes_{<I>}\).
  Post-composing with the Segal map of \cref{Segal_map} 
  we get our desired $\ul{G/H}$-object \(\ul{G/H} \to \ul\C^\otimes_{<I>} \xto{\sim} \prod\limits_I \ul\C \ultimes \ul{U} \), which by abuse of notation we also denote by $\sigma_x \colon  \ul{G/H} \to \prod\limits_I \ul\C \ultimes \ul{U}$. 
\end{rem}

Unpacking the construction of the Segal maps \eqref{bare_Segal_map} in \cref{def:G_SM_cat} gives the following fiberwise characterization of $G$-symmetric monoidal categories, which is easier to verify. 
\begin{lem} \label{G_SM_cat_Fiberwise_char}
  A coCartesian fibration $\ul\C^\otimes \fib \GFin_*$ is a $G$-symmetric monoidal category (\cref{def:G_SM_cat}) if and only if for each finite pointed $G$-set \( J=(V \to G/K) \in \GFin_* \) the functor 
    \[ \ul\C^\otimes_J \to \prod_{W\in\orb(\psi^* U)} \ul\C_{[W]}\]
    is an equivalence of $\infty$-categories,  
    where $\ul\C^\otimes_J$ is the fiber of \( \ul\C^\otimes \fib \GFin_* \) over $J=(V \to G/K)$
    and the above functor is the product of  \( \ul\C^\otimes_J \to \ul\C_{[W]} \) associated to the $\GFin_*$ edges
    \begin{align*}
     \forall W\in \orb(V):\quad 
     \xymatrix{
       V \ar[d] & \ar[l] W \ar[r]^{=}  \ar[d]^{=} & W \ar[d]^{=} \\
       G/K & \ar[l] W \ar[r]^{=} & W .
     }
    \end{align*}
\end{lem}
\begin{proof}
  The Segal condition of $G$-symmetric monoidal $G$-categories states that the Segal map \(\phi_{<I>} \) is a parametrized equivalence, i.e for each $(G/H \xfrom{\psi} G/K ) \in \ul{G/H}$, the Segal map \(\phi_{<I>} \) induces an equivalence between the fibers
  \[ (\ul\C^\otimes_{<I>})_{[\psi]} \to  \ul{\Fun}_{\ul{G/H}}(\ul{U}, \ul{G/H} \ultimes \ul\C)_{[\psi]} . \] 
  The fiber of $\ul\C^\otimes_{<I>}$ over $\psi$ is the fiber of \( \ul\C^\otimes \fib \GFin_* \) over the finite pointed $G$-set \( J:= (\psi^* U \to G/K ) \).
  The fiber of $ \ul{\Fun}_{\ul{G/H}}(\ul{U}, \ul{G/H} \ultimes \ul\C)$ over $\psi$ is the $\infty$-category of $G$-functors \( \Fun_{\OGop}(\ul{\psi^* U}, \ul\C) \). 
  Decomposing the finite $G$-set \( \psi^* U = \coprod_{W\in\orb(\psi^* U)} W \) into orbits we have 
  \[ \Fun_{\OGop}(\ul{\psi^* U}, \ul\C) \= \Fun_{\OGop}(\coprod \ul{W}, \ul\C) \simeq \prod_W \Fun_{\OGop}( \ul{W}, \ul\C) \simeq \prod_{W\in\orb(\psi^* U)} \ul\C_{[W]}. \]
  Since both sides depend only on \( J=(\psi^* U \to G/K)\in\GFin_* \) the result follows. 
\end{proof}

We end this appendix with the definition of parametrized tensor product functors in a $G$-symmetric monoidal category.
\begin{mydef} \label{def:param_tensor_prod}
  Let $\ul\C^\otimes \fib \ulFin_*^G$ be a $G$-symmetric monoidal category.
  Let $I = (U\to G/H), \, J = (V \to G/H) \in \ulFin_*^G$ be two object over the orbit $G/H$, and \( f \colon I \to J \) a morphism in  $(\ulFin_*^G)_{[G/H]}$, given by 
  \begin{align*}
   \xymatrix{
     U \ar[d] & \ar[l]_{=} U \ar[d] \ar[r]^{f} & V \ar[d] \\
     G/H & \ar[l]_{=} G/H \ar[r]^{=} & G/H .
   }
  \end{align*}
  The  morphism $f$ corresponds to a functor $\Delta^1 \to (\ulFin_*^G)_{[G/H]}$, or equivalently to a $G$-functor \( \sigma_{<f>} \colon  \ul{G/H} \times \Delta^1 \to \ulFin_*^G \), which restricts to $\sigma_{<I>}$ over $\ul{G/H}\times\set{0}$ and to $\sigma_{<J>}$ over $\ul{G/H}\times\set{1}$.

  Pulling back $\ul\C^\otimes \fib \ulFin_*^G$ along $\sigma_{<f>}$ we get a coCartesian fibration \( \ul\C^\otimes_{<f>} \fib \ul{G/H} \times \Delta^1 \) which restricts to $\ul\C^\otimes_{<I>}$ over $\ul{G/H}\times\set{0}$ and to $\ul\C^\otimes_{<J>}$ over $\ul{G/H}\times\set{1}$.
  Therefore this coCartesian fibration classifies a $\ul{G/H}$-functor  
  \begin{align*} 
    \xymatrix@=2ex{
       \otimes_f  \colon  \ul\C^\otimes_{<I>} \ar@{->>}[dr] \ar[rr] & & \ul\C^\otimes_{<J>}  \ar@{->>}[dl] \\
      & \ul{G/H}
    }
  \end{align*}
  which we refer to as the tensor product over $f$.
  Composing with the Segal maps of \cref{Segal_map},
  we can rewrite the tensor product over $f$ as 
  \begin{align*} 
    \xymatrix@=2ex{
       \otimes_f  \colon  \prod\limits_I \ul\C \ultimes \ul{U} \ar@{->>}[dr] \ar[rr] & & \prod\limits_J \ul\C \ultimes \ul{V}  \ar@{->>}[dl] \\
      & \ul{G/H}.
    }
  \end{align*}

\end{mydef}

\section{Mapping spaces in over-categories} \label{sec:maping_spaces}

We prove some simple properties of mapping spaces in over categories.

\begin{lem} \label{lem1:over_mappin_spaces}
  Consider the over category $\C_{/b}$ for $b$ an object in an $\infty$-category $\C$.
  Let $x\to b,y_1 \to b, y_2 \to b$ be objects in $\C_{/b}$, and a morphism $\varphi$ in $\C_{/b}$ from $y_1 \to b$ to $y_2 \to b$.
  Then 
  \begin{align*}
    \xymatrix{
      \Map_{\C_{/b}}( x \to b , y_1 \to b ) \ar[r] \ar[d]^{\varphi_*} & \Map_{\C} (x,y_1) \ar[d]^{\varphi_*} \\
      \Map_{\C_{/b}}( x \to b, y_2 \to b) \ar[r] & \Map_{\C}(x,y_2)
    }
  \end{align*}
  is a homotopy pullback square.
\end{lem}
\begin{proof}
  The mapping space $\Map_{\C_{/b}}( x \to b, y \to b)$ is the homotopy fiber of the postcomposition map $ \Map_{\C} (x,y) \to \Map_{\C} (x,b) $.
  Therefore the lower square and outer rectangle in the diagram
  \begin{align*}
    \xymatrix{
      \Map_{\C_{/b}}( x \to b , y_1 \to b ) \ar[r] \ar[d]^{\varphi_*} & \Map_{\C} (x,y_1) \ar[d]^{\varphi_*} \\
      \Map_{\C_{/b}}( x \to b, y_2 \to b) \ar[r] \ar[d] & \Map_{\C} (x,y_2) \ar[d]^{(y_2 \to b)_*} \\
      \ast \ar[r]^{x\to b} & \Map_{/\C} (x,b) 
    } 
  \end{align*}
  are homotopy pullback diagram.
  It follows that the top square is a homotopy pullback square.
\end{proof}

\begin{lem} \label{lem2:over_mappin_spaces}
  Let $f \colon b \to b'$ be a morphism in an $\infty$-category $\C$, and consider the postcomposition functor \( f_* \colon  \C_{/b} \to \C_{/b'} \).
  Let $x\to b,y_1 \to b, y_2 \to b$ be objects in $\C_{/b}$, and a morphism $\varphi$ in $\C_{/b}$ from $y_1 \to b$ to $y_2 \to b$.
  Then
  \begin{align*}
    \xymatrix{
      \Map_{\C_{/b}}( x \to b , y_1 \to b ) \ar[r]^-{f_*} \ar[d]^{\varphi_*} & \Map_{\C_{/b'}}( x \to b \xto{f} b', y_1 \to b \xto{f} b') \ar[d]^{\varphi_*}  \\
      \Map_{\C_{/b}}( x \to b, y_2 \to b) \ar[r]^-{f_*} & \Map_{\C_{/b'}}( x \to b \xto{f} b', y_2 \to b \xto{f} b') 
    }
  \end{align*}
  is a homotopy pullback square.
\end{lem}
\begin{proof}
  Consider the commutative diagram
  \begin{align*}
    \xymatrix{
      \Map_{\C_{/b}}( x \to b , y_1 \to b ) \ar[r]^-{f_*} \ar[d]^{\varphi_*} & \Map_{\C_{/b'}}( x \to b \xto{f} b', y_1 \to b \xto{f} b') \ar[r] \ar[d]^{\varphi_*} & \Map_{\C} (x,y_1) \ar[d]^{\varphi_*} \\
      \Map_{\C_{/b}}( x \to b, y_2 \to b) \ar[r]^-{f_*} & \Map_{\C_{/b'}}( x \to b \xto{f} b', y_2 \to b \xto{f} b') \ar[r] & \Map_{\C}(x,y_2) .
    }
  \end{align*}
  By \cref{lem1:over_mappin_spaces} the right square and the outer rectangle are homotopy pullback squares, hence the left square is a homotopy pullback square.
\end{proof}

Next, let $f \colon b \to b'$ be a morphism in an $\infty$-category $\C$ as before, and $T \colon \M \to \C_{/b'}$ a functor of $\infty$-categories. 
Define an $\infty$-category $\M_T$ as the pullback 
\begin{align*}
  \xymatrix{
    \M_T \ar[d]^{u} \ar[r] & \C_{/b} \ar[d]^{f_*} \\
    \M \ar[r]^T & \C_{/b'} .
  }
\end{align*}

\begin{lem} \label{lem3:over_mappin_spaces}
  Let $\C, f \colon b \to b', T \colon \M \to \C_{/b'}$ and $\M_T$ be as above.
  Let $X,Y_1,Y_2$ be objects in $\M_T$ and $\varPhi \colon Y_1 \to Y_2 $ be morphism in $\M_T$.
  Then
  \begin{align*}
    \xymatrix{
      \Map_{\M_T} (X,Y_1) \ar[r]^{\varPhi_*} \ar[d] & \Map_{\M_T} (X,Y_2) \ar[d] \\
      \Map_{\M} ( u (X) , u (Y_1) ) \ar[r]^{u(\varPhi)_*} & \Map_{\M} ( u (X), u (Y_2) )
    }
  \end{align*}
  is a homotopy pullback square.
\end{lem}
\begin{proof}
  Denote the images of $X,Y_1,Y_2 \in \M_F$ in $\C_{/b}$ by $x \to b,\, x \to y_1, x \to y_2$, and the image of $\varPhi$ by $\varphi$.
  Using the equivalences  \( u(X) \simeq (x\to b \xto{f} b'),  u(Y_1) \simeq (y_1 \to b \xto{f} b'),  u(Y_2) \simeq (y_2 \to b \xto{f} b') \) in $\C_{/b'}$ 
  we can identify the mapping spaces 
  \begin{align*}
    \Map_{\C_{/b'}} ( T u (X), T u (Y_i) ) \simeq \Map_{\C_{/b'}}(x \to b \xto{f} b', y_i \xto{f} b') , \quad i=1,2 .
  \end{align*}
  Under these identifications we have a commutative diagram
  \[ 
    \begin{tikzcd}
      \Map_{\M_T} (X,Y_1) \ar{r} \ar{d} & \Map_{\M} ( u (X) , u (Y_1) ) \ar{d}{T} \\
      \Map_{\C_{/b}}(x \to b, y_1 \to b) \ar{d}{\varphi_*} \ar{r}{f_*} & \Map_{\C_{/b'}}(x \to b \xto{f} b', y_1 \to b \xto{f} b')  \ar{d}{\varphi_*} \\
      \Map_{\C_{/b}}(x \to b, y_2 \to b) \ar{r}{f_*} & \Map_{\C_{/b'}}(x \to b \xto{f} b', y_2 \xto{f} b') .
    \end{tikzcd}
  \]
   The top square is a homotopy pullback square by the definition of $\M_T$ as a pullback, and the bottom square is a homotopy pullback square by \cref{lem2:over_mappin_spaces}.
   Therefore the outer rectangle is a homotopy pullback square. 
   On the other hand, this is also the outer rectangle in the diagram
  \begin{align*}
    \xymatrix{
      \Map_{\M_T} (X,Y_1) \ar[r]^{\varPhi_*} \ar[d] & \Map_{\M_T} (X,Y_2) \ar[d] \ar[r] & \Map_{\C_{/b}}(x \to b, y_2 \to b) \ar[d]^{f_*} \\
      \Map_{\M} ( u (X) , u (Y_1) ) \ar[r]^{u(\varPhi)_*} & \Map_{\M} ( u (X), u (Y_2) ) \ar[r]^-{T} & \Map_{\C_{/b'}}(x \to b \xto{f} b', y_2 \xto{f} b') .
    }
   \end{align*}
   By definition of $\M_T$ is a pullback the right square is a homotopy pullback square, hence the left square is a homotopy pullback square, as claimed. 
\end{proof}

\end{appendices}

\addcontentsline{toc}{section}{References}
\bibliographystyle{alpha}
\bibliography{references}
\end{document}